\documentclass[12pt,reqno]{amsart}
\usepackage{amsmath,amssymb,latexsym, amsthm, amscd, mathrsfs, stmaryrd} %, txfonts}
\usepackage[linktocpage=true]{hyperref}
\hypersetup{colorlinks,linkcolor=blue,urlcolor=cyan,citecolor=blue}

\usepackage[all]{xy}
\usepackage{hyperref}
\usepackage{color}
\newcommand{\nc}{\newcommand}
\nc{\browntext}[1]{\textcolor{brown}{#1}}
\nc{\greentext}[1]{\textcolor{green}{#1}}
\nc{\redtext}[1]{\textcolor{red}{#1}}
\nc{\bluetext}[1]{\textcolor{blue}{#1}}
\nc{\brown}[1]{\browntext{ #1}}
\nc{\green}[1]{\greentext{ #1}}
\nc{\red}[1]{\redtext{ #1}}
\nc{\blue}[1]{\bluetext{ #1}}
\nc{\zb}[1]{\redtext{From zb: #1}}

\newtheorem{theorem}{Theorem}  [section]
\newtheorem{corollary}[theorem]{Corollary}
\newtheorem{lemma}[theorem]{Lemma}

\newtheorem{proposition}[theorem]{Proposition}
\newtheorem{definition}[theorem]{Definition}

\theoremstyle{remark}
\newtheorem{remark}[theorem]{Remark}

\numberwithin{equation}{section}

\def\haTh{\widehat{\Theta}}
\def \haH{\widehat{H}}

\def \bt{\mathbf t}
\def \bn{\mathbf n}

\newcommand{\cc}{{\mathcal C}}

\def \dB{\Theta}

\def \co{\mathcal O}

\def \ch{{\mathcal H}}
\def \cm{{\mathcal M}}
\def \ct{{\mathcal T}}
\def \bS{\mathbf S}

\def \haB{\widehat{B}}
\newcommand{\tMHX}{{}^\imath\widetilde{\ch}(\X_\bfk)}
\newcommand{\tCMH}{{}^\imath\widetilde{\cc}(\bfk Q)}
\newcommand{\tCMHX}{{}^\imath\widetilde{\cc}(\X_\bfk)}
\renewcommand{\mod}{\operatorname{mod}\nolimits}
\newcommand{\tCMHC}{{}^\imath\widetilde{\cc}(\bfk C_n)}

\numberwithin{equation}{section}

\newtheorem{innercustomthm}{{\bf Theorem}}
\newenvironment{customthm}[1]
{\renewcommand\theinnercustomthm{#1}\innercustomthm}
{\endinnercustomthm}

\renewcommand{\ker}{\operatorname{Ker}\nolimits}

\newcommand{\Hom}{\operatorname{Hom}\nolimits}
\newcommand{\RHom}{\operatorname{RHom}\nolimits}

\newcommand{\Aut}{\operatorname{Aut}\nolimits}
\newcommand{\Id}{\operatorname{Id}\nolimits}

\newcommand{\Ext}{\operatorname{Ext}\nolimits}

\newcommand{\add}{\operatorname{add}\nolimits}

\def \y{{B}}
\def \haB{\widehat{B}}

\newcommand{\mbf}{\mathbf}

\newcommand{\mrm}{\mathrm}

\newcommand{\End}{\mrm{End}}
\newcommand{\rank}{\mrm{rank}}
\newcommand{\de}{\delta}
\def \C{{\mathbb C}}

\newcommand{\N}{\mathbb N}
\newcommand{\bbZ}{\mathbb Z}

\newcommand{\qbinom}[2]{\begin{bmatrix} #1\\#2 \end{bmatrix} }
\newcommand{\Q}{\mathbb Q}
\newcommand{\sll}{\mathfrak{sl}}
\newcommand{\T}{\texttt{\rm T}}

\newcommand{\U}{\mbf U}
\newcommand{\K}{\mathbb K}
\newcommand{\F}{\mathbb F}

\newcommand{\bs}{\mathbf s}
\newcommand{\arxiv}[1]{\href{http://arxiv.org/abs/#1}{\tt arXiv:\nolinkurl{#1}}}

\newcommand{\Ui}{{\mbf U}^\imath}

\newcommand{\Z}{\mathbb Z}

\newcommand{\TT}{\mathbf T}

\def \X{\mathbb X}
\def \fg{\mathfrak{g}}

\def \I{\mathbb{I}}

\def \cv{\mathcal V}
\def \cu{\mathcal U}

\newcommand{\tUiD}{{}^{\text{Dr}}\tUi}
\def \nua{a}

\newcommand{\tUi}{\widetilde{{\mathbf U}}^\imath}
\newcommand{\sqq}{{\bf v}}

\newcommand{\coker}{\operatorname{Coker}\nolimits}

\newcommand{\tU}{\widetilde{\mathbf U}}

\def \cl{L}

\def \K{\mathbb{K}}
\def \R{\mathbb{R}}
\def \SS{\mathbb{S}}

\def \BF{\digamma}

\def \cc{\mathcal C}
\def\ca{\mathcal A}
\newcommand{\Iso}{\operatorname{Iso}\nolimits}
\renewcommand{\Im}{\operatorname{Im}\nolimits}
\newcommand{\res}{\operatorname{res}\nolimits}
\newcommand{\iH}{{}^\imath\widetilde{\ch}}
\newcommand{\Mod}{\operatorname{Mod}\nolimits}
\newcommand{\coh}{\operatorname{coh}\nolimits}
\newcommand{\rep}{\operatorname{rep}\nolimits}
\newcommand{\Ker}{\operatorname{Ker}\nolimits}

\def \PL{\mathbb{P}^1_{\bfk}}

\def \scrf{\mathscr F}
\def \scrt{\mathscr T}
\def \P{\mathbb P}
\def \cI{\mathcal I}

\def \cs{{\mathcal{S}}}
\def \ch{\mathcal H}
\def \cd{\mathcal D}
\def\bfk{\mathbf{k}}
\def \ck{\mathcal K}
\def \bp{\mathbf p}
\def \ul{\underline}
\def \fp{\mathfrak p}

\def \II{\I_0}
\def \Lg{L\fg}

\def \bvs{{\boldsymbol{\varsigma}}}

\def \cn{\mathcal N}

\newcommand{\calc}{{\mathcal C}}
\newcommand{\haT}{\widehat{\Theta}}

\newcommand{\wt}{\text{wt}}
\def \cR{\mathcal R}
\def \ce{\mathcal E}

\def \bla{\boldsymbol{\lambda}}
\def \blx{x}

\allowdisplaybreaks

\begin{document}

	%%%%%
	\title[$\imath$Hall algebras of weighted projective lines and quantum symmetric pairs]
	{$\imath$Hall algebras of weighted projective lines and quantum symmetric pairs}
	
	\author[Ming Lu]{Ming Lu}
	\address{Department of Mathematics, Sichuan University, Chengdu 610064, P.R.China}
	\email{luming@scu.edu.cn}

	\author[Shiquan Ruan]{Shiquan Ruan}
	\address{ School of Mathematical Sciences,
		Xiamen University, Xiamen 361005, P.R.China}
	\email{sqruan@xmu.edu.cn}

	%\author[Weiqiang Wang]{Weiqiang Wang}
	%\address{Department of Mathematics, University of Virginia, Charlottesville, VA 22904, USA}
	%\email{ww9c@virginia.edu}

	
	
	\begin{abstract}
		The $\imath$Hall algebra of a weighted projective line is defined to be the semi-derived Ringel-Hall algebra of the category of $1$-periodic complexes of coherent sheaves on the weighted projective line over a finite field. We show that this Hall algebra provides a realization of the $\imath$quantum loop algebra, which is a generalization of the $\imath$quantum group arising from the quantum symmetric pair of split affine type ADE in its Drinfeld type presentation. The $\imath$Hall algebra of the $\imath$quiver algebra of split affine type A was known earlier to realize the same algebra in its Serre presentation. We then establish a derived equivalence which induces an isomorphism of these two $\imath$Hall algebras, explaining the isomorphism of the $\imath$quantum group of split affine type A under the two presentations.
	\end{abstract}
	
	\maketitle
	\setcounter{tocdepth}{1}
	\tableofcontents

	%\newpage
	%%%%%%%
	%%%%%%%
	
	%\newpage
	%%%%%%%
	%%%%%%%
	\section{Introduction}
	
	\subsection{Backgrounds}
	
	\subsubsection{Quantum groups and Hall algebras}
	Ringel \cite{Rin90} in 1990s used the Hall algebra associated to a Dynkin quiver $Q$ over a finite field $\mathbb F_q$ to realize half a quantum group $\U^+ =\U^+_v(\fg)$; see Green \cite{Gr95} for an extension to arbitrary acyclic quivers. Later, inspired by Ringel's construction, Lusztig \cite{Lus90} gave a geometric realization of $\U^+$, and additionally constructed its canonical basis \cite{Lus90}. These constructions can be regarded as the earliest examples of categorifications of half parts of quantum groups.
	
	The Hall algebra realization of the whole quantum group was achieved by Bridgeland \cite{Br13} who in fact provided a Hall algebra construction of the Drinfeld double $\tU$ of a quantum group $\U$. %; see Ringel \cite{Rin90}, Lusztig \cite{Lus90} and Green \cite{Gr95} on Hall algebra constructions of half a quantum group $\U^+$.
	As a further generalization of Bridgeland's construction, the first author and Peng \cite{LP16} formulated the {\em semi-derived Ringel-Hall algebras} for arbitrary hereditary abelian categories via $2$-periodic complexes (also called $\Z_2$-graded complexes); compare with the semi-derived Hall algebras defined by Gorsky for {\em Frobenius categories} \cite{Gor18}. The construction of semi-derived Ringel-Hall algebras was further extended to 1-Gorenstein algebras in \cite[Appendix A]{LW19a}. %; compare with the  {\em semi-derived Hall algebras} introduced by Gorsky \cite{Gor18} for exact cateRingwelgories with enough projective objects.

	There has been a current realization of the affine quantum groups formulated by Drinfeld \cite{Dr88, Be94, Da15}, which plays a crucial role on (algebraic and geometric) representation theory. The Hall algebra of the projective line was studied in a visionary paper by Kapranov \cite{Ka97} and then extended by Baumann-Kassel \cite{BKa01} to realize the current half of quantum affine $\sll_2$. The Hall algebra of a weighted projective line was developed in \cite{Sch04} to realize half an affine quantum group of ADE type, which was then upgraded to the whole quantum group via Drinfeld double techniques \cite{DJX12, BS13} (in fact, they considered more general simply-laced quantum loop algebras of star-shaped graphs). The semi-derived Ringel-Hall algebras formulated in \cite{LP16} are isomorphic to the Drinfeld double Ringel-Hall algebras, and so a geometric realization of quantum loop algebras in Drinfeld presentation was achieved in some sense based on \cite{DJX12} by considering the categories of $2$-periodic complexes over weighted projective lines; see \cite[Remark 4.16]{LP16}.

	\subsubsection{$\imath$Quantum groups and $\imath$Hall algebras}
	The $\imath$quantum groups $\Ui$ arising from quantum symmetric pairs $(\U, \Ui)$ associated to Satake diagrams \cite{Let99} (see \cite{Ko14}) can be viewed as a vast generalization of Drinfeld-Jimbo quantum groups associated to Dynkin diagrams; see the survey \cite{W22} and references therein. 
	
	The universal $\imath$quantum group $\tUi =\langle B_i, \K_i^{\pm1} \mid i\in \I \rangle$ is by definition a subalgebra of $\tU$, and $(\tU, \tUi)$ forms a quantum symmetric pair; cf. \cite{LW19a}. %A novelty of $\tUi$ is that it admits a Cartan subalgebra generated by $\tk_i$ $(i\in \I =Q_0)$, which contains various central elements.
	A central reduction of $\tUi$ reproduces the $\imath$quantum group $\Ui =\Ui_\bvs$ (with parameters $\bvs$) which arises from the construction of quantum symmetric pairs $(\U, \Ui)$ by Letzter \cite{Let99}. %Note $\tUi$ is naturally $\Z \I$-graded while $\Ui$ (with its inhomogeneous $\imath$Serre relations) is not.
	Recently, Lu and Wang \cite{LW20b} established a Drinfeld type presentation for the $\imath$quantum groups arising from quantum symmetric pairs of split affine ADE type.

	In recent years, the first author and Wang \cite{LW19a, LW20a} have developed $\imath$Hall algebras of $\imath$quivers to realize the (universal) quasi-split $\imath$quantum groups of Kac-Moody type. %A universal $\imath$quantum group admits a Serre type presentation and contains various central generators, which replace the parameters in an $\imath$quantum group arising from quantum symmetric pairs \`a la G.~Letzter \cite{Let99, Ko14}.
	The $\imath$Hall algebras are constructed in the framework of semi-derived Ringel-Hall algebras of 1-Gorenstein algebras (\cite[Appendix A]{LW19a}, \cite{LW20a}). 
	In particular, Bridgeland's Hall algebras realization of quantum groups is reformulated as $\imath$Hall algebras of $\imath$quivers of diagonal type; see \cite{LW19a}.
	
	Together with Wang, we  \cite{LRW20a} used the $\imath$Hall algebra of the projective line to give a geometric realization of the $q$-Onsager algebra (i.e., $\imath$quantum group of split affine $\mathfrak{sl}_2$) in its Drinfeld type presentation, a generalization of the geometric realization of the current half of quantum affine $\mathfrak{sl}_2$ via the Hall algebra of the projective line \cite{Ka97,BKa01}.
	
	We view Letzter's $\imath$quantum groups and our universal $\imath$quantum groups as a vast generalization of the Drinfeld-Jimbo quantum groups. The $\imath$-program as outlined by  Bao and Wang \cite{BW18} aims at generalizing various fundamental constructions from quantum groups to $\imath$quantum groups; in addition to the works mentioned above, see \cite{BW18, BK19, BW18b} for generalizations of (quasi) R-matrix and canonical bases. %, and also see \cite{BKLW,LW19c} (and \cite{Li19}) for geometric realizations and
	%\cite{BSWW} for KLR type categorification of a class of (modified) $\Ui$.

	\subsection{Goal}
	
	This is a sequel of \cite{LRW20a} in our program devoted to developing a geometric realization of the $\imath$quantum groups in Drinfeld type presentations, a vast generalization of \cite{Ka97, BKa01,Sch04,DJX12,BS13}. Explicitly, we use the $\imath$Hall algebras of weighted projective lines to give a geometric realization of the {so-called} simply-laced $\imath$quantum loop algebra of a star-shaped graph in Drinfeld type presentation (a generalization of $\imath$quantum groups arising from quantum symmetric pairs of split affine ADE type considered in \cite{LW20b}).
	We further show that the isomorphism of $\imath$quantum groups of split affine type A in two (Serre vs Drinfeld) presentations is induced from a derived equivalence of the categories underlying the two (quivers vs weighted projective lines) $\imath$Hall algebra realizations.
	
	The construction in this paper should be viewed as a geometric counterpart of the $\imath$Hall algebra construction in \cite{LW19a,LW20a}. 
	It is our hope that this work will stimulate further connection between the theory of weighted projective lines and the theory of  quantum symmetric pairs.
	%It is our hope that this work will stimulate further interactions between communities on weighted projective lines and on quantum symmetric pairs.
	
	%\red{TBA}
	
	\subsection{Main results}
	
	\subsubsection{}
	Let $\Gamma:=\T_{p_1,\dots,p_\bt}$ be a star-shaped graph with its vertex set $\II$ as in \eqref{star-shaped}, $\fg$ be its associated Lie algebra, $\tUiD$ be its associated $\imath$quantum loop algebra \cite{LW20b} defined as in Definition \ref{def:iDR}. Note that if $\Gamma$ is of ADE type, then $\tUiD$ is the Drinfeld type presentation of the $\imath$quantum group arising from quantum symmetric pairs of split affine ADE type; see Theorem \ref{thm:ADE}. Recall that $\tUiD$ is generated by $B_{\star,l},  B_{[i,j],l}, H_{\star,m},H_{[i,j],m}, \Theta_{\star,m},\Theta_{[i,j],m}$ and central elements $\K_\star^{\pm1}, \K_{[i,j]}^{\pm}, C^{\pm1}$ for $l\in\Z,m>0$, and $1\leq i\leq \bt$, $1\leq j\leq p_i-1$, subject to \eqref{iDR1a}--\eqref{exp h}.

	Let $\X_\bfk$ be the weighted projective line over a finite field $\bfk=\F_q$ associated to $\Gamma$. Let $\sqq=\sqrt{q}$. The $\imath$Hall algebra $\iH(\X_\bfk)$ considered in this paper is the twisted semi-derived Ringel-Hall algebra of the category of $1$-periodic complexes over $\coh(\X_\bfk)$; cf. \cite{LP16,LW19a}.

	\begin{customthm}{\bf A} [Theorem~\ref{thm:morphi}]
		There exists a $\Q(\sqq)$-algebra homomorphism
		\begin{align}
			\Omega: \tUiD_{ |v=\sqq}\longrightarrow \iH(\X_\bfk),
		\end{align}
		which sends
		\begin{align}
			\label{eq:THA1}
			&\K_{\star}\mapsto [K_{\co}], \qquad \K_{[i,j]}\mapsto [K_{S_{ij}}], \qquad C\mapsto [K_\de];&
			\\
			\label{eq:THA2}
			&B_{\star,l}\mapsto \frac{-1}{q-1}[\co(l\vec{c})],\qquad\Theta_{\star,r} \mapsto \widehat{\Theta}_{\star,r}, \qquad H_{\star,r} \mapsto \widehat{H}_{\star,r};\\
			\label{eq:THA3}
			&\y_{[i,j],l}\mapsto \frac{-1}{q-1}\haB_{[i,j],l}, \quad \Theta_{[i,j],r}\mapsto\widehat{\Theta}_{[i,j],r}, \quad H_{[i,j],r}\mapsto \widehat{H}_{[i,j],r}
		\end{align}
		for any $[i,j]\in\II-\{\star\}$, $l\in\Z$, $r>0$. %Moreover, $\Omega$ is an embedding if $\fg$ is of finite or affine type.
	\end{customthm}
	
	Let us explain the counterparts of the Drinfeld type generators in $\iH(\X_\bfk)$. First, $\widehat{\Theta}_{\star,r}$ and $\haH_{\star,r}$ are defined by using the corresponding ones $\widehat{\Theta}_r$, $\haH_{r}$ defined in $\iH(\P^1_\bfk)$ \cite{LRW20a} via the natural embedding (see \eqref{the embedding functor F on algebra})
	$$F_{\X,\P^1}:\iH(\P^1_\bfk) \longrightarrow \iH(\X_\bfk).$$
	Second, for $1\leq i\leq \bt$, there is an algebra embedding $$\widetilde{\psi}_{C_{p_i}}: \tUi_\sqq(\widehat{\mathfrak{sl}}_{p_i})\longrightarrow \iH(\bfk C_{p_i})$$ by \cite[Theorem 9.6]{LW20a}, where $C_{p_i}$ is the cyclic quiver with $p_i$ vertices. As a consequence, each branch of $\Gamma$  corresponds to a subalgebra of $\iH(\X_\bfk)$ isomorphic to
	$\tUiD_v(\widehat{\mathfrak{sl}}_{p_i})$ by the composition $\Omega_{C_{p_i}}$ of $\widetilde{\psi}_{C_{p_i}}$ and the isomorphism of two presentations \cite{LW20b} (see \eqref{the map Psi A}) $$\Phi:\tUiD_v(\widehat{\mathfrak{sl}}_{p_i})\stackrel{\cong}{\longrightarrow} \tUi_v(\widehat{\mathfrak{sl}}_{p_i}).$$
	On the other hand, the subcategory $\scrt_{\bla_i}$ of $\coh(\X_\bfk)$, consisting of torsion sheaves supported at the distinguished point $\bla_i$, is isomorphic to $\rep^{\rm nil}_\bfk( C_{p_i})$, the category of finite-dimensional nilpotent representations of  $C_{p_i}$ over $\bfk$. This isomorphism induces an algebra embedding (see \eqref{eq:embeddingx})
	$$\iota_i: \iH(\bfk C_{p_i})\longrightarrow\iH(\X_\bfk).$$
	So we define $\haB_{[i,j],l},\widehat{\Theta}_{[i,j],r},\widehat{H}_{[i,j],r}$ to be the images of the Drinfeld generators of $\tUiD_v(\widehat{\mathfrak{sl}}_{p_i})$ under the composition $\iota_i\circ \Omega_{C_{p_i}}$; see \eqref{def:haBThH}.
	To show that $\Omega: \tUiD_{_{ |v=\sqq}} \rightarrow \tMHX$ is a homomorphism, we must verify the Drinfeld type relations \eqref{iDR1a}--\eqref{exp h} for the $\imath$quantum group $\tUiD$.
	
	The counterparts in $\tMHX$ of the relations at the star point $\star$ follow from \cite[Theorem A]{LRW20a}, and the counterparts in $\tMHX$ of the relations involving all the vertices $[i,j]$ follow from the above definitions and the fact any two torsion sheaves supported at distinct points have zero Hom and $\Ext^1$-spaces.
	
	In order to check the remaining Drinfeld type relations in $\iH(\X_\bfk)$, the key part is to verify the relations between $\star$ and $[i,j]$ for $1\leq i\leq \bt$, especially the ones between $\star$ and $[i,1]$, which consists of plenty of highly non-trivial computations in \S\ref{subsec:relationsstarjneq1}--\S\ref{sec:Relationsstari1 II}. For these relations we have to determine some of the Drinfeld generators in $\iH(\X_\bfk)$ clearly, especially the ones at $[i,1]$; see \S\ref{sec:rooti1}.

	For each $r>0$ and $1\leq i\leq \bt$, $\cm_{r\delta\pm\alpha_{i1}}$ is defined in \eqref{eq:Mrde+alpha11}, and $\cm_{i,r\de}$  is defined in \eqref{def:Mr}.
	\begin{customthm}{\bf B}[Proposition~\ref{prop:realroot}, Proposition~\ref{prop:imageroot}]
		%\label{prop:realroot}
		For any $r>0$ and $1\leq i\leq \bt$,  we have
		\begin{align}
			%\label{realroot}
			\haB_{[i,1],r}= &(q-1)\sum\limits_{[M]\in\cm_{r\delta+\alpha_{i1}}}\bn(\ell(M)-1)[\![M]\!]; \\
			\haB_{[i,1],-r}=&(1-q)\sum\limits_{[M]\in\cm_{r\delta-\alpha_{i1}}}\bn(\ell(M)-1)[\![M]\!]*[K_{-r\de+\alpha_{i1}}];
			\\
			\haTh_{[i,1],r}=
			&\frac{\sqq}{q-1} \sum\limits_{|\lambda|=r}\bn(\ell(\lambda))  [\![S_{i,0}^{(\lambda)}]\!]
			+\sqq^{-1} \sum\limits_{[M]\in\cm_{i,r\de}}\bn(\ell(M)-1) [\![M]\!].
		\end{align}
	\end{customthm}
	
	It is very helpful and interesting to describe all the root vectors at any $[i,j]$ in the $\imath$Hall algebra $\iH(\X_\bfk)$.
	In fact, all the root vectors can be described in $\imath$Hall algebra with the quiver being suitably oriented, and then in the $\imath$Hall algebra of the cyclic quiver $C_n$ by using Fourier transformations; see Proposition \ref{rem:rootvectors}.

	\subsubsection{}
	Let $\X$ be a weighted projective line of domestic type, and $Q$ be the affine quiver of the associated type.
	Then there exists a tilting object in $\coh(\X_\bfk)$ which gives a derived equivalence (see \cite{GL87})
	$$\cd^b(\coh(\X_\bfk))\stackrel{\simeq}{\longrightarrow} \cd^b(\mod(\bfk Q^{op})).$$
	%Let $\La^\imath$ be the $\imath$quiver algebra of $Q$ of split type, i.e., $$\La^\imath:=\bfk Q\otimes_\bfk \bfk[\varepsilon]/(\varepsilon^2).$$
	%There is a derived equivalence (see Remark \ref{prop:tilt})
	%$$\cd^b(\cc_1(\coh(\X))) \stackrel{\simeq}{\longrightarrow} \cd^b(\mod(\La^{\imath,op})).$$
	In the case $\X_\bfk$ being of domestic type A, we have $\bfk Q\cong \bfk C_{p_1,p_2}$; see \eqref{figure: quasi-split i}.
	As a special case of the main result in \cite{LW20a}, there is a realization of the universal $\imath$quantum group $\tUi$ (in its Serre type presentation) via the $\imath$Hall algebra $\iH(\bfk C_{p_1,p_2}^{op})$ of the quiver $C_{p_1,p_2}$, that is, we have an injective homomorphism (see Theorem \ref{thm:main})
	$$\widetilde\psi: \tUi_{ |v=\sqq} \longrightarrow \iH(\bfk C_{p_1,p_2}^{op}).$$
	The aforementioned derived equivalence induces an isomorphism of $\imath$Hall algebras
	$$\BF: \tMHX\stackrel{\cong}{\longrightarrow} \iH(\bfk C_{p_1,p_2}^{op}),$$ 
	providing a realization of the algebra isomorphism \cite{LW20b}
	$$\Phi:\tUiD \stackrel{\cong}{\longrightarrow} \tUi$$
	of the $\imath$quantum group of split affine type A in two (Serre and Drinfeld type) presentations. We summarize this result as follows.

	\begin{customthm}{\bf C} [Proposition~\ref{prop:F}, Theorem~\ref{main thm2}]
		We have the following commutative diagram % of algebra homomorphisms
		\[
		\xymatrix{ \ar[r]^{\Phi}   {}^{\text{Dr}}\tUi_{ |v=\sqq} \ar[d]^{\Omega} & \tUi_{ |v=\sqq} \ar[d]^{\widetilde{\psi}}
			\\
			\tMHX\ar[r]^-{\BF} &  \iH(\bfk C_{(p_1,p_2)}^{op}).  }\]
		%	where $\Phi, \BF$ are isomorphisms and $\widetilde{\psi},\Omega$ are monomorphisms.
		In particular, the homomorphism $\Omega$ is injective.
	\end{customthm}

	This result generalizes the corresponding result in \cite{LRW20a} for affine $\mathfrak{sl}_2$.  For $\fg$ of type DE and $Q$ the associated affine quiver, we also have a derived equivalence $\cd^b(\coh(\X_\bfk))\simeq \cd^b(\bfk Q^{op})$, which induces  an algebra isomorphism $\BF:\tMHX\stackrel{\cong}{\rightarrow} \iH(\bfk Q^{op})$. However, the morphism $\Phi$ can not be the one induced by $\BF$, in fact it can not be induced by any derived equivalent functors by the same argument as in \cite[\S9.2]{DJX12}.
	
	%\subsubsection{}
	
	%\red{As an analog of \cite{Sch02,Hu05,Hu}, we expect that the whole  $\imath$Hall algebra can be used to realize the $\imath$quantum group of split affine ${\mathfrak{gl}}_{n}$, which we hope to return elsewhere.}
	
	% Comparing with \cite{BKa01,Sch04,DJX12,BS13}, we also describe some additional root vectors clearly, for example, the counterpart of $\haB_{[i,1],r}$ (see  \eqref{realroot}) in %quantum groups is only known when $r\geq 0$ in quantum group (see \cite[\S6.5]{Sch04}), while the counterpart of $\haTh_{[1,1],r}$ is not known.

	\subsection{Comparison with previous works}
	The $\imath$Hall algebras used in this paper to give a realization of the whole $\imath$quantum loop algebras are categorical and intrinsic \cite{LP16, LW19a}, while the Drinfeld double Ringel-Hall algebras to give a realization of the whole quantum loop algebras are constructed from two copies of Ringel-Hall algebras by algebraic technique, and have much unnatural flexibility, such as Grothendieck groups, Hopf pairings, etc. For weighted projective lines, their Ringel-Hall algebras are only {\em topological bialgebras}, and we need to compute coproducts in completed spaces to obtain products of elements, even though their Drinfeld double Ringel-Hall algebras could be defined without considering completed spaces (it is also not easy to prove this fact); see \cite{BS12a,Cr10}.

	Compared with the geometric construction of quantum loop algebras in \cite{Ka97,BKa01,Sch04,DJX12,BS13}, for $\imath$quantum groups, our approach of proving the Drinfeld type relations in $\iH(\X_\bfk)$ is completely different. In fact, we can not use the same method as in \cite{DJX12} to define the corresponding root vectors in Hall algebras. 
	In addition, the relations \eqref{iDR2} (equivalently, \eqref{eq:hB1}), \eqref{iDR3b} and \eqref{iDR5} contain terms involving $\K_\de$ which do not arise in the computations {\em loc. cit.} of similar identities in the (Drinfeld double) Hall algebra of the (weighted) projective line; some new homological computations are needed to determine these $\K_\de$ terms.
	In particular, the complexity of \eqref{iDR2} makes our computations much more challenging than \cite{Sch04,DJX12}.
	
	The root vectors of $\imath$quantum loop algebras at $[i,1]$ are also described more explicitly. In fact, the counterpart of $\haB_{[i,1],l}$ (see  \eqref{realroot}) for quantum loop algebras is only known when $l\geq 0$ (see \cite[\S6.5]{Sch04}), while the counterpart of $\haTh_{[i,1],r}$ is not known. To our knowledge, it is  the first attempt to describe the root vectors of quantum and $\imath$quantum loop algebras by using the technique of Fourier transformations as in Proposition \ref{rem:rootvectors}.

	%For the injectivity of $\Omega$, the corresponding result of Theorem {\bf C} for the whole quantum loop algebras has not been proved except for the special case affine $\mathfrak{sl}_2$ \cite{BS12}; see also \cite{DJX12}. It is proved in \cite{Sch04} that a ``positive part'' of the quantum loop algebras can be realized by the Hall algebra of $\X_\bfk$ if $\fg$ is of finite or affine type. However, it can not induce the injectivity of the morphism from quantum loop algebras to Drinfeld double Hall algebras directly by applying Drinfeld double technique. 
	
	The corresponding result of Theorem {\bf C} for the whole quantum group is not known except for the special case affine $\mathfrak{sl}_2$; see \cite{BS12}. In fact, for affine $\mathfrak{sl}_2$, the Drinfeld-Beck isomorphism $\Phi$ can be interpreted by Beilinson's derived equivalence \cite{BS12}.

	%\red{The corresponding result of Theorem {\bf E} for the whole quantum group is not known \blue{except \blue{for} the special case affine $\mathfrak{sl}_2$} (only the situation of affine $\mathfrak{sl}_2$ has been considered in \cite{BS12}).} In particular, for affine $\mathfrak{sl}_2$, the Drinfeld-Beck isomorphism \blue{$\Phi$} can be interpreted by Beilinson's derived equivalence \cite{BS12}. \green{delete? As mentioned above, the injectivity of the morphism from quantum groups to Drinfeld double Hall algebras of weighted projective lines is also only proved for type affine $\mathfrak{sl}_2$; see \cite{BS12}.}

	%In a sequel of this work, we shall use the semi-derived Hall algebras \cite{LP16} of weighted projective lines to give a whole categorical realization of quantum loop algebras, which consolidates \cite[Remark 4.16]{LP16} and strengthens the results of \cite{Ka97,BKa01,Sch04,DJX12,BS12,BS13}.
	
	\subsection{Further works}
	This work opens up further research directions. In a sequel of this paper, we shall describe the composition algebra of the $\imath$Hall algebra of a weighted projective line, prove that the morphism $\Omega: \tUiD\rightarrow \iH(\X_\bfk)$ is injective if $\fg$ is of ADE type,
	and give a PBW basis for the $\imath$quantum group via coherent sheaves. For a cyclic quiver $C_n$, its $\imath$Hall algebra can be used to realize the $\imath$quantum groups of split affine ${\mathfrak{sl}}_{n}$ and ${\mathfrak{gl}}_{n}$.
	
	The $\imath$quantum loop algebra defined in Definition \ref{def:iDR} is a generalization of the Drinfeld type presentation of $\imath$quantum group of split affine ADE type. It will be interesting to study its algebraic structures. In particular, for $\fg$ of affine type $D_4, E_6, E_7, E_8$, this $\imath$quantum loop algebra is the $\imath$quantum group of elliptic (or $2$-toroidal) algebras of types $D_4^{(1,1)}$, $E_6^{(1,1)}$, $E_7^{(1,1)}$ and $E_8^{(1,1)}$ (in this case, $\X_\bfk$ is of tubular type), whose algebraic structure and geometric realizations should be studied more deeply. We expect the morphism $\Omega: \tUiD_{|v=\sqq}\rightarrow \iH(\X_\bfk)$ to be injective for arbitrary Kac-Moody algebra $\fg$.
	It will also be interesting to study the $\imath$Hall algebras of higher genus curves, in particular, of elliptic curves, which shall be studied deeply in the future.

	A Drinfeld type presentation of $\imath$quantum groups arising from quasi-split quantum symmetric pairs of affine ADE type is in preparation.
	We shall develop connections between $\imath$Hall algebras of weighted projective lines and the $\imath$quantum groups of quasi-split affine ADE type in Drinfeld type presentation, where a new kind of categories (different from the categories of $1$-periodic complexes) over weighted projective lines shall be introduced. 
	
	As quantum groups are $\imath$quantum groups arising from quantum symmetric pairs of diagonal type (which can be viewed as $\imath$quantum groups of quasi-split type), we shall reformulate the geometric realization \cite{DJX12} (see also \cite{Sch04,BS13})  of quantum groups in Drinfeld presentation via the twisted semi-derived Hall algebras of weighted projective lines, which consolidates \cite[Remark 4.16]{LP16} and strengthens the results of \cite{Ka97,BKa01,Sch04,DJX12,BS12,BS13}.

	\subsection{Organization}

	The paper is organized as follows.
	In Section~\ref{sec:main}, we review the materials on affine $\imath$quantum groups of split ADE type, $\imath$quantum loop algebras and their Drinfeld type presentations.
	The $\imath$Hall algebra on the categories of $1$-periodic complexes and the  $\imath$quiver algebras of split type  are summarized in Section~\ref{sec:Hall}.
	In Section~\ref{sec:WPL}, we recalled weighted projective lines and their categories of coherent sheaves. The $\imath$canonical algebras, Grothendieck groups and Euler forms are also reviewed.

	In Section~\ref{sec:cyclic}, we studied the $\imath$Hall algebra of a cyclic quiver. The map $\Omega:\tUiD_{|v=\sqq}\rightarrow \iH(\X_\bfk)$ in Theorem~{\bf A} is defined in Section~\ref{sec:hom}, and the proof of $\Omega$  being an algebra homomorphism occupies Sections~\ref{sec:Relationtube}--\ref{sec:Relationsstari1 II}. In particular, in Subsection~\ref{sec:rooti1}, we describe the root vectors at $[i,1]$ in $\iH(\X_\bfk)$ as stated in Theorem {\bf B}.
	%In Section \ref{sec:Hall cyclic}, we use the whole $\imath$Hall algebra of the cyclic quiver $C_n$ to realize the $\imath$quantum loop algebra of affine $\mathfrak{gl}_n$, and Theorem {\bf C} is proved.
	%Section \ref{sec:injectivity} is devoted to proving the injectivity of $\Omega$ when $\fg$ is of finite or affine type as stated in Theorem {\bf D}.
	
	A derived equivalence leading to the isomorphism of $\imath$Hall algebras is established and then Theorem~{\bf C} is proved in Section~\ref{sec:derived}.
	
	In order to make the paper readable and the idea of proofs clear, we put a dozen proofs in Appendix.

	\subsection{Acknowledgments}
	ML deeply thanks Weiqiang Wang for guiding him to study the $\imath$quantum groups, and also his continuing encouragement.
	ML thanks University of Virginia for hospitality and support. We thank Liangang Peng and Jie Xiao for helpful discussions on Ringel-Hall algebras of weighted projective lines. We thank the anonymous referee for very helpful suggestions and comments.
	
	ML is partially supported by the National Natural Science Foundation of China (No. 12171333).
	SR is partially supported by the National Natural Science Foundation of China (No. 11801473) and the Fundamental Research Funds for Central Universities of China (No. 20720210006).
	%
	%

	%%%%%%%%
	%%%%%%%%
	\section{A Drinfeld type presentation of $\imath$quantum groups }
	\label{sec:main}

	In this section, we collect  materials on (affine) Lie algebras, $\imath$quantum groups and their Drinfeld presentations.

	\subsection{The loop algebras and affine Lie algebras}
	
	Let $\II$ be an index set. Let $A=(a_{ij})_{i,j\in \II}$ be an irreducible generalized Cartan matrix (GCM) of a Kac-Moody algebra $\fg$. Let $\cR_0$ be the set of roots for $\fg$, and fix a set $\cR^+_0$ of positive roots  with simple roots $\alpha_i$ $(i\in \II)$. Let $\Lg$ be the {loop algebra} of $\fg$, and $\widehat{\fg}:=\fg[t,t^{-1}]\oplus \C c$. If $\fg$ is a complex simple Lie algebra, then $\Lg\cong \widehat{\fg}$, which is an affine Kac-Moody algebra \cite{Ga80}. In general $\Lg$ and $\widehat{\fg}$ are not Kac-Moody algebras, and there is a surjective (not injective) homomorphism $\Lg \rightarrow \widehat{\fg}$. %; moreover, if $\fg$ is a complex affine Lie algebra, then $\Lg$ is a toroidal algebra \cite{MRY90}.

	Let $\widehat{\fg}$ be the (untwisted) affine Kac-Moody algebra with affine Cartan matrix denoted by $(a_{ij})_{i,j\in\I}$, where $\I=\{0\} \cup \II$ with the affine node $0$. Let $\alpha_i$ $(i\in \I)$ be the simple roots of $\widehat{\fg}$, and $\alpha_0=\de -\theta$, where $\de$ denotes the basic imaginary root, and $\theta$ is the highest root of $\fg$. The root system $\cR$  for $\widehat{\fg}$ is defined to be
	\begin{align}
		\cR &=\{\pm (\beta + k \delta) \mid \beta \in \cR_0^+, k  \in \Z\}  \cup \{m \delta \mid m \in \Z\backslash \{0\} \}.
		% \quad \cR_{\text{im}}^+ = \{m \delta \mid m \ge 1\},
		\label{eq:roots} % \\
		%\cR_> &= \{k \delta +\beta \mid \beta \in \cR_0^+, k  \ge 0\}, \quad
		%\cR_< = \{k \delta -\beta \mid \beta \in \cR_0^+, k > 0\},
		%  \label{eq:rootsR} 
		%\\
		%\cR^+ &= \{k \delta +\beta \mid \beta \in \cR_0^+, k  \ge 0\}
		%\cup  \{k \delta -\beta \mid \beta \in \cR_0^+, k > 0\}
		%\cup \{m \delta \mid m \ge 1\}.
		%\cR_> \cup \cR_{\text{im}}^+ \cup \cR_<,  \quad
		%\cR^+[\II] = \cR_> \cup (\cR_{\text{im}}^+ \times \II) \cup \cR_<.
		% \label{eq:roots+}
	\end{align}
	% The set $\cR^+[\II]$ stands for positive roots with multiplicity.
	For $\gamma =\sum_{i\in \I} n_i \alpha_i \in \N \I$, the height $\text{ht} (\gamma)$ is defined as $\text{ht} (\gamma) =\sum_{i\in \I} n_i$.
	
	Let $P$ and $Q$ denote the weight and root lattices of $\fg$, respectively. Let $\omega_i \in P$ $(i\in \II)$ be the fundamental weights of $\fg$. Note $\alpha_i =\sum_{j\in \II} a_{ij}\omega_j$. We define a bilinear pairing $\langle \cdot, \cdot \rangle : P\times Q \rightarrow \Z$ such that $\langle \omega_i, \alpha_j \rangle =\delta_{i,j}$, for $i,j \in \II$, and thus $\langle \alpha_i, \alpha_j \rangle = a_{ij}$.
	
	%The coweight lattice $P^\vee$ of $\fg$ is a lattice over $\Z$ with a basis given by $\omega_i^\vee$, $1\leq i\leq \bt$. The coroot lattice $Q^\vee$ of $\fg$ is a sublattice of $P^\vee$ with a basis given by $\alpha_i^\vee$, $1\leq i\leq \bt$, where $\alpha_i^\vee=\sum_{1\leq j\leq n} a_{ij}\omega_j^\vee$.
	%The root lattice $Q=\Hom(P^\vee,\Z)$ with a basis given by $\alpha_i$ such that $\langle \alpha_i,\omega_j^\vee \rangle=\delta_{ij}$ for any $1\leq i,j\leq n$.
	
	The Weyl group $W_0$ of $\fg$ is generated by the simple reflection $s_i$, for $i \in \II$. It acts on $P$ by
	$s_i(x)=x-\langle x, \alpha_i \rangle\alpha_i$ for $x\in P$. The extended affine Weyl group $\widetilde{W}$ is the semi-direct product $W_0 \ltimes P$, which contains the affine Weyl group $W:=W_0 \ltimes Q =\langle s_i \mid i \in \I \rangle$ as a subgroup; we denote
	\[
	t_\omega =(1, \omega) \in \widetilde W, \quad \text{ for } \omega \in P.
	\]
	We identify $P/Q$ with a finite group $\Omega$ of Dynkin diagram automorphisms, and so $\widetilde{W} =\Omega \cdot W$. %There is a length function $\ell(\cdot)$ on $\widetilde{W}$ such that $\ell(s_i)=1$, for $i\in \I$, and each element in $\Omega$ has length 0.
	
	For $i\in \II$, as in \cite{Be94}, we define
	\begin{equation}
		\label{eq:tomega}
		%\ell(t_{\omega_i}) =\ell(\omega_i')+1, \qquad
		%\text{ where }
		\omega_i':= t_{\omega_i} s_i \in \widetilde{W}.
	\end{equation}
	
	\subsection{Affine $\imath$quantum groups of split type ADE}

	For $n\in \Z, r\in \N$, denote by
	\[
	[n] =\frac{v^n -v^{-n}}{v-v^{-1}},\qquad
	\qbinom{n}{r} =\frac{[n][n-1]\ldots [n-r+1]}{[r]!}.
	\]
	For $A, B$ in a $\Q(v)$-algebra, we shall denote $[A,B]_{v^a} =AB -v^aBA$, and $[A,B] =AB - BA$. We also use the convenient notation throughout the paper
	\begin{align}
		[A,B,C]_{v^a}=\big[[A,B]_{v^a},C\big]_{v^a}.
	\end{align}

	Recall the Cartan matrix $A=(a_{ij})_{i,j\in \I}$ of affine type ADE, for $\I = \II \cup \{0\}$ with the affine node $0$. %In particular, $a_{ij}=0, -1$, for all $i\neq j$, unless it is of the affine $A_1$ type; in that case $a_{01}=a_{10}=-2$.
	The notion of (quasi-split) universal $\imath$quantum groups $\tUi$ was formulated in \cite{LW19a}.
	
	The {\em universal affine $\imath$quantum group of split type ADE} is the $\Q(v)$-algebra $\tUi =\tUi(\widehat{\fg})$ with generators $B_i$, $\K_i^{\pm 1}$ $(i\in \I)$, subject to the following relations, for $i, j\in \I$:
	\begin{eqnarray}
		\label{eq:KK}
		\K_i\K_i^{-1} =\K_i^{-1}\K_i=1, \qquad\qquad \quad \K_i  \text{ is central},&
		\\
		B_iB_j -B_j B_i=0, \qquad\qquad\qquad\qquad\qquad &&\text{ if } a_{ij}=0,
		\label{eq:S1} \\
		B_i^2 B_j -[2] B_i B_j B_i +B_j B_i^2 = - v^{-1}  B_j \K_i,  \qquad\qquad\qquad &&\text{ if }a_{ij}=-1,
		\label{eq:S2} \\
		\sum_{r=0}^3 (-1)^r \qbinom{3}{r} B_i^{3-r} B_j B_i^{r} = -v^{-1} [2]^2 (B_iB_j-B_jB_i) \K_i,  & & \text{ if }a_{ij}=-2.
		\label{eq:S3}
	\end{eqnarray}
	%The universal split $\imath$quantum group of rank 1 is also known as the $q$-Onsager algebra, and this is the only case where the relation \eqref{eq:S3} is needed; see Definition~\ref{def:Onsager}.
	
	For $\mu = \mu' +\mu''  \in \Z \I := \oplus_{i\in \I} \Z \alpha_i$,  define $\K_\mu$ such that
	\begin{align}
		\K_{\alpha_i} =\K_i, \quad
		\K_{-\alpha_i} =\K_i^{-1}, \quad
		\K_{\mu} =\K_{\mu'} \K_{\mu''},
		\quad  \K_\delta =\K_0 \K_\theta.
	\end{align}
	%The algebra $\tUi$ is endowed with a filtered algebra structure
	%\begin{align}  \label{eq:filt1}
	%\widetilde{\U}^{\imath,0} \subset \widetilde{\U}^{\imath,1} \subset \cdots \subset \widetilde{\U}^{\imath,m} \subset \cdots
	%\end{align}
	%by setting %$|B_i| =1, |\K_i|=0$, for $i\in \I$, and then $|B_{i_1}\ldots B_{i_m} \K_\mu |=m$, for any $m$ and $\mu$.
	%\begin{align}  \label{eq:filt}
	%\widetilde{\U}^{\imath,m} =\Q(v)\text{-span} \{ B_{i_1} B_{i_2} \ldots B_{i_r} \K_\mu \mid \mu \in \N\I, i_1, \ldots, i_r \in \I, r\le m \}.
	%\end{align}
	%Note that
	%\begin{align}  \label{eq:UiCartan}
	%\widetilde{\U}^{\imath,0} =\bigoplus_{\mu \in \N\I} \Q(v) \K_\mu,
	%\end{align}
	%is the $\Q(v)$-subalgebra generated by $\K_i$ for $i\in \I$.
	%Then, according to a basic result of Letzter and Kolb on quantum symmetric pairs adapted to our setting of $\tUi$ (cf. \cite{Let02, Ko14}), the associated graded $\gr \tUi$ with respect to \eqref{eq:filt1}--\eqref{eq:filt} can be identified with
	%\begin{align}   \label{eq:filter}
	%\gr \tUi \cong \U^- \otimes \Q(v)[\K_i^\pm | i\in \I],
	%\qquad \overline{B_i}\mapsto F_i,  \quad
	%\overline{\K}_i \mapsto \K_i \; (i\in \I).
	%\end{align}
	
	\begin{remark}
		\label{rem:Kk}
		The generator $\K_i$ here, which corresponds to the acyclic complex $K_{S_i}$  in the $\imath$Hall algebra, is related to the generator $\tilde{k}_i$ used in \cite{LW19a, LW20a} (which is natural from the viewpoint of Drinfeld doubles) by a rescaling:
		$\K_i = -v^2\tilde{k}_i$.
		%A Serre presentation of a quasi-split (including split) universal $\imath$quantum group $\tUi$ of Kac-Moody type can be given uniformly via $\imath$divided powers, see \cite[Theorem 4.7]{LW20a}.
		The precise relation between the algebra $\tUi$ and the $\imath$quantum group $\Ui$ arising from quantum symmetric pairs \cite{Ko14} is explained {\em loc. cit.} %; also see \S\ref{subsec:parameter} below.
	\end{remark}
	
	\begin{remark}
		The $\Q(v)$-algebra $\tUi$ is $\Z \I$-graded by letting
		\begin{align}
			\label{eq:deg}
			\wt (B_i) =\alpha_i, \quad \wt (\K_i) =2\alpha_i, \quad \text{ for } i \in \I.
		\end{align}
		A variant of $\tUi$, in which $\K_i$ is not assumed to be invertible, is $\N\I$-graded by \eqref{eq:deg}. \end{remark}
	
	\begin{lemma} [\cite{LW21b}; \text{also cf. \cite{KP11, BK20}}]
		\label{lem:Ti}
		For $i\in \I$, there exists an automorphism $\TT_i$ of the $\Q(v)$-algebra $\tUi$ such that
		$\TT_i(\K_\mu) =\K_{s_i\mu}$ for $\mu\in \Z\I$, and
		\[
		\TT_i(B_j)= \begin{cases}
			\K_i^{-1} B_i,  &\text{ if }j=i,\\
			B_j,  &\text{ if } a_{ij}=0, \\
			B_jB_i-vB_iB_j,  & \text{ if }a_{ij}=-1, \\
			{[}2]^{-1} \big(B_jB_i^{2} -v[2] B_i B_jB_i +v^2 B_i^{2} B_j \big) + B_j\K_i,  & \text{ if }a_{ij}=-2,
		\end{cases}
		\]
		for $j\in \I$.
		Moreover,  $\TT_i$ $(i\in \I)$ satisfy the braid group relations, i.e., $\TT_i \TT_j =\TT_j \TT_i$ if $a_{ij}=0$, and $\TT_i \TT_j \TT_i =\TT_j \TT_i \TT_j$ if $a_{ij}=-1$.
	\end{lemma}

	%\begin{lemma}
	%\label{lem:anti}
	%There exists a $\Q(v)$-algebra anti-involution $\sigma:\tUi\rightarrow \tUi$ such that
	%$$\sigma(B_i)=B_i, \quad \sigma(\K_i)= \K_{ i},
	%\quad \forall i\in \I.$$
	%\end{lemma}
	
	%\begin{proof}
	%It is easy to check that $\sigma$ preserves \eqref{eq:KK}--\eqref{eq:S3}.
	%\end{proof}
	
	\subsection{A Drinfeld type presentation}
	\label{subsec:Dr2}
	
	Let $A=(a_{ij})_{i,j\in \II}$ be a generalized Cartan matrix (GCM) of a {\em simply-laced} Kac-Moody algebra $\fg$. The $\imath$quantum loop algebra $\tUiD=\tUiD(\Lg)$ of split type is the $\Q(v)$-algebra  generated by $\K_{i}^{\pm1}$, $C^{\pm1}$, $H_{i,m}$, $\Theta_{i,l}$ and $\y_{i,l}$, where  $i\in \II$, $m \in \Z_{+}$, $l\in\Z$, subject to some relations.
	In order to give the explicit definition of $\tUiD$,  we introduce some shorthand notations below.
	
	% such that $a_{ij}\in\{-1,0,2\}$.
	Let $k_1, k_2, l\in \Z$ and $i,j \in \II$. Set
	\begin{align}
		\begin{split}
			S(k_1,k_2\mid l;i,j)
			&=  B_{i,k_1} B_{i,k_2} B_{j,l} -[2] B_{i,k_1} B_{j,l} B_{i,k_2} + B_{j,l} B_{i,k_1} B_{i,k_2},
			%&=\Sym_{k_1,k_2}\sum_{t=0}^2(-1)^t \qbinom{2}{t}B_{i,k_1}\cdots B_{i,k_t}B_{j,l} B_{i,k_{t+1}} \cdots B_{i,k_2}\notag \\
			%&=\sum_{t=0}^2(-1)^t \qbinom{2}{t}B_{i,k_1}\cdots B_{i,k_t}B_{j,l} B_{i,k_{t+1}} \cdots B_{i,k_2},
			\\
			\SS(k_1,k_2\mid l;i,j)
			&= S(k_1,k_2\mid l;i,j)  + \{k_1 \leftrightarrow k_2 \}.
			\label{eq:Skk}
		\end{split}
	\end{align}
	Here and below, $\{k_1 \leftrightarrow k_2 \}$ stands for repeating the previous summand with $k_1, k_2$ switched if $k_1\neq k_2$, so the sums over $k_1, k_2$ are symmetric.
	We also denote
	\begin{align}
		\begin{split}
			R(k_1,k_2\mid l; i,j)
			&=   \K_i  C^{k_1}
			\Big(-\sum_{p\geq0} v^{2p}  [2] [\Theta _{i,k_2-k_1-2p-1},\y_{j,l-1}]_{v^{-2}}C^{p+1}
			\label{eq:Rkk} \\
			&\quad -\sum_{p\geq 1} v^{2p-1}  [2] [\y_{j,l},\Theta _{i,k_2-k_1-2p}]_{v^{-2}} C^{p}
			- [\y_{j,l}, \Theta _{i,k_2-k_1}]_{v^{-2}} \Big),
			\\
			\R(k_1,k_2\mid l; i,j) &= R(k_1,k_2\mid l;i,j) + \{k_1 \leftrightarrow k_2\}.
		\end{split}
	\end{align}
	Sometimes, it is convenient to rewrite part of the summands in \eqref{eq:Rkk} as
	\begin{align*}
		&-\sum_{p\geq 1} v^{2p-1}  [2] [\y_{j,l},\Theta _{i,k_2-k_1-2p}]_{v^{-2}} C^{p}
		- [\y_{j,l}, \Theta _{i,k_2-k_1}]_{v^{-2}}\\
		&=
		-\sum_{p\geq 0} v^{2p-1}  [2] [\y_{j,l},\Theta _{i,k_2-k_1-2p}]_{v^{-2}} C^{p}
		+v^{-2}[\y_{j,l}, \Theta _{i,k_2-k_1}]_{v^{-2}}.
	\end{align*}
	
	%We shall often omit $i,j$ and write $S(k_1,k_2|l)= S(k_1,k_2|l;i,j)$, $\SS(k_1,k_2|l)= \SS(k_1,k_2|l;i,j)$, and similarly for $R$ and $\R$, whenever $i,j$ are clear from the context.
	
	\begin{definition}[$\imath$quantum loop algebras, \cite{LW20b}]
		\label{def:iDR}
		Let $A=(a_{ij})_{i,j\in \II}$ be a generalized Cartan matrix (GCM) of a simply-laced Kac-Moody algebra $\fg$. % such that $a_{ij}\in\{-1,0,2\}$.
		The $\imath$quantum loop algebra $\tUiD$ of split type is the $\Q(v)$-algebra  generated by $\K_{i}^{\pm1}$, $C^{\pm1}$, $H_{i,m}$ and $\y_{i,l}$, where  $i\in \II$, $m >0$, $l\in\Z$, subject to the following relations (for $m,n >0$ and $k,l\in \Z$):
		\begin{align}
			%1
			& \K_i, C \text{ are central, } \quad \K_i\K_i^{-1}=1, \;\; C C^{-1}=1,
			\label{iDR1a}
			\\%2
			%\label{iDR1b}
			&[H_{i,m},H_{j,n}]=0,\label{iDR1b}\\
			&[H_{i,m},\y_{j,l}]=\frac{[ma_{ij}]}{m} \y_{j,l+m}-\frac{[ma_{ij}]}{m} \y_{j,l-m}C^m,
			\label{iDR2}
			\\%3
			&[\y_{i,k} ,\y_{j,l}]=0,   \text{ if }a_{ij}=0,  \label{iDR4}
			\\%4
			&[\y_{i,k}, \y_{j,l+1}]_{v^{-a_{ij}}}  -v^{-a_{ij}} [\y_{i,k+1}, \y_{j,l}]_{v^{a_{ij}}}=0, \text{ if }i\neq j,
			\label{iDR3a}
			\\ %5
			&[\y_{i,k}, \y_{i,l+1}]_{v^{-2}}  -v^{-2} [\y_{i,k+1}, \y_{i,l}]_{v^{2}}
			\label{iDR3b} \\
			&=v^{-2}\Theta_{i,l-k+1} C^k \K_i-v^{-4}\Theta_{i,l-k-1} C^{k+1} \K_i
			+v^{-2}\Theta_{i,k-l+1} C^l \K_i-v^{-4}\Theta_{i,k-l-1} C^{l+1} \K_i, \notag
			\\
			\label{iDR5}
			&    \SS(k_1,k_2\mid l; i,j) =  \R(k_1,k_2\mid l; i,j), \text{ if }a_{ij}=-1.
		\end{align}
		Here we set
		\begin{align}  \label{Hm0}
			{\Theta}_{i,0} =(v-v^{-1})^{-1}, \qquad {\Theta}_{i,m} =0, \; \text{ for }m<0;
		\end{align}
		and $\Theta_{i,m}$ ($m\geq1$)  are related to  $H_{i,m}$ by the following equation:
		\begin{align}
			\label{exp h}
			1+ \sum_{m\geq 1} (v-v^{-1})\Theta_{i,m} u^m  = \exp\Big( (v-v^{-1}) \sum_{m\geq 1} H_{i,m} u^m \Big).
		\end{align}
	\end{definition}
	Let us mention that in spite of its appearance, the RHS of \eqref{iDR3b} typically has two non-zero terms, thanks to the convention \eqref{Hm0}.

	\begin{lemma}
		[\cite{LW20b}]
		\label{lem:equiv}
		\quad
		\begin{enumerate}
			\item
			The relation~\eqref{iDR1b} is equivalent to
			\begin{align}
				\label{eq:hh1}
				[\Theta_{i,m},\Theta_{j,n}] = 0 \quad (m, n\geq 1).
			\end{align}
			\item
			The relation \eqref{iDR2} is equivalent to
			\begin{align}
				\label{eq:hB1}
				&[\Theta_{i,m}, \y_{j,l}]+[\Theta_{i,m-2},\y_{j,l}]C\\
				&=v^{a_{ij}}[\Theta_{i,m-1},\y_{j,l+1}]_{v^{-2a_{ij}}}+v^{-a_{ij}}[\Theta_{i,m-1},\y_{j,l-1}]_{v^{2a_{ij}}}C \quad (m\geq 1, l\in\Z).
				\notag
			\end{align}
		\end{enumerate}
	\end{lemma}
	
	Let $A$ be a GCM of an affine Lie algebra $\fg$ of type ADE.
	Define a sign function
	\begin{align}
		\label{eq:sign}
		o(\cdot): \,\II \longrightarrow \{\pm 1\}
	\end{align}
	such that $o(i) o(j)=-1$ whenever $a_{ij} <0$ (there are clearly and exactly 2 such functions).
	We define the following elements in $\tUi$, for $i\in \II$, $k\in \Z$ and $m\ge 1$: %(compare with the rank 1 formulas \eqref{eq:B1n}--\eqref{eq:dB1} and \eqref{eq:dBB}, where $\TT_{\omega_1} =\dag \TT_1$):
	\begin{align}
		&B_{i,k} = o(i)^k \TT_{\omega_i}^{-k} (B_i),
		\label{Bik} \\\notag
		&\acute{\Theta}_{i,m} =  o(i)^m \Big(-B_{i,m-1} \TT_{\omega_i'} (B_i) +v^{2} \TT_{\omega_i'} (B_i) B_{i,m-1}+ (v^{2}-1)\sum_{p=0}^{m-2} B_{i,p} B_{i,m-p-2}  \K_{i}^{-1}\K_{\de} \Big),\\
		&\dB_{i,m} =\acute{\Theta}_{i,m} - \sum\limits_{a=1}^{\lfloor\frac{m-1}{2}\rfloor}(v^2-1) v^{-2a} \acute{\Theta}_{i,m-2a}\K_{a\de} -\de_{m,ev}v^{1-m} \K_{\frac{m}{2}\de}.
		\label{Thim}
	\end{align}
	and $H_{i,m}$ is defined via \eqref{exp h}.
	Here, $\de_{m,ev}=1$ if $m$ is even and $\de_{m,ev}=0$ otherwise. In particular, $B_{i,0}=B_i$. %imilarly, we define
	
	\begin{theorem}[\cite{LW20b}]
		\label{thm:ADE}
		If $A=(a_{ij})_{i,j\in\II}$ is a GCM of a Lie algebra $\fg$ of type ADE, then there is a $\Q(v)$-algebra isomorphism ${\Phi}: \tUiD(\Lg) \rightarrow\tUi(\widehat{\fg})$, which sends
		\begin{align}   \label{eq:map}
			\y_{i,k}\mapsto B_{i,k},  \quad H_{i,m}\mapsto H_{i,m}, \quad \Theta_{i,m}\mapsto \Theta_{i,m}, \quad
			\K_i\mapsto \K_i, %-v^{2}\tk_i,
			\quad C\mapsto \K_\de,
		\end{align}
		for $i\in \II, k\in \Z$, and $m \ge 1$.
	\end{theorem}
	
	\subsection{Star-shaped graph}
	
	For $\bp=(p_1,\dots,p_\bt)$ with $p_i\geq2$, $\bt\geq1$, let us consider the following star-shaped graph $\Gamma=\T_{p_1,\dots,p_\bt}$:
	
	\begin{center}\setlength{\unitlength}{0.8mm}
		\begin{equation}
			\label{star-shaped}
			\begin{picture}(110,30)(0,35)
				\put(0,40){\circle*{1.4}}
				\put(2,42){\line(1,1){16}}
				\put(20,60){\circle*{1.4}}
				\put(23,60){\line(1,0){13}}
				\put(40,60){\circle*{1.4}}
				\put(43,60){\line(1,0){13}}
				\put(60,58.5){\large$\cdots$}
				\put(70,60){\line(1,0){13}}
				\put(88,60){\circle*{1.4}}
				%\put(90,58){\line(1,-1){16}}
				%\put(108,40){\circle*{1.4}}

				\put(3,41){\line(4,1){13}}
				\put(20,45){\circle*{1.4}}
				\put(23,45){\line(1,0){13}}
				\put(40,45){\circle*{1.4}}
				\put(43,45){\line(1,0){13}}
				\put(60,43.5){\large$\cdots$}
				\put(70,45){\line(1,0){13}}
				\put(88,45){\circle*{1.4}}
				%\put(91,44){\line(4,-1){13}}
				
				\put(19,30){\Large$\vdots$}
				
				\put(39,30){\Large$\vdots$}
				
				\put(87,30){\Large$\vdots$}
				
				\put(2,38){\line(1,-1){16}}
				\put(20,20){\circle*{1.4}}
				\put(23,20){\line(1,0){13}}
				\put(40,20){\circle*{1.4}}
				\put(43,20){\line(1,0){13}}
				\put(60,18.5){\large$\cdots$}
				\put(70,20){\line(1,0){13}}
				\put(88,20){\circle*{1.4}}
				%\put(90,22){\line(1,1){16}}
				
				\put(-4.5,39){$\star$}
				
				\put(16,62){\tiny$[1,1]$}
				
				\put(36,62){\tiny$[1,2]$}
				\put(81,62){\tiny$[1,p_1-1]$}

				\put(16.5,47){\tiny$[2,1]$}
				\put(36.5,47){\tiny$[2,2]$}
				\put(81,47){\tiny$[2,p_2-1]$}

				\put(16.5,16){\tiny$[\bt,1]$}
				\put(36.5,16){\tiny$[\bt,2]$}
				\put(81,16){\tiny$[\bt,p_\bt-1]$}

				%\put(109.5,39.5){\tiny$c$}
				
			\end{picture}
		\end{equation}
		\vspace{1cm}
	\end{center}
	The set of vertices is denoted by $\II$. As marked in the graph, the central vertex is denoted by $\star$. Let $J_1,\dots,J_\bt$ be the branches, which are subdiagrams of type $A_{p_1-1},\dots,A_{p_\bt-1}$ respectively. Denote by $[i,j]$ the $j$-th vertex in the $i$-th branch. These examples include all finite-type Dynkin diagrams as well as the affine Dynkin diagrams of types $D_4^{(1)},E_6^{(1)},E_7^{(1)}$, and $E_8^{(1)}$.
	
	Let $\fg$ be the Kac-Moody algebra corresponding to $\Gamma$, and $\cl \fg$ be its loop algebra. Then we have the $\imath$quantum loop algebra $\tUiD$.
	In view of the graph $\Gamma$, for the root system $\cR$ of $\cl\fg$, its simple roots are denoted by $\alpha_\star$ and $\alpha_{ij}$, for $1\leq i\leq \bt$ and $1\leq j\leq p_i-1$.

	%%%%%%%%
	%%%%%%%%
	
	\begin{lemma}
		\label{lem:reduced generators}
		For any $a\in\Z$,
		$\tUiD$ is  generated by $\y_{\star,a}$, $\y_{\star,a+1}$, $\y_{i,0}$, $\K_{i}^{\pm1}$ and $C^{\pm1}$ for all  $i\in \II-\{\star\}$.
	\end{lemma}
	
	\begin{proof}
		Let $U$ be the subalgebra of $\tUiD$ generated by $\y_{\star,a}$, $\y_{\star,a+1}$, $\y_{i,0}$, $\K_{i}^{\pm1}$ and $C^{\pm1}$ for $i\in\II-\{\star\}$.
		First, we have $\Theta_{\star,1}\in U$ by taking $i=\star$, $k=a=l$ in \eqref{iDR3b} since $\Theta_{\star,m}=0$ for $m<0$. From \eqref{iDR2} by taking $i=j=\star$,  since $H_{\star,1}=\Theta_{\star,1}$, we have
		$\y_{\star,a-1}\in U$ by taking $m=1$ and $l=a$; and $\y_{\star,a+2}\in U$ by taking $m=1$ and $l=a+1$. Inductively, we have
		$\y_{\star,l}\in U$ for any $l\in\Z$.

		From \eqref{iDR3b}, we have $\Theta_{\star,2}\in U$ by taking $k=0$, $l=1$ since $\Theta_{\star,0}=\frac{1}{v-v^{-1}}$. Inductively, one can prove that $\Theta_{\star,m}\in U$ for any $m\geq1$. It follows from \eqref{exp h} that $H_{\star,m}\in U$ for any $m\geq1$.

		In order to prove that $U=\tUiD$, it is enough to prove that $\Theta_{j,m},\y_{j,l}\in U$ for any $m\geq1$, $l\in\Z$, provided  $\Theta_{i,m},\y_{i,l}\in U$ and $a_{ij}=-1$.

		Assume $\Theta_{i,m},\y_{i,l}\in U$ for any $m\geq1$, $l\in\Z$ and $a_{ij}=-1$.
		Considering $\mathbb{S}(k_1,k_1+1|0;i,j)$, it follows from \eqref{iDR5} that $B_{j,-1}\in U$. Together with \eqref{iDR2}, we have $B_{j,1}\in U$; using \eqref{iDR3b} by taking $k=0=l$, we have $\Theta_{j,1}\in U$.
		We prove $B_{j,\pm m}$ and $\Theta_{j,m}$ are in $U$ by induction on $m\geq0$. By considering $\mathbb{S}(k_1,k_1+1|m;i,j)$, it follows from \eqref{iDR5} that $B_{j,-m-1}\in U$; by \eqref{iDR2}, we have $B_{j,m+1}\in U$;
		by \eqref{iDR3b}, we have $\Theta_{j,m+1}\in U$.
	\end{proof}
	
	%%%%%%%%%%%
	\section{$\imath$Hall algebras}
	\label{sec:Hall}
	
	For an additive category $\ce$ and $M\in \ce$, we denote
	
	$\triangleright$ $\add M$ -- subcategory of $\ce$ whose objects are the direct summands of finite direct sums of copies of $M$,
	
	%$\triangleright$ $\Ind (\ca)$ --  set of the isoclasses of indecomposable objects in $\ca$,
	
	$\triangleright$ $\Iso(\ce)$ -- set of the isoclasses of objects in $\ce$,
	
	$\triangleright$ $K_0(\ce)$ --  Grothendieck group of $\ce$ if $\ce$ is an exact category,
	
	$\triangleright$ $\widehat{M}$ -- the class of $M \in \ce$ in $K_0(\ce)$. 
	
	In this paper, we take the field $\bfk=\mathbb F_q$, a finite field of $q$ elements. Let $\ca$ be a hereditary abelian category over $\bfk$. We define the $\imath$Hall algebra $\iH(\ca)$ as a twisted semi-derived Ringel-Hall algebra for the category $\cc_1(\ca)$ of $1$-periodic complexes over $\ca$.

	\subsection{Hall algebras}
	
	Let $\ce$ be an essentially small exact category in the sense of Quillen, linear over  $\bfk$.
	%For the basics on exact categories, we refer to \cite{Buh} and references therein.
	Assume that $\ce$ has finite morphism and extension spaces, i.e.,
	\[
	|\Hom(M,N)|<\infty,\quad |\Ext^1(M,N)|<\infty,\,\,\forall M,N\in\ce.
	\]
	
	Given objects $M,N,L\in\ce$, define $\Ext^1(M,N)_L\subseteq \Ext^1(M,N)$ as the subset parameterizing extensions whose middle term is isomorphic to $L$. We define the {\em Ringel-Hall algebra} $\ch(\ce)$ (or {\em Hall algebra} for short) to be the $\Q$-vector space whose basis is formed by the isoclasses $[M]$ of objects $M$ in $\ce$, with the multiplication defined by (see \cite{Br13})
	\begin{align}
		\label{eq:mult}
		[M]\diamond [N]=\sum_{[L]\in \Iso(\ce)}\frac{|\Ext^1(M,N)_L|}{|\Hom(M,N)|}[L].
	\end{align}
	
	For any three objects $L,M,N$, let
	\begin{align}
		\label{eq:Fxyz}
		F_{MN}^L:= \big |\{X\subseteq L \mid X \cong N,  L/X\cong M\} \big |.
	\end{align}
	The Riedtmann-Peng formula states that
	\[
	F_{MN}^L= \frac{|\Ext^1(M,N)_L|}{|\Hom(M,N)|} \cdot \frac{|\Aut(L)|}{|\Aut(M)| |\Aut(N)|}.
	\]
	For any object $M$, %let $\ell(M)$ denote the number of indecomposable direct summand of $M$, and define
	let
	\begin{align}
		\label{eq:doublebrackets}
		[\![M]\!]:=\frac{[M]}{|\Aut(M)|}.
	\end{align}
	Then the Hall multiplication \eqref{eq:mult} can be reformulated to be
	\begin{align}
		[\![M]\!]\diamond [\![N]\!]=\sum_{[\![L]\!]}F_{M,N}^L[\![L]\!],
	\end{align}
	which is the original version of Hall multiplication used in \cite{Rin90}.

	\subsection{The category of $1$-periodic complexes}
	%%%
	
	\label{subsec:periodic}
	
	Let $\ca$ be a hereditary abelian category  which is essentially small with finite-dimensional homomorphism and extension spaces. %In this paper, we shall take  $\ca$ to be either $\coh(\PL)$ or $\rep^{\rm nil}_\bfk(Q)$, where $Q$ is the Kronecker quiver or the Jordan quiver.
	
	A $1$-periodic complex $X^\bullet$ in $\ca$ is a pair $(X,d)$ with $X\in\ca$ and a differential $d:X\rightarrow X$. A morphism $(X,d) \rightarrow (Y,e)$ is given by a morphism $f:X\rightarrow Y$ in $\ca$ satisfying $f\circ d=e\circ f$. Let $\cc_1(\ca)$ be the category of all $1$-periodic complexes in $\ca$. Then $\cc_1(\ca)$ is an abelian category.
	A $1$-periodic complex $X^\bullet=(X,d)$ is called acyclic if $\ker (d)=\Im (d)$. We denote by $\cc_{1,ac}(\ca)$ the full subcategory of $\cc_1(\ca)$ consisting of acyclic complexes.
	Denote by $H(X^\bullet)\in\ca$ the cohomology group of $X^\bullet$, i.e., $H(X^\bullet)=\ker (d)/\Im (d)$, where $d$ is the differential of $X^\bullet$.
	
	The category $\cc_1(\ca)$ is Frobenius with respect to the degree-wise split exact structure.
	The $1$-periodic homotopy category $\ck_1(\ca)$ is obtained as the stabilization
	of $\cc_1(\ca)$, and the $1$-periodic derived category $\cd_1(\ca)$ is the localization of the homotopy category $\ck_1(\ca)$ with respect to quasi-isomorphisms. Both $\ck_1(\ca)$ and $\cd_1(\ca)$ are triangulated categories.% ; see e.g. \cite{St17}. %; see \cite{PX97}.
	%Let \red{$\ck_1(\ca)$} be the homotopy category of $\cc_1(\ca)$, and let $\cd_{1}(\ca)$ be the \red{$\Z/1$-graded?} derived category of $\ca$, i.e., the localization of the homotopy category of \red{$\ck_1(\ca)$} with respect to quasi-isomorphisms.
	
	Let $\cc^b(\ca)$ be the category of bounded complexes over $\ca$ and $\cd^b(\ca)$ be the corresponding derived category with the shift functor $\Sigma$.  Then there is a covering functor $\pi:\cc^b(\ca)\rightarrow \cc_1(\ca)$, induces a covering functor $\pi:\cd^b(\ca)\rightarrow \cd_1(\ca)$ which is dense (see, e.g.,
	\cite[Lemma 5.1]{St17}). The orbit category $\cd^b(\ca)/\Sigma$ is a triangulated category \cite{Ke2}, and we have
	\begin{align}
		\label{dereq}
		\cd_1(\ca)\simeq \cd^b(\ca)/\Sigma.
	\end{align}
	
	For any $X\in\ca$, denote the stalk complex by
	\[
	C_X =(X,0)
	\]
	(or just by $X$ when there is no confusion), and denote by $K_X$ the following acyclic complex:
	\[
	K_X:=(X\oplus X, d),
	\qquad \text{ where }
	d=\left(\begin{array}{cc} 0&\Id \\ 0&0\end{array}\right).
	\]
	
	\begin{lemma}
		[\text{\cite[Lemma 2.2]{LRW20a}}]
		\label{lem:pd acyclic}
		For any acyclic complex $K^\bullet$ and $p \ge 2$, we have
		\begin{align}
			\label{Ext2vanish}
			\Ext^p_{\cc_1(\ca)}(K^\bullet,-)=0=\Ext^p_{\cc_1(\ca)}(-,K^\bullet).
		\end{align}
	\end{lemma}

	For any $K^\bullet\in\cc_{1,ac}(\ca)$ and $M^\bullet\in\cc_{1}(\ca)$, by \cite[Corollary 2.4]{LRW20a}, define
	\begin{align*}
		\langle K^\bullet, M^\bullet\rangle=\dim_{\bfk} \Hom_{\cc_{1}(\ca)}( K^\bullet, M^\bullet)-\dim_{\bfk}\Ext^1_{\cc_{1}(\ca)}(K^\bullet, M^\bullet),
		\\
		\langle  M^\bullet,K^\bullet\rangle=\dim_{\bfk} \Hom_{\cc_{1}(\ca)}( M^\bullet, K^\bullet)-\dim_{\bfk}\Ext^1_{\cc_{1}(\ca)}( M^\bullet,K^\bullet).
	\end{align*}
	These formulas give rise to well-defined bilinear forms (called {\em Euler forms}), again denoted by $\langle \cdot, \cdot \rangle$, on the Grothendieck groups $K_0(\cc_{1,ac}(\ca))$
	and $K_0(\cc_{1}(\ca))$.
	
	Denote by $\langle \cdot,\cdot\rangle_\ca$ (also denoted by $\langle \cdot,\cdot\rangle$ if no confusion arises) the Euler form of $\ca$,  i.e.,
	\[
	\langle M,N\rangle_\ca=\dim_\bfk\Hom_\ca(M,N)-\dim_\bfk\Ext_\ca^1(M,N).
	\]
	Let $\res: \cc_1(\ca)\rightarrow\ca$ be the restriction functor by forgetting differentials.
	
	\begin{lemma}
		[\text{\cite[Lemma 2.7]{LRW20a}}]
		\label{lemma compatible of Euler form}
		We have
		\begin{itemize}
			\item[(1)]
			$\langle K_X, M^\bullet\rangle = \langle X,\res (M^\bullet) \rangle_\ca$,\; $\langle M^\bullet,K_X\rangle =\langle \res(M^\bullet), X \rangle_\ca$, for $X\in \ca$, $M^\bullet\in\cc_1(\ca)$;
			\item[(2)] $\langle M^\bullet,N^\bullet\rangle=\frac{1}{2}\langle \res(M^\bullet),\res(N^\bullet)\rangle_\ca$, for  $M^\bullet,N^\bullet\in\cc_{1,ac}(\ca)$.
		\end{itemize}
	\end{lemma}
	
	\begin{lemma}[\text{\cite[Lemma 2.3]{LRW20a}}]
		\label{lem: iso in singularity}
		For any $X^\bullet,Y^\bullet\in\cc_1(\ca)$, we have $H(X^\bullet)\cong H(Y^\bullet)$ if and only if
		there exist two short exact sequences
		\begin{align*}
			0\longrightarrow U_1^\bullet\longrightarrow Z^\bullet \longrightarrow X^\bullet\longrightarrow0,
			\qquad 0\longrightarrow U_2^\bullet\longrightarrow Z^\bullet\longrightarrow Y^\bullet\longrightarrow0
		\end{align*}
		with $U^\bullet_1,U^\bullet_2\in \cc_{1,ac}(\ca)$. %, $Z\in\mod(\Lambda^\imath)$.
	\end{lemma}

	%%%%%
	\subsection{Semi-derived Hall algebras and $\imath$Hall algebras}
	%%%%%%%
	
	Define
	\[
	\sqq :=\sqrt{q}.
	\]
	
	We continue to work with a hereditary abelian category $\ca$ as in \S\ref{subsec:periodic}.
	Let $\ch(\cc_1(\ca))$ be the Ringel-Hall algebra of $\cc_1(\ca)$ over $\Q(\sqq)$, i.e., $\ch(\cc_1(\ca))=\bigoplus_{[M^\bullet]\in \Iso(\cc_1(\ca))} \Q(\sqq)[M^\bullet]$, with multiplication defined by
	\begin{align*}
		[M^\bullet]\diamond[N^\bullet]=\sum_{[L^\bullet]\in\Iso(\cc_1(\ca)) } \frac{|\Ext^1(M^\bullet,N^\bullet)_{L^\bullet}|}{|\Hom(M^\bullet,N^\bullet)|}[L^\bullet].
	\end{align*}

	%Let $I$ be the ideal of $\ch(\cc_1(\ca))$ generated by
	%all differences $[L^\bullet]-[K^\bullet\oplus M^\bullet]$ if there is a short exact sequence $K^\bullet\rightarrowtail L^\bullet\twoheadrightarrow M^\bullet$ in $\cc_{1}(\ca)$ with $K$ acyclic.
	For any $M\in\ca$, we denote by $\widehat{M}$ its class in the Grothendieck group $K_0(\ca)$ of $\ca$. Following \cite{LP16,LW19a,LW20a}, we consider the ideal $\cI$ of $\ch(\cc_1(\ca))$ generated by
	\begin{align}
		\label{eq:ideal}
		&\Big\{ [M^\bullet]-[N^\bullet]\mid H(M^\bullet)\cong H(N^\bullet), \quad {\widehat{\Im (d_{M^\bullet})}=\widehat{\Im (d_{N^\bullet})} \Big\}.}
	\end{align}
	We denote
	\[
	\cs:=\{ a[K^\bullet] \in \ch(\cc_1(\ca))/\cI \mid a\in \Q(\sqq)^\times, K^\bullet\in \cc_1(\ca) \text{ acyclic}\},
	\]
	a multiplicatively closed subset of $\ch(\cc_1(\ca))/ \cI$ with the identity $[0]$.
	
	\begin{lemma}
		[\text{\cite[Proposition A.5]{LW19a}}]
		The multiplicatively closed subset $\cs$ is a right Ore, right reversible subset of $\ch(\cc_1(\ca))/\cI$. Equivalently, there exists the right localization of
		$\ch(\cc_1(\ca))/\cI$ with respect to $\cs$, denoted by $(\ch(\cc_1(\ca))/\cI)[\cs^{-1}]$.
	\end{lemma}
	
	\begin{proof}
		It follows by noting that $[K^\bullet]\diamond [M^\bullet]=q^{\langle M^\bullet,K^\bullet\rangle- \langle K^\bullet,M^\bullet\rangle}[M^\bullet]\diamond [K^\bullet]$ in $\ch(\cc_1(\ca))/\cI$ for any acyclic $K^\bullet$.
	\end{proof}

	The algebra $(\ch(\cc_1(\ca))/\cI)[\cs^{-1}]$ is the {\em semi-derived Ringel-Hall algebra} of $\cc_1(\ca)$ in the sense of \cite{LP16, LW19a} (also cf. \cite{Gor18}), and will be denoted by $\cs\cd\ch(\cc_1(\ca))$.

	\begin{definition}  \label{def:iH}
		The $\imath$Hall algebra of a hereditary abelian category $\ca$, denoted by $\iH(\ca)$,  is defined to be the twisted semi-derived Ringel-Hall algebra of $\cc_1(\ca)$, that is, the $\Q(\sqq)$-algebra on the same vector space as $\cs\cd\ch(\cc_1(\ca)) =(\ch(\cc_1(\ca))/\cI)[\cs^{-1}]$ equipped with the following modified multiplication (twisted via the restriction functor $\res: \cc_1(\ca)\rightarrow\ca$)
		\begin{align}
			\label{eq:tH}
			[M^\bullet]* [N^\bullet] =\sqq^{\langle \res(M^\bullet),\res(N^\bullet)\rangle_\ca} [M^\bullet]\diamond[N^\bullet].
		\end{align}
	\end{definition}
	For any complex $M^\bullet$ and acyclic complex $K^\bullet$, we have
	\[
	[K^\bullet]*[M^\bullet]=[K^\bullet\oplus M^\bullet]=[M^\bullet]*[ K^\bullet].
	\]
	
	For any $\alpha\in K_0(\ca)$,  there exist $X,Y\in\ca$ such that $\alpha=\widehat{X}-\widehat{Y}$. Define $[K_\alpha]:=[K_X]* [K_Y]^{-1}$. This is well defined, see e.g.,
	\cite[\S 3.2]{LP16}.
	It follows that $[K_\alpha]\; (\alpha\in K_0(\ca))$ are central in the algebra $\iH(\ca)$.
	
	%It follows that $[K_\alpha]$ are central in the algebra $\cs\cd\widetilde{\ch}(\ca)$ for any $\alpha\in K_0(\ca)$.
	
	The {\em quantum torus} $\widetilde{\ct}(\ca)$ is defined to be the subalgebra of $\iH(\ca)$ generated by $[K_\alpha]$, for $\alpha\in K_0(\ca)$. %Then $\widetilde{\ct}(\ca)$ is the group algebra of $K_0(\ca)$.
	
	\begin{proposition}
		[\text{\cite[Proposition 2.9]{LRW20a}}]
		\label{prop:hallbasis}
		The following hold in $\iH(\ca)$:
		\begin{enumerate}
			\item
			The quantum torus $\widetilde{\ct}(\ca)$ is a central subalgebra of $\iH(\ca)$.
			\item
			The algebra $\widetilde{\ct}(\ca)$ is isomorphic to the group algebra of the abelian group $K_0(\ca)$.
			\item
			$\iH(\ca)$ has an ($\imath$Hall) basis given by
			\begin{align*}
				\{[M]*[K_\alpha]\mid [M]\in\Iso(\ca), \alpha\in K_0(\ca)\}.
			\end{align*}
		\end{enumerate}
	\end{proposition}

	For any $f:X\rightarrow Y$ in $\ca$, we denote by
	\begin{align*}
		C_f:=\Big( Y\oplus X,  \begin{pmatrix} 0&f \\0&0  \end{pmatrix}  \Big)\in\cc_1(\ca).
	\end{align*}
	
	\begin{lemma}
		[\text{\cite[Lemma 2.10]{LRW20a}}]
		\label{lem:Cf}
		For any $M^\bullet=(M,d)$, we have $[M^\bullet]=[H(M^\bullet)]*[K_{\Im (d)}]$ in $\iH(\ca)$. In particular, for any $f:X\rightarrow Y$, we have
		\begin{align*}
			[C_f]= [\ker (f)\oplus \coker(f)]*[K_{\Im (f)}].
		\end{align*}
	\end{lemma}

	\subsection{Realization of $\imath$quantum groups via quivers}
	\label{subsec:iquiver algebra}
	%Let us recall from  \cite[\S 2--3]{LW19a} (see also \cite{LW20a}) the $\imath$quiver algebra of split type since we only consider the $\imath$quantum groups of split type.
	
	Let $Q=(Q_0,Q_1)$ be a quiver (not necessarily acyclic), and we sometimes write $\I =Q_0$.
	%Let $\bfk=\F_{q}$ be a finite field of $q$ elements. 
	For a quiver $Q$, we denote
	
	$\triangleright$ $\rep_\bfk(Q)$ -- category of finite-dimensional representations of $Q$ over $\bfk$,
	
	$\triangleright$ $\rep^{\rm nil}_\bfk(Q)$ -- subcategory of $\rep_\bfk(Q)$ consisting of nilpotent representations of $Q$.
	
	%by $\rep_\bfk(Q)$ the category of finite-dimensional representations over the field $\bfk$, and denote by $\rep^{\rm nil}_\bfk(Q)$ the subcategory of $\rep_\bfk(Q)$ consisting of  nilpotent representations.
	\noindent Note that $\rep_{\bfk}(Q)= \rep^{\rm nil}_\bfk(Q)$ if $Q$ is acyclic.
	
	%	For a quiver with relations $(Q,I)$, let $\Lambda=\bfk Q/(I)$ be its (not necessarily finite-dimensional) quiver algebra. We define $\rep_\bfk(Q,I)$ and $\rep_\bfk^{\rm nil}(Q,I)$ similarly. We also denote  $\rep(\Lambda)\stackrel{\rm def}{=}\rep_\bfk(Q,I)$ and $\rep^{\rm nil}(\Lambda)\stackrel{\rm def}{=}\rep^{\rm nil}_\bfk(Q,I)$. Note that $\rep(\Lambda)=\rep^{\rm nil}(\Lambda)$ if $\Lambda$ is finite-dimensional.

	Let $S_i$ be the simple $\bfk Q$-module associated to $i\in\I$. 
	We denote by $\iH(\bfk Q)$ the twisted semi-derived Ringel-Hall algebra of $\cc_1(\rep_\bfk^{\rm nil}(Q))$. Then we have the following result.

	\begin{theorem}
		[\text{\cite{LW20a}}]
		\label{thm:main}
		Let $Q$ be an arbitrary quiver. Then there exists a $\Q(\sqq)$-algebra embedding
		\begin{align*}
			\widetilde{\psi}_Q: \tUi_{|v= \sqq} &\longrightarrow \iH(\bfk Q),
		\end{align*}
		which sends
		\begin{align}
			B_i \mapsto \frac{-1}{q-1}[S_{i}],
			&\qquad
			\K_i \mapsto [K_{S_i}], \qquad\text{ for }i \in \I.
			\label{eq:split}
		\end{align}
	\end{theorem}
	Let $\tCMH$  be the composition subalgebra of $\iH(\bfk Q)$ generated by $[S_i]$ and $[K_{S_i}]^{\pm1}$ for $i\in\I$.
	Then \eqref{eq:split} gives a $\Q(\sqq)$-algebra isomorphism
	\begin{align}
		\label{eq:isocomp}
		\widetilde{\psi}_Q: \tUi_{|v= \sqq} &\stackrel{\cong}{\longrightarrow} \tCMH.
	\end{align}

	\subsection{An $\imath$Hall multiplication formula}
	
	%Let $\ce$ be an exact category. For any short exact sequence $0\rightarrow A \xrightarrow{f} B\xrightarrow{g}C\rightarrow0$ in $\ce$, we denote by $\ov{(f,g)}$ the corresponding element in $\Ext^1_\ce(C,A)$.
	The following multiplication formulas in the $\imath$Hall algebra are useful and applicable.
	\begin{proposition}
		[\text{\cite[Proposition 3.10]{LW20a}}]
		\label{prop:iHallmult}
		Let $\ca$ be a hereditary abelian category over $\bfk$.
		For any $A,B\in\ca\subset \cc_1(\ca)$, the following formulas hold in $\iH(\ca)$.
		%\begin{align}
		%\label{Hallmult}
		%[A]\diamond[B]=&%\sum_{\alpha\in K_0(\ca)}
		%\sum_{[M]\in\Iso(\ca)}\sum_{[L],[N]\in\Iso(\ca)} q^{ \langle M,A-N\rangle} q^{\langle N,B\rangle-\langle N,A\rangle +\langle N,N\rangle -\langle A,B\rangle}\frac{|\Ext^1(N, %L)_{M}|}{|\Hom(N,L)| }
		%\\
		%\notag
		%&\cdot |\{s\in\Hom(A,B)\mid \Ker s\cong N, \coker s\cong L%\widehat{\Im d}=\alpha
		%\}|[M]\diamond[K_{\widehat{A}-\widehat{N}}].
		%\end{align}
		\begin{align}
			\label{Hallmult1}
			[A]*[B]=&%\sum_{\alpha\in K_0(\ca)}
			\sum_{[M]\in\Iso(\ca)}\sum_{[L],[N]\in\Iso(\ca)} \sqq^{-\langle A,B\rangle}  q^{\langle N,L\rangle}\frac{|\Ext^1(N, L)_{M}|}{|\Hom(N,L)| }
			\\
			\notag
			&\cdot |\{f\in\Hom(A,B)\mid \ker (f)\cong N, \coker(f)\cong L%\widehat{\Im d}=\alpha
			\}|[M]*[K_{\widehat{A}-\widehat{N}}],
			\\
			\label{Hallmult2}
			[A]*[B]=&%\sum_{\alpha\in K_0(\ca)}
			\sum_{[M]\in\Iso(\ca)}\sum_{[L],[N]\in\Iso(\ca)} \sqq^{-\langle A,B\rangle}  \frac{|\Ext^1(N, L)_{M}|}{|\Ext^1(N,L)| }
			\\
			\notag
			&\cdot |\{f\in\Hom(A,B)\mid \ker (f)\cong N, \coker(f)\cong L%\widehat{\Im d}=\alpha
			\}|[M]*[K_{\widehat{A}-\widehat{N}}]
		\end{align}
		%$\widehat{A}+\widehat{B}-\widehat{M}}{2}$
	\end{proposition}
	
	\begin{proof}
		The proof of \eqref{Hallmult1} is completely same to \cite[Proposition 3.10]{LW20a}, hence omitted here. \eqref{Hallmult2} follows from \eqref{Hallmult1} by a direct computation.
	\end{proof}
	
	\begin{corollary}
		\label{ihall formula}
		%Let $\ca$ be a hereditary abelian category over $\bfk$. 
		For any $A,B\in\ca$, the following formulas hold in $\iH(\ca)$.
		
		(1) If any non-zero map in $\Hom(A,B)$ is injective, then
		\begin{align*}
			[\![A]\!]*[\![B]\!]=&\sqq^{\langle A,B\rangle}
			\sum_{[M]\in\Iso(\ca)}F_{A,B}^M[\![M]\!]
			\\
			&+\sqq^{{-}\langle A,B\rangle}
			\sum_{[M]\in\Iso(\ca)}\frac{|\Ext^1(M, A)_{B}|}{|\Hom(M,A)|\cdot |\Aut(A)| }[\![M]\!]*[K_{\widehat{A}}].
		\end{align*}
		(2) If any non-zero map in $\Hom(A,B)$ is surjective, then
		\begin{align*}
			[\![A]\!]*[\![B]\!]=&\sqq^{\langle A,B\rangle}
			\sum_{[M]\in\Iso(\ca)}F_{A,B}^M[\![M]\!]
			\\
			&+\sqq^{{-}\langle A,B\rangle}
			\sum_{[M]\in\Iso(\ca)}\frac{|\Ext^1(B,M)_{A}|}{|\Hom(B,M)|\cdot |\Aut(B)| }[\![M]\!]*[K_{\widehat{B}}].
		\end{align*}
	\end{corollary}
	
	\begin{proof}
		We only prove (1) since (2) is dual. For any non-zero map $f:A\rightarrow B$, by assumption $f$ is injective, i.e., $\Ker (f)=0$. So
		it follows from \eqref{Hallmult1} that
		\begin{align*}
			[A]*[B]=&\sum_{[M]\in\Iso(\ca)}\sqq^{\langle A,B\rangle} \frac{|\Ext^1(A,B)_M|}{|\Hom(A,B)|}[M]\\
			&+ \sum_{[M]\in\Iso(\ca)}\sqq^{-\langle A,B\rangle} |\{0\neq f\in\Hom(A,B)\mid  \coker(f)\cong M%\widehat{\Im d}=\alpha
			\}|\cdot [M]*[K_{\widehat{A}}]
			\\
			=&\sum_{[M]\in\Iso(\ca)}\Big(\sqq^{\langle A,B\rangle} \frac{|\Ext^1(A,B)_M|}{|\Hom(A,B)|}[M]+ \sqq^{-\langle A,B\rangle}F_{M,A}^B|\Aut(A)|\cdot [M]*[K_{\widehat{A}}]\Big).
		\end{align*}

		Then the desired formula follows from the Riedtmann-Peng formula and \eqref{eq:doublebrackets}.
	\end{proof}

	%%%%%%%%%%%%%%%%%%
	\section{Weighted projective lines and coherent sheaves}
	\label{sec:WPL}
	%%%%%%%%%%%%%%%
	
	In this section, we recall some basic facts on the category of coherent sheaves on a weighted projective line given by Geigle and Lenzing \cite{GL87}.
	%Throughout this section, we assume $\bfk=\mathbb{F}_{q}$ be a finite field of $q$ elements.
	%\red{$\rep(\Lambda)$ vs $\mod(\Lambda)$? }

	\subsection{Weighted projective lines}
	Recall $\bfk=\F_q$.
	Fix a positive integer $\bs$ such that $2\leq \bs\leq q$ in the following.
	Let $\bp=(p_1,p_2,\dots,p_\bs)$ with $p_i\geq1$. Without loss of generality, we assume $p_i\geq2$ for $1\leq i\leq \bt$, and $p_i=1$ for $\bt+1\leq i\leq \bs$; compare with the star-shaped graph $\Gamma=\T_{p_1,\dots,p_\bt}$.
	
	Let $L(\bp)$
	be the rank one abelian group on generators $\vec{x}_1$, $\vec{x}_2$, $\cdots$, $\vec{x}_{\bs}$ with relations
	$p_1\vec{x}_1=p_2\vec{x}_2=\cdots=p_{\bs}\vec{x}_{\bs}$. We call $\vec{c}:=p_i\vec{x}_i$ the \emph{canonical element} of $L(\bp)$.
	Obviously, the polynomial ring $\bfk[X_1,\dots,X_{\bs}]$ is a $L(\bp)$-graded algebra by setting $\deg X_i=\vec{x}_i$, which is denoted by $\bS(\bp)$.
	
	Let $\ul{\bla}=\{\bla_1,\dots,\bla_{\mathbf{s}}\}$ be a collection of distinguished closed points (of degree one) on the projective line $\PL$, normalized such that $\bla_1=\infty$, $\bla_2=0$, $\bla_3=1$.
	Let $I(\bp,\ul{\bla})$ be the $L(\bp)$-graded ideal of $\bS(\bp)$ generated by
	$$\{X_i^{p_i}-(X_2^{p_2}-\bla_iX_1^{p_1})\mid 3\leq i\leq \bs \}.$$
	Then $\bS(\bp,\ul{\bla}):=\bS(\bp)/I(\bp,\ul{\bla})$ is an $L(\bp)$-graded algebra.
	For any $1\leq i\leq \bs$, denote by $x_i$ the image of $X_i$ in $\bS(\bp,\ul{\bla})$.
	
	The \emph{weighted projective line} $\X:=\X_{\bp,\ul{\bla}}$ is the set of all non-maximal prime homogeneous ideals of $\bS:=\bS(\bp,\ul{\bla})$, which is also denoted by $\X_\bfk$ to emphasis the base field $\bfk$.
	For any homogeneous element $f\in \bS$, let $V_f:=\{\fp\in\X\mid f\in\fp\}$ and $D_f:=\X\setminus V_f$. The set $\{D_f\mid f\in \bS\}$ forms a basis of the Zariski topology on $\X$.
	The pair $(\bp,\underline{\bla})$ (or just $\bp$ if $\underline{\bla}$ is obvious) is called the weight type of $\X_{\bp,\ul{\bla}}$, and $p_i$ is the weight of the point $\bla_i$.  Roughly speaking, $\X$ is the projective line $\PL$ with a finite number of marked points $\bla_1,\dots,\bla_\bs$ with attached weights $p_1,\dots,p_\bs$.
	In particular, the weighted projective line $\X_\bfk$ with $\bp=(1,1)$ is $\P^1_\bfk$.

	The classification of the closed points in $\X$ is provided in  \cite[Proposition 1.3]{GL87}. First, each $\bla_i$ corresponds to the prime ideal generated by $x_i$, called the \emph{exceptional point}. Second, any other homogeneous prime ideal is of the form $(f(x_1^{p_1},x_2^{p_2}))$, where $f\in \bfk [y_1,y_2]$ is an irreducible homogeneous polynomial in $y_1,y_2$, which is different from $y_1$ and $y_2$, called the \emph{ordinary point}.
	
	For any closed point $x$ in $\X$, define
	\begin{align}
		\label{def:dx}
		d_x:=\begin{cases} 1, & \text{ if }x=\bla_i \text{ (exceptional)},
			\\
			\deg(f), & \text{ if }x=(f(x_1^{p_1},x_2^{p_2})) \text{ (ordinary)};
		\end{cases}
	\end{align}
	compare with the degree of $x$ defined in \cite{Sch04,DJX12,BS13}.

	\subsection{Coherent sheaves}
	
	%For any homogeneous element $f\in S$, let $V_f:=\{\fp\in\X\mid f\in\fp\}$ and $D_f:=\X\setminus V_f$.
	The structure sheaf $\co:=\co_\X$ is the sheaf of $L(\bp)$-graded algebras over $\X$ associated to the presheaf $D_f\mapsto \bS_f$, where $\bS_f:=\{gf^{-l}\mid g\in \bS,l>0\}$ is the localization of $\bS$. Let $\Mod(\co_\X)$ be the category of sheaves of $L(\bp)$-graded $\co_\X$-modules.
	For any $\cm\in\Mod(\co_\X)$, we denote by $\cm(\vec{x})$ the twisting of $\cm$ by $\vec{x}\in L(\bp)$.
	
	A sheaf $\cm\in\Mod(\co_\X)$ is called coherent if for any ${\blx}\in\X$, there exists a neighbourhood $U$ of ${\blx}$ and an exact sequence
	\begin{align}
		\bigoplus_{j=1}^n\co_\X(\vec{l}_j)|_{U}\rightarrow\bigoplus_{i=1}^m\co_\X(\vec{k}_i)|_{U}\rightarrow \cm|_{U}\rightarrow 0.
	\end{align}
	We denote by $\coh(\X)$ the full subcategory of $\Mod(\co_\X)$ consisting of all coherent sheaves on $\X$. Then $\coh(\X)$ is a $\bfk$-linear hereditary, Hom-finite and Ext-finite abelian category; see \cite[Subsection 2.2]{GL87}. Moreover, by \cite[Proposition 1.8]{GL87}, $\coh(\X)$ is equivalent to the Serre quotient $\mod^{L(\bp)}(\bS)/ \mod^{L(\bp)}_0(\bS)$ of the category of finitely generated $L(\bp)$-graded
	$\bS$-modules modulo the subcategory of finite-length modules.

	A coherent sheaf $\cm\in\coh(\X)$ is called a \emph{vector bundle} (or \emph{locally free sheaf})   of rank $n$ if there is an open covering $\{U_i\}$ of $\X$, and an isomorphism
	$\bigoplus_{j=1}^n \co(\vec{l}_j)|_{U_i}\xrightarrow{\simeq} \cm|_{U_i}$ for some suitable $\vec{l}_j$ for each $U_i$; and is called a \emph{torsion sheaf} if it is a finite-length object in $\coh(\X)$.
	Let $\scrf$ be the full subcategory of $\coh(\X)$ consisting of all locally free sheaves, and $\scrt$ be the full subcategory consisting of all torsion sheaves. Then $\scrf$ and $\scrt$ are extension-closed. Moreover, $\scrt$ is a hereditary length abelian category.
	
	\begin{lemma}[\cite{GL87}]
		(1) For any sheaf $\cm\in\coh(\X)$, it can be decomposed as $\cm_t\oplus \cm_f$ where $\cm_f\in\scrf$ and $\cm_t\in \scrt$;
		
		(2) $\Hom(\cm_t,\cm_f)=0=\Ext^1(\cm_f,\cm_t)$ for any $\cm_f\in\scrf$ and $\cm_t\in\scrt$.
	\end{lemma}

	\subsection{Torsion sheaves}
	
	In order to describe the category $\scrt$, we shall introduce the representation theory of cyclic quivers.
	We consider the oriented cyclic quiver $C_n$ for $n\geq2$ with the vertex set $\Z/n\Z=\{0,1,2,\dots,n-2,n-1\}$:
	%\[\xymatrix{ &&0\red{=n} \ar[drr]&&  \\
	%1\ar[urr]&2 \ar[l]&\cdots\ar[l]&n-2\ar[l]&n-1 \ar[l]}\]
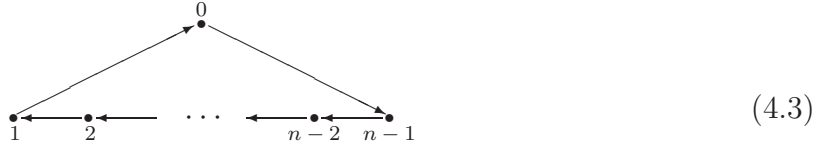
	\begin{center}\setlength{\unitlength}{0.5mm}
		%\vspace{-2cm}
		\begin{equation}
			\label{fig:Cn}
			\begin{picture}(100,30)
				\put(2,0){\circle*{2}}
				\put(22,0){\circle*{2}}
				
				\put(82,0){\circle*{2}}
				\put(102,0){\circle*{2}}
				\put(52,25){\circle*{2}}
				
				\put(20,0){\vector(-1,0){16}}
				\put(40,0){\vector(-1,0){16}}
				\put(47.5,-2){$\cdots$}
				\put(80,0){\vector(-1,0){16}}
				\put(100,0){\vector(-1,0){16}}
				
				\put(54,24.5){\vector(2,-1){47}}
				\put(3,1){\vector(2,1){47}}
				
				\put(50.5,27){\tiny $0$}
				\put(1,-6){\tiny $1$}
				\put(21,-6){\tiny $2$}
				%\put(41,-6){\tiny $3$}
				\put(75,-6){\tiny $n-2$}
				\put(95,-6){\tiny $n-1$}
			\end{picture}
		\end{equation}
		\vspace{-0.2cm}
	\end{center}
	Let $C_1$ be the Jordan quiver, i.e., the quiver with only one vertex and one loop arrow.
	
	%Denote by $\rep^{\rm nil}_\bfk(C_n)$ the category of finite-dimensional nilpotent representations of $C_n$ over the field $\bfk$.

	Then the structure of $\scrt$ is described in the following.
	\begin{lemma}[\cite{GL87}; also see \cite{Sch04}]
		\label{lem:isoclasses Tor}
		(1) The category $\scrt$ decomposes as a coproduct $\scrt=\coprod_{x\in\X} \scrt_{x}$, where $\scrt_{x}$ is the subcategory of torsion sheaves with support at $x$.
		
		(2) For any ordinary point $x$ of degree $d$, let $\bfk_{x}$ denote the residue field at $x$, i.e., $[\bfk_x:\bfk]=d$. Then $\scrt_{x}$ is equivalent to the category $\rep^{\rm nil}_{\bfk_{x}}(C_1)$.
		
		(3) For any exceptional point $\bla_i$ ($1\leq i\leq \bs$), the category $\scrt_{\bla_i}$ is equivalent to $\rep^{\rm nil}_{\bfk}(C_{p_i})$.
	\end{lemma}

	For any ordinary point $\blx$ of degree $d$, let $\pi_{\blx}$ be the prime homogeneous polynomial corresponding to $\blx$. The multiplication by $\pi_{\blx}$ gives the exact sequence
	$$0\longrightarrow \co\stackrel{\pi_{\blx}}{\longrightarrow} \co(d\vec{c})\longrightarrow S_{\blx}\rightarrow0,$$
	where $S_{\blx}$ is the unique (up to isomorphism) simple sheaf in the category $\scrt_{\blx}$. Then $S_{\blx}(\vec{l})=S_{\blx}$ for any $\vec{l}\in L(\bp)$.
	
	For any exceptional point $\bla_i$, multiplication by $x_i$ yields the short exact sequence
	$$0\longrightarrow \co((j-1)\vec{x}_i)\stackrel{x_i}{\longrightarrow} \co(j\vec{x}_i)\longrightarrow S_{ij}\rightarrow0,\text{ for } 1\le j\le p_i;$$
	where $\{S_{ij}\mid  j\in\Z_{p_i}\}$ is a complete set of pairwise non-isomorphic simple sheaves in the category $\scrt_{\bla_i}$ for any $1\leq i\leq \bs$. Moreover, $S_{ij}(\vec{l})=S_{i,j+l_i}$ for any $\vec{l}=\sum_{1\leq i\leq \bt} l_i\vec{x}_i$.

	In order to describe the indecomposable objects in $\scrt$, we need the following well-known results.
	%\begin{itemize}
	% \item[(1)] {\red{For $\rep^{\rm nil}_\bfk(C_1)$, the set of isoclasses of indecomposable objects is $\{S^{(a)}\mid a\in\N\}$, where $S=S^{(1)}$ is the only simple representation, and $S^{(a)}$ is the unique indecomposable representation of length $a$. Define $S^{(\mu)}=\bigoplus_{i} S^{(\mu_i)}$ for any partition $\mu=(\mu_1\geq \cdots\geq \mu_l)$. Then any object in $\rep^{\rm nil}_\bfk(C_1)$ is isomorphic to $S^{(\mu)}$ for some partition $\mu$.}}
	%\item[(2)]  For $\rep^{\rm nil}_\bfk(C_n)$, we have $n$ simple representations $S_j=S_j^{(1)}$ ($1\leq j\leq n$). Denote by $S_j^{(a)}$ the unique indecomposable representation with top $S_j$ and length $a$. Then $\{S_j^{(a)}\mid 1\le j\le n, a\in\N\}$ is  the set of isoclasses of indecomposable objects of $\rep^{\rm nil}_\bfk(C_n)$. For any $1\le j\le n$, define $S_j^{(\mu)}=\bigoplus_{i} S_j^{(\mu_i)}$ for any partition $\mu=(\mu_1\geq \cdots\geq \mu_l)$. %Then any object in $\rep^{\rm nil}_\bfk(C_1)$ is isomorphic to $S^{(\mu)}$ for some partition $\mu$.
	%\end{itemize}
	For $\rep^{\rm nil}_\bfk(C_n)$, we have $n$ simple representations $S_j=S_j^{(1)}$ ($1\leq j\leq n$). Denote by $S_j^{(a)}$ the unique indecomposable representation with top $S_j$ and length $a$. Then $\{S_j^{(a)}\mid 1\le j\le n, a\in\N\}$ is  the set of isoclasses of indecomposable objects of $\rep^{\rm nil}_\bfk(C_n)$. For any $1\le j\le n$, define $S_j^{(\mu)}=\bigoplus_{i} S_j^{(\mu_i)}$ for any partition $\mu=(\mu_1\geq \cdots\geq \mu_l)$. %Then any object in $\rep^{\rm nil}_\bfk(C_1)$ is isomorphic to $S^{(\mu)}$ for some partition $\mu$.
	Combining with Lemma \ref{lem:isoclasses Tor}, we can give a classification of indecomposable objects in $\scrt$.
	
	We denote by $S_{\blx}^{(a)}$ the unique indecomposable object of length $a$ in $\scrt_{\blx}$ for any ordinary point ${\blx}$; and denote by $S_{ij}^{(a)}$ the unique indecomposable object with top $S_{ij}$ and length $a$ in $\scrt_{\bla_i}$.

	\subsection{Grothendieck groups and Euler forms}
	%Let $\X$ be a weighted projective line.
	The Grothendieck group $K_0(\coh(\X))$ of $\coh(\X)$ satisfies
	\begin{align}
		K_0(\coh(\X))\cong \Big(\Z\widehat{\co}\oplus \Z\widehat{\co(\vec{c})}\oplus \bigoplus_{i,j}\Z\widehat{S_{ij}}\Big)/I,
	\end{align}
	where $I$ is the subgroup generated by $\{\sum_{j=1}^{p_i} \widehat{S_{ij}} +\widehat{\co} -\widehat{\co(\vec{c})}\mid i=1,\dots,\bs\}$; see \cite{GL87}.
	Let $\de:=\widehat{\co(\vec{c})}-\widehat{\co}= \sum_{j=1}^{p_i} \widehat{S_{ij}}$. Note that $\widehat{S_{i1}}=\de$ if $p_i=1$.
	The Euler form $\langle-,-\rangle$ on $K_0(\coh(\X))$ is described in \cite{Sch04}, which is given by
	\begin{align}
		&\langle \widehat{\co},\widehat{\co}\rangle=1,\quad \langle \widehat{\co},\delta\rangle=1,\quad  \langle\delta,\widehat{\co}\rangle=-1,
		\\
		&\langle \delta,\delta\rangle=0,\quad \langle \widehat{S_{ij}},\delta\rangle=0,\quad  \langle \delta, \widehat{S_{ij}}\rangle=0,
		\\
		&\langle \widehat{\co}, \widehat{S_{ij}}\rangle =\begin{cases} 1& \text{ if } j=p_i, \\ 0& \text{ if }j\neq p_i , \end{cases}
		\qquad \langle \widehat{S_{ij}}, \widehat{\co}\rangle =\begin{cases} -1& \text{ if } j=1, \\ 0& \text{ if }j\neq 1 ,\end{cases}
		\\
		&\langle\widehat{ S_{ij}}, \widehat{S_{i',j'}}\rangle =\begin{cases} 1& \text{ if } i=i', j=j',
			\\ -1& \text{ if }i=i', j\equiv j'+1 (\mod p_i),
			\\
			0& \text{ otherwise. }\end{cases}
	\end{align}
	%The symmetric Euler form $(-,-)$ on $K_0(\coh(\X))$ is defined by $(\alpha,\beta)=\langle \alpha,\beta\rangle+\langle \beta,\alpha\rangle$.
	
	Recall the root system $\cR$ of $L\fg$ of the star-shaped graph $\Gamma=\T_{p_1,\dots,p_\bt}$. Then there is a natural isomorphism of $\Z$-modules $K_0(\coh(\X))\cong \cR$ given as below (for $1\leq i\leq \bt,\;1\le j\le p_i-1, $ $r\in\N$ and $l\in\Z$):
	\begin{align*}
		&\widehat{S_{ij}} \mapsto \alpha_{ij}, %\text{ for }1\leq i\leq \bt,\;1\le j\le p_i-1,
		&\widehat{S_{i,0}}\mapsto \de-\sum_{j=1}^{p_i-1}\alpha_{ij},\quad% \text{ for }1\leq i\leq \bt,
		&\widehat{S_{i,0}^{(p_i-1)}}\mapsto \de-\alpha_{i1}, %\text{ for }1\leq i\leq \bt,
		&\widehat{S_{i,0}^{(rp_i)}}\mapsto r\de,\quad
		&\widehat{\co(l\vec{c})}\mapsto \alpha_\star+l\de.% \text{ for }k\in\Z.
	\end{align*}
	Under this isomorphism, the symmetric Euler form on $K_0(\coh(\X))$ coincides with the Cartan form on $\cR$. So we always identify $K_0(\coh(\X))$ with $\cR$ via this isomorphism in the following.
	
	Let $p=\mathrm{l.c.m}(p_1,\dots,p_\bt)$. We can define the degree and rank for each coherent sheaf, which induces two functions on $K_0(\coh(\X))$ denoted by $\deg$ and $\rank$ respectively. In fact, they are uniquely determined by
	\begin{align}
		\deg(\co)=0, \quad\deg(S_{ij})=\frac{p}{p_i};&&\rank(\co)=1,\quad \rank(S_{ij})=0;
	\end{align}
	for any $1\leq i\leq \bt, 0\leq j\leq p_i-1$.
	
	%%%%%%%%%%%%%
	\subsection{$\imath$Hall algebras of weighted projective lines}
	%%%%%%%%%%%%

	By \cite{GL87}, we know that $\coh(\X)$ is a hereditary abelian category, so we can define the $\imath$Hall algebra (a.k.a. the twisted semi-derived Ringel-Hall algebra) of $\X$, which is denoted by $\iH(\X_\bfk)$.
	
	%\begin{lemma}
	%Twisting with $\vec{c}$ gives an isomorphism of $\iH(\X_\bfk)$.
	%\end{lemma}
	
	%\begin{proof}
	%It is well know that twisting with $\vec{c}$ induces an equivalence of $\coh(\X)$, and then an equivalence of $\cc_1(\coh(\X))$. So it gives an isomorphism of $\widetilde{\ch}(\X_\bfk)$, and then an isomorphism of $\iH(\X_\bfk)$.
	%\end{proof}

	A direct computation in $\iH(\X_\bfk)$ shows that for any $l\in\mathbb{Z}$ and $1\leq i\leq \bt$:
	%\begin{align} \label{co and si1}
	%[\co(l\vec{c})]*[S_{i,1}]=& [ S_{i,1}\oplus\co(l\vec{c})],
	%\\\label{si1 and co}
	%[S_{i,1}]* [\co(l\vec{c})]=& \sqq^{-1}[ S_{i,1}\oplus\co(l\vec{c})]+(\sqq-\sqq^{-1})[\co(l\vec{c}+\vec{x}_i)],
	%\\\label{si0 and co}
	%[S_{i,0}^{(p_i-1)}]* [\co(l\vec{c})]=& [S_{i,0}^{(p_i-1)}\oplus\co(l\vec{c})], %+\sqq^{-1}(q-1)[\co(l\vec{c}+\vec{x}_i)],
	%\\\label{co and si0}
	%[\co(l\vec{c})]*[S_{i,0}^{(p_i-1)}]=& \sqq^{-1}[S_{i,0}^{(p_i-1)}\oplus\co(l\vec{c})]+(\sqq-\sqq^{-1})[\co((l-1)\vec{c}+\vec{x}_i)]*[K_{ \de-\alpha_{11}}]. \end{align}
	
	\begin{align} \begin{split}\label{Si1 and co double direction}
			[\co(l\vec{c})]*[S_{i,1}]=& [ S_{i,1}\oplus\co(l\vec{c})],\\
			[S_{i,1}]* [\co(l\vec{c})]=& \sqq^{-1}[ S_{i,1}\oplus\co(l\vec{c})]+(\sqq-\sqq^{-1})[\co(l\vec{c}+\vec{x}_i)];
	\end{split}\end{align}
	and
	\begin{align} \begin{split}\label{Si0 and co double direction}
			[S_{i,0}^{(p_i-1)}]* [\co(l\vec{c})]=& [S_{i,0}^{(p_i-1)}\oplus\co(l\vec{c})],\\
			[\co(l\vec{c})]*[S_{i,0}^{(p_i-1)}]=& \sqq^{-1}[S_{i,0}^{(p_i-1)}\oplus\co(l\vec{c})]+(\sqq-\sqq^{-1})[\co((l-1)\vec{c}+\vec{x}_i)]*[K_{ \de-\alpha_{i}}].
	\end{split}\end{align}
	
	The following results will be used in the proof of Lemma \ref{level one pi ij and star} and Proposition \ref{prop:coBil}. %, and their proofs are in Appendix~\ref{sec:proofSecWPL}.
	\begin{lemma}
		\label{middle and ending terms same}
		For any $1\leq i\leq \bt$ and $1\leq j< k\leq p_i$, the following equations hold in $\iH(\X_\bfk)$:
		\begin{align}
			\label{eq:OSij}
			\big[[\co(l\vec{c})], [S_{ij}^{(k)}]\big]=&\big[[\co(l\vec{c})], [S_{ij}^{(j)}\oplus S_{i,0}^{(k-j)}]\big],\quad\text{ for any }l\in\Z;
			\\
			\label{eq:Sipiij}
			\big[[S_{i,0}^{(rp_i)}], [S_{ij}^{(k)}]\big]=&\big[[S_{i,0}^{(rp_i)}], [S_{ij}^{(j)}\oplus S_{i,0}^{(k-j)}]\big],\quad\text{ for any }r\geq 1.
		\end{align}
	\end{lemma}

	\begin{proof}%[Proof of Lemma \ref{middle and ending terms same}]
		First, we prove \eqref{eq:OSij}. %The proof for the second one is similar.
		Observe that we have the following exact sequence
		\begin{equation}
			\label{exact sequence for sij}
			0\rightarrow S_{i,0}^{(k-j)}\longrightarrow S_{ij}^{(k)} \longrightarrow S_{ij}^{(j)}\rightarrow0.
		\end{equation}
		The following computations are in $\cc_1(\coh(\X))$.
		
		By applying $\Hom_{\cc_1(\coh(\X))}(\co(l\vec{c}),-)$ and $\Hom_{\cc_1(\coh(\X))}(-,\co(l\vec{c}))$ to (\ref{exact sequence for sij}), we know that
		\begin{align*}
			\Hom(\co(l\vec{c}),S_{ij}^{(k)})\cong \Hom(\co(l\vec{c}),S_{i,0}^{(k-j)})\cong \bfk,\quad
			\Ext^1(S_{ij}^{(k)}, \co(l\vec{c}))\cong \Ext^1(S_{ij}^{(j)}, \co(l\vec{c}))\cong \bfk.
		\end{align*}
		Let $X\in \{S_{ij}^{(k)}, S_{ij}^{(j)}\oplus S_{i,0}^{(k-j)}\}$.
		For any non-zero morphism $f:\co(l\vec{c})\rightarrow X$, there are
		$$\Im(f)=S_{i,0}^{(k-j)}, \quad\ker(f)=\co(l\vec{c}-(k-j)\vec{x}_i),\quad \coker(f)=S_{ij}^{(j)},$$ while the non-trivial extension $\Ext^1(X, \co(l\vec{c}))$ is given by
		$$0\rightarrow \co(l\vec{c})\longrightarrow \co(l\vec{c}+j\vec{x}_i)\oplus S_{i,0}^{(k-j)} \longrightarrow X\rightarrow0.$$
		Therefore, by \eqref{Hallmult1} we have
		\begin{align*}
			&\big[[\co(l\vec{c})], [S_{ij}^{(k)}]\big]\\
			&=[\co(l\vec{c})]\ast [S_{ij}^{(k)}]-[S_{ij}^{(k)}]\ast[\co(l\vec{c})]\\
			&=\sqq^{\langle \widehat{\co(l\vec{c})}, \widehat{S_{ij}^{(k)}}\rangle}\cdot\Big(\frac{1}{q}[\co(l\vec{c})\oplus S_{ij}^{(k)}]+\frac{q-1}{q}[\co(l\vec{c}-(k-j)\vec{x}_i)\oplus S_{ij}^{(j)}]\ast [K_{S_{i,0}^{(k-j)}}]\Big)\\
			&\quad -\sqq^{\langle \widehat{S_{ij}^{(k)}},\widehat{\co(l\vec{c})}\rangle}\cdot\Big([\co(l\vec{c})\oplus S_{ij}^{(k)}]+(q-1)[\co(l\vec{c}+j\vec{x}_i)\oplus S_{i,0}^{(k-j)}]\Big)\\
			&=(\sqq-\sqq^{-1})\Big([\co(l\vec{c}-(k-j)\vec{x}_i)\oplus S_{ij}^{(j)}]* [K_{S_{i,0}^{(k-j)}}]-[\co(l\vec{c}+j\vec{x}_i)\oplus S_{i,0}^{(k-j)}]\Big).
		\end{align*}
		Similarly,
		%On the other hand, for any non-zero morphism $g:\co(l\vec{c})\rightarrow S_{ij}^{(j)}\oplus S_{i,0}^{(k-j)}$, there are
		%$$\Im(g)=S_{i,0}^{(k-j)}, \quad\ker(g)=\co(l\vec{c}-(k-j)\vec{x}_i),\quad \coker(g)=S_{ij}^{(j)};$$ while the non-trivial extension $\Ext^1(S_{ij}^{(j)}\oplus S_{i,0}^{(k-j)}, \co(l\vec{c}))$ is given by
		%$$0\rightarrow \co(l\vec{c})\longrightarrow \co(l\vec{c}+j\vec{x}_i)\oplus S_{i,0}^{(k-j)} \longrightarrow S_{ij}^{(j)}\oplus S_{i,0}^{(k-j)}\rightarrow0.$$
		%It follows that
		\begin{align*}
			&\big[[\co(l\vec{c})], [S_{ij}^{(j)}\oplus S_{i,0}^{(k-j)}]\big]\\
			&=(\sqq-\sqq^{-1})\Big([\co(l\vec{c}-(k-j)\vec{x}_i)\oplus S_{ij}^{(j)}]\ast [K_{S_{i,0}^{(k-j)}}]-[\co(l\vec{c}+j\vec{x}_i)\oplus S_{i,0}^{(k-j)}]\Big),
			%=&\Big[[\co(l\vec{c})], [S_{ij}^{(k)}]\Big].
		\end{align*}
		and then \eqref{eq:OSij} follows.
		
		For \eqref{eq:Sipiij}, observe that the composition factors of $S_{ij}^{(j)}$ are $S_{i,1}, S_{i,2},\cdots, S_{ij}$ (from socle to top), while the composition factors of $S_{i,0}^{(k-j)}$ are $S_{i,p_i+1-k+j},\cdots, S_{i,p_i-1}, S_{i,0}$. Hence
		\begin{align*}
			\Hom_\X(S_{i,0}^{(rp_i)},S_{ij}^{(j)})=0= \Ext^1_\X(S_{i,0}^{(rp_i)},S_{ij}^{(j)}),
			\quad
			\Hom_\X(S_{ij}^{(k-j)}, S_{i,0}^{(rp_i)})=0=\Ext^1_\X(S_{ij}^{(k-j)}, S_{i,0}^{(rp_i)}).
		\end{align*}
		By applying $\Hom_{\X}(S_{i,0}^{(rp_i)},-)$ and $\Hom_{\X}(-, S_{i,0}^{(rp_i)})$ to (\ref{exact sequence for sij}), we obtain
		\begin{align*}
			&\Hom_\X(S_{i,0}^{(rp_i)},S_{ij}^{(k)})\cong\Hom(S_{i,0}^{(rp_i)},S_{i,0}^{(k-j)}),
			\qquad
			\Ext^1_\X(S_{i,0}^{(rp_i)},S_{ij}^{(k)})\cong\Ext^1(S_{i,0}^{(rp_i)},S_{i,0}^{(k-j)})),
			\\
			&\Hom_\X(S_{ij}^{(k)}, S_{i,0}^{(rp_i)})\cong \Hom_\X(S_{i,0}^{(k-j)}, S_{i,0}^{(rp_i)}),
			\qquad
			\Ext^1_\X(S_{ij}^{(k)}, S_{i,0}^{(rp_i)})\cong \Ext^1_\X(S_{ij}^{(j)}, S_{i,0}^{(rp_i)}),
		\end{align*}
		which are of dimension $1$.
		Then \eqref{eq:Sipiij} holds by a similar argument as \eqref{eq:OSij}.
	\end{proof}

	\begin{lemma}
		\label{pi+1 and co}
		For any $1\leq i\leq \bt$ and $l\in\mathbb{Z}$,
		\begin{align*}
			\big[[S_{i,1}^{(p_i+1)}]-[S_{i,1}\oplus S_{i,0}^{(p_i)}],[\co(l\vec{c})]\big]_{\sqq^{-1}}
			=\sqq\big[[S_{i,1}],[\co((l+1)\vec{c})]\big]_\sqq.
		\end{align*}
	\end{lemma}

	\begin{proof}%[Proof of Lemma \ref{pi+1 and co}]
		By (\ref{Si1 and co double direction}), we have
		\begin{align*}
			\big[[S_{i,1}],[\co((l+1)\vec{c})]\big]_\sqq=(\sqq-\sqq^{-1})\Big([\co((l+1)\vec{c}+\vec{x}_i)]-[\co((l+1)\vec{c})\oplus S_{i,1}]\Big).
		\end{align*}
		On the other hand, by definition we have
		\begin{align*}
			&\big[[S_{i,1}^{(p_i+1)}],[\co(l\vec{c})]\big]_{\sqq^{-1}}\\
			&=[S_{i,1}^{(p_i+1)}]*[\co(l\vec{c})]-\sqq^{-1} [\co(l\vec{c})]*[S_{i,1}^{(p_i+1)}]\\
			&=\sqq^{-2}\Big([S_{i,1}^{(p_i+1)}\oplus \co(l\vec{c})]+(q-1) [\co(l\vec{c}+\vec{x}_i)\oplus S_{i,0}^{(p_i)}]+(q^2-q)[\co((l+1)\vec{c}+\vec{x}_i)]\Big)\\
			&\quad -\sqq^{-1}\cdot \frac{\sqq}{q} \Big([S_{i,1}^{(p_i+1)}\oplus \co(l\vec{c})]+(q-1)[S_{i,1}\oplus \co((l-1)\vec{c})]*[K_\de]\Big)\\
			&=(q-1)[\co((l+1)\vec{c}+\vec{x}_i)]+\frac{q-1}{q}[\co(l\vec{c}+\vec{x}_i)\oplus S_{i,0}^{(p_i)}]-\frac{q-1}{q}[S_{i,1}\oplus \co((l-1)\vec{c})]*[K_\de],
		\end{align*}
		and
		\begin{align*}
			&\big[[S_{i,1}\oplus S_{i,0}^{(p_i)}],[\co(l\vec{c})]\big]_{\sqq^{-1}}\\
			&=[S_{i,1}\oplus S_{i,0}^{(p_i)}]*[\co(l\vec{c})]-\sqq^{-1} [\co(l\vec{c})]*[S_{i,1}\oplus S_{i,0}^{(p_i)}]\\
			&=\sqq^{-2}\Big([S_{i,1}\oplus S_{i,0}^{(p_i)}\oplus \co(l\vec{c})]+(q-1) [\co(l\vec{c}+\vec{x}_i)\oplus S_{i,0}^{(p_i)}]+(q^2-q)[\co((l+1)\vec{c})\oplus S_{i,1}]\Big)\\
			&\quad -\sqq^{-1}\cdot \frac{\sqq}{q} \Big([S_{i,1}\oplus S_{i,0}^{(p_i)}\oplus \co(l\vec{c})]+(q-1)[S_{i,1}\oplus \co((l-1)\vec{c})]*[K_\de]\Big)\\
			&=(q-1)[\co((l+1)\vec{c})\oplus S_{i,1}]+\frac{q-1}{q}[\co(l\vec{c}+\vec{x}_i)\oplus S_{i,0}^{(p_i)}]-\frac{q-1}{q}[S_{i,1}\oplus \co((l-1)\vec{c})]*[K_\de].
		\end{align*}
		Hence, by combining the above two formulas we obtain
		\begin{align*}
			&\big[[S_{i,1}^{(p_i+1)}]-[S_{i,1}\oplus S_{i,0}^{(p_i)}],[\co(l\vec{c})]\big]_{\sqq^{-1}}\\
			&=(q-1)\Big([\co((l+1)\vec{c}+\vec{x}_i)]-[\co((l+1)\vec{c})\oplus S_{i,1}]\Big)\\
			&=\sqq\big[[S_{i,1}],[\co((l+1)\vec{c})]\big]_\sqq.
		\end{align*}
	\end{proof}

	%%%%%%%%%%%
	\section{$\imath$Hall algebras of cyclic quivers}
	\label{sec:cyclic}

	Let $C_n$ be the cyclic quiver with $n$ vertices; see \eqref{fig:Cn}.
	Recall that $\rep^{\rm nil}_\bfk(C_n)$ is the category of finite-dimensional nilpotent representations of $C_n$. For the Grothendieck group $K_0(\rep^{\rm nil}_\bfk(C_n))$, we denote $\alpha_j=\widehat{S_j}$ and $\de=\sum_{j=1}^n\widehat{S_j}$ by abusing the notations. For any positive real root $\beta$ of $\widehat{\mathfrak{sl}}_{n}$, by Gabriel-Kac Theorem, we denote by $M(\beta)$ the unique
	(up to isomorphism) $\bfk C_n$-module with its class $\beta$ in $K_0(\rep^{\rm nil}_\bfk(C_n))$.
	
	For a finite-dimensional representation $M$ of $C_n$, we denote by $\mathrm{soc}(M)$ and $\mathrm{top}(M)$ the socle and top of $M$ respectively.
	
	%%%%%%%%%%%%%%%%%
	\subsection{Root vectors}
	
	Let $\tUi:=\tUi_v(\widehat{\mathfrak{sl}_n})$ be the universal affine $\imath$quantum group of type $A_{n-1}$, and $\tUiD$ be its Drinfeld presentation throughout this section.
	By Theorem \ref{thm:ADE}, there exists an algebra isomorphism $\Phi:\tUiD\rightarrow\tUi$. We choose the sign function $o(\cdot)$ such that $o(j)=(-1)^j$ for any $1\leq j\leq n-1$ throughout this paper.
	%\begin{align}
	%\Phi:{}^{Dr}\tUi(\widehat{\mathfrak{sl}}_{n+1})\rightarrow\tUi; \qquad
	%\red{\K_j\mapsto \K_{j},\quad C\mapsto C;\quad
	%B_{j,0}\mapsto B_{j}, \qquad B_{j,k}\mapsto (-1)^{jk}T_{\omega_j}^{-k}(B_{j}), \,\,\forall 1\leq j\leq n}.
	%\end{align}
	On the other hand, for $n\geq2$, \cite[Theorem 9.6]{LW20a} gives an algebra monomorphism:
	%$\widetilde{\psi}: \tUi\rightarrow \tMH$ which sends
	\begin{align}
		\label{eq:isoCnsln}
		\widetilde{\psi}_{C_n}: \tUi_{|v=\sqq}\longrightarrow &\iH(\bfk C_n)
		\\
		B_j \mapsto& \frac{-1}{q-1}[S_{j}],
		\quad  \K_j \mapsto [K_{S_j}], \text{ for }0\leq j\leq n-1. 
	\end{align}
	%If $n=2$, one can check that the above algebra monomorphism is also valid; see \cite[Remark 9.7]{LW20a} and a complete proof will be given elsewhere.
	
	Define
	\begin{align}\label{the map Psi A}
		\Omega_{C_n}:=\widetilde{\psi}_{C_n}\circ \Phi:\tUiD_{|v=\sqq}\longrightarrow \iH(\bfk C_n).
	\end{align}
	Then $\Omega_{C_n}$ sends:
	\begin{align}
		\label{eq:HaDrA1}\K_j\mapsto [K_{S_j}],\qquad C\mapsto [K_\de],\qquad
		B_{j,0}\mapsto \frac{-1}{q-1}[S_j], \qquad \forall \,\,1\leq j\leq n-1.
	\end{align}
	For any $1\leq j\leq n-1$, $l\in\Z$ and $r\geq 1$, we define
	\begin{align}
		\label{def:haB}
		\haB_{j,l}:=(1-q)\Omega_{C_n}(B_{j,l}),\qquad
		\widehat{\Theta}_{j,r}:= \Omega_{C_n}(\Theta_{j,r}),\qquad
		\widehat{H}_{j,r}:= \Omega_{C_n}(H_{j,r}).
	\end{align}
	In particular, $\haB_{j,0}=[S_j]$ for any $1\leq j\leq n-1 $.

	For $1\leq j\leq {n}$ and any $\alpha\in K_0(\rep^{\rm nil}_\bfk (C_n))$, set
	\begin{align}
		\label{def:Mjalpha}
		\cm_{j,\alpha}:=\{[M]\mid \widehat{M}=\alpha,{\rm soc}(M)\subseteq S_1\oplus  \cdots \oplus S_j\}.
	\end{align}
	%$\cm_{j,\alpha}$ be the set of all isoclasses of modules $M$ such that $\underline{\dim} M=\alpha$ and ${\rm soc}(M)\in\add(S_1\oplus  \cdots \oplus S_j)$.
	Define
	%$$\eta_{j}
	%=\frac{-\sqq^{-j}}{\sqq-\sqq^{-1}}\sum\limits_{\cm_{j+1,\delta-e_{j}}}(-1)^{\dim\End(M)}[M];$$
	\begin{align}
		\pi_{j,1}=\frac{-\sqq^{-j}}{\sqq-\sqq^{-1}}\sum\limits_{[M]\in\cm_{j,\delta}}(-1)^{\dim\End(M)}[M].
	\end{align}
	Set $\pi_{0,1}=0$.
	
	\begin{proposition}\label{prop:DrGenA}
		For any $1\leq j\leq n-1$, we have
		$$\haB_{j,-1}=\sqq^{-j+1}\sum\limits_{[M]\in\cm_{j+1,\delta-e_{j}}}(-1)^{\dim\End(M)}[M]*[K_{\de-\alpha_j}]^{-1}.$$
	\end{proposition}

	In order to prove Proposition \ref{prop:DrGenA}, we need to recall the isomorphisms of $\imath$Hall algebras induced by the BGP type reflection functors introduced in \cite{LW19b,LW21b}.
	
	Let $Q=(Q_0,Q_1)$ be a general quiver. For any sink $\ell \in Q_0$, define the quiver $\bs_\ell^+ Q$ by reversing all the arrows of $Q$ ending at $\ell $. %Denote by
	%$\bs_{\ell }^+\Lambda^{\imath}:=\bfk \bs_\ell^+Q\otimes_\bfk R_1$. The quiver $\ov{Q'}$ of $\bs^+_\ell \Lambda^{\imath}$ can be constructed from $\ov{Q}$ by reversing all the arrows of $Q$ ending at $\ell $.
	Associated to a sink $\ell \in Q_0$, the BGP reflection functor induces a {\rm reflection functor} (see \cite[\S3.2]{LW19b}):
	\begin{align}  \label{eq:Fl}
		F_\ell ^+: \cc_1(\rep_\bfk(Q))  \longrightarrow \cc_1(\rep_\bfk(\bs^+_\ell Q)).
	\end{align}
	The functor $F_\ell^+$ induces an isomorphism $\Gamma_\ell: \iH(\bfk Q)  \stackrel{\sim}{\rightarrow} \iH(\bfk \bs_\ell^+ Q)$ by \cite[Theorem 4.3]{LW19b},
	and then a commutative diagram
	\begin{equation}
		\label{eq:comm}
		\xymatrix{ \tUi_{ |v={\sqq}} \ar[r]^{ \TT_\ell} \ar[d]^{\widetilde{\psi}_{Q}} & \tUi_{ |v={\sqq}} \ar[d]^{\widetilde{\psi}_{\bs_\ell Q}}
			\\
			\iH(\bfk Q) \ar[r]^{\Gamma_\ell}  &\iH(\bfk \bs_\ell Q),}
	\end{equation}
	where $\TT_\ell$ is defined in Lemma \ref{lem:Ti}.
	%Then we have the following commutative diagram
	
	Dually, we can define $\bs^{-}_\ell Q$ for any source $\ell$ by reversing all the arrows of $Q$ starting from $\ell $, and a reflection functor $F_\ell^-$ associated to a source $\ell\in Q_0$ is also defined.
	This functor $F_\ell^-$ also induces an isomorphism $\Gamma_\ell^-: \iH(\bfk Q)  \stackrel{\sim}{\rightarrow} \iH(\bfk \bs_\ell^- Q)$.
	In particular, the action of $F_\ell ^+$ (and $F_\ell^-$) on $\rep_\bfk( Q)\subseteq \cc_1(\rep_\bfk( Q))$ is the same as that of the classic BGP reflection functor.
	
	A sequence $i_1,\dots,i_t$ in $Q_0$ is called {\em $\imath$-admissible} if $i_r$ is a sink in $\bs_{i_{r-1}}\cdots \bs_{i_1}(Q)$, for each $1\leq r\leq t$.

	%\subsection{The proof of Proposition \ref{prop:DrGenA} }
	
	\begin{proof}[Proof of Proposition \ref{prop:DrGenA}]
		Fix $1\leq j\leq n-1$. We assume $Q(j)$ is the following quiver: %, and let ${}^j\Lambda^\imath$ be its split $\imath$quiver algebra.
		\begin{center}\setlength{\unitlength}{0.5mm}
			%\vspace{-2cm}
			\begin{equation}
				\label{quiver:affineA}
				\begin{picture}(100,30)
					\put(2,0){\circle*{2}}
					\put(51,0){\circle*{2}}
					
					%\put(82,0){\circle*{2}}
					\put(102,0){\circle*{2}}
					\put(52,25){\circle*{2}}
					
					\put(3.5,0){\vector(1,0){16}}
					\put(21,-2){$\cdots$}
					
					\put(32,0){\vector(1,0){16}}
					
					\put(71.5,-2){$\cdots$}
					\put(70,0){\vector(-1,0){16}}
					\put(100,0){\vector(-1,0){16}}
					
					\put(54,24.5){\vector(2,-1){47}}
					\put(50,24.5){\vector(-2,-1){47}}
					
					\put(50.5,27){\tiny $0$}
					\put(1,-6){\tiny $1$}
					%\put(21,-6){\tiny $2$}
					%\put(41,-6){\tiny $3$}
					\put(50,-6){\tiny $j$}
					\put(95,-6){\tiny $n-1$}
				\end{picture}
			\end{equation}
			\vspace{-0.2cm}
		\end{center}

		Recall that
		\begin{align}
			\label{eq:wj}
			\omega_j=\sigma^j (s_{n-j}\cdots s_1)(s_{n-j+1}\cdots s_2) \cdots (s_{n-1}\cdots s_j),
			%,\quad \omega'_i:= \omega_is_i^{-1}.
		\end{align}
		where $\sigma$ is the automorphism of the underlying graph of $C_n$ such that $\sigma(i)=i+1$; see e.g. \cite[\S6.2]{DJX12}.
		Note that $\bs_{\omega_j}(Q(j))=Q(j)$, and $j,\dots, n-1,\dots, 2,\dots, n-j+1,1,\dots, n-j$ is an $\imath$-admissible sequence. Then \cite[Corollary 3.10, Proposition 4.4]{LW19b} show that
		\begin{align}
			\label{eq:M-de}
			\Gamma_{\omega_j}([S_j])=[M(\de-\alpha_j)]*[K_{\de-\alpha_j}]^{-1},
		\end{align}
		where $M(\de-\alpha_j)$ is the unique (up to isomorphism) $\bfk Q(j)$-module with $\de-\alpha_j$ as its class in $K_0(\bfk Q(j))$.
		
		Let $N$ be the unique (up to isomorphism) $\bfk Q(j)$-module with $\widehat{S_0}+\widehat{S_{n-1}}+\cdots + \widehat{S_{j+1}}$ as its class in $K_0(\bfk Q(j))$.
		A direct calculation in $\iH(\bfk Q(j))$ shows that
		\begin{align*}
			[M(\de-\alpha_j)]=\frac{(-1)^{j-1}}{(q-1)^{j-1}}\big[[S_{j-1}],[S_{j-2}],\cdots,[S_{1}],[N]\big]_{\sqq}.
		\end{align*}
		Then the commutative diagram \eqref{eq:comm} yields that
		\begin{align}
			\label{eq:QJBJ}
			\widetilde{\psi}_{Q(j)}((-1)^{j}\TT_{\omega_j}(B_j))
			=&(-1)^{j}\cdot\frac{-1}{q-1}\Gamma_{\omega_j}([S_j])
			\\\notag
			=&\frac{1}{(q-1)^{j}}\big[[S_{j-1}],[S_{j-2}],\cdots,[S_{1}],[N]\big]_{\sqq}*[K_{\de-\alpha_j}]^{-1}.
		\end{align}
		%Recall $\Lambda^\imath$ is the $\imath$quiver algebra of the cyclic quiver $C_n$. 
		It is important to note that $N$ can be viewed as a $\bfk C_n$-module.
		Together with \cite[Lemma 6.3]{LW19b}, $\widetilde{\psi}_{C_n}((-1)^{j}\TT_{\omega_j}(B_j))$ has the same expression as in \eqref{eq:QJBJ}.
		Therefore, by \eqref{Bik}, a direct computation in $\iH(\bfk C_n)$ shows that
		\begin{align*}
			\widehat{B}_{j,-1}=&(1-q)\widetilde{\psi}_{C_n}((-1)^{j}\TT_{\omega_j}(B_j))
			\\
			=&\frac{-1}{(q-1)^{j-1}}\big[[S_{j-1}],[S_{j-2}],\cdots,[S_{1}],[N]\big]_{\sqq}*[K_{\de-\alpha_j}]^{-1}
			\\
			=&\sqq^{-j+1}\sum\limits_{\cm_{j+1,\delta-e_{j}}}(-1)^{\dim\End(M)}[M]*[K_{\de-\alpha_j}]^{-1}.
		\end{align*}
	\end{proof}
	
	%\begin{remark}
	%From the proof of Proposition \ref{prop:DrGenA}, we get that
	%\begin{align*}
	%\Gamma_{\omega_j}^k([S_j])=[M(k\de-\alpha_j)]*[K_{k\de-\alpha_j}]^{-1}
	%\end{align*}
	%for any $k>0$ in $\iH(\bfk Q(j))$.
	%\end{remark}
	
	\begin{corollary}
		\label{cor:DrGenA}
		For any $1\leq j\leq n-1$, we have
		$$\widehat{H}_{j,1}=\widehat{\Theta}_{j,1}=\pi_{j+1,1}-(\sqq+\sqq^{-1})\pi_{j,1}+\pi_{j-1,1}.$$
	\end{corollary}
	
	\begin{proof}%[Proof of Corollary \ref{cor:DrGenA}]
		
		Let
		\begin{align*}
			\cm_{j,\de}^j:=&\{[M]\in \cm_{j,\de}\mid S_j\text{ is a direct summand of }{\rm soc}(M)\},
			\\
			\cn_{j+1,\de}^j:=&\{[M]\in \cm_{j+1,\de}\mid S_j\text{ is not a direct summand of }{\rm soc}(M)\}.
		\end{align*}
		Then we have the following disjoint unions
		\begin{align}
			\label{eq:MMN}
			\cm_{j,\de}=&\cm_{j,\de}^j\sqcup \cm_{j-1,\de},\\
			\label{eq:MMNM}
			\cm_{j+1,\de}^{j+1}=& \cn_{j+1,\de}^j\sqcup  \{[M\oplus S_j]\mid [M]\in \cm_{j+1,\de-e_j}\}.
		\end{align}
		
		Recall that $\Phi(\Theta_{j,1})=(-B_{j}B_{j,-1}+v^2B_{j,-1}B_j)\K_{\de-\alpha_j}$ in $\tUi$.
		By applying $\Omega_{C_n}$ and using Proposition \ref{prop:DrGenA}, we have
		\begin{align*}
			\widehat{\Theta}_{j,1}
			=&\frac{1}{(q-1)^2}\big(-\haB_{j,0}*\haB_{j,-1}+\sqq^2\haB_{j,-1}*\haB_{j,0}\big)*[K_{\de-\alpha_j}]
			\\
			=&\frac{-\sqq^{-j+1}}{(q-1)^2}\sum\limits_{[M]\in\cm_{j+1,\delta-e_{j}}}(-1)^{\dim\End(M)}\Big[[S_j],[M]\Big]_{\sqq^2}
			\\
			=&\frac{\sqq^{-j}}{(q-1)^2}\sum_{[M]\in\cm_{j+1,\de-e_j}}(-1)^{\dim\End(M\oplus S_j)}[M\oplus S_j]-\frac{\sqq^{-j}}{q-1}\sum_{[N]\in\cn_{j+1,\de}^j}(-1)^{\dim\End(N)}[N]
			\\
			&-\frac{\sqq^{-j+2}}{(q-1)^2}\sum_{[M]\in\cm_{j+1,\de-e_j}}(-1)^{\dim\End(M\oplus S_j)}[M\oplus S_j]+\frac{\sqq^{-j+2}}{q-1}\sum_{[L]\in\cm_{j,\de}^j}(-1)^{\dim\End(L)}[L]
			\\
			=&-\frac{\sqq^{-j-1}}{\sqq-\sqq^{-1}}\sum_{[M]\in\cm_{j+1,\de-e_j}}(-1)^{\dim\End(M\oplus S_j)}[M\oplus S_j]
			\\
			&-\frac{\sqq^{-j-1}}{\sqq-\sqq^{-1}}\sum_{[N]\in\cn_{j+1,\de}^j}(-1)^{\dim\End(N)}[N]
			%-\frac{\sqq^{-j+2}}{(q-1)^2}\sum_{[M]\in\cm_{j+1,\de-e_j}}(-1)^{\dim\End(M\oplus S_j)-1}[M\oplus S_j]
			+\frac{\sqq^{-j+1}}{\sqq-\sqq^{-1}}\sum_{[L]\in\cm_{j,\de}^j}(-1)^{\dim\End(L)}[L]
			\\
			\stackrel{\eqref{eq:MMNM}}{=}&-\frac{\sqq^{-j-1}}{\sqq-\sqq^{-1}}\sum_{[M]\in\cm_{j+1,\de}^{j+1}}(-1)^{\dim\End(M)}[M]
			%-\frac{\sqq^{-j+2}}{(q-1)^2}\sum_{[M]\in\cm_{j+1,\de-e_j}}(-1)^{\dim\End(M\oplus S_j)-1}[M\oplus S_j]
			+\frac{\sqq^{-j+1}}{\sqq-\sqq^{-1}}\sum_{[L]\in\cm_{j,\de}^j}(-1)^{\dim\End(L)}[L]
			\\
			\stackrel{\eqref{eq:MMN}}{=}& -\frac{\sqq^{-j-1}}{\sqq-\sqq^{-1}}\big(\sum_{[M]\in\cm_{j+1,\de}}(-1)^{\dim\End(M)}[M]-\sum_{[M]\in\cm_{j,\de}}(-1)^{\dim\End(M)}[M]\big)
			%-\frac{\sqq^{-j+2}}{(q-1)^2}\sum_{[M]\in\cm_{j+1,\de-e_j}}(-1)^{\dim\End(M\oplus S_j)-1}[M\oplus S_j]
			\\
			&+\frac{\sqq^{-j+1}}{\sqq-\sqq^{-1}}\big(\sum_{[L]\in\cm_{j,\de}}(-1)^{\dim\End(L)}[L]-\sum_{[L]\in\cm_{j-1,\de}}(-1)^{\dim\End(L)}[L]\big)
			\\
			=&\pi_{j+1,1}-(\sqq+\sqq^{-1})\pi_{j,1} +\pi_{j-1,1}.
		\end{align*}
		%where the fifth equality follows by \eqref{eq:MMN}, and the last one follows by \eqref{eq:MMNM}.
		%$$\pi_{j,1}=\frac{-\sqq^{-j}}{\sqq-\sqq^{-1}}\sum\limits_{[M]\in\cm_{j,\delta}}(-1)^{\dim\End(M)}[M].$$
	\end{proof}

	\begin{corollary}
		\label{cor for level one of i1}
		We have
		\begin{align}
			\label{eq:Theta11}
			&\widehat{\Theta}_{1,1}=\frac{1}{q-1}\Big(\sqq^{-1}[S_{1}^{(n)}]-\sqq^{-1}[S_{0}^{(n-1)}\oplus S_{1}]-\sqq[S_{0}^{(n)}]\Big),
			\\
			\label{eq:B1-1}
			&\haB_{1,-1}=-[S_{0}^{(n-1)}]*[K_{\alpha_1-\de}],
			\\\label{eq:B11}
			&\haB_{1,1}=\frac{1}{q}[S_1^{(n+1)}] -\frac{1}{q}[S_1\oplus S_{0}^{(n)}],
			\\\label{eq:B1-2}
			&\haB_{1,-2}=(-\frac{1}{q}[S_0^{(2n-1)}] +\frac{1}{q}[S_{0}^{(n)}\oplus S_{0}^{(n-1)}])\ast[K_{\alpha_1-2\de}].
		\end{align}
	\end{corollary}

	\begin{proof}
		\eqref{eq:Theta11}--\eqref{eq:B1-1} follow from Corollary \ref{cor:DrGenA} and Proposition \ref{prop:DrGenA}. %, we have
		%$$\haB_{1,-1}=-[S_{0}^{(n-1)}]*[K_{S_{0}^{(n-1)}}]^{-1}$$
		%and $$\widehat{\Theta}_{1,1}=\pi_{2,1}-(\sqq+\sqq^{-1})\pi_{1,1}=\frac{1}{q-1}\Big(\sqq^{-1}[S_{1}^{(n)}]-\sqq^{-1}[S_{0}^{(n-1)}\oplus S_{1}]-\sqq[S_{0}^{(n)}]\Big).$$
		
		For \eqref{eq:B11}, a direct computation shows that
		\begin{align*}
			&\big[[S_{0}^{(n-1)}\oplus S_{1}],[S_1]\big]=\frac{q-1}{q}[S_{0}^{(n)}\oplus S_1]-\frac{q-1}{q}[S_{1}^{(n)}\oplus S_1];\\\notag
			&\big[[S_{1}^{(n)}],[S_1]\big]=\frac{1-q}{q}[S_{1}^{(n)}\oplus S_1]+\frac{q-1}{q}[S_{1}^{(n+1)}]+\frac{q(q-1)}{q}[S_{0}^{(n-1)}]*[K_{S_1}];
			\\\notag
			&\big[[S_{0}^{(n)}],[S_1]\big]=\frac{q-1}{q}[S_{0}^{(n)}\oplus S_1]-\frac{q-1}{q}[S_{1}^{(n+1)}]-\frac{q(q-1)}{q}[S_{0}^{(n-1)}]*[K_{S_1}].
		\end{align*}
		It follows that
		\begin{align*}
			\notag
			\big[\widehat{\Theta}_{1,1}, \haB_{1,0}\big]&=\frac{1}{q-1}\big[\sqq^{-1}[S_{1}^{(n)}]-\sqq^{-1}[S_{0}^{(n-1)}\oplus S_{1}]-\sqq[S_{0}^{(n)}], [S_1]\big]\\
			&=\frac{[2]_\sqq}{q}[S_1^{(n+1)}] -\frac{[2]_\sqq}{q}[S_1\oplus S_{0}^{(n)}] +[2]_\sqq[S_{0}^{(n-1)}]*[K_{S_1}].
		\end{align*}
		By (\ref{iDR2}) we have
		\begin{align}\notag
			\haB_{1,1}&=\frac{1}{[2]_\sqq}\big[\widehat{\Theta}_{1,1}, \haB_{1,0}\big]+\haB_{1,-1}*[K_\de]\\\notag
			&=\frac{1}{[2]_\sqq}\cdot\Big(\frac{[2]_\sqq}{q}[S_1^{(n+1)}] -\frac{[2]_\sqq}{q}[S_1\oplus S_{0}^{(n)}] +[2]_\sqq[S_{0}^{(n-1)}]*[K_{S_1}]\Big)-[S_{0}^{(n-1)}]*[K_{S_1}]\\\notag
			&=\frac{1}{q}[S_1^{(n+1)}] -\frac{1}{q}[S_1\oplus S_{0}^{(n)}].
		\end{align}
		
		For \eqref{eq:B1-2},
		a direct computation shows that
		\begin{align*}
			&\big[[S_{0}^{(n-1)}\oplus S_{1}],[S_0^{(n-1)}]\big]=\frac{q-1}{q}[S_{1}^{(n)}\oplus S_0^{(n-1)}]-\frac{q-1}{q}[S_{0}^{(n)}\oplus S_0^{(n-1)}];\\\notag
			&\big[[S_{1}^{(n)}],[S_0^{(n-1)}]\big]=\frac{q-1}{q}[S_{1}^{(n)}\oplus S_0^{(n-1)}]-\frac{q-1}{q}[S_{0}^{(2n-1)}]-\frac{q(q-1)}{q}[S_{1}]*[K_{S_0^{(n-1)}}];
			\\\notag
			&\big[[S_{0}^{(n)}],[S_0^{(n-1)}]\big]=\frac{1-q}{q}[S_{0}^{(n)}\oplus S_0^{(n-1)}]+\frac{q-1}{q}[S_{0}^{(2n-1)}]+\frac{q(q-1)}{q}[S_1]*[K_{S_0^{(n-1)}}].
		\end{align*}
		It follows that
		\begin{align}\notag
			\big[\widehat{\Theta}_{1,1}, \haB_{1,-1}\big]%&=\frac{1}{q-1}\Big[\sqq^{-1}[S_{1}^{(n)}]-\sqq^{-1}[S_{0}^{(n-1)}\oplus S_{1}]-\sqq[S_{0}^{(n)}], [S_1]\big]\\
			&=\frac{[2]_\sqq}{q}[S_0^{(2n-1)}]*[K_{\alpha_1-\de}] -\frac{[2]_\sqq}{q}[S_0^{(n)}\oplus S_{0}^{(n-1)}]*[K_{\alpha_1-\de}]  +[2]_\sqq[S_1].
		\end{align}
		Then by (\ref{iDR2}) we have
		\begin{align}\notag
			\haB_{1,-2}&=-\frac{1}{[2]_\sqq}\big[\widehat{\Theta}_{1,1}, \haB_{1,-1}\big]*[K_{-\de}]+\haB_{1,0}*[K_{-\de}]\\\notag
			%&=\frac{1}{[2]_\sqq}\cdot\Big(\frac{[2]_\sqq}{q}[S_1^{(n+1)}] -\frac{[2]_\sqq}{q}[S_1\oplus S_{0}^{(n)}] +[2]_\sqq[S_{0}^{(n-1)}]*[K_{S_1}]\Big)-[S_{0}^{(n-1)}]*[K_{S_1}]\\\notag
			&=\big(-\frac{1}{q}[S_0^{(2n-1)}] +\frac{1}{q}[S_{0}^{(n)}\oplus S_{0}^{(n-1)}]\big)\ast[K_{\alpha_1-2\de}].
		\end{align}
	\end{proof}
	
	%The proofs of Corollaries \ref{cor:DrGenA}--\ref{cor for level one of i1} are in Appendix~\ref{sec:proofscylic}.

	The following lemma will be used in the proof of Lemma \ref{lem:HBstar11-1}. It is understood that the terms involving $S^{(r)}$ vanish if $r<0$.
	\begin{lemma}
		\label{theta and B 4 term in a tube}
		For any $r\geq 0$, we have
		\begin{align}
			&\big[[S_{0}^{(rn)}], \widehat{B}_{1,-1}\big]+q\big[[S_{0}^{((r-2)n)}],\widehat{B}_{1,-1}\big]*[K_{\delta}]
			\\\notag
			&=
			\big[[S_{0}^{((r-1)n)}],\widehat{B}_{1,0}\big]_{\sqq^{2}}
			+q\big[[S_{0}^{((r-1)n)}],\widehat{B}_{1,-2}\big]_{\sqq^{-2}}*[K_{\delta}].
		\end{align}
	\end{lemma}

	\begin{proof}%[Proof of Lemma \ref{theta and B 4 term in a tube}]
		By Corollary \ref{cor for level one of i1}, we only need to check that
		\begin{align}\label{equiv formula}
			&\big[[S_{0}^{(rn)}],[S_{0}^{(n-1)}]\big]
			+q\big[[S_{0}^{((r-2)n)}],[S_{0}^{(n-1)}]\big]*[K_{\delta}]
			\\\notag
			&=
			q\big[[S_{1}],[S_{0}^{((r-1)n)}]\big]_{\sqq^{-2}}*[K_{\delta-\alpha_{11}}]
			+\big[[S_{0}^{((r-1)n)}], [S_{0}^{(2n-1)}]-
			[S_{0}^{(n)}\oplus S_{0}^{(n-1)}]
			\big]_{\sqq^{-2}}.
		\end{align}
		
		For any $r\geq 0$, an easy computation shows that
		\begin{align*}%\label{S0 n-1 and rn}
			\begin{split}
				[S_{0}^{(n-1)}]*[S_{0}^{(rn)}]
				=&[S_{0}^{(rn)}\oplus S_{0}^{(n-1)}],\\
				[S_{0}^{(rn)}]*[S_{0}^{(n-1)}]=&\frac{1}{q}[S_{0}^{(rn)}\oplus S_{0}^{(n-1)}]+\frac{q-1}{q}[S_{0}^{((r+1)n-1)}]+(q-1)[S_{1}^{((r-1)n+1)}]*[K_{\de-\alpha_1}];
				%=&\frac{1}{q}\Big([S_{0}^{(rn)}\oplus S_{0}^{(n-1)}]+(q-1)[S_{0}^{((r+1)n-1)}]+q(q-1)[S_{1}^{((r-1)n+1)}]*[K_{\de-\alpha_1}]\Big).
			\end{split}
			\\
			%\label{S0 rn and S1}
			\begin{split}
				[S_{0}^{(rn)}]*[S_{1}]
				=&[S_{0}^{(rn)}\oplus S_{1}],\\
				[S_{1}]*[S_{0}^{(rn)}]=&\frac{1}{q}[S_{0}^{(rn)}\oplus S_{1}]+\frac{q-1}{q}[S_{1}^{(rn+1)}]+(q-1)[S_{0}^{(rn-1)}]*[K_{\alpha_{1}}].
		\end{split}\end{align*}
		Then we obtain
		%Hence
		\begin{align*}%\label{S1and S0^rn}
			\big[[S_{1}], [S_{0}^{(rn)}]\big]_{\sqq^{-2}}
			=\frac{q-1}{q}[S_{1}^{(rn+1)}]+(q-1)[S_{0}^{(rn-1)}]*[K_{\alpha_{1}}],
		\end{align*}
		and \begin{align*}%\label{S0^{n-1}and S0^rn}
			\big[[S_{0}^{(rn)}],[S_{0}^{(n-1)}]\big]
			=\frac{q-1}{q}\Big([S_{0}^{((r+1)n-1)}]-[S_{0}^{(rn)}\oplus S_{0}^{(n-1)}]+q[S_{1}^{((r-1)n+1)}]*[K_{\de-\alpha_1}]\Big).
		\end{align*}
		%Similarly, we can obtain
		On the other hand, for $r\geq 1$, by Proposition \ref{prop:iHallmult}, we have
		% (for $\blue{r\geq 1}$, \red{not only for $r\geq 2$})
		\begin{align*}
			&[S_{0}^{(2n-1)}]*[S_{0}^{(rn)}]\\
			&=\frac{1}{q}[S_{0}^{(rn)}\oplus S_{0}^{(2n-1)}]+\frac{q-1}{q}[S_{0}^{((r+1)n)}\oplus S_{0}^{(n-1)}]
			+\frac{q(q-1)}{q}[S_{0}^{((r-1)n)}\oplus S_{0}^{(n-1)}]*[K_\de],
		\end{align*}
		\begin{align*}
			&[S_{0}^{(rn)}]*[S_{0}^{(2n-1)}]\\
			&=\frac{1}{q^2}[S_{0}^{(rn)}\oplus S_{0}^{(2n-1)}]+\frac{q-1}{q^2}[S_{0}^{((r+1)n-1)}\oplus S_{0}^{(n)}]+\frac{q(q-1)}{q^2}[S_{0}^{((r+2)n-1)}]\\
			&\quad +\frac{q-1}{q}[S_{1}^{((r-1)n+1)}\oplus S_{0}^{(n)}]*[K_{S_{0}^{(n-1)}}]+
			\frac{(q-1)^2}{q}[S_{1}^{(rn+1)}]*[K_{S_{0}^{(n-1)}}]\\
			&\quad +q(q-1)[S_{1}^{((r-2)n+1)}]*[K_{S_{0}^{(2n-1)}}];
		\end{align*}
		and
		\begin{align*}
			&[S_{0}^{(n)}\oplus S_{0}^{(n-1)}]*[S_{0}^{(rn)}]\\
			&=\frac{1}{q}[S_{0}^{(n)}\oplus S_{0}^{(n-1)}\oplus S_{0}^{(rn)}]
			+\frac{q-1}{q}[S_{0}^{((r+1)n)}\oplus S_{0}^{(n-1)}]
			+\frac{q(q-1)}{q}[S_{0}^{((r-1)n)}\oplus S_{0}^{(n-1)}]*[K_\de],
		\end{align*}
		\begin{align*}
			&[S_{0}^{(rn)}]*[S_{0}^{(n)}\oplus S_{0}^{(n-1)}]\\
			&=\frac{1}{q^2}[S_{0}^{(rn)}\oplus S_{0}^{(n)}\oplus S_{0}^{(n-1)}]+\frac{q-1}{q^2}[S_{0}^{((r+1)n-1)}\oplus S_{0}^{(n)}]
			+\frac{q(q-1)}{q^2}[S_{0}^{((r+1)n)}\oplus S_{0}^{(n-1)}]
			\\
			&\quad +\frac{q-1}{q}[S_{1}^{((r-1)n+1)}\oplus S_{0}^{(n)}]*[K_{S_{0}^{(n-1)}}]
			+
			\frac{(q-1)^2}{q}[S_{1}^{(rn+1)}]*[K_{S_{0}^{(n-1)}}]
			\\
			&\quad +\frac{q(q-1)}{q}[S_{0}^{((r-1)n)}\oplus S_{0}^{(n-1)}]*[K_\de] +\frac{q(q-1)^2}{q}[S_{0}^{(rn-1)}]*[K_\de].
		\end{align*}
		By combining the above formulas, we obtain
		\begin{align*}%\label{S0^{2n-1}and S0^rn}
			&\big[[S_{0}^{(rn)}],[S_{0}^{(2n-1)}]\big]_{\sqq^{-2}}-\big[[S_{0}^{(rn)}],[S_{0}^{(n)}\oplus S_{0}^{(n-1)}]\big]_{\sqq^{-2}}
			\\\notag
			&=\frac{q-1}{q}[S_{0}^{((r+2)n-1)}]-
			\frac{q-1}{q}[S_{0}^{((r+1)n)}\oplus S_{0}^{(n-1)}]+q(q-1)[S_{1}^{((r-2)n+1)}]*[K_{S_{0}^{(2n-1)}}]
			\\\notag
			&\quad-(q-1)[S_{0}^{((r-1)n)}\oplus S_{0}^{(n-1)}]*[K_\de]-\frac{q(q-1)^2}{q}[S_{0}^{(rn-1)}]*[K_\de].
		\end{align*}
		In particular, for $r=0$, we have
		
		\begin{align*}
			&\Big[[S_{0}^{(rn)}],[S_{0}^{(2n-1)}]-[S_{0}^{(n)}\oplus S_{0}^{(n-1)}]\Big]_{\sqq^{-2}}=(1-\sqq^{-2})\Big([S_{0}^{(2n-1)}]-[S_{0}^{(n)}\oplus S_{0}^{(n-1)}]\Big)
		\end{align*}
		since $S_{0}^{(0)}=0$.
		Then (\ref{equiv formula}) follows since both sides are equal to
		\begin{align*}
			&\frac{q-1}{q}\Big([S_{0}^{((r+1)n-1)}]-[S_{0}^{(rn)}\oplus S_{0}^{(n-1)}]+q[S_{1}^{((r-1)n+1)}]*[K_{S_{0}^{(n-1)}}]\Big)\\
			&+(q-1)\Big([S_{0}^{((r-1)n-1)}]-[S_{0}^{((r-2)n)}\oplus S_{0}^{(n-1)}]+q[S_{1}^{((r-3)n+1)}]*[K_{S_{0}^{(n-1)}}]\Big)*[K_{\delta}].
		\end{align*}
	\end{proof}

	\subsection{Another description of root vectors}
	It is helpful and interesting to compute all root vectors %$\haB_{i,l},\haTh_{i,k}$ ($l\in\Z,k>0$, $1\leq i\leq \bt-1$)
	clearly in $\iH(\bfk C_n)$. In particular, the root vectors at the vertex $1$ are described in Subsection \ref{sec:rooti1}. There is another way to interpret the root vectors in $\imath$Hall algebras, which we give in the following.
	
	Let $Q$ be a quiver with its underlying graph the same as $C_n$. Then there exists an isomorphism $\texttt{FT}_{C_n,Q}:\tCMH\stackrel{\cong}{\longrightarrow} \tCMHC$ %, called a {\em Frourier transformation},
	such that $\texttt{FT}_{C_n,Q}([S_i])=[S_i]$, $\texttt{FT}_{C_n,Q}([K_{S_i}])=[K_{S_i}]$ for any $i\in\I$; cf. \eqref{eq:isocomp}. For any positive real root $\beta$ of $\widehat{\mathfrak{sl}}_n$, we denote by $M(\beta)$ the unique (up to isomorphism) $\bfk Q$-module with its class $\beta$ in $K_0(\mod(\bfk Q))$.

	For any hereditary abelian category $\ca$ and a full subcategory $\cs$, recall from \cite{GL87} that the right (resp. left) perpendicular category $\cs{}^\bot$ (resp. ${}^\bot \cs$) is defined as follows:
	\begin{align}
		\label{def:perpcat right}
		\cs{}^\bot:=\{X\in\ca\mid \Hom_\ca(C,X)=0=\Ext_\ca^1(C,X)\text{ for any }C\in\cs\};\\
		\label{def:perpcat}
		{}^\bot\cs:=\{X\in\ca\mid \Hom_\ca(X,C)=0=\Ext_\ca^1(X,C)\text{ for any }C\in\cs\}.
	\end{align}
	Then the category $\cs{}^\bot$ (resp. ${}^\bot \cs$) is a hereditary abelian category, which is extension-closed in $\ca$; see \cite[Proposition 3.2]{BS13}.
	If $\cs=\{M\}$ for an object $M\in\ca$, then we also denote $M{}^\bot:=\cs{}^\bot$ and ${}^\bot M:={}^\bot\cs$.

	Recall $Q(j)$ in \eqref{quiver:affineA}.
	\begin{proposition}
		\label{rem:rootvectors}
		For any $l>0$, we have
		\begin{align}
			\label{eq:haB-}
			\haB_{j,-l}=& (-1)^{jl}\texttt{FT}_{C_n,Q(j)}([M(l\de-\alpha_j)])*[K_{l\de-\alpha_j}]^{-1},
			\\
			\label{eq:haB+}
			\haB_{j,l}= &(-1)^{jl}\texttt{FT}_{C_n,Q(j)}([M(l\de+\alpha_j)]),
			\\
			\label{eq:haTh+}
			\widehat{\Theta}_{j,l}=&\frac{(-1)^{jl}}{(q-1)^2\sqq^{l-1}}\texttt{FT}_{C_n,Q(j)}\Big(\sum_{0\neq f:S_j\rightarrow M(l\de+\alpha_j) } [\coker(f)]\Big).
		\end{align}
	\end{proposition}
	
	\begin{proof}
		From the proof of Proposition \ref{prop:DrGenA}, in $\iH(\bfk Q(j))$ we have
		\begin{align*}
			\Gamma_{\omega_j}^l([S_j])=[M(l\de-\alpha_j)]*[K_{l\de-\alpha_j}]^{-1}.
		\end{align*}
		Then $B_{j,-l}$ satisfies that $\Phi(B_{j,-l})= (-1)^{jl} \TT_{\omega_j}^{l} (B_j)$ and so
		$$\widetilde{\psi}_{Q(j)}\circ\Phi(B_{j,-l})=(-1)^{jl}\cdot\frac{-1}{q-1}[M(l\de-\alpha_j)]*[K_{l\de-\alpha_j}]^{-1} $$
		in $\iH(\bfk Q(j))$. Similarly, one can see that
		$$\widetilde{\psi}_{Q(j)}\circ\Phi(B_{j,l})=(-1)^{jl}\cdot\frac{-1}{q-1}[M(l\de+\alpha_j)];$$
		cf. \cite[Lemma 5.5]{LRW20a}.
		So we have \eqref{eq:haB-}--\eqref{eq:haB+}.
		%\begin{align}
		%\label{eq:haB-}
		%\haB_{j,-l}=& (-1)^{jl}\texttt{FT}_{C_n,Q(j)}([M(l\de-\alpha_j)])*[K_{l\de-\alpha_j}]^{-1},
		%\\
		%\label{eq:haB+}
		% \haB_{j,l}= &(-1)^{jl}\texttt{FT}_{C_n,Q(j)}([M(l\de+\alpha_j)]).
		%\end{align}
		%for any $l>0$.
		
		For \eqref{eq:haTh+}, let $\cs:=\{S_i\in\mod(\bfk Q(j))\mid i\neq 0,j\}$. Then the subcategory ${}^\bot\cs$ is equivalent to $\mod(\bfk Q_{\texttt{Kr}})$, where $Q_{\texttt{Kr}}$ is the Kronecker quiver $\xymatrix{0\ar@<0.5ex>[r] \ar@<-0.5ex>[r] & 1}$. So there exists an embedding functor $\F: \mod(\bfk Q_{\texttt{Kr}})\simeq {}^\bot\cs\rightarrow \mod(\bfk Q(j))$. In particular, $\F(S_1)=S_j$, and $\F(M(l\alpha_0+(l+1)\alpha_1))=M(l\de+\alpha_j)$ for any $l\geq0$.
		So $\F$ induces an embedding of algebras $\widetilde{\ch}(\cc_1(\mod(\bfk Q_{\texttt{Kr}})))\rightarrow \widetilde{\ch}(\cc_1(\mod(\bfk Q(j))))$, and then an embedding of algebras $F:\iH(\bfk Q_{\texttt{Kr}})\rightarrow \iH(\bfk Q (j))$.
		Using \cite[Theorem 5.11]{LRW20a}, in order to compute $\haTh_{j,l}$ in $\iH(\bfk Q_{\texttt{Kr}})$, it is equivalent to compute them in $\iH(\P^1_\bfk)$. By  \cite[(4.2)]{LRW20a}, we have
		$$\haTh_{1,l}=\frac{1}{(q-1)^2\sqq^{l-1}}\sum_{0\neq f:S_1\rightarrow M(l\alpha_0+(l+1)\alpha_1) } [\coker(f)]$$
		in $\iH(\bfk Q_{\texttt{Kr}})$. By applying $F$, we have
		\begin{align}
			\label{haTh}
			\haTh_{j,l}=\frac{(-1)^{jl}}{(q-1)^2\sqq^{l-1}}\sum_{0\neq f:S_j\rightarrow M(l\de+\alpha_j) } [\coker(f)]
		\end{align}
		in $\iH(\bfk Q (j))$. Here the scalar $(-1)^{jl}$ comes from \eqref{eq:haB-}--\eqref{eq:haB+}. Then \eqref{eq:haTh+} follows.
	\end{proof}

	\subsection{Two useful formulas}
	\label{subsec:twoformu}
	
	Return to the cyclic quiver $C_n$. Let $\mathcal{S}=\{S_{j}\mid 2\leq j\leq n-1\}$. Observe that any indecomposable object in ${}^{\perp}\mathcal{S}$ has the form $S_{0}^{(an)}$, $S_{1}^{(an)}$, $S_{0}^{(an-1)}$ or $S_{1}^{(bp+1)}$, for some  $a\geq 1$ and $b\geq 0$. In particular, ${}^{\perp}\mathcal{S}$ is a rank-two tube with two simples $S_1$ and $S_{0}^{(n-1)}$.
	
	For any nilpotent $\bfk C_n$-module $N$, let $\ell(N)$ be the number of its indecomposable direct summands. In particular, $\ell(N)=1$ if and only if $N$ is indecomposable.
	
	The following two lemmas will be used to prove Lemma \ref{co times Theta 11}; see Appendix \ref{subsec:proofOtimesTH11}.

	\begin{lemma}
		\label{sum of l(N)=l(M)+1 for p-1}
		Assume $a\geq1$ and $M\in {}^{\perp}\mathcal{S}$. We have
		\begin{align}\label{Ext with fix middle length 1}
			\sum_{[N];\ell(N)=\ell(M)+1}|\Ext^1(M, S_{0}^{(an-1)})_{N}|=\frac{|\Hom(M, S_{0}^{(an-1)})|}{|\Hom(M, S_{0}^{(n-1)})|}.
		\end{align}
	\end{lemma}
	
	\begin{proof}
		We prove \eqref{Ext with fix middle length 1} by induction on $a$.
		If $a=1$, then by $M\in{}^{\perp}\mathcal{S}$ and $S_{0}^{(n-1)}$ is simple in ${}^{\perp}\mathcal{S}$, we see that any extension $0\rightarrow S_{0}^{(n-1)}\rightarrow N \rightarrow M \rightarrow 0$ with $\ell(N)=\ell(M)+1$ must be split. Then the result holds by noting
		\begin{align*}
			\sum_{[N];\ell(N)=\ell(M)+1}|\Ext^1(M,S_{0}^{(n-1)})_{N}|=1.
		\end{align*}
		
		For $a>1$, we fix a short exact sequence
		$0\rightarrow S_{1}^{(n)}\xrightarrow{f_1} S_{0}^{(an-1)} \xrightarrow{f_2} S_{0}^{((a-1)n-1)}\rightarrow0$ in $\rep^{\rm nil}_\bfk(C_n)$.
		By applying $\Hom(M,-)$, we have the following long exact sequence
		\begin{align}
			\label{eqn:long exact of S_1^p}
			0\longrightarrow &\Hom(M,S_{1}^{(n)})\longrightarrow \Hom(M,S_{0}^{(an-1)}) \longrightarrow \Hom(M,S_{0}^{((a-1)n-1)})\longrightarrow
			\\
			\notag &\Ext^1(M,S_{1}^{(n)}) \longrightarrow \Ext^1(M,S_{0}^{(an-1)}) \stackrel{\phi}{\longrightarrow} \Ext^1(M,S_{0}^{((a-1)n-1)})\longrightarrow0,
		\end{align}
		where the action of $\phi$ is induced by $f_2$ via the following pushout diagram
		\[\xymatrix{ S_{1}^{(n)}\ar@{=}[d] \ar[r]^{f_1} & S_{0}^{(an-1)} \ar[r]^{f_2} \ar[d]^{h_1} & S_{0}^{((a-1)n-1)} \ar[d]^{g_1} \\
			S_{1}^{(n)} \ar[r]^{l_1} & N  \ar[r]^{l_2} \ar[d]^{h_2} & L \ar[d]^{g_2} \\
			& M \ar@{=}[r] & M .}\]
		We claim that $\ell(N)=\ell(L)$, or equivalently, $\ell(N)\neq \ell(L)+1$. Suppose for a contradiction that $\ell(N)=\ell(L)+1$. In this case, the short exact sequence in the second row can not be split by observing the one in the first row is not split. Together with the assumption $\ell(N)=\ell(L)+1$, $l_1$ should have a component $S_{1}^{(n)}\twoheadrightarrow X$, for some indecomposable direct summand $X$ of $N$. Recall that ${}^{\perp}\mathcal{S}$ is a rank-two tube with two simples $S_1$ and $S_{0}^{(n-1)}$.
		It follows from $N\in{}^{\perp}\mathcal{S}$ that $X= S_1$. By observing $\Hom(S_{0}^{(an-1)}, S_{1})=0$, the above commutative diagram does not exist, giving a contradiction. So the claim holds.
		
		Therefore,
		\begin{align*}
			&\sum_{[N];\ell(N)=\ell(M)+1}|\Ext^1(M,S_0^{(an-1)})_{N}|\\
			&=\sum_{[L];\ell(L)=\ell(M)+1}|\phi^{-1}(\Ext^1(M,S_0^{((a-1)n-1)})_{L}) |\\
			&=\sum_{[L];\ell(L)=\ell(M)+1}|\Ext^1(M,S_0^{((a-1)n-1)})_{L}|\cdot|\ker(\phi)|.
		\end{align*}
		By \eqref{eqn:long exact of S_1^p} we obtain that
		\begin{align*}
			|\ker(\phi)|=&\frac{|\Ext^1(M,S_{1}^{(n)})|\cdot|\Hom(M,S_{0}^{(an-1)})|}{|\Hom(M,S_{1}^{(n)})|\cdot |\Hom(M,S_{0}^{((a-1)n-1)})|}
			\\
			=&\frac{|\Hom(M,S_{0}^{(an-1)})|}{|\Hom(M,S_{0}^{((a-1)n-1)})|}.
		\end{align*}
		Now by induction, we have
		\begin{align*}
			&\sum_{[N];\ell(N)=\ell(M)+1}|\Ext^1(M,S_0^{(an-1)})_{N}|\\
			&=\frac{|\Hom(M, S_{0}^{((a-1)n-1)})|}{|\Hom(M, S_{0}^{(n-1)})|}\cdot \frac{|\Hom(M,S_{0}^{(an-1)})|}{|\Hom(M,S_{0}^{((a-1)n-1)})|}=\frac{|\Hom(M,S_{0}^{(an-1)})|}{|\Hom(M,S_{0}^{(n-1)})|}.
		\end{align*}
	\end{proof}
	
	\begin{lemma}\label{sum of l(N)=l(M)+1 for p}
		Assume $a\geq1$ and $M\in {}^{\perp}\mathcal{S}$. We have
		\begin{align}
			\sum_{[N];\ell(N)=\ell(M)+1}|\Ext^1(M, S_{0}^{(an)})_{N}|=\frac{|\Hom(M, S_{0}^{(an)})|}{|\Ext^1(M, S_{0}^{(n-1)})|}.
		\end{align}
	\end{lemma}
	
	\begin{proof}
		We fix a short exact sequence
		$0\rightarrow S_{1}\xrightarrow{f_1} S_{0}^{(an)} \xrightarrow{f_2} S_{0}^{(an-1)}\rightarrow0$ in $\rep^{\rm nil}_\bfk(C_n)$.
		By applying $\Hom(M,-)$, we have the following long exact sequence
		\begin{align}
			\label{eqn:long exact of S_1}
			0\longrightarrow &\Hom(M,S_{1})\longrightarrow \Hom(M,S_{0}^{(an)}) \longrightarrow \Hom(M,S_{0}^{(an-1)})\longrightarrow
			\\
			\notag &\Ext^1(M,S_{1}) \longrightarrow \Ext^1(M,S_{0}^{(an)}) \stackrel{\phi}{\longrightarrow} \Ext^1(M,S_{0}^{(an-1)})\longrightarrow0,
		\end{align}
		where the action of $\phi$ is induced by $f_2$ via the following pushout diagram
		\[\xymatrix{ S_{1}\ar@{=}[d] \ar[r]^{f_1} & S_{0}^{(an)} \ar[r]^{f_2} \ar[d]^{h_1} & S_{0}^{(an-1)} \ar[d]^{g_1} \\
			S_{1} \ar[r]^{l_1} & N  \ar[r]^{l_2} \ar[d]^{h_2} & L \ar[d]^{g_2} \\
			& M \ar@{=}[r] & M. }\]
		Since $S_1$ is simple, we obtain that $\ell(N)= \ell(L)+1$ if and only if the short exact sequence in the second row is split, contradicting to the non-splitness of the one in the first row. Hence $\ell(N)\neq \ell(L)+1$, or equivalently, $\ell(N)=\ell(L)$.
		Therefore,
		\begin{align*}
			\sum_{[N];\ell(N)=\ell(M)+1}|\Ext^1(M,S^{(an)})_{N}|=&\sum_{[L]:\ell(L)=\ell(M)+1}|\phi^{-1}(\Ext^1(M,S^{(an-1)})_{L}) |\\
			=&\sum_{[L];\ell(L)=\ell(M)+1}|\Ext^1(M,S^{(an-1)})_{L}|\cdot|\ker(\phi)|.
		\end{align*}
		By \eqref{eqn:long exact of S_1} we obtain that
		\begin{align*}
			|\ker(\phi)|=&\frac{|\Ext^1(M,S_{1})|\cdot|\Hom(M,S_{0}^{(an)})|}{|\Hom(M,S_{1})|\cdot |\Hom(M,S_{0}^{(an-1)})|}
			\\
			=&q^{-\langle M, S_{1}\rangle}\cdot\frac{|\Hom(M,S_{0}^{(an)})|}{|\Hom(M,S_{0}^{(an-1)})|}.
		\end{align*}
		Now by Lemma \ref{sum of l(N)=l(M)+1 for p-1}, we have
		\begin{align*}
			&\sum_{[N];\ell(N)=\ell(M)+1}|\Ext^1(M,S_0^{(an)})_{N}|\\
			&=\frac{|\Hom(M,S_{0}^{(an-1)})|}{|\Hom(M,S_{0}^{(n-1)})|}\cdot q^{-\langle M, S_{1}\rangle}\cdot\frac{|\Hom(M,S_{0}^{(an)})|}{|\Hom(M,S_{0}^{(an-1)})|}
			\\
			&=q^{-\langle M, S_{1}\rangle}\cdot\frac{|\Hom(M,S_{0}^{(an)})|}{|\Hom(M,S_{0}^{(n-1)})|}
			= \frac{|\Hom(M, S_{0}^{(an)})|}{|\Ext^1(M, S_{0}^{(n-1)})|}.
		\end{align*}
		Here the last equality uses \eqref{eqn:long exact of S_1} with $a=1$.
	\end{proof}

	%%%%%%%%%%%%%%%
	\section{$\imath$Hall algebras and $\imath$quantum loop algebras}
	\label{sec:hom}
	
	In this section, we shall formulate the main result of this paper. 
	
	Let us fix some notations. 
	For a partition $\lambda=(1^{l_1},\dots,n^{l_n})$, we denote
	$|\lambda|:=\sum\limits_{i=1}^n il_i$,
	and $\ell(\lambda):=\sum\limits_{i=1}^n l_i$.

	\subsection{Embedding of $\imath$Hall algebras I}
	\label{subsec:embeddingtube}
	Let $\X$ be a weighted projective line of weight type $(\bp,\ul{\bla})$.
	Recall that $\bla_i$ is the exceptional closed point of $\X$ of weight $p_i$ for any $1\leq i\leq \bt$, and $\scrt_{\bla_i}$ is the Serre subcategory of $\coh(\X)$ consisting of torsion sheaves supported at $\bla_i$. There is an equivalence $\scrt_{\bla_i}\cong\rep^{\rm nil}_\bfk(C_{p_i})$, which induces a canonical embedding of Ringel-Hall algebras
	$\ch(\cc_1(\rep_\bfk^{\rm nil} (C_{p_i})))\rightarrow \ch(\cc_1(\coh(\X_\bfk)))$, and then an embedding of $\imath$Hall algebras:
	\begin{equation}
		\label{eq:embeddingx}
		\iota_i: \iH(\bfk C_{p_i})\longrightarrow\iH(\X_\bfk).
	\end{equation}
	
	Inspired by \eqref{def:haB}, we define
	\begin{align}
		\label{def:haBThH}
		\haB_{[i,j],l}:= \iota_i (\haB_{j,l}),\quad \widehat{\Theta}_{[i,j],r}:=\iota_i(\widehat{\Theta}_{j,r}),\quad \widehat{H}_{[i,j],r}:=\iota_{i}(\widehat{H}_{j,r}),
	\end{align}
	for any $1\leq j\leq p_i-1$, $l\in\Z$ and $r>0$.
	%the images of $\haB_{j,l},\;
	%\widehat{\Theta}_{j,r},\;
	%\widehat{H}_{j,r}$ under $\Omega_{\lambda_i}$ by $\haB_{[i,j],l},\;
	%\widehat{\Theta}_{[i,j],r},\;
	%\widehat{H}_{[i,j],r}$
	%respectively.
	By \eqref{exp h} we have
	\begin{align}
		\label{exp haH}
		1+ \sum_{m\geq 1} (\sqq-\sqq^{-1})\haTh_{[i,j],m} u^m  = \exp\Big( (\sqq-\sqq^{-1}) \sum_{m\geq 1} \widehat{H}_{[i,j],m} u^m \Big).
	\end{align}
	
	For convenience, we set $\haTh_{[i,j],0}= \frac{1}{\sqq-\sqq^{-1}}$ for any $1\leq i\leq \bt,1\le j\le p_i-1$.

	\subsection{Embedding of $\imath$Hall algebras II}%{\red{Embedding of coherent sheaf categories from projective line to weighted projective line}}{\red{Embedding of coherent sheaf categories between weighted projective lines}}
	\label{Embedding from projective line to weighted projective line}
	Let $\mathcal{C}$ be the Serre subcategory of $\coh(\X)$ generated by those simple sheaves $S$ satisfying $\Hom(\co, S)=0$. Recall from \cite{BS13} that the Serre quotient $\coh(\X)/\mathcal{C}$ is equivalent to the category $\coh(\P^1)$ and the canonical functor $\coh(\X)\rightarrow\coh(\X)/\mathcal{C}$ has an exact fully faithful right adjoint functor
	\begin{equation}
		\label{the embedding functor F}
		\mathbb{F}_{\X,\P^1}: \coh(\P^1)\rightarrow\coh(\X),
	\end{equation}
	which sends $$\co_{\P^1}(l)\mapsto\co(l\vec{c}),\quad S_{\bla_i}^{(r)}\mapsto S_{i,0}^{(rp_i)},\quad S_{x}^{(r)}\mapsto S_{x}^{(r)}$$ % and $\haT_{m}\mapsto \haT_{\star,m}$, $\haH_m\mapsto \haH_{\star,m}$
	for any $l\in\Z, r\geq 1$, $1\leq i\leq \bt$ and $x\in\PL\setminus\{\bla_1,\cdots, \bla_\bt\}$. Then $\mathbb{F}_{\X,\P^1}$ induces an exact fully faithful functor $\cc_1(\coh(\P^1))\rightarrow \cc_1(\coh(\X))$, which is also denoted by $\mathbb{F}_{\X,\P^1}$. This functor $\mathbb{F}_{\X,\P^1}$ induces a canonical embedding of Hall algebras $\ch\big(\cc_1(\coh(\P^1))\big)\rightarrow \ch\big(\cc_1(\coh(\X))\big)$, and then an embedding of $\imath$Hall algebras:
	\begin{equation}
		\label{the embedding functor F on algebra}
		F_{\X,\P^1}:\iH(\P^1_\bfk) \longrightarrow \iH(\X_\bfk).
	\end{equation}
	
	Recall from \cite[(4.2), Proposition 6.3]{LRW20a} that for any $m\geq0$, in $\iH(\P^1_\bfk)$
	$$\haT_{m}= \frac{1}{(q-1)^2\sqq^{m-1}}\sum_{0\neq f:\co(s)\rightarrow \co(m+s) } [\coker(f)],
	$$
	which is independent of $s\in\Z$, and
	\begin{align*}
		%\label{formula for Hm}
		\haH_{m}
		=&\sum_{x,d_x|m} \frac{[m]_\sqq}{m} d_x \sum_{|\lambda|=\frac{m}{{d_x}}} \bn_x(\ell(\lambda)-1)\frac{[S_x^{(\lambda)}]}{\big|\Aut(S_x^{(\lambda)})\big|}  -\de_{m,ev}\frac{[m]_\sqq}{m} [K_{\frac{m}{2}\de}],
	\end{align*}
	where $d_x$ is defined in \eqref{def:dx}, and
	\[\bn_x(l)=\prod_{i=1}^l(1-\sqq^{2id_x}), \,\,\forall l\geq1.\]

	Hence, we define
	\begin{align}
		\label{def:Theta star}
		\widehat{\Theta}_{\star,m}:=F_{\X,\P^1}(\haT_m)=\frac{1}{(q-1)^2\sqq^{m-1}}\sum_{0\neq f:\co(s\vec{c})\rightarrow \co((m+s)\vec{c}) } [\coker(f)],
	\end{align}
	and
	\begin{align}
		\label{formula for Hm}
		\haH_{\star,m}:=F_{\X,\P^1}(\haH_{m})
		=&\sum_{x,d_x|m} \frac{[m]_\sqq}{m} d_x \sum_{|\lambda|=\frac{m}{{d_x}}} \bn_x(\ell(\lambda)-1)\frac{F_{\X,\P^1}([S_x^{(\lambda)}])}{\big|\Aut(S_x^{(\lambda)})\big|}  -\de_{m,ev}\frac{[m]_\sqq}{m} [K_{\frac{m}{2}\de}].
	\end{align}
	Here, for any partition $\lambda=(\lambda_1, \lambda_2,\cdots)$, $F_{\X,\P^1}([S_x^{(\lambda)}])=[S_x^{(\lambda)}]$ for $x\in\PL\setminus\{\bla_1,\cdots, \bla_\bt\}$, while $F_{\X,\P^1}([S_{\bla_i}^{(\lambda)}])=S_{i,0}^{(\lambda)}:=\bigoplus_{i}S_{i,0}^{(\lambda_ip_i)}$ for $1\leq i\leq \bt$.
	Note that $\haTh_{\star,m}=0$ if $m<0$ and $\haTh_{\star,0}=\frac{1}{\sqq-\sqq^{-1}}$.

	\begin{lemma}
		\label{lem:Hxm}	
		For any $m\geq 1$ and $x\in\P^1_\bfk$ such that $d_x| m$, we define
		\begin{align}
			\label{def:Hxm}
			\widehat{H}_{x,m}:= \frac{[m]_\sqq}{m} d_x \sum_{|\lambda|=\frac{m}{{d_x}}} \bn_x(\ell(\lambda)-1)\frac{F_{\X,\P^1}([S_x^{(\lambda)}])}{\big|\Aut(S_x^{(\lambda)})\big|}  -\de_{\frac{m}{d_x},ev}\cdot d_x\sqq^{-\frac{m}{2}}\frac{[\frac{m}{2}]_\sqq }{m} [K_{\frac{m}{2}\de}].
		\end{align}
		Then
		\begin{align}
			\widehat{H}_{\star,m}=\sum_{x,d_x|m} \widehat{H}_{x,m}.
		\end{align}
	\end{lemma}
	
	\begin{proof}
		Note that for $m\geq 1$,
		$$\sum_{x\in \P^1_\bfk;\, d_x |m} d_x =q^m+1.$$
		A direct computation shows that
		\begin{align*}
			&\sum_{x, d_x|m}(\sqq_x-\sqq_x^{-1})\delta_{\frac{m}{d_x}, ev}\sqq^{-\frac{m}{2}} \frac{[m/(2d_x)]_{\sqq_x}}{m/d_x} [K_{\frac{m}{2}\de}]\\
			&=\sum_{x, d_x|m} \delta_{\frac{m}{d_x}, ev} \sqq^{-\frac{m}{2}} \frac{\sqq^{\frac{m}{2}}-\sqq^{-\frac{m}{2}}}{m} d_x [K_{\frac{m}{2}\de}]
			= \delta_{m, ev} \frac{1-\sqq^{-m}}{m} \sum_{x, d_x|\frac{m}{2}}d_x [K_{\frac{m}{2}\de}]\\
			&= \delta_{m, ev} \frac{1-\sqq^{-m}}{m} (q^{\frac{m}{2}}+1) [K_{\frac{m}{2}\de}]
			%=& \delta_{m, ev} \frac{\sqq^m-\sqq^{-m}}{m} [K_{\frac{m}{2}\de}]\\
			=(\sqq-\sqq^{-1})\de_{m,ev}\frac{[m]_\sqq}{m} [K_{\frac{m}{2}\de}].
		\end{align*}	
		Then the desired result follows from %{\red{(when $\frac{m}{d_x}$ even)}}
		\begin{align*}
			\frac{\sqq_x-\sqq_x^{-1}}{\sqq-\sqq^{-1}} [\frac{m}{2d_x}]_{\sqq_x}=[\frac{m}{2}]_\sqq.
		\end{align*}
	\end{proof}

	\subsection{The homomorphism $\Omega$}
	\label{subsec:homo}
	
	Recall the star-shaped graph $\Gamma=T_{p_1,\dots,p_\bt}$ in \eqref{star-shaped}. The following is one of the main results of this paper.
	
	\begin{theorem}
		\label{thm:morphi}
		For any star-shaped graph $\Gamma$, let $\fg$ be the Kac-Moody algebra and $\X$ be the weighted projective line associated to $\Gamma$. Then there exists a $\Q(\sqq)$-algebra homomorphism
		\begin{align}
			\Omega: \tUiD_{ |v=\sqq}\longrightarrow \iH(\X_\bfk),
		\end{align}
		which sends
		\begin{align}
			\label{eq:mor1}
			&\K_{\star}\mapsto [K_{\co}], \qquad \K_{[i,j]}\mapsto [K_{S_{ij}}], \qquad C\mapsto [K_\de];&
			\\
			\label{eq:mor2}
			&{B_{\star,l}\mapsto \frac{-1}{q-1}[\co(l\vec{c})]},\qquad\Theta_{\star,r} \mapsto {\widehat{\Theta}_{\star,r}}, \qquad H_{\star,r} \mapsto {\widehat{H}_{\star,r}};\\
			\label{eq:mor3}
			&\y_{[i,j],l}\mapsto {\frac{-1}{q-1}}\haB_{[i,j],l}, \quad \Theta_{[i,j],r}\mapsto\widehat{\Theta}_{[i,j],r}, \quad H_{[i,j],r}\mapsto \widehat{H}_{[i,j],r};
		\end{align}
		for any $[i,j]\in\II-\{\star\}$, $l\in\Z$, $r>0$.
		%extends to a unique algebra homomorphism $\Omega: \tUiD\rightarrow \tM(\X)$, which sends
		%$$\y_{[i,j],l}\mapsto \red{\frac{-1}{q-1}}\haB_{[i,j],l}; \quad \Theta_{j,r}\mapsto\widehat{\Theta}_{j,r}; \quad H_{[i,j],r}\mapsto \widehat{H}_{[i,j],r}.$$
		%Moreover, \red{if $\Gamma$ is finite or affine,} if $\X$ is domestic type or tubular type, then $\Omega$ is an isomorphism.
	\end{theorem}
	
	We shall verify that the relations~\eqref{iDR1a}--\eqref{iDR5} are preserved by $\Omega$, and the proof of this theorem consists of the next 3 sections. %thanks to Lemma~\ref{lem:equiv}.
	%(Later in Theorem~\ref{thm:injectivity}, we shall strengthen Theorem ~\ref{thm:morphi} by showing that $\Omega$ is injective if $\fg$ is of finite or affine type.)

	The composition subalgebra $\tCMHX$ is the subalgebra of $\tMHX$ generated by the elements $[\co(l\vec{c})]$ for $l\in\Z$, $\widehat{\Theta}_{\star,r}$ for $r\geq1$, $[S_{ij}]$ for $1\leq i\leq \bt$ and $1\leq j\leq p_i-1$, together with $[K_{\co}]^{\pm1}, [K_{\delta}]^{\pm1}$ and $[K_{S_{ij}}]^{\pm1}$.
	
	\begin{corollary}
		\label{cor:epimorphism}
		For any star-shaped graph $\Gamma$, let $\fg$ be the Kac-Moody algebra and $\X$ be the weighted projective line associated to $\Gamma$. Then there exists a $\Q(\sqq)$-algebra epimorphism
		\begin{align}
			\Omega: \tUiD_{ |v=\sqq}\longrightarrow \tCMHX,
		\end{align}
		defined as in \eqref{eq:mor1}--\eqref{eq:mor3}.
	\end{corollary}
	
	\begin{proof}
		Using Lemma \ref{lem:reduced generators}, it follows from Theorem \ref{thm:morphi}.
	\end{proof}

	%\begin{remark}
	%The structure of $\tCMHX$ is rich, which will be studied in depth in a separate publication; compare with \cite{BS13}. Moreover, a PBW basis will be given for the $\imath$quantum group via coherent sheaves if $\X_\bfk$ is of domestic type.
	%\end{remark}
	
	%\red{uniqueness follows from Lemma \ref{lem:reduced generators}.}
	
	%\green{The images of all generators shall be put in \S\ref{sec:cyclic}.}
	%

	%\begin{corollary}
	%Assume $\Gamma$ is of finite type. Then there exists an algebra monomorphism $\Omega: \tUiD\stackrel{\cong}{\rightarrow} \tCMHX$
	%which is uniquely determined by \eqref{eq:mor1}--\eqref{eq:mor3}.
	%\end{corollary}

	%%%%%%%%%%%%%%%%%
	\section{Relations in tubes and commutative relations}
	\label{sec:Relationtube}
	
	In this section, we shall verify the relations \eqref{iDR1a}--\eqref{iDR5} for the indexes $\mu,\nu\in\II$, where $\{\mu,\nu\}\neq\{\star, [i,1]\}$ for any $1\leq i\leq \bt$.
	More precisely, we will consider the following three cases: $\mu,\nu\in\II-\{\star\}$, or $\mu=\star=\nu$, or $\{\mu,\nu\}=\{\star, [i,j]\}$ for $2\leq j\leq p_i-1$ in \S\ref{subsec:Relations tor}, \S\ref{subsec:Relationstar} and \S\ref{subsec:relationsstarjneq1} respectively.
	
	%In the following \S\ref{subsec:Relations tor}, we shall verify the relations \eqref{iDR1a}--\eqref{iDR5} for $[i,j],[i',j']\in\II-\{\star\}$; in \S\ref{subsec:Relationstar}, we shall verify the relations \eqref{iDR1a}--\eqref{iDR5} for the point $\star$; in \S\ref{subsec:relationsstarjneq1}, we shall verify the relations  \eqref{iDR1a}--\eqref{iDR5} between $\star$ and $[i,j]$ for $j\neq 1$. %in $\iH(\scrt)$.

	%Recall that $\scrt$, which consists of torsion sheaves, is a hereditary abelian category. Let $\iH(\scrt)$ be the $\imath$Hall algebra of $\scrt$. Then $\iH(\scrt)$ is a subalgebra of $\iH(\X_\bfk)$.
	%For simplicity, we denote the vertices $[i,\nu]$'s by $\mu,\nu$
	Let $k_1, k_2, l\in \Z$ and $\mu,\nu \in \II$. Inspiring by \eqref{eq:Skk}--\eqref{eq:Rkk}, we introduce shorthand notations:
	\begin{align}
		%\begin{split}
		\notag
		\widehat{S}(k_1,k_2\mid l;\mu,\nu)
		&:=  \haB_{\mu,k_1} * \haB_{\mu,k_2} *\haB_{\nu,l} -[2]_\sqq \haB_{\mu,k_1} *\haB_{\nu,l} * \haB_{\mu,k_2} + \haB_{\nu,l} *\haB_{\mu,k_1} * \haB_{\mu,k_2},
		%&=\Sym_{k_1,k_2}\sum_{t=0}^2(-1)^t \qbinom{2}{t}B_{\mu,k_1}\cdots B_{\mu,k_t}B_{\nu,l} B_{\mu,k_{t+1}} \cdots B_{\mu,k_2}\notag \\
		%&=\sum_{t=0}^2(-1)^t \qbinom{2}{t}B_{\mu,k_1}\cdots B_{\mu,k_t}B_{\nu,l} B_{\mu,k_{t+1}} \cdots B_{\mu,k_2},
		\\
		\widehat{\SS}(k_1,k_2\mid l;\mu,\nu)
		&:= \widehat{S}(k_1,k_2\mid l;\mu,\nu)  + \{k_1 \leftrightarrow k_2 \}.
		\label{eq:haSkk}
		% \end{split}
	\end{align}
	%Here and below, $\{k_1 \leftrightarrow k_2 \}$ stands for repeating the previous summand with $k_1, k_2$ switched, so the sums over $k_1, k_2$ are symmetrized.
	We also denote
	\begin{align}
		\notag
		%\begin{split}
		\widehat{R}(k_1,k_2\mid l; \mu,\nu)
		&:=  [K_{k_1\de+\alpha_\mu}]* %\K_\mu  C^{k_1}
		\Big(-\sum_{p\geq0} \sqq^{2p}  [2]_\sqq \cdot [\haTh _{\mu,k_2-k_1-2p-1},\haB_{\nu,l-1}]_{\sqq^{-2}}*[K_\de]^{p+1}
		\\\notag
		&\qquad -\sum_{p\geq 1} \sqq^{2p-1}  [2]_\sqq \cdot [\haB_{\nu,l},\haTh _{\mu,k_2-k_1-2p}]_{\sqq^{-2}} * [K_\de]^{p}
		- [\widehat{B}_{\nu,l}, \haTh _{\mu,k_2-k_1}]_{\sqq^{-2}} \Big),
		\\
		\widehat{ \R}(k_1,k_2\mid l; \mu,\nu) &:=\widehat{ R}(k_1,k_2\mid l;\mu,\nu) + \{k_1 \leftrightarrow k_2\}.\label{eq:haRkk}
		%\end{split}
	\end{align}
	
	By direct computations, we have %({\red{the following two are special cases of the third one!}})
	\begin{align}
		\label{eqnSSkk1}
		&\widehat{ \R}(k,k\mid l; \mu,\nu)=- [\widehat{B}_{\nu,l}, \haTh _{\mu,0}]_{\sqq^{-2}}*[K_{k\de+\alpha_\mu}]=-\sqq^{-1}\widehat{B}_{\nu,l} *[K_{k\de+\alpha_\mu}],\\
		\label{eqnSSkk+1}
		&\widehat{ \R}(k,k+1\mid l; \mu,\nu)=\Big(-[2]_\sqq \cdot [\haTh _{\mu,0},\haB_{\nu,l-1}]_{\sqq^{-2}}*[K_\de]- [\widehat{B}_{\nu,l}, \haTh _{\mu,1}]_{\sqq^{-2}}\Big)*[K_{k\de+\alpha_\mu}],
	\end{align}
	and for $k_2\geq k_1+2$,
	\begin{align}
		\label{Theta and B k_2>k_1+1}
		&\widehat{ \R}(k_1,k_2\mid l; \mu,\nu)-q\widehat{ \R}(k_1,k_2-2\mid l; \mu,\nu)*[K_\de]
		\\\notag&=\Big(-[2]_\sqq \cdot [\haTh _{\mu,k_2-k_1-1},\haB_{\nu,l-1}]_{\sqq^{-2}}*[K_\de] -  [\haB_{\nu,l},\haTh _{\mu,k_2-k_1-2}]_{\sqq^{-2}}* [K_\de]
		- [\widehat{B}_{\nu,l}, \haTh _{\mu,k_2-k_1}]_{\sqq^{-2}} \Big)\\\notag
		&\quad *[K_{k_1\de+\alpha_\mu}].
	\end{align}

	\subsection{Relations for vertices in $\II-\{\star\}$}
	\label{subsec:Relations tor}
	
	\begin{proposition}
		For any $\mu=[i,j],\nu=[i',j']\in\II-\{\star\}$, the following relations hold in $\iH(\X_\bfk)$, where $k,l\in\Z$, $m,n\geq 1:$
		\begin{align}
			\label{HHmn}
			&[\widehat{H}_{\mu,m},\widehat{H}_{\nu,n}]=0,\\
			&[\widehat{H}_{\mu,m},\widehat{\y}_{\nu,l}]=\frac{[ma_{\mu\nu}]_\sqq}{m} \widehat{\y}_{\nu,l+m}-\frac{[ma_{\mu\nu}]_\sqq}{m} \widehat{\y}_{\nu,l-m}*[K_{m\de}],\label{humbvl}
			\\%3
			&[\widehat{\y}_{\mu,k} ,\widehat{\y}_{\nu,l}]=0,   \text{ if }a_{\mu\nu}=0,
			\\%4
			&[\widehat{\y}_{\mu,k}, \widehat{\y}_{\nu,l+1}]_{\sqq^{-a_{\mu\nu}}}  -\sqq^{-a_{\mu\nu}} [\widehat{\y}_{\mu,k+1}, \widehat{\y}_{\nu,l}]_{\sqq^{a_{\mu\nu}}}=0, \text{ if } a_{\mu\nu}=-1,
			\\ %5
			\label{Bmuk Bmul}
			&[\widehat{\y}_{\mu,k}, \widehat{\y}_{\mu,l+1}]_{\sqq^{-2}}  -\sqq^{-2} [\widehat{\y}_{\mu,k+1}, \widehat{\y}_{\mu,l}]_{\sqq^{2}}
			=(1-q)^2
			\Big(\sqq^{-2}\widehat{\Theta}_{\mu,l-k+1} *[K_{k\de+\alpha_{ij}}]
			\\
			&\qquad-\sqq^{-4}\widehat{\Theta}_{\mu,l-k-1} * [K_{(k+1)\de+\alpha_{ij}}]
			+\sqq^{-2}\widehat{\Theta}_{\mu,k-l+1}* [K_{l\de+\alpha_{ij}}]-\sqq^{-4}\widehat{\Theta}_{\mu,k-l-1}* [K_{(l+1)\de+\alpha_{ij}}]\Big), \notag
			\\
			\label{SerrehaBB}
			&    \widehat{\SS}(k_1,k_2\mid l; \mu,\nu) = (1-q)^2 \widehat{\R}(k_1,k_2\mid l; \mu,\nu), \text{ if } a_{\mu\nu}=-1.
		\end{align}
	\end{proposition}
	
	\begin{proof}
		%Let $\bla_i$ be the exceptional closed point of $\X$ of weight $p_i$. %, and let $\scrt_{\bla_i}$ be the Serre subcategory of $\coh(\X)$ consisting of torsion sheaves supported at $\bla_i$. It is well known that $\scrt_{\bla_i}\simeq \rep_{\bfk}^{\rm nil}(C_{p_i})$.
		% , where $Q_x$ is a cyclic quiver consisting of $p_x$ many vertices. This induces a canonical embedding of Hall algebras:
		%$\widetilde{\ch}\big(\calc_1(\rep_{\bfk_x}^{\rm nil}(Q_x))\big) \rightarrow \widetilde{\ch}\big(\calc_1(\coh(\PL))\big)$ and then an embedding
		Recall from  \eqref{eq:embeddingx} and \eqref{the map Psi A} that we have the following algebra embeddings
		\[
		\tUiD_\sqq(\widehat{\mathfrak{sl}}_{p_i})\stackrel{\Omega_{C_{p_i}}\,\,}{\longrightarrow}  {}^\imath\widetilde{\ch}(\bfk C_{p_i}) \stackrel{\iota_i}{\longrightarrow} \iH(\X_\bfk).
		\]
		%induced by $\scrt_{\bla_i}\simeq \rep_{\bfk}^{\rm nil}(C_{p_i})$.
		%Recall the algebra embedding $\Omega_{C_{p_i}}:\tUiD(\widehat{\mathfrak{sl}}_{p_i})\longrightarrow \iH(\bfk C_{p_i})$ defined in (\ref{the map Psi A}) and the algebra embedding
		%that there exists an algebra monomorphism $\Omega_{C_n}: {}^{Dr}\tUi(\widehat{\mathfrak{sl}}_{p_x+1})\rightarrow {}^\imath\widetilde{\ch}(\bfk_xQ_x)$. Hence we get an embedding
		%$$\iota_i\circ \Omega_{C_{p_i}}:\tUiD(\widehat{\mathfrak{sl}}_{p_i})\longrightarrow\iH(\X_\bfk).$$
		Now it follows from \eqref{def:haB} and \eqref{def:haBThH} that all the relations \eqref{HHmn}--\eqref{SerrehaBB} hold for $i'=i$.
		
		Moreover, for any torsion sheaves $X,Y$ supported at distinct points, we have $\Ext^1_\X(X,Y)=0=\Hom_\X(X,Y)$. It is well known that $\Ext^1_{\cc_1(\coh(\X))}(X,Y)\cong \Ext^1_\X(X,Y)\oplus \Hom_\X(X,Y)$, and $\Hom_{\cc_1(\coh(\X))}(X,Y)=\Hom_\X(X,Y)$; see e.g. \cite[Lemma 2.6]{LRW20a}. Hence $[X]*[Y]=[X\oplus Y]=[Y]*[X]$ in $\widetilde{\ch}(\cc_1(\coh(\X)))$ and in $\iH(\X_\bfk)$. Therefore, all the relations \eqref{HHmn}--\eqref{SerrehaBB} hold for $i'\neq i$.
	\end{proof}
	
	\subsection{Relations for $\star$}
	\label{subsec:Relationstar}

	The following proposition shows that (\ref{iDR1b}) and (\ref{iDR3b}) hold  in $\iH(\X_\bfk)$ for $\mu=\nu=\star$.
	\begin{proposition}
		\label{prop:OO}
		For any $m,n\geq1$, $k,l\in\Z$, the following relations hold in $\iH(\X_\bfk)$:
		\begin{align}
			\label{eq:HaDr1}
			[\haH_{\star, m},\haH_{\star,n}]&=0,\\
			\label{eq:HaDr2}
			\big[\haH_{\star,m}, [\co(l\vec{c})]\big] &=\frac{[2m]_\sqq}{m} [\co((l+m)\vec{c})]-\frac{[2m]_\sqq}{m} [\co((l-m)\vec{c})]*[K_{m\de}],
			%\big[\haT_{\star,m},&[\co(l)] \big] + \big[\haT_{\star,m-2},[\co(l)] \big]*[K_\de]
			%\\\notag
			%&\qquad\qquad\quad =\sqq^{2} \big[\haT_{\star,m-1},[\co(l+1)] \big]_{\sqq^{-4}} +\sqq^{-2} \big[\haT_{\star,m-1},[\co(l-1)] \big]_{\sqq^{4}}*[K_\de].
		\end{align}
		\begin{align}
			\label{eq:HaDr3}
			&\big[[\co(k\vec{c})], [\co((l+1)\vec{c})]\big]_{\sqq^{-2}}  -\sqq^{-2} \big[[\co((k+1)\vec{c})], [\co(l\vec{c})]\big]_{\sqq^{2}}
			\\\notag
			&= (1-q)^2\Big(\sqq^{-2}\haTh_{\star,l-k+1} *[K_\de]^k* [K_{\co}]-\sqq^{-4} \haTh_{\star,l-k-1}* [K_\de]^{k+1} *[K_{\co}]
			\\\notag
			&\quad +\sqq^{-2}\haTh_{\star,k-l+1}*[K_\de]^l *[K_{\co}]-\sqq^{-4} \haTh_{\star,k-l-1}*[ K_\de]^{l+1} *[K_\co]\Big).
		\end{align}
		%\red{Consequently, under the map $\Omega: \tUiD\rightarrow \tM(\X)$, all the relations in Definition \ref{def:iDR} are satisfied in $\tM(\X)$ for $i=j=\star$.}
	\end{proposition}
	
	\begin{proof}
		Recall from \S \ref{Embedding from projective line to weighted projective line}
		that there is an embedding \begin{equation*}
			F_{\X,\P^1}:\iH(\P^1_\bfk) \longrightarrow \iH(\X_\bfk),
		\end{equation*} which sends
		$[\co_{\P^1}(l)]\mapsto[\co(l\vec{c})]$, $[S_{\bla_i}^{(r)}]\mapsto [S_{i,0}^{(rp_i)}]$, $\haT_{r}\mapsto \haT_{\star,r}$ and $\haH_r\mapsto \haH_{\star,r}$
		for any $l\in\Z, r\geq 1$.
		Then the formulas follow from \cite[Theorem 4.2]{LRW20a} and the Drinfeld type presentation of $\tUi_v(\widehat{\mathfrak{sl}}_2)$ therein.
	\end{proof}

	\subsection{Relations between $\star$ and $[i,j]$ with $j\neq 1$}
	\label{subsec:relationsstarjneq1}

	In this subsection, we verify the relations \eqref{iDR1b}--\eqref{iDR4} between $\star$ and $[i,j]$ with $j\neq 1$.
	First we give some preliminary results.
	
	The set $\cm_{j,\alpha}$ defined in \eqref{def:Mjalpha} can be viewed as a set of isoclasses of $\coh(\X)$ via $\iota_i$; see \S\ref{subsec:embeddingtube}.
	Recall from Corollary \ref{cor:DrGenA} that
	$$\widehat{\Theta}_{[i,j],1}=\widehat{H}_{[i,j],1}=\pi_{[i,j+1],1}-(\sqq+\sqq^{-1})\pi_{[i,j],1}+\pi_{[i,j-1],1},$$
	where \begin{align*}
		\pi_{[i,j],1}=\frac{-\sqq^{-j}}{\sqq-\sqq^{-1}}\sum\limits_{\cm_{j,\delta}}(-1)^{\dim\End(M)}[M].
	\end{align*}
	
	\begin{lemma}
		\label{level one pi ij and star}
		For any $1\leq i\leq \bt$, $1\leq j\leq p_i-1$, and any $l\in\mathbb{Z}$, $r\geq 1$, we have
		\begin{align}\label{pi and o}
			\big[\pi_{[i,j],1}, [\co(l\vec{c})]\big]=&\frac{\sqq^{-j}}{\sqq-\sqq^{-1}}\big[[S_{i,0}^{(p_i)}], [\co(l\vec{c})]\big];\\
			\label{pi and S(rp)}
			\big[\pi_{[i,j],1}, [S_{i,0}^{(rp_i)}]\big]=&0.
		\end{align}
	\end{lemma}
	
	\begin{proof}
		
		For any $1\leq j',k\leq p_i$, the torsion sheaf $S_{i,j'}^{(k)}$ is uniserial, and $S_{i,1}$ is a composition factor of $S_{i,j'}^{(k)}$ if and only if $k\geq j'$. Hence, for any $[M]\in\cm_{j,\delta}$, there exists a unique direct summand $S_{i,j'}^{(k)}$ of $M$ with $k\geq j'$, namely, $M\cong S_{i,j'}^{(k)}\oplus M'$ for some $M'$.
		We emphasize that $$[S_{i,j'}^{(k)}\oplus M']\in\cm_{j,\delta}\quad \text{if\; and\; only\; if} \quad[S_{i,j'}^{(j')}\oplus S_{i,0}^{(k-j')}\oplus M']\in\cm_{j,\delta}.$$
		For any $k>j'$, we have $$\dim\End(S_{i,j'}^{(j')}\oplus S_{i,0}^{(k-j')}\oplus M')=\dim\End(S_{i,j'}^{(k)}\oplus M')+1;$$ moreover, there are no homomorphisms and extensions between $\co(l\vec{c})$ and $M'$ since both $S_{i,0}$ and $S_{i,1}$ are not composition factors of $M'$, hence by (\ref{eq:OSij}) we get
		$$\big[[S_{i,j'}^{(k)}\oplus M'], [\co(l\vec{c})]\big]=\big[[S_{i,j'}^{(j')}\oplus S_{i,0}^{(k-j')}\oplus M'], [\co(l\vec{c})]\big].$$
		Therefore, only for the case $k=j'=p_i$, i.e., the term $[S_{i,0}^{(p_i)}]$ in $\pi_{[i,j],1}$, has non-trivial contribution to $\big[\pi_{[i,j],1}, [\co(l\vec{c})]\big]$.
		Then (\ref{pi and o}) follows immediately.
		
		For \eqref{pi and S(rp)}, by similar arguments as above and using (\ref{eq:Sipiij}), we obtain that
		\begin{align*}
			\big[\pi_{[i,j],1}, [S_{i,0}^{(rp_i)}]\big]=\frac{\sqq^{-j}}{\sqq-\sqq^{-1}}\big[[S_{i,0}^{(p_i)}], [S_{i,0}^{(rp_i)}]\big]=0,
		\end{align*}
		where the last equality follows from \cite[Lemma 4.5]{LRW20a} and the embedding (\ref{the embedding functor F on algebra}). %Then (\ref{h_ij^1 and s(rp) commutes}) follows from (\ref{pi and S(rp)}).
		%The proof is completed.
	\end{proof}
	
	Consequently, we have the following result.

	\begin{lemma}\label{level one hij and star}
		For any $1\leq i\leq \bt$ and $2\leq j\leq p_i-1$, we have
		\begin{itemize}
			\item[(1)] $\big[\widehat{\Theta}_{[i,j],1}, [\co(l\vec{c})]\big]=0$ for any $l\in\Z$;
			\item[(2)] $\big[\widehat{\Theta}_{[i,j],1}, \widehat{\Theta}_{\star,m}\big]=0$ for any $m\geq 1$.
		\end{itemize}
	\end{lemma}
	
	\begin{proof} For the statement (1), by (\ref{pi and o}) we have
		\begin{align*}
			&\big[\widehat{\Theta}_{[i,j],1}, [\co(l\vec{c})]\big]\\
			&=\big[\pi_{[i,j+1],1}-(\sqq+\sqq^{-1})\pi_{[i,j],1}+\pi_{[i,j-1],1},[\co(l\vec{c})]\big]\\
			&=\big(\frac{\sqq^{-(j+1)}}{\sqq-\sqq^{-1}}-(\sqq+\sqq^{-1})\cdot\frac{\sqq^{-j}}{\sqq-\sqq^{-1}}+\frac{\sqq^{-(j-1)}}{\sqq-\sqq^{-1}}\big)\big[[S_{i,0}^{(p_i)}], [\co(l\vec{c})]\big]\\
			&=0.
		\end{align*}
		For the statement (2), observe that there are no non-zero homomorphisms and extensions between distinct tubes, then by the definition of $\widehat{\Theta}_{\star,m}$, it suffices to show that
		\begin{align*}\big[\widehat{\Theta}_{[i,j],1}, [S_{i,0}^{(rp_i)}]\big]=0, \quad\forall\; r\geq 1,
		\end{align*}
		which then follows from (\ref{pi and S(rp)}) immediately.
	\end{proof}

	\begin{lemma}\label{level one eta and star}
		For any $1\leq i\leq \bt$ and $2\leq j\leq p_i-1$, we have
		\begin{itemize}
			\item[(1)] $\big[\widehat{B}_{[i,j],-1}, [\co(l\vec{c})]\big]=0$ for any $l\in\Z$;
			\item[(2)] $\big[\widehat{B}_{[i,j],-1}, \widehat{\Theta}_{\star,m}\big]=0$ for any $m\geq 1$.
		\end{itemize}
	\end{lemma}
	
	\begin{proof}
		Recall from Proposition \ref{prop:DrGenA} and (\ref{eq:embeddingx})
		that $$\haB_{[i,j],-1}=\sqq^{-j+1}\sum\limits_{\cm_{j+1,\delta-e_{j}}}(-1)^{\dim\End(M)}[M]*[K_{\de-\alpha_j}]^{-1}.$$

		For any term $[M]$ in $\cm_{j+1,\delta-e_{j}}$, we know that $S_{i,1}$ is a composition factor of $M$ with multiplicity one. Hence there exists a unique direct summand $S_{i,j'}^{(k)}$ of $M$ with $k\geq j'$. Then by similar arguments as in the proof of Lemma \ref{level one pi ij and star} (observing $[S_{i,0}^{(rp_i)}]\notin\cm_{j+1,\delta-e_{j}}$), we obtain $\big[\widehat{B}_{[i,j],-1}, [\co(l\vec{c})]\big]=0$ for any $l\in\bbZ$, and $\big[\widehat{B}_{[i,j],-1}, [S_{i,0}^{(rp_i)}]\big]=0$ for any $r\geq 1$. Then the result follows immediately.
	\end{proof}
	
	\subsubsection{Relation (\ref{iDR4}) between $\star$ and $[i,j]$ with $j\neq 1$}
	
	\begin{proposition}
		\label{co kc and Bijl}
		For any $1\leq i\leq \bt, \;2\leq j\leq p_i-1$, and any $k,l\in\Z$, we have
		\begin{align}\label{co kc and Bij l}
			\big[[\co(k\vec{c})],\widehat{\y}_{[i,j],l}\big]=0.
		\end{align}
	\end{proposition}
	
	\begin{proof}
		%\emph{\blue{Method 1: short proof via simplify relations on algebras}}
		
		%By using similar arguments as in Lemma \ref{bij and bjl reduced to 01}, we only need to prove (\ref{co kc and Bij l}) for $l=0$ and -1.
		
		%Recall that $\widehat{\y}_{ij,0}=[S_{ij}]$, and there are no non-zero homomorphisms and extensions between $\co(k\vec{c})$ and $S_{ij}$ for $2\leq j\leq p_i-1$. Hence %$\big[[\co(k\vec{c})],[S_{ij}]\big]=0$. Moreover, by Lemma \ref{level one eta and star} (1), we have
		%$\big[[\co(k\vec{c})], \widehat{\y}_{ij,-1}\big]=0$. The proof is completed.
		
		%\emph{\red{Method 2: complete calculation in category}}
		%\red{Delete Lemma \ref{bij and bjl reduced to 01} }
		For $2\leq j\leq p_i-1$, by (\ref{humbvl}) we have
		$$[\widehat{H}_{[i,j],1},\widehat{\y}_{[i,j],l}]=[2]_\sqq \widehat{\y}_{[i,j],l+1}-[2]_\sqq \widehat{\y}_{[i,j],l-1}*[K_{\delta}].$$
		Hence \begin{equation}\label{l and l pm 1}
			\big[[\co(k\vec{c})],[\widehat{H}_{[i,j],1},\widehat{\y}_{[i,j],l}]\big]=[2]_\sqq \big[[\co(k\vec{c})],\widehat{\y}_{[i,j],l+1}\big]-[2]_\sqq \big[[\co(k\vec{c})],\widehat{\y}_{[i,j],l-1}\big]*[K_{\delta}].\end{equation}
		By Lemma \ref{level one hij and star} (1), we have $\big[[\co(k\vec{c})],\widehat{H}_{[i,j],1}\big]=0$. Hence
		$$\big[[\co(k\vec{c})],[\widehat{H}_{[i,j],1},\widehat{\y}_{[i,j],l}]\big]=\big[\widehat{H}_{[i,j],1},\big[[\co(k\vec{c})],\widehat{\y}_{[i,j],l}\big]\big].$$
		%\begin{align*}
		%&[H_{\star, m},[H_{ij,1},\y_{ij,l}]]\\
		%=&[[H_{\star, m},H_{ij,1}],\y_{ij,l}]+[H_{ij,1},[H_{\star, m},\y_{ij,l}]]\\
		%=&[H_{ij,1},[H_{\star, m},\y_{ij,l}]].
		%\end{align*}
		If $\big[[\co( k\vec{c})],\y_{[i,j],l}\big]=0$, then LHS of $\eqref{l and l pm 1}=0$, hence $\big[[\co( k\vec{c})],\haB_{[i,j],l+1}\big]=0$ if and only if $\big[[\co( k\vec{c})],\haB_{[i,j],l-1}\big]=0$. Therefore, by induction, we only need to prove (\ref{co kc and Bij l}) for $l=0$ and $-1$.
		
		Recall that $\widehat{\y}_{[i,j],0}=[S_{ij}]$, and there are no non-zero homomorphisms and extensions between $\co(k\vec{c})$ and $S_{ij}$ for $2\leq j\leq p_i-1$. Hence $\big[[\co(k\vec{c})],[S_{ij}]\big]=0$. Moreover, by Lemma \ref{level one eta and star} (1), we have
		$\big[[\co(k\vec{c})], \widehat{\y}_{[i,j],-1}\big]=0$. %\green{Add details, can not follow}
		The proof is completed.
	\end{proof}

	\subsubsection{Relation (\ref{iDR2}) between $\star$ and $[i,j]$ with $j\neq 1$}

	\begin{proposition}\label{star and ij for j neq 1}
		For any $1\leq i\leq \bt, \;2\leq j\leq p_i-1$, and any $m\geq 1, l\in\Z$, we have
		\begin{itemize}
			\item[(1)] $\big[\widehat{H}_{\star, m},\haB_{[i,j],l}\big]=0$;
			\item[(2)] $\big[\haH_{[i,j], m}, [\co(l\vec{c})]\big]=0.$
		\end{itemize}
	\end{proposition}
	
	\begin{proof}
		(1) It is equivalent to prove that $\big[\widehat{\Theta}_{\star, m},\haB_{[i,j],l}\big]=0$.
		By using similar arguments as in Proposition \ref{co kc and Bijl}, thanks to (\ref{humbvl}) and Lemma \ref{level one hij and star} (2), we only need to prove the formula for $l=0$ and $-1$. For $l=0$, $\Big[\widehat{\Theta}_{\star, m},[S_{ij}]\Big]=0$ is obvious since $\big[ [S_{i,0}^{(p_i)}], [S_{ij}] \big]=0$ for $2\leq j\leq p_i-1$; while for $l=-1$, it follows from Lemma \ref{level one eta and star} (2).

		(2)
		By (\ref{Bmuk Bmul}) we have
		\begin{align}
			\label{bb to hij for n=1}
			&\big[\haB_{[i,j],l+1}, \haB_{[i,j],l+1}\big]_{\sqq^{-2}}  -\sqq^{-2} \big[\haB_{[i,j],l+2}, \haB_{[i,j],l}\big]_{\sqq^{2}}
			\\\notag
			&=(1-q)^2\big(\sqq^{-3} [K_{(l+1)\de+\alpha_{ij}}]+\sqq^{-2}\widehat{\Theta}_{[i,j],2} *[K_{l\de+\alpha_{ij}}]\big). %C^l \K_{ij}.
		\end{align}
		Apply the operator $\big[-, [\co(l\vec{c})]\big]$ to (\ref{bb to hij for n=1}). Then it follows from Proposition \ref{co kc and Bijl} and Lemma \ref{level one hij and star} (1) that $$\big[\widehat{\Theta}_{[i,j],1}, [\co(l\vec{c})]\big]=0=\big[\widehat{\Theta}_{[i,j],2}, [\co(l\vec{c})]\big].$$
		For $n\geq 2$, by (\ref{Bmuk Bmul}) again we obtain
		\begin{align}
			\label{bb to hij for general n}
			&\big[\widehat{\y}_{[i,j],l+n}, \haB_{[i,j],l+1}\big]_{\sqq^{-2}}  -\sqq^{-2} \big[\widehat{\y}_{[i,j],l+n+1}, \haB_{[i,j],l}\big]_{\sqq^{2}}
			\\\notag
			&=(1-q)^2\big(\sqq^{-2}\widehat{\Theta}_{[i,j],n+1} * [K_{l\de+\alpha_{ij}}] -\sqq^{-4}\widehat{\Theta}_{[i,j],n-1} *[K_{(l+1)\de+\alpha_{ij}}]\big). %C^{l+1} \K_{ij}.
		\end{align}
		By applying $\big[-, [\co(l\vec{c})]\big]$ to (\ref{bb to hij for general n}), we obtain $[\widehat{\Theta}_{[i,j],m}, [\co(l\vec{c})]]=0$ for any $m\geq 1$ inductively. Then (2) follows from (\ref{exp h}).
	\end{proof}
	
	\subsubsection{Relation (\ref{iDR1b}) between $\star$ and $[i,j]$ with $j\neq 1$}
	
	\begin{proposition}
		For any $1\leq i\leq \bt, \;2\leq j\leq p_i-1$, and any $m,n\geq 1,$
		\begin{align}\label{H star m and Hij r}
			[\widehat{H}_{\star, m}, \widehat{H}_{[i,j],n}]=0.
		\end{align}
	\end{proposition}
	
	\begin{proof}
		The proof is similar to that of Proposition \ref{star and ij for j neq 1} (2).  Apply the operator $[\widehat{\Theta}_{\star, m}, -]$ to (\ref{bb to hij for n=1}). By Proposition \ref{co kc and Bijl} and Lemma \ref{level one hij and star} (2), we have
		$$[\widehat{\Theta}_{\star, m}, \widehat{\Theta}_{[i,j],1}]=0=[\widehat{\Theta}_{\star, m}, \widehat{\Theta}_{[i,j],2}].$$
		Now by applying $[\widehat{\Theta}_{\star, m}, -]$ to (\ref{bb to hij for general n}), we obtain
		$[\widehat{\Theta}_{\star, m}, \widehat{\Theta}_{[i,j],n}]=0$ for any $n\geq 1$. Then (\ref{H star m and Hij r}) follows from (\ref{exp h}).
	\end{proof}

	%%%%%%%%%%%%%%%%%%
	\section{Relations between star point and $[i,1]$, I}
	\label{sec:Relationsstari1 I}
	
	In this section and Section \ref{sec:Relationsstari1 II}, we shall verify the relations \eqref{iDR1a}--\eqref{iDR5} between $\star$ and $[i,1]$. In this section, we shall prove the following formulas in $\iH(\X_\bfk)$:
	\begin{align}
		\label{coBi1}
		&\big[[\co(k\vec{c})], \widehat{\y}_{[i,1],l+1}\big]_{\sqq}  -\sqq \big[[\co((k+1)\vec{c})], \widehat{\y}_{[i,1],l}\big]_{\sqq^{-1}} =0,
		\\
		\label{Serrestari1}
		&\widehat{\SS}(k_1,k_2\mid l;\star, [i,1])=(1-q)^2\widehat{\R}(k_1,k_2\mid l;\star, [i,1]),
		\\
		\label{humstarbvl}
		&[\widehat{H}_{\star,m},\widehat{\y}_{[i,1],l}]=-\frac{[m]_\sqq}{m} \widehat{\y}_{[i,1],l+m}+\frac{[m]_\sqq}{m} \widehat{\y}_{[i,1],l-m}*[K_{m\de}],
		\\
		\label{eq:HHstari1}
		&[\widehat{H}_{\star,m},\widehat{H}_{[i,1],r}]=0,
	\end{align}
	for any $k, k_1,k_2,l\in\Z$, $m,r>0$ and $1\leq i\leq \bt$.
	
	Without loss of generality, we assume $i=1$ throughout this section.
	
	\subsection{The relation \eqref{coBi1}}
	
	\begin{proposition}
		\label{prop:coBil}
		We have \eqref{coBi1} holds for any %$1\leq i\leq \bt$ and
		$k,l\in\Z$.
		%\begin{align}\label{relation BB}
		% \Big[[\co(k\vec{c})], \widehat{\y}_{[i,1],l+1}\Big]_{\sqq}  -\sqq \Big[[\co((k+1)\vec{c})], \widehat{\y}_{[i,1],l}\Big]_{\sqq^{-1}} =0.
		%\end{align}
	\end{proposition}

	\begin{proof}
		Obviously, \eqref{coBi1} is equivalent to
		\begin{align}
			\label{coBi1equiv}
			\big[[\co(k\vec{c})], \widehat{\y}_{[1,1],l+1}\big]_{\sqq}  + \big[ \widehat{\y}_{[1,1],l},[\co((k+1)\vec{c})]\big]_{\sqq} =0,\qquad\forall\,\, k,l\in\Z.
		\end{align}
		%\blue{By using similar arguments as in Lemma \ref{lem:append}, and thanks to Lemma \ref{level one hij and star} (1) and (\ref{humbvl}),
		First, we prove \eqref{coBi1equiv} for the case when $l=0$ and $-1$.
		Recall from Corollary \ref{cor for level one of i1} that
		$$\haB_{[1,1],0}=[S_{1,1}],\quad \haB_{[1,1],-1}=-[S_{1,0}^{(p_1-1)}]*[K_{\alpha_{11}-\de}], \quad \haB_{[1,1],1}=\frac{1}{q}[S_{1,1}^{(p_1+1)}] -\frac{1}{q}[S_{1,1}\oplus S_{1,0}^{(p_1)}].$$
		%Observe that $\Ext^1_{\cc_1(\coh(\X))}(S_{i,0}^{(p_i-1)}, \co(k\vec{c}+\vec{x}_i))\cong\bfk\cong\Ext^1_{\cc_1(\coh(\X))}(S_{i,1},\co(k\vec{c})),$
		%and the following exact sequences:
		%$$0\longrightarrow \co(k\vec{c}+\vec{x}_i)\longrightarrow \co((k+1)\vec{c})\longrightarrow S_{i,0}^{(p_i-1)}\longrightarrow 0,$$
		%$$0\longrightarrow \co(k\vec{c})\longrightarrow \co(k\vec{c}+\vec{x}_i)\longrightarrow S_{i,1}\longrightarrow 0.$$
		%\green{delete Lemma \ref{A2 case}}
		%Then by Lemma \ref{A2 case}, we have
		By (\ref{Si1 and co double direction}) and (\ref{Si0 and co double direction}) we obtain
		$$\big[[\co(k\vec{c})], [S_{1,1}]\big]_{\sqq}=-(q-1)[\co(k\vec{c}+\vec{x}_1)]=\big[[S_{1,0}^{(p_1-1)}]*[K_{\alpha_{11}-\de}], [\co((k+1)\vec{c})]\big]_{\sqq}.$$
		%A direct computation shows that
		%\begin{align*}
		%[S_{i,0}^{(p_i-1)}]* [\co((k+1)\vec{c})]=& [S_{i,0}^{(p_i-1)}\oplus\co((k+1)\vec{c})], %+\sqq^{-1}(q-1)[\co(k\vec{c}+\vec{x}_i)],
		%\\
		%[\co((k+1)\vec{c})]*[S_{i,0}^{(p_i-1)}]=& \sqq^{-1}[S_{i,0}^{(p_i-1)}\oplus\co((k+1)\vec{c})]+\sqq^{-1}(q-1)[\co(k\vec{c}+\vec{x}_i)]*[K_{ S_{i,0}^{(p_i-1)}}],
		%\\
		%[S_{i,1}]* [\co(k\vec{c})]=& \sqq^{-1}[ S_{i,1}\oplus\co(k\vec{c})]+\sqq^{-1}(q-1)[\co(k\vec{c}+\vec{x}_i)],
		%\\
		%[\co(k\vec{c})]*[S_{i,1}]=& [ S_{i,1}\oplus\co(k\vec{c})].
		%\end{align*}
		Hence \eqref{coBi1equiv} holds for $l=-1$.
		For $l=0$, \eqref{coBi1equiv} follows from Lemma \ref{pi+1 and co}.

		For any $l\in\Z$, from \eqref{humbvl}, we have
		\begin{align*}
			&\big[\widehat{\Theta}_{[1,1],1}, \haB_{[1,1],l}\big]=[2]_\sqq\haB_{[1,1],l+1}-[2]_\sqq\haB_{[1,1],l-1}\ast[K_{\delta}].
			%\\
			%&\big[\widehat{\Theta}_{[i,1],1}, [\co(l\vec{c})]\big]=-[\co((l+1)\vec{c})]+[\co((l-1)\vec{c})]\ast[K_{\delta}]
		\end{align*}
		Moreover, by Corollary \ref{cor:DrGenA} and Lemma \ref{level one pi ij and star},
		\begin{align}
			\label{eq:Theta111Ol}
			&\big[\widehat{\Theta}_{[1,1],1}, [\co(k\vec{c})]\big]\\\notag
			&=\big[\pi_{[1,2],1}-(\sqq+\sqq^{-1})\pi_{[1,1],1},[\co(k\vec{c})]\big]\\\notag
			&=\big(\frac{\sqq^{-2}}{\sqq-\sqq^{-1}}-\frac{\sqq^{-1}(\sqq+\sqq^{-1})}{\sqq-\sqq^{-1}}\big)\big[[S_{1,0}^{(p_1)}], [\co(k\vec{c})]\big]\\\notag
			&=-\frac{1}{\sqq-\sqq^{-1}}\cdot \sqq^{-1}(q-1)\big([\co((k+1)\vec{c})]-[\co((k-1)\vec{c})]\ast[K_{\delta}]\big)\\\notag
			&=-[\co((k+1)\vec{c})]+[\co((k-1)\vec{c})]\ast[K_{\delta}].
		\end{align}
		Therefore, we have
		\begin{align*}
			&\Big[\haTh_{[1,1],1}, \big[[\co(k\vec{c})], \widehat{\y}_{[1,1],l}\big]_{\sqq}  + \big[ \widehat{\y}_{[1,1],l-1},[\co((k+1)\vec{c})]\big]_{\sqq}\Big]
			\\
			&=\Big[\haTh_{[1,1],1}, \big[[\co(k\vec{c})], \haB_{[1,1],l}\big]_{\sqq}\Big]  + \Big[\haTh_{[1,1],1}, \big[\haB_{[1,1],l-1},[\co((k+1)\vec{c})] \big]_{\sqq}\Big]
			\\
			&=\Big[\big[\haTh_{[1,1],1}, [\co(k\vec{c})]\big], \haB_{[1,1],l}\Big]_{\sqq}+\Big[[\co(k\vec{c})], \big[\haTh_{[1,1],1},  \haB_{[1,1],l}\big]\Big]_{\sqq}
			\\
			&\quad +\Big[\big[\haTh_{[1,1],1}, \haB_{[1,1],l-1}\big], [\co((k+1)\vec{c})]\Big]_{\sqq}+\Big[\haB_{[1,1],l-1}, \big[\haTh_{[1,1],1},  [\co((k+1)\vec{c})]\big]\Big]_{\sqq}\\
			&=\Big[- [\co((k+1)\vec{c})]+ [\co((k-1)\vec{c})]*[K_\de], \haB_{[1,1],l}\Big]_\sqq\\
			&\quad +[2]_\sqq\Big[[\co(k\vec{c})], \haB_{[1,1],l+1}-\haB_{[1,1],l-1}*[K_\de]\Big]_{\sqq} \\
			&\quad+[2]_\sqq\Big[ \haB_{[1,1],l}-\haB_{[1,1],l-2}*[K_\de], [\co((k+1)\vec{c})]\Big]_{\sqq}\\
			&\quad+\Big[\haB_{[1,1],l-1}, - [\co((k+2)\vec{c})]+ [\co(k\vec{c})]*[K_\de]\Big]_\sqq\\
			&= -\Big(\big[[\co((k+1)\vec{c})],\haB_{[1,1],l}\big]_{\sqq}+\big[\haB_{[1,1],l-1},[\co((k+2)\vec{c})]\big]_{\sqq}\Big)
			\\
			&\quad +\Big(\big[[\co((k-1)\vec{c})],\haB_{[1,1],l}\big]_{\sqq}+\big[\haB_{[1,1],l-1},[\co(k\vec{c})]\big]_{\sqq}\Big)*[K_\de] \\
			&\quad +[2]_\sqq\Big(\big[[\co(k\vec{c})],\haB_{[1,1],l+1}\big]_{\sqq}+\big[\haB_{[1,1],l},[\co((k+1)\vec{c})]\big]_{\sqq}\Big)
			\\
			&\quad -[2]_\sqq\Big(\big[[\co(k\vec{c})],\haB_{[1,1],l-1}\big]_{\sqq}+\big[\haB_{[1,1],l-2},[\co((k+1)\vec{c})]\big]_{\sqq}\Big) *[K_{\delta}].
		\end{align*}
		Observe that twisting with $\vec{c}$ induces an automorphism on $\iH(\X_\bfk)$. Then \eqref{coBi1equiv} holds for any $l\in\Z$ by induction.
	\end{proof}

	\subsection{The relation \eqref{Serrestari1}}
	\label{sec:intrelationsI}

	We shall prove \eqref{Serrestari1} by induction on $l$. First we check \eqref{Serrestari1} for $l=0,-1$.
	
	\begin{lemma}
		\label{lem:SerreO11}
		For any $k_1,k_2\in\Z$, we have
		\begin{align}
			\label{SerreO11}
			\widehat{\SS}(k_1,k_2\mid 0;\star,[1,1])=&(1-q)^2\widehat{\R}(k_1,k_2\mid 0;\star,[1,1]).
		\end{align}
	\end{lemma}
	
	The proof of Lemma \ref{lem:SerreO11} is provided in Appendix \ref{proofSerre011}
	
	\begin{lemma}
		\label{lem:SerreO111}
		For any $k_1,k_2\in\Z$, we have
		\begin{align}
			\label{SerreO11-1}
			\widehat{\SS}(k_1,k_2\mid -1;\star,[1,1])=&(1-q)^2\widehat{\R}(k_1,k_2\mid -1;\star,[1,1]).
		\end{align}
	\end{lemma}
	
	The proof of Lemma \ref{lem:SerreO111} is provided in  Appendix \ref{subsec:proofSerre0111}.

	\begin{proposition}
		\label{prop:SerreO11l}
		For any $k_1,k_2,l\in\Z$, we have
		\begin{align}
			\label{SerreO11l}
			\widehat{\SS}(k_1,k_2\mid l;\star,[1,1])=&(1-q)^2\widehat{\R}(k_1,k_2\mid l;\star,[1,1]).
		\end{align}
	\end{proposition}
	
	\begin{proof}
		We only need to consider the case $l\ge0$, since the other case $l<0$ is similar.
		By Lemma \ref{lem:SerreO11}, we do induction on $l$, and assume that \eqref{SerreO11l} holds for $0\le l'\leq l$.
		%$$\widehat{\SS}(k_1,k_2\mid l;\star,[1,1])=(1-q)^2\widehat{\R}(k_1,k_2\mid l;\star,[1,1]).$$
		
		Using \eqref{humbvl} and \eqref{eq:Theta111Ol}, a direct computation shows
		\begin{align}
			\label{eq:THO}
			&\big[\haTh_{[1,1],1},\widehat{\SS}(k_1,k_2\mid l;\star,[1,1])\big]
			\\\notag
			&=-\widehat{\SS}(k_1+1,k_2\mid l;\star,[1,1])+\widehat{\SS}(k_1-1,k_2\mid l;\star,[1,1])*[K_\de]
			\\\notag
			&\quad -\widehat{\SS}(k_1,k_2+1\mid l;\star,[1,1])+\widehat{\SS}(k_1,k_2-1\mid l;\star,[1,1])*[K_\de]
			\\\notag
			&\quad +[2]_\sqq\widehat{\SS}(k_1,k_2\mid l+1;\star,[1,1])-[2]_\sqq\widehat{\SS}(k_1,k_2\mid l-1;\star,[1,1])*[K_\de]
			\\\notag
			&=(1-q)^2\Big(-\widehat{\R}(k_1+1,k_2\mid l;\star,[1,1])+\widehat{\R}(k_1-1,k_2\mid l;\star,[1,1])*[K_\de]
			\\\notag
			&\quad -\widehat{\R}(k_1,k_2+1\mid l;\star,[1,1])+\widehat{\R}(k_1,k_2-1\mid l;\star,[1,1])*[K_\de]\Big)
			\\\notag
			&\quad +[2]_\sqq\widehat{\SS}(k_1,k_2\mid l+1;\star,[1,1])-[2]_\sqq\widehat{\SS}(k_1,k_2\mid l-1;\star,[1,1])*[K_\de]
			\\\notag
			&=[2]_\sqq\widehat{\SS}(k_1,k_2\mid l+1;\star,[1,1])-[2]_\sqq\widehat{\SS}(k_1,k_2\mid l-1;\star,[1,1])*[K_\de].
		\end{align}
		Here the second equality uses the inductive assumptions, and the third equality uses the following equality obtained from \eqref{eq:haRkk}:
		\begin{align*}
			\widehat{\R}(k_1,k_2\pm 1\mid l;\star,[1,1])=\widehat{\R}(k_1\mp1,k_2\mid l;\star,[1,1])*[K_\de]^{\pm1}.
		\end{align*}
		
		On the other hand, using the same argument as in Lemma \ref{level one hij and star}~(2), one can obtain that
		$\big[\widehat{\Theta}_{[1,1],1}, \widehat{\Theta}_{\star,m}\big]=0$ for any $m>0$. Together with \eqref{humbvl}, we have
		\begin{align}
			\label{eq:THR}
			&[\haTh_{[1,1],1},\widehat{\R}(k_1,k_2\mid l;\star,[1,1])]
			\\\notag
			&= [2]_\sqq\widehat{\R}(k_1,k_2\mid l+1;\star,[1,1])-[2]_\sqq \widehat{\R}(k_1,k_2\mid l-1;\star,[1,1])*[K_\de].
		\end{align}
		Comparing \eqref{eq:THO} and \eqref{eq:THR}, using the inductive assumption, one can obtain
		$$\widehat{\SS}(k_1,k_2\mid l+1;\star,[1,1])=(1-q)^2\widehat{\R}(k_1,k_2\mid l+1;\star,[1,1]).$$
		The proof is completed.
	\end{proof}

	%%%%%%%%%%%%
	\subsection{The relations \eqref{humstarbvl} and \eqref{eq:HHstari1}}
	
	\begin{proposition}
		\label{prop:humstarbvl}
		Relation \eqref{humstarbvl} holds in $\iH(\X_\bfk)$ for any $m>0$ and $l\in\Z$.
	\end{proposition}
	
	\begin{proof}
		Denote
		\begin{align*}
			%\label{eq:Y}
			\widehat{Y}_{m|l} :=&[\haB_{[1,1],l+1}, \haTh _{\star,m-1}]_{\sqq^{-2}} - [\haTh _{\star,m-1},\haB_{[1,1],l-1}]_{\sqq^{-2}}*[K_\de]
			\\\notag
			&+\sqq^{-1} [\haTh_{\star,m},\haB_{[1,1],l}] + \sqq^{-1} [\haTh_{\star,m-2},\haB_{[1,1],l}] *[K_\de].
		\end{align*}
		Using Lemma \ref{lem:equiv}~(2), it is equivalent to prove that $\widehat{Y}_{m|l}=0$ for any $m\geq0$.
		Clearly, $\widehat{Y}_{m|l}=0$ if $m\le 0$. By induction, we assume that $\widehat{Y}_{m|l}=0$ for any $m\leq k$.

		%First, we prove \eqref{humstarbvl}.
		%For any $k_1,k_2\in\Z$ with $k_1\leq k_2$, let $k=k_2-k_1\geq0$.

		We also denote
		for $n\in \Z$,
		\begin{align*}
			\notag
			\widehat{X}_{n |l}=&\sum_{p\ge 0} \sqq^{2p+1}[\haB_{[1,1],l+1},\haTh_{\star,n-2p-1}]_{\sqq^{-2}} *[K_\de]^{p+1}
			+\sum_{p\ge 1} \sqq^{2p-1} [2]_\sqq [\haTh_{\star,n-2p},\haB_{[1,1],l}] *[K_\de]^{p+1}
			\\
			& -\sum_{p\ge 0} \sqq^{2p+1}[\haTh_{\star,n-2p-1},\haB_{[1,1],l-1}]_{\sqq^{-2}}*[K_\de]^{p+2}
			+ [\haTh_{\star,n},\haB_{[1,1],l}]*[K_\de].
			%\label{eq:X}
		\end{align*}
		Clearly, $\widehat{X}_{n |l}=0$, for $n\leq 0$. By using the same argument as in \cite[Lemma 4.10]{LW20b}, we have
		$\widehat{X}_{n|l}-q\widehat{X}_{n-2|l}=\sqq\widehat{Y}_{n|l}$. Then
		$\widehat{X}_{n |l}=0$ for any $n\leq k$ by using the inductive assumption.

		Using the same argument as in \cite[Lemma 4.14]{LW20b}, we have
		\begin{align}
			\label{RRR}
			& \widehat{\R}(k_1,k_2+1 \mid l;\star,[1,1]) + \widehat{\R}(k_1+1,k_2\mid l;\star,[1,1]) -[2]_\sqq \widehat{\R}(k_1,k_2\mid l+1;\star,[1,1])
			\\\notag
			&= [2]_\sqq^2 \widehat{X}_{k-1|l} *[K_\de]^{k_1}*[K_\co] + [2]_\sqq \widehat{Y}_{k+1|l} *[K_\de]^{k_1}* [K_\co]
			\\\notag
			&\quad +  \big( -[\haTh_{\star, k_2-k_1+1}, \haB_{[1,1],l}]_{\sqq^{-2}} *[K_\de]^{k_1} +\sqq^{-2} [\haTh_{\star, k_2-k_1-1}, \haB_{[1,1],l}]_{\sqq^{-2}} *[K_\de]^{k_1+1} \big) *[K_{\co}]
			\\\notag
			&\quad + \{k_1 \leftrightarrow k_2\}.% + \{k_1 \leftrightarrow k_2\}.
		\end{align}

		With the help of Proposition \ref{prop:coBil} and \eqref{eq:HaDr3}, by using the same argument as in \cite[Lemma 4.13]{LW20b},  we have
		\begin{align}
			\label{eq:SSSS}
			& \widehat{\SS}(k_1,k_2+1 \mid l;\star,[1,1]) + \widehat{\SS}(k_1+1,k_2\mid l;\star,[1,1]) -[2]_\sqq \widehat{\SS}(k_1,k_2\mid l+1;\star,[1,1])
			\\\notag
			&=(1-q)^2\Big( -\big[\haTh_{\star, k_2-k_1+1}, \haB_{[1,1],l}\big]_{\sqq^{-2}}* [K_\de]^{k_1} +\sqq^{-2} \big[\haTh_{\star, k_2-k_1-1}, \haB_{[1,1],l}\big]_{\sqq^{-2}} *[K_\de]^{k_1+1} \Big) \\\notag&\quad*[K_\co]
			+ \{k_1 \leftrightarrow k_2\}.
		\end{align}
		
		From Proposition \ref{prop:SerreO11l}, we have
		\begin{align}
			\label{eq:SSS}
			&\widehat{\SS}(k_1,k_2+1 \mid l;\star,[1,1]) + \widehat{\SS}(k_1+1,k_2\mid l;\star,[1,1]) -[2]_\sqq \widehat{\SS}(k_1,k_2\mid l+1;\star,[1,1])
			\\\notag
			&=(1-q)^2\Big( \widehat{\R}(k_1,k_2+1 \mid l;\star,[1,1]) + \widehat{\R}(k_1+1,k_2\mid l;\star,[1,1]) -[2]_\sqq \widehat{\R}(k_1,k_2\mid l+1;\star,[1,1])\Big).
		\end{align}
		
		By combining \eqref{RRR}, \eqref{eq:SSSS} and \eqref{eq:SSS}, it follows that
		$\widehat{Y}_{k+1|l}=0$
		since $\widehat{X}_{k-1 |l}=0$ by the inductive assumption. The proof is completed.
	\end{proof}

	\begin{proposition}
		Relation \eqref{eq:HHstari1} holds in $\iH(\X_\bfk)$ for any $m,r>0$.
	\end{proposition}
	
	\begin{proof}
		With the help of \eqref{humstarbvl} proved in Proposition \ref{prop:humstarbvl} and
		\eqref{Bmuk Bmul} for $\mu=[1,1]$, we can apply the same argument as in \cite[Proposition 4.8]{LW20b} to prove that \eqref{eq:HHstari1} holds in $\iH(\X_\bfk)$ for any $m,r>0$.
	\end{proof}
	
	%%%%%%%%%%%%
	%%%%%%%%%%%%
	\section{Relations between star point and $[i,1]$, II}
	\label{sec:Relationsstari1 II}
	
	In this section, we shall prove the following formulas in $\iH(\X_\bfk)$:
	\begin{align}
		\label{humbvlstar}
		&\big[\widehat{H}_{[i,1],m},[\co(l\vec{c})]\big]=-\frac{[m]_\sqq}{m} [\co((l+m)\vec{c})]+\frac{[m]_\sqq}{m} [\co((l-m)\vec{c})]*[K_{m\de}],
		\\
		\label{Serrei1star}
		&\widehat{\SS}(k_1,k_2\mid l;[i,1],\star)=(1-q)^2\widehat{\R}(k_1,k_2\mid l; [i,1],\star),
	\end{align}
	for any $k_1,k_2,l\in\Z$, $m>0$, $1\leq i\leq \bt$.
	%%%%%%%%%%%%%%%

	\subsection{Real and imaginary roots at $[i,1]$}
	\label{sec:rooti1}
	%%%%%%%%%%

	In this subsection we give explicit descriptions of the real and imaginary roots at $[i,1]$ in $\iH(\X_\bfk)$.
	
	Without loss of generality, we assume $i=1$ in this whole section.

	%%%%%%%%%%%%%%%%%%
	\subsubsection{Real roots at $[1,1]$}
	\label{subsec:realroot}
	% through this section. %

	For any $r>0$, denote by
	\begin{align}
		\label{eq:realroots11set}
		\begin{split}
			&\cm_{r\delta+\alpha_{11}}=\{[S_{1,1}^{(bp_1+1)}\oplus S_{1,0}^{(\nu)}]\mid b\geq 0, \;|\nu|+b=r\};\\
			&\cm_{r\delta-\alpha_{11}}=\{[S_{1,0}^{(ap_1-1)}\oplus S_{1,0}^{(\nu)}]\mid a\geq 1, \;|\nu|+a=r\};
		\end{split}
	\end{align}
	and set
	\begin{align}
		\label{eq:Mrde+alpha11}
		[M_{r\delta\pm\alpha_{11}}]:=\sum\limits_{[M]\in\cm_{r\delta\pm\alpha_{11}}}\bn(\ell(M)-1)[\![M]\!],
	\end{align}
	where
	\[[\![M]\!]=\frac{[M]}{|\Aut(M)|},\qquad\qquad\bn(l)=\prod_{i=1}^l(1-\sqq^{2i}),\,\, \forall l\geq1.\]

	\begin{proposition}
		\label{prop:realroot}
		For any $r>0$, we have
		\begin{align}
			\label{realroot}
			\haB_{[1,1],r}= (q-1)[M_{r\de+\alpha_{11}}], \qquad \haB_{[1,1],-r}=(1-q) [M_{r\de-\alpha_{11}}]*[K_{-r\de+\alpha_{11}}].
		\end{align}
	\end{proposition}

	The proof of Proposition \ref{prop:realroot} is provided in  Appendix \ref{subsec:proofrealroot}.
	
	%%%%%%%%%%
	
	%%%%%%%%%%%%
	\subsubsection{Imaginary roots at $[1,1]$}
	\label{subsec:imageroot}

	For any $r>0$, set
	\begin{align}
		\label{def:Mr}
		\cm_{i,r\de}:=&\big\{[S_{i,1}^{(\nua  p_1)}\oplus S_{i,0}^{(\nu)}]\mid a\geq 1, a+|\nu|=r\big\}
		\\\notag
		&\bigcup\big\{[S_{i,0}^{(\nua  p_i-1)}\oplus S_{i,1}^{(b  p_i+1)}\oplus S_{i,0}^{(\nu)}]\mid  a+b+|\nu|=r\big\}.
	\end{align}
	Since we focus on $i=1$ in the following, we also denote
	$$\cm_{r\de}:=\cm_{1,r\de}.$$
	\begin{proposition}
		\label{prop:imageroot}
		For any $r>0$, we have
		\begin{align}
			\label{eq:TH11r}
			\haTh_{[1,1],r}=
			&\frac{\sqq}{q-1} \sum\limits_{|\lambda|=r}\bn(\ell(\lambda))  [\![S_{1,0}^{(\lambda)}]\!]
			+\sqq^{-1} \sum\limits_{[M]\in\cm_{r\de}}\bn(\ell(M)-1) [\![M]\!].
		\end{align}
	\end{proposition}
	
	The proof of Proposition \ref{prop:imageroot} is provided in  Appendix \ref{subsec:proofimroot}.

	\subsection{The relation \eqref{humbvlstar}}
	
	%\begin{align}
	%\label{iDR2-reform}
	% [ & \haTh_{[1,1],r},[\co]]+[\haTh_{[1,1],r-2},[\co]]\K_\de
	% \\
	%&= v^{-1}[\haTh_{[1,1],r-1},[\co(\vec{c})]]_{v^2}+v[\haTh_{[1,1],r-1},[\co(-\vec{c})]]_{v^{-2}}\K_\de, \quad \text{ for } r>0.
	%\notag
	%\end{align}
	In this subsection, we shall prove \eqref{humbvlstar}. Note that twisting with $\vec{c}$ induces an  automorphism on $\iH(\X_\bfk)$. By using Lemma \ref{lem:equiv}~(2), it is equivalent to prove that for any $r>0$,
	\begin{align}
		\label{humbvlstar-reform}
		&\big[ \haTh_{[1,1],r},[\co]\big]+\big[\haTh_{[1,1],r-2},[\co]\big]*[K_\de]
		\\
		&= \sqq^{-1}\big[\haTh_{[1,1],r-1},[\co(\vec{c})]\big]_{\sqq^2}+\sqq\big[\haTh_{[1,1],r-1},[\co(-\vec{c})]\big]_{\sqq^{-2}}*[K_\de].
		\notag
	\end{align}

	%Recall $\cm_{r\de}$ defined in \eqref{def:Mr}.
	The proof of \eqref{humbvlstar-reform} needs the following Lemmas \ref{Theta 11 times co} and \ref{co times Theta 11}, whose proofs are technical, and we put them in Appendix \ref{subsec:proofTH11timesO}--\ref{subsec:proofOtimesTH11}.
	\begin{lemma}
		\label{Theta 11 times co}
		For any $r>0$, we have
		\begin{align*}
			\haTh_{[1,1],r}\ast [\![\co]\!]=&\frac{\sqq^{r+1}}{q-1}\sum\limits_{|\mu|=r}\bn(\ell(\mu))[\![\co\oplus S_{1,0}^{(\mu)} ]\!]+\sqq^{r-1}\sum\limits_{[N]\in\cm_{r\de}}\bn(\ell(N)-1)[\![\co\oplus N]\!]\\&
			-\sum\limits_{k>0}
			\sum\limits_{|\mu|=r-k}\sqq^{r-2k+1}\cdot \bn(\ell(\mu))[\![\co(k\vec{c})\oplus S_{1,0}^{(\mu)} ]\!]\\
			&+
			\sum\limits_{k>0}\sum\limits_{[N]\in\cm_{(r-k)\de}} \sqq^{r-2k-1}(1-q)\cdot \bn(\ell(N)-1)[\![\co(k\vec{c})\oplus N]\!].
		\end{align*}
	\end{lemma}

	\begin{lemma}
		\label{co times Theta 11}
		For any $r>0$, we have
		\begin{align*}
			[\![\co]\!]\ast \haTh_{[1,1],r}=&\frac{\sqq^{r+1}}{q-1}\sum\limits_{|\mu|={r}}   \bn(\ell(\mu))[\![\co\oplus S_{1,0}^{(\mu)}]\!]+\sqq^{r-1}\sum\limits_{[N]\in\cm_{r\de}}  \bn(\ell(N)-1)[\![\co\oplus N]\!]\\
			&- \sum\limits_{k>0}\sum\limits_{|\mu|={r-k}}  \sqq^{r-2k+1}\cdot \bn(\ell(\mu))[\![\co(-k\vec{c})\oplus S_{1,0}^{(\mu)}]\!]\ast [K_{k\delta}]
			\\
			&+\sum\limits_{k> 0}\sum\limits_{[N]\in\cm_{(r-k)\de}}  \sqq^{r-2k-1}(1-q)\cdot {\bf n}(\ell(N)-1)[\![\co(-k\vec{c})\oplus N]\!]\ast [K_{k\delta}].
		\end{align*}
	\end{lemma}

	\begin{proposition}
		\label{prop:humbvlstar-reform}
		Relation \eqref{humbvlstar-reform} holds in $\iH(\X_\bfk)$ for any $r>0$.
	\end{proposition}
	
	\begin{proof}
		By using Lemmas \ref{Theta 11 times co} and \ref{co times Theta 11}, (\ref{humbvlstar-reform}) follows by easy (but tedious) cancellations, which we omit here.
	\end{proof}

	\subsection{The relation \eqref{Serrei1star}}
	%%%%%%%%%%%%
	
	In this subsection, we shall prove the relation \eqref{Serrei1star}. First, let us give a lemma.
	
	\begin{lemma}
		\label{lem:serrel=0}
		The following formula holds in $\iH(\X_\bfk)$ for any $l\in\Z$:
		\begin{align*}
			%&%[S_{1,1}]*[S_{1,1}]*[\co(l\vec{c})]-[2]_\sqq [S_{1,1}]*[\co(l\vec{c})]*[S_{1,1}]+[\co(l\vec{c})]*[S_{1,1}] *[S_{1,1}]
			%\\
			\widehat{\SS}(0,0 \mid l;[1,1],\star)=(1-q)^2\widehat{\R}(0,0 \mid l;[1,1],\star).
		\end{align*}
	\end{lemma}
	
	\begin{proof}
		Using \eqref{Si1 and co double direction}, we obtain
		\begin{align*}
			&\widehat{\SS}(0,0\mid l;[1,1],\star)
			\\
			&=[S_{1,1}]*[S_{1,1}]*[ \co(l\vec{c})]-[2]_\sqq[S_{1,1}]*[ \co(l\vec{c})]*[S_{1,1}]+[\co(l\vec{c})]*[S_{1,1}]*[S_{1,1}]\\
			&=[S_{1,1}]*\big(\sqq^{-1}[\co(l\vec{c})]*[S_{1,1}]+(\sqq-\sqq^{-1}) [\co(l\vec{c}+\vec{x}_1)]\big)-[2]_\sqq[S_{1,1}]*[ \co(l\vec{c})]*[S_{1,1}]
			\\
			&\quad +\big(\sqq[S_{1,1}]*[\co(l\vec{c})] -(q-1)[\co(l\vec{c}+\vec{x}_1)]  \big)*[S_{1,1}]
			\\
			&=(\sqq-\sqq^{-1})\big[[S_{1,1}], [\co(l\vec{c}+\vec{x}_1)]\big]_\sqq
			\\
			&=-\sqq^{-1}(q-1)^2[\co(l\vec{c})]*[K_{\alpha_{11}}]
			\\
			&=(1-q)^2\widehat{\R}(0,0 \mid l;[1,1],\star).
		\end{align*}
	\end{proof}
	
	\begin{proposition}
		Relation \eqref{Serrei1star} holds in $\iH(\X_\bfk)$ for any $k_1,k_2,l\in\Z$.
	\end{proposition}

	\begin{proof}
		%First, we prove \eqref{humstarbvl}.
		We can assume  that $k_1\leq k_2$.
		%Let $k=k_2-k_1\geq0$.
		With the help of Proposition \ref{prop:coBil} and \eqref{Bmuk Bmul} for $\mu=[1,1]$, by using the same argument as in the proof of \cite[Lemmas 4.9 and 4.13]{LW20b},  we have
		\begin{align*}
			&\widehat{\SS}(k_1,k_2+1 \mid l;[1,1],\star) + \widehat{\SS}(k_1+1,k_2\mid l;[1,1],\star) -[2]_\sqq \widehat{\SS}(k_1+1,k_2+1 \mid l-1;[1,1],\star)
			\\
			&=  {(1-q)^2}\Big( -\big[[\co(l\vec{c})], \haTh_{[1,1], k_2-k_1+1}\big]_{\sqq^{-2}} *[K_\de]^{k_1} \\
			&\quad+\sqq^{-2} \big[[\co(l\vec{c})], \haTh_{[1,1], k_2-k_1-1}\big]_{\sqq^{-2}} *[K_\de]^{k_1+1} \Big)* [K_{\alpha_{11}}]+ \{k_1 \leftrightarrow k_2\},
		\end{align*}
		and
		\begin{align*}
			& \widehat{\SS}(k_1,k_2+1 \mid l;[1,1],\star) + \widehat{\SS}(k_1+1,k_2\mid l;[1,1],\star) -[2]_\sqq \widehat{\SS}(k_1,k_2\mid l+1;[1,1],\star)
			\\
			&= (1-q)^2\Big( -\big[\haTh_{[1,1], k_2-k_1+1}, [\co(l\vec{c})]\big]_{\sqq^{-2}} *[K_\de]^{k_1} \\
			&\quad+\sqq^{-2} \big[\haTh_{[1,1], k_2-k_1-1}, [\co(l\vec{c})]\big]_{\sqq^{-2}} *[K_\de]^{k_1+1} \Big) *[K_{\alpha_{11}}] + \{k_1 \leftrightarrow k_2\}.
		\end{align*}

		Denote for $m\in\Z$,
		\begin{align*}  %\label{Y}
			\widehat{Y}_{m|l} =&\big[[\co((l+1)\vec{c})], \haTh _{[1,1],m-1}\big]_{\sqq^{-2}} - \big[\haTh _{[1,1],m-1},[\co((l-1)\vec{c})]\big]_{\sqq^{-2}}*[K_\de]
			\\\notag
			&+\sqq^{-1} \big[\haTh_{[1,1],m},[\co(l\vec{c})]\big]+ \sqq^{-1} \big[\haTh_{[1,1],m-2},[\co(l\vec{c})]\big] *[K_\de].
		\end{align*}
		Clearly, $\widehat{Y}_{m|l}=0$ if $m\le 0$.
		Using  Proposition \ref{prop:humbvlstar-reform}, we have $\widehat{Y}_{m|l}=0$ for any $m$.
		
		We also denote
		for $n\in \Z$,
		\begin{align*}
			% \label{X}
			\widehat{X}_{n |l}=&\Big(\sum_{p\ge 0} \sqq^{2p+1}\big[[\co((l+1)\vec{c})],\haTh_{[1,1],n-2p-1}\big]_{\sqq^{-2}}
			+\sum_{p\ge 1} \sqq^{2p-1} [2]_\sqq \big[\haTh_{[1,1],n-2p},[\co(l\vec{c})]\big]
			\\
			&-\sum_{p\ge 0} \sqq^{2p+1}\big[\haTh_{[1,1],n-2p-1},[\co((l-1)\vec{c})]\big]_{\sqq^{-2}}*[K_\de]\Big)*[K_\de]^{p+1}
			+ \big[\haTh_{[1,1],n},[\co(l\vec{c})]\big]*[K_\de].
			\notag
		\end{align*}
		Clearly, $\widehat{X}_{n |l}=0$ for $n\leq 0$. By using the same argument as in \cite[Lemma 4.10]{LW20b} and Proposition
		\ref{prop:humbvlstar-reform}, we have
		$\widehat{X}_{n |l}=0$ for any $n$.
		
		Using the same argument as in \cite[Lemmas 4.11 and 4.14]{LW20b}, we have
		\begin{align*}
			&\widehat{\R}(k_1,k_2+1 \mid l;[1,1],\star) + \widehat{\R}(k_1+1,k_2\mid l;[1,1],\star)  -[2]_\sqq \widehat{\R}(k_1+1,k_2+1\mid l-1;[1,1],\star)
			\\
			&= - [2]_\sqq^2 \widehat{X}_{k_2-k_1| l-1} *[K_\de]^{k_1} *[K_{\alpha_{11}}]&
			\\&\quad- \Big( \big[[\co(l\vec{c})], \haTh_{[1,1], k_2-k_1+1}\big]_{\sqq^{-2}} * -\sqq^{-2} \big[[\co(l\vec{c})], \haTh_{[1,1], k_2-k_1-1}\big]_{\sqq^{-2}} *[K_\de] \Big)* [K_{k_1\delta+\alpha_{11}}]
			\\
			&\quad+ \{k_1 \leftrightarrow k_2\},\\
			& \widehat{\R}(k_1,k_2+1 \mid l;[1,1],\star) + \widehat{\R}(k_1+1,k_2\mid l;[1,1],\star)  -[2]_\sqq \widehat{\R}(k_1,k_2\mid l+1;[1,1],\star)
			\\
			&= [2]_\sqq^2 \widehat{X}_{k_2-k_1-1 |l} * [K_\de]^{k_1}*[K_{\alpha_{11}}] + [2]_\sqq \widehat{Y}_{k_2-k_1+1\mid l} [K_\de]^{k_1} *[K_{\alpha_{11}}]
			\\
			&\quad-\Big( \big[\haTh_{[1,1], k_2-k_1+1}, [\co(l\vec{c})]\big]_{\sqq^{-2}}  -\sqq^{-2} \big[\haTh_{[1,1], k_2-k_1-1}, [\co(l\vec{c})]\big]_{\sqq^{-2}} *[K_\de] \Big) *[K_{k_1\delta+\alpha_{11}}]\\
			&\quad + \{k_1 \leftrightarrow k_2\}.
		\end{align*}
		
		Combining the above computations, for any $k_1,k_2,l\in\Z$ we have
		\begin{align}
			\label{serreind1}
			&\widehat{\SS}(k_1,k_2+1 \mid l;[1,1],\star) + \widehat{\SS}(k_1+1,k_2\mid l;[1,1],\star) -[2]_\sqq \widehat{\SS}(k_1+1,k_2+1\mid l-1;[1,1],\star)
			\\\notag
			&=(1-q)^2\Big( \widehat{\R}(k_1,k_2+1 \mid l;[1,1],\star) + \widehat{\R}(k_1+1,k_2\mid l;[1,1],\star)
			\\\notag
			&\quad -[2]_\sqq \widehat{\R}(k_1+1,k_2+1\mid l-1;[1,1],\star)\Big),
			\\
			\label{serreind2}
			&\widehat{\SS}(k_1,k_2+1 \mid l;[1,1],\star) + \widehat{\SS}(k_1+1,k_2\mid l;[1,1],\star) -[2]_\sqq \widehat{\SS}(k_1,k_2\mid l+1;[1,1],\star)
			\\\notag
			&= (1-q)^2\Big(\widehat{\R}(k_1,k_2+1 \mid l;[1,1],\star) + \widehat{\R}(k_1+1,k_2\mid l;[1,1],\star)
			\\\notag
			& \quad-[2]_\sqq\widehat{\R}(k_1,k_2\mid l+1;[1,1],\star)\Big).
		\end{align}
		By Lemma \ref{lem:serrel=0}, we have
		\begin{align*}
			\widehat{\SS}(0,0 \mid l;[1,1],\star)=(1-q)^2\widehat{\R}(0,0 \mid l;[1,1],\star).
		\end{align*}
		By induction and using \eqref{serreind1}--\eqref{serreind2}, one can obtain that
		\begin{align*}
			\widehat{\SS}(k_1,k_2 \mid l;[1,1],\star)=(1-q)^2\widehat{\R}(k_1,k_2 \mid l;[1,1],\star)
		\end{align*}
		for any $k_1,k_2,l\in\Z$.
	\end{proof}

	%%%%%%%%%
	\section{Derived equivalence and two presentations of affine type A}
	\label{sec:derived}
	%%%%%%%%%%%%%%%
	
	Let $\X$ be a weighted projective line of domestic type. Let $Q$ be an acyclic quiver of the associated type. In this section, we establish a derived equivalence between $\cd^b(\cc_{1}(\coh(\X)))$ and $\cd^b(\cc_1(\rep_\bfk^{\rm nil}(Q^{op})))$, which induces an isomorphism of $\iH(\X_\bfk)$ and $\iH(\bfk Q^{op})$. Moreover, if $\X$ is of weight type $(p_1,p_2)$, then it induces the isomorphism of the $\imath$quantum group in the Serre  and Drinfeld type presentations. 
	
	%%%%%%%
	\subsection{A derived equivalence}

	Let $\X$ be a weighted projective line of domestic type. Then there exists a tilting sheaf $T$ such that $\End_\X(T)$ is hereditary of affine type ADE. Let $Q$ be the quiver such that $\End_\X(T)\cong \bfk Q$.
	Then there is a derived equivalence
	\begin{align}  \label{Bei}
		\RHom_{\X}(T,-):\cd^b(\coh(\X))\stackrel{\simeq}{\longrightarrow} \cd^b(\rep_\bfk^{\rm nil}( Q^{op})).
	\end{align}

	Using the proof of \cite[Proposition 5.9]{LRW20a}, one can obtain the following remark.
	\begin{remark}  \label{prop:tilt}
		Let $T$ be the tilting sheaf as above. Then $K_T=(T\oplus T, d)$ is a tilting object in $\cd^b(\calc_1(\coh(\X_\bfk)))$, which gives rise to an equivalence:
		$$\RHom (K_T,-): \cd^b(\calc_1(\coh(\X_\bfk))) \stackrel{\simeq}{\longrightarrow} \cd^b(\cc_1(\rep_\bfk^{\rm nil}(Q^{op}))).$$
	\end{remark}
	
	Let $\cu:=\{M\in\coh(\X)\mid \Ext^1_\X(T,M)=0\}$ and $\cv:=\{M\in\coh(\X)\mid \Hom_\X(T,M)=0\}$.
	%\begin{align*}
	%\cu:=\{M\in\coh(\X)\mid \Ext^1_\X(T,M)=0\},\qquad	\cv:=\{M\in\coh(\X)\mid \Hom_\X(T,M)=0\}.
	%\end{align*}
	Then $(\cu,\cv)$ is a {\em torsion pair} of $\coh(\X)$.
	Let $\cc_1(\cu)$ and $\cc_1(\cv)$ be the categories of all $1$-periodic complexes in $\cu$ and $\cv$ respectively.
	
	\begin{lemma}
		[\text{\cite[Lemma 5.8]{LRW20a}}]
		$(\cc_1(\cu),\cc_1(\cv))$ is a torsion pair of $\cc_1(\coh(\X))$. In particular, any $M\in\cc_1(\coh(\X))$ admits a short exact sequence of the form
		\begin{align}
			\label{T-resol}
			0 \longrightarrow M \longrightarrow X_M \longrightarrow T_M \longrightarrow 0,
		\end{align}
		where $X_M\in \cc_1(\cu)$ and $T_M\in\add K_T$.
	\end{lemma}

	\begin{proposition}
		\label{prop:F}
		Let $F= \Hom(K_T,-)$. Then there exists an isomorphism
		\begin{align*}
			\BF:\tMHX&\stackrel{\sim}{\longrightarrow} \iH(\bfk Q^{op})\\\notag
			[M]&\mapsto [F(T_M)]^{-1}* [F(X_M)],
		\end{align*}
		where $X_M\in \cc_1(\cu)$ and $T_M\in\add K_T$, are defined in the short exact sequence \eqref{T-resol}.
	\end{proposition}
	
	\begin{proof}
		The proof is the same as in \cite[Proposition 5.10]{LRW20a}, hence omitted here.
	\end{proof}
	%\begin{proof}
	%Follows from \cite[Theorem A. 22]{LW19a} with the help of Lemma~\ref{prop:tilt}.
	%\end{proof}

	%%
	\subsection{A commutative diagram}

	In this subsection, we always assume that $\X$ is of weight type $(p_1,p_2)$. Let
	$T=\bigoplus_{0\leq \vec{x}\leq \vec{c}} \co_\X(\vec{x})$ be the canonical tilting sheaf in $\coh(\X)$.
	Then $\End_\X(T)=\bfk Q$ is the canonical algebra \cite{Rin84} with $Q=C_{p_1,p_2}$ as shown below:

	\begin{center}\setlength{\unitlength}{0.75mm}
		\begin{equation}
			%	\begin{figure}
			\label{figure: quasi-split i}
			\begin{picture}(80,20)(0,20)
				\put(-20,20){$\circ$}
				\put(-17,22){\vector(2,1){17}}
				\put(-17,21){\vector(2,-1){17}}
				\put(-22,17.5){$\star$}
				\put(0,10){$\circ$}
				\put(0,30){$\circ$}
				\put(50,10){$\circ$}
				\put(50,30){$\circ$}
				\put(72,10){$\circ$}
				\put(72,30){$\circ$}
				\put(92,20){$\circ$}
				\put(-3,6){\tiny$[2,1]$}
				\put(-3,34){\tiny${[1,1]}$}
				\put(44,6){\tiny $[2,p_2-2]$}
				\put(44,34){\tiny $[1,p_1-2]$}
				\put(67,6){\tiny $[2,p_2-1]$}
				\put(67,34){\tiny $[1,p_1-1]$}
				\put(94.5,17.5){$c$.}
				
				\put(3,11.5){\vector(1,0){16}}
				\put(3,31.5){\vector(1,0){16}}
				\put(23,10){$\cdots$}
				\put(23,30){$\cdots$}
				\put(33.5,11.5){\vector(1,0){16}}
				\put(33.5,31.5){\vector(1,0){16}}
				\put(53,11.5){\vector(1,0){18.5}}
				\put(53,31.5){\vector(1,0){18.5}}
				
				\put(75,12){\vector(2,1){17}}
				\put(75,31){\vector(2,-1){17}}

			\end{picture}
			%
			%		
			%	\end{figure}
		\end{equation}
	\end{center}
	
	\vspace{1.1cm}
	
	Let $\tUi=\tUi(\widehat{\mathfrak{sl}}_{n})$ for $n=p_1+p_2$.
	Let $\widetilde{\psi}: \tUi_{ |v={\sqq}}\rightarrow \iH(\bfk C_{(p_1,p_2)}^{op})$ be the algebra embedding obtained in  Theorem \ref{thm:main}. The following result is a generalization of
	\cite[Theorem 5.11]{LRW20a} for $\tUi(\widehat{\mathfrak{sl}}_{2})$.
	
	\begin{theorem}
		\label{main thm2}
		We have the following commutative diagram of algebra homomorphisms
		%\[
		%\xymatrix{\tUi_{ |v={\sqq}} \ar[d]^{\widetilde{\psi}} \ar[r]^{\Phi^{-1}}   & {}^{\text{Dr}}\tUi_{ |v=\sqq} \ar[d]^{\Omega}
		%\\
		%\tMHL \ar[r]^{\BF^{-1}} & \tMHX  }\]
		\[
		\xymatrix{ \ar[r]^{\Phi}   {}^{\text{Dr}}\tUi_{ |v=\sqq} \ar[d]^{\Omega} & \tUi_{ |v={\sqq}} \ar[d]^{\widetilde{\psi}}
			\\
			\tMHX\ar[r]^-{\BF} &  \iH(\bfk C_{(p_1,p_2)}^{op}).  }\]
		%where $\Phi, \BF$ are isomorphisms and $\widetilde{\psi},\Omega$ are monomorphisms.
		In particular, the homomorphism $\Omega$ is injective.
	\end{theorem}

	\begin{proof}
		A direct computation shows that
		\begin{align*}
			\widetilde{\psi}\circ\Phi(B_{\star,0})=&\widetilde{\psi}(B_\star)=\frac{1}{1-q}[S_\star],
			\\
			\BF\circ \Omega(B_{\star,0})=&\frac{1}{1-q}\BF([\co])=\frac{1}{1-q}[S_\star],
		\end{align*}
		and then $\widetilde{\psi}\circ\Phi(B_{\star,0})=\BF\circ \Omega(B_{\star,0})$.

		Denote by $M(l\de-\alpha_\star)$ the indecomposable $\bfk C_{(p_1,p_2)}^{op}$-module with dimension vector $l\de-\alpha_\star$. Similar to \eqref{eq:M-de}, we have
		\begin{align*}
			\Gamma_{\omega_\star}([S_\star])=[M(\de-\alpha_\star)]*[K_{\de-\alpha_\star}]^{-1},
		\end{align*}
		and then
		\begin{align*}
			\widetilde{\psi}\circ\Phi(B_{\star,-1})
			= \widetilde{\psi}( \TT_{\omega_\star}(B_\star))
			=\frac{-1}{q-1}\Gamma_{\omega_\star}([S_\star])
			=\frac{-1}{q-1}[M(\de-\alpha_\star)]*[K_{\de-\alpha_\star}]^{-1}.
		\end{align*}
		On the other hand, there exists a short exact sequence in $\coh(\X)$
		\begin{align*}
			0\longrightarrow \co(-\vec{c}) \longrightarrow \co\oplus \co\stackrel{(g_1,g_2)}{\longrightarrow} \co(\vec{c})\longrightarrow0,
		\end{align*}
		and then a short exact sequence in $\cc_1(\coh(\X))$
		\begin{align*}
			0\longrightarrow \co(-\vec{c})\longrightarrow N^\bullet:=(\co\oplus \co\oplus \co(\vec{c}),d)\longrightarrow K_{\co(\vec{c})}\longrightarrow0,
		\end{align*}
		where $d= \begin{bmatrix} 0&0&0\\ 0&0&0\\g_1&g_2&0 \end{bmatrix}$. Note that
		$$\BF([N^\bullet])=[(S_\star\oplus S_\star\oplus P_c,F(d))]= [M(\de-\alpha_\star)]*[K_{2\alpha_\star}].$$
		Here $P_c$ is the projective $\bfk C_{(p_1,p_2)}^{op}$-module corresponding to the vertex $c$. Then
		\begin{align*}
			\BF([\co(-\vec{c})])= &\BF([N^\bullet])* [F(K_{\co(\vec{c})})]^{-1}
			\\
			=&[M(\de-\alpha_\star)]*[K_{2\alpha_\star}]*[K_{P_c}]^{-1}
			\\
			=&[M(\de-\alpha_\star)]*[K_{ \de-\alpha_\star}]^{-1}
		\end{align*}
		by Proposition \ref{prop:F}.
		Therefore,
		\begin{align*}
			\BF\circ \Omega(B_{\star,-1})=\frac{-1}{q-1}\BF([\co(-\vec{c})])= \frac{-1}{q-1}[M(\de-\alpha_\star)]*[K_{ \de-\alpha_\star}]^{-1}= \widetilde{\psi}\circ\Phi(B_{\star,-1}).
		\end{align*}
		
		We consider $S_{ij}\in\coh(X)$ with $j\neq0$ in the following. Then we have
		a short exact sequence
		$$0\longrightarrow \co((j-1)\vec{x}_i)\stackrel{x_i}{\longrightarrow} \co(j\vec{x}_i)\longrightarrow S_{ij}\rightarrow0.$$
		Similarly to the above, we have
		\begin{align*}
			\BF([S_{ij}])=[S_{[i,j]}],
		\end{align*}
		where $S_{[i,j]}$ is the simple $\bfk C_{p_1,p_2}^{op}$-module corresponding to $[i,j]$.
		Then
		\begin{align*}
			\BF\circ \Omega(B_{[i,j],0})=\frac{-1}{q-1}\BF([S_{ij}])= \frac{-1}{q-1}[S_{[i,j]}]= \widetilde{\psi}(B_{[i,j]})= \widetilde{\psi}\circ\Phi(B_{[i,j],0}).
		\end{align*}

		Finally, we verify that
		\begin{align*}
			\BF\circ\Omega(\K_{[i,j]})=&\BF([K_{S_{ij}}])=  [K_{S_{[i,j]}}]= \widetilde{\psi}(\K_{[i,j]})= \widetilde{\psi}\circ\Phi(\K_{[i,j]}),
			\\
			\BF\circ\Omega(\K_{\star})=&\BF([K_{\co}])=  [K_{S_{\star}}]= \widetilde{\psi}(\K_{\star})= \widetilde{\psi}\circ\Phi(\K_{\star}),
			\\
			\BF\circ\Omega(C)=& \BF([K_{\co(\vec{c})}]*[K_{\co}]^{-1})=[K_{P_c}]*[K_{S_\star}]^{-1}=[K_\de]=\widetilde{\psi}\circ\Phi(C).
		\end{align*}
		Summarizing, by using Lemma \ref{lem:reduced generators}, we have proved $\BF\circ\Omega=  \widetilde{\psi}\circ\Phi$.
		
		The injectivity of $\Omega$ follows by the injectivity of $\widetilde{\psi}$ and the commutative diagram.
	\end{proof}

	\begin{remark}
		For $\fg$ of type DE and $Q$ the associated affine quiver, we also have a derived equivalence $\cd^b(\coh(\X_\bfk))\simeq \cd^b(\rep_\bfk (Q^{op}))$, which induces a derived equivalence on the categories of $1$-periodic complexes. Similar to Proposition \ref{prop:F}, there is an algebra isomorphism $\BF:\tMHX\stackrel{\cong}{\rightarrow} \iH(\bfk Q^{op})$. However,  the morphism $\Phi$ can not be the one induced by $\BF$, in fact it can not be induced by any derived equivalent functors by the same argument as in \cite[\S9.2]{DJX12}.
	\end{remark}

	\begin{remark}
		For the commutative diagram obtained in Theorem \ref{main thm2}, the corresponding result for quantum groups is only known for  quantum affine $\sll_2$ \cite{BS12} by using the approach of Drinfeld double Ringel-Hall algebras. In fact,  Burban-Schiffmann \cite{BS12} interpreted the Drinfeld-Beck isomorphism of quantum affine $\sll_2$ via Beilinsion's derived equivalence between the Kronecker quiver and the projective line, which gives an isomorphism of their Drinfeld double Ringel-Hall algebras \cite{Cr10}.
		
		We also expect the corresponding result of Theorem \ref{main thm2} holds for quantum groups of affine type $A$, and it will be considered elsewhere by using semi-derived Ringel-Hall algebras.
	\end{remark}
	
	\appendix

	%%%%%

	%%%%%%%%%

	%%%%%%%%%%
	\section{Proofs of results in Section \ref{sec:Relationsstari1 I} }
	%%%%%%%%%%%%%%

	\subsection{Proof of Lemma \ref{lem:SerreO11}}
	\label{proofSerre011}

	We first introduce the following notations.
	Denote by $\mathfrak{S}(a,b)$ the set of pairs of coprime homogeneous polynomials $(J, L)\in\mathbf{k}[X,Y]$ of degree $a$ and $b$ respectively, and denote by $\mathfrak{T}(a,b)$ its subset consisting of those $(J,L)$'s such that {$X$} does not divide $J$. The cardinalities of $\mathfrak{S}(a,b)$ and $\mathfrak{T}(a,b)$ are denoted by $\mathfrak{s}(a,b)$ and $\mathfrak{t}(a,b)$ respectively.
	
	\begin{lemma}
		\label{the formula of psi(a,b)}
		For any $a,b\geq0$,
		$$\mathfrak{t}(a,b)=\begin{cases}
			(q-1)(q^{b+1}-1),  &\text{ if }a=0,\\
			(q-1)^2q^{a+b},  &\text{ if } a\neq 0.
		\end{cases}$$
		Consequently, $\mathfrak{t}(a,b)=\mathfrak{t}(a-k,b+k)$ for $1\leq k\leq a-1$.
	\end{lemma}

	\begin{proof}
		Recall from \cite[Lemma 9]{BKa01} that
		$$\mathfrak{s}(a,b)=\begin{cases}
			(q-1)(q^{a+b+1}-1),  &\text{ if }ab=0,\\
			(q-1)(q^2-1)q^{a+b-1},  &\text{ if } ab\neq 0.
		\end{cases}$$
		By definition of $\mathfrak{S}(a,b)$ and $\mathfrak{T}(a,b)$, we know that $\mathfrak{t}(0,b)=\mathfrak{s}(0,b)=(q-1)(q^{b+1}-1)$ for $b\geq 1$, and
		$\mathfrak{t}(a,0)=\mathfrak{s}(a,0)-\mathfrak{s}(a-1,0)=(q-1)(q^{a+1}-q^{a})=(q-1)^2q^{a}$ for $a\geq 1$.
		Moreover, for any $a,b\geq 1$, according to \cite[(6.10)]{Sch04}, we have
		$$\mathfrak{t}(a,b)=\mathfrak{s}(a,b)-\mathfrak{t}(b,a-1).$$
		%Let us prove that $\mathfrak{t}(a,b)=(q-1)^2q^{a+b}$ by induction on $a+b$.
		%By inductive assumption,
		Then by induction on $a+b$, we get $$\mathfrak{t}(a,b)
		%=\mathfrak{s}(a,b)-\mathfrak{t}(b,a-1)
		=(q-1)(q^2-1)q^{a+b-1}-(q-1)^2q^{a+b-1}=(q-1)^2q^{a+b}.$$
		The proof is completed.
	\end{proof}

	\begin{corollary}
		\label{extension between line bundles}
		For any $l\geq1$, we have
		\begin{align}
			\label{co(l) and x1}
			&[\co(l\vec{c})]* [\co(\vec{x}_1)]-\sqq^{-1}[\co(l\vec{c}+\vec{x}_1)]*[\co]
			\\\notag&=
			\sqq^{-l-1}\Big([\co(l\vec{c})\oplus\co(\vec{x}_1)]-[\co(l\vec{c}+\vec{x}_1)\oplus\co]\Big),
			\\\label{co(l+x1) and c}
			&[\co(l\vec{c}+\vec{x}_1)]*[\co(\vec{c})]-\sqq^{-1}[\co((l+1)\vec{c})]* [\co(\vec{x}_1)]
			\\\notag
			&=
			\sqq^{-l-1}\Big([\co(l\vec{c}+\vec{x}_1)\oplus\co(\vec{c})]-[\co((l+1)\vec{c})\oplus\co(\vec{x}_1)]\Big).
		\end{align}
	\end{corollary}
	
	\begin{proof}
		For any $1\leq k\leq l-1$, set $M_k:=\co((l-k)\vec{c})\oplus \co(k\vec{c}+\vec{x}_1)$. According to \cite[Lemma 6.5]{Sch04}, we have $$F_{\co(l\vec{c}), \co(\vec{x}_1)}^{M_k}=\frac{\mathfrak{t}(k,l-k-1)}{q-1}, \quad F_{\co(l\vec{c}+\vec{x}_1), \co}^{M_k}=\frac{\mathfrak{t}(l-k,k)}{q-1}.$$
		Then it follows from Riedtmann-Peng formula that
		\begin{align*}
			|\Ext^1(\co(l\vec{c}), \co(\vec{x}_1))_{M_k}|&=\frac{(q-1)\cdot\mathfrak{t}(k,l-k-1)}{|\Aut(M_k)|},\\
			|\Ext^1(\co(l\vec{c}+\vec{x}_1), \co)_{M_k}|&=\frac{(q-1)\cdot\mathfrak{t}(l-k,k)}{|\Aut(M_k)|}.
		\end{align*}
		%Comparing with \cite[Lemma 6.5]{Sch04}, we have
		So
		\begin{align}
			\label{eq:OOx1}
			[\co(l\vec{c})]* [\co(\vec{x}_1)]=&
			\sqq^{1-l}\Big([\co(l\vec{c})\oplus\co(\vec{x}_1)]+\sum\limits_{k=1}^{l-1}\frac{(q-1)\cdot\mathfrak{t}(k,l-k-1)}{|\Aut(M_k)|}[M_k]\Big),
			\\
			\label{eq:Ox1O}
			[\co(l\vec{c}+\vec{x}_1)]* [\co]=&
			\sqq^{-l}\Big([\co(l\vec{c}+\vec{x}_1)\oplus\co]+(q-1)[\co(l\vec{c})\oplus\co(\vec{x}_1)]
			\\\notag
			&\quad+
			\sum\limits_{k=1}^{l-1}\frac{(q-1)\cdot\mathfrak{t}(k,l-k)}{|\Aut(M_k)|}[M_k]\Big).
		\end{align}
		By Lemma \ref{the formula of psi(a,b)}, we have $\mathfrak{t}(k,l-k)=q\cdot\mathfrak{t}(k,l-k-1)$ for $1\leq k\leq l-1$,
		and then (\ref{co(l) and x1}) follows by combining \eqref{eq:OOx1}--\eqref{eq:Ox1O}.
		
		For (\ref{co(l+x1) and c}), we define
		\begin{align}
			\label{def:S}
			\mathcal{S}:=\{S_{1,j}\mid 2\leq j\leq p_1-1\}\cup \{S_{ij}\mid 2\leq i\leq \bt, 1\leq j\leq p_i-1\}.
		\end{align}
		%$\mathcal{S}$ be the subset of exceptional simple sheaves in $\coh(\X_{\bfk})$ different from $S_{1,1}$ and $S_{i,0}$ for any $1\leq i\leq \bt$.
		Then the left perpendicular category ${}^\perp\mathcal{S}$ is equivalent to {$\coh(\mathbb{Y}_{\bfk})$}, where $\mathbb{Y}_{\bfk}$ is a weighted projective line of weight type $\bp=(2,1)$. By similar arguments as in \S\ref{Embedding from projective line to weighted projective line}, we have an embedding
		\begin{equation}
			\label{the embedding functor from Y to X}
			F_{\mathbb{X},\mathbb{Y}}:\iH(\mathbb{Y}_{\bfk}) \longrightarrow \iH(\X_\bfk).
		\end{equation}
		Observe that each term in (\ref{co(l) and x1}) and (\ref{co(l+x1) and c}) belongs to the image of $F_{\mathbb{X},\mathbb{Y}}$, hence we can restrict to $\mathbb{Y}_{\bfk}$. In this case, $2\vec{x}_1=\vec{c}$. Then (\ref{co(l+x1) and c}) follows from (\ref{co(l) and x1}) by twisting with $\vec{x}_1$.
		%The proof of  the second formula is similar, hence omitted here.
	\end{proof}
	
	For any $m\geq 0$,  define
	\begin{align}
		\label{eq:THstarm+}
		\haTh_{\star,m}^+:=&\frac{1}{(q-1)^2\sqq^{m-1}}\sum_{0\neq f:\co(s\vec{c})\rightarrow \co((m+s)\vec{c}+\vec{x}_1) } [\coker(f)],
		%=\frac{1}{(q-1)\sqq^{m-1}}\sum\limits_{||\bn||=m}[S_\bn],
		\\
		\label{eq:THstarm-}
		\haTh_{\star,m}^-:=&\frac{1}{(q-1)^2\sqq^{m-1}}\sum_{0\neq f:\co(s\vec{c}+\vec{x}_1)\rightarrow \co((m+s)\vec{c}) } [\coker(f)],
		%=\frac{1}{(q-1)\sqq^{m-1}}\sum\limits_{||\bn||=m}[S_\bn],
	\end{align}
	which are independent of $s\in\Z$.
	%Observe that $\haTh_{\red{\star},0}^+=\frac{1}{\sqq-\sqq^{-1}}[S_{1,1}]$,
	For convenience we set $\haTh_{{\star},m}^\pm=0$ for $m<0$.

	%\section{}
	%\label{sec:Append A}

	\begin{lemma}
		For any $l\geq 1$,
		\begin{align}
			\label{eq:commcolcox1}
			&\big[ [\co(l\vec{c})], [\co(\vec{x}_1)] \big]_{\sqq^{-1}} +\big[ [\co],[\co(l\vec{c}+\vec{x}_1)]\big]_{\sqq^{-1}}
			\\\notag
			&=
			(\sqq-\sqq^{-1})^2 \big( \haTh_{\star,l}^+*[K_\co]-\haTh_{\star,l}^-*[K_{\co(\vec{x}_1)}]\big),	
			\\
			\label{eq:commcolx1co}
			&\big[ [\co(l\vec{c}+\vec{x}_1)], [\co(\vec{c})] \big]_{\sqq^{-1}} +\big[ [\co(\vec{x}_1)],[\co((l+1)\vec{c})]\big]_{\sqq^{-1}}
			\\\notag
			&=
			(\sqq-\sqq^{-1})^2 \big(\sqq \haTh_{\star,l+1}^-*[K_{\co(\vec{x}_1)}]-\sqq^{-1}\haTh_{\star,l-1}^+*[K_{\co(\vec{c})}]\big),
		\end{align}
	\end{lemma}
	
	\begin{proof}
		Observe that
		\begin{align*}
			&[\co(\vec{x}_1)]*[\co(l\vec{c})]\\
			&=\sqq^{-l}[\co(\vec{x}_1)\oplus \co(l\vec{c})]
			+\sqq^{-l}\sum_{0\neq f:\co(\vec{x}_1)\rightarrow \co(l\vec{c})} [\coker(f)]*[K_{\co(\vec{x}_1)}]\\
			&=\sqq^{-l}[\co(\vec{x}_1)\oplus \co(l\vec{c})]+\sqq^{-1}(q-1)^2\haTh_{\star,l}^-*[K_{\co(\vec{x}_1)}],
		\end{align*}
		and
		\begin{align*}
			[\co]*[\co(l\vec{c}+\vec{x}_1)]
			=&\sqq^{-l-1}[\co\oplus \co(l\vec{c}+\vec{x}_1)]+\sqq^{-l-1}\sum_{0\neq g:\co\rightarrow \co(l\vec{c}+\vec{x}_1)} [\coker (g)]*[K_{\co}]\\
			=&\sqq^{-l-1}[\co\oplus \co(l\vec{c}+\vec{x}_1)]+\sqq^{-2}(q-1)^2\haTh_{\star,l}^+*[K_{\co}].
		\end{align*}
		
		Moreover, by \eqref{co(l) and x1} we have
		\begin{align*}
			&\sqq^{-1} [\co(l\vec{c}+\vec{x}_1)]*[\co]-[\co(l\vec{c})]* [\co(\vec{x}_1)]\\
			&=\sqq^{-l-1}\big([\co(l\vec{c}+\vec{x}_1)\oplus\co]
			-[\co(l\vec{c})\oplus\co(\vec{x}_1)]\big).
		\end{align*}
		Combining the above formulas, \eqref{eq:commcolcox1} follows.
		
		For \eqref{eq:commcolx1co}, similar to the proof of \eqref{co(l+x1) and c}, we can restrict to consider $\X_\bfk$ of weight type $\bp=(2,1)$. In this case, by definition we have
		\begin{align*}
			\haTh_{\star,l}^+(\vec{x}_1)=\sqq\haTh_{\star,l+1}^-,\qquad \haTh_{\star,l}^-(\vec{x}_1)
			=\sqq^{-1} \haTh_{\star,l-1}^+.
		\end{align*}
		Then \eqref{eq:commcolx1co} follows from \eqref{eq:commcolcox1} by twisting with  $\vec{x}_1$.
	\end{proof}
	
	\begin{lemma}
		\label{lem:A1}
		For any $m\in\Z$,
		\begin{align}\label{Theta+}
			\haTh_{\star,m}^+=&
			q\haTh_{\star, m-2}^+*[K_{\delta}]+\frac{q}{q-1}\big[[S_{1,1}],\haTh_{\star, m}\big]_{\sqq^{-2}}-\frac{q\sqq}{q-1}\big[\haTh_{\star, m-1}, [S_{1,0}^{(p_1-1)}]\big]_{\sqq^{-2}}*[K_{\alpha_{11}}],
			\\
			\label{Theta-}
			\haTh_{\star,m}^-=&
			q\haTh_{\star, m-2}^-*[K_{\delta}]+\frac{\sqq}{q-1}\big[\haTh_{\star, m-1}, [S_{1,0}^{(p_1-1)}]\big]_{\sqq^{-2}}-\frac{q}{q-1}\big[[S_{1,1}],\haTh_{\star, m-2}\big]_{\sqq^{-2}}*[K_{\de-\alpha_{11}}].
		\end{align}
	\end{lemma}
	
	\begin{proof}
		For $m\leq0$, the above equalities hold obviously by definition. We assume $m>0$ in the following.
		
		For any $r\geq 0$, an easy computation shows that
		\begin{align}
			\label{S10S10p}
			[S_{1,0}^{(p_1-1)}]*[S_{1,0}^{(rp_1)}]
			=&[S_{1,0}^{(rp_1)}\oplus S_{1,0}^{(p_1-1)}],\\
			\label{S10pS10}
			[S_{1,0}^{(rp_1)}]*[S_{1,0}^{(p_1-1)}]=&\frac{1}{q}\Big([S_{1,0}^{(rp_1)}\oplus S_{1,0}^{(p_1-1)}]+(q-1)[S_{1,0}^{((r+1)p_1-1)}]\\
			&+q(q-1)[S_{1,1}^{((r-1)p_1+1)}]*[K_{\de-\alpha_{11}}]\Big).\notag
		\end{align}
		Here and in the following, it is understood that the terms involving $S_{1,1}^{((r-1)p_1+1)}$ do not exist if $r=0$. %\blue{Fix in the Notation?}}
		Hence it follows by combining \eqref{S10S10p} and \eqref{S10pS10} that
		\begin{align}
			\label{la}\big[[S_{1,0}^{(rp_1)}],[S_{1,0}^{(p_1-1)}]\big]_{\sqq^{-2}}
			=\frac{q-1}{q}[S_{1,0}^{((r+1)p_1-1)}]+(q-1)[S_{1,1}^{((r-1)p_1+1)}]*[K_{\de-\alpha_{11}}].
		\end{align}

		For any non-zero map $f:\co\to\co(k\vec{c})$, assume $f=x_1^{rp_1}\cdot g$ for some $r\geq 0$ and $x_1\nmid g$.
		So $\coker(f)\cong S_{1,0}^{(rp_1)}\oplus M$, where $M$ has no direct summand supported at ${\bla_1}=\infty$.
		Note that any two torsion sheaves supported at distinct points have zero Hom and $\Ext^1$-spaces. Hence $[M]*[X]=[M\oplus X]=[X]*[M]$ for any $X\in\scrt_{\bla_1}$.
		Therefore,
		\begin{align*}
			&\big[[\coker(f)],[S_{1,0}^{(p_1-1)}]\big]_{\sqq^{-2}}=[M]*\big[[S_{1,0}^{(rp_1)}],[S_{1,0}^{(p_1-1)}]\big]_{\sqq^{-2}}\\
			&\stackrel{(\ref{la})}{=} [M]*\Big(\frac{q-1}{q}[S_{1,0}^{((r+1)p_1-1)}]+(q-1)[S_{1,1}^{((r-1)p_1+1)}]*[K_{\de-\alpha_{11}}]\Big)\\
			&=\frac{q-1}{q}[M\oplus S_{1,0}^{((r+1)p_1-1)}]+(q-1)[M\oplus S_{1,1}^{((r-1)p_1+1)}]*[K_{\de-\alpha_{11}}]\\
			&=\frac{q-1}{q}[\coker(f_1)]+(q-1) [\coker(f_2)]*[K_{\de-\alpha_{11}}],
		\end{align*}
		where $f_1=x_1^{(r+1)p_1-1}\cdot g:\co(\vec{x}_1)\rightarrow \co((k+1)\vec{c})$ and $f_2=x_1^{(r-1)p_1+1}\cdot g:\co(\vec{c})\rightarrow \co(k\vec{c}+\vec{x}_1)$. Here it is understood that $f_2$ and then the term involving $[\coker (f_2)]$ does not exist if $r=0$. 
		Hence,  \begin{align*}
			&\sum_{0\neq f:\co\rightarrow \co(m\vec{c})} \big[[\coker(f)],[S_{1,0}^{(p_1-1)}]\big]_{\sqq^{-2}}\\
			&=\frac{q-1}{q}\sum_{0\neq f_1:\co(\vec{x}_1)\rightarrow \co((m+1)\vec{c})} [\coker (f_1)]+(q-1)\sum_{0\neq f_2:\co(\vec{c})\rightarrow \co(m\vec{c}+\vec{x}_1)} [\coker (f_2)]*[K_{\de-\alpha_{11}}].
		\end{align*}
		Then it follows from \eqref{def:Theta star}, \eqref{eq:THstarm+} and \eqref{eq:THstarm-} that
		\begin{align}\label{S_10andTheta}
			\big[\haTh_{\star, m},[S_{1,0}^{(p_1-1)}]\big]_{\sqq^{-2}}
			=(\sqq-\sqq^{-1})\haTh_{\star, m+1}^-+(\sqq-\sqq^{-1})\haTh_{\star, m-1}^+*[K_{\de-\alpha_{11}}].\end{align}

		Similarly, for any $r\geq 0$ we have
		\begin{align*}
			\big[[S_{1,1}], [S_{1,0}^{(rp_1)}]\big]_{\sqq^{-2}}
			=\frac{q-1}{q}[S_{1,1}^{(rp_1+1)}]+(q-1)[S_{1,0}^{(rp_1-1)}]*[K_{\alpha_{11}}],
		\end{align*}
		which implies
		\begin{align*}
			&\big[[S_{1,1}],\sum_{0\neq f:\co\rightarrow \co(m\vec{c})} [\coker(f)]\big]_{\sqq^{-2}}\\
			&=\frac{q-1}{q}\sum_{0\neq f_1:\co\rightarrow \co(m\vec{c}+\vec{x}_1)} [\coker (f_1)]
			+(q-1)\sum_{0\neq f_2:\co(\vec{x}_1)\rightarrow \co(m\vec{c})} [\coker (f_2)]*[K_{\alpha_{11}}].
		\end{align*}
		Then by definition we have
		\begin{align}\label{S_11andTheta}
			\big[[S_{1,1}],\haTh_{\star, m}\big]_{\sqq^{-2}}=\frac{q-1}{q}\haTh_{\star,m}^++(q-1)\haTh_{\star,m}^-*[K_{\alpha_{11}}].
		\end{align}
		Combining (\ref{S_10andTheta}) and (\ref{S_11andTheta}), we obtain (\ref{Theta+}) and (\ref{Theta-}).
	\end{proof}

	Now we can give the proof of Lemma \ref{lem:SerreO11}.
	\begin{proof}[Proof of Lemma \ref{lem:SerreO11}]
		For any $k\in\Z$,
		by (\ref{Si1 and co double direction}), we have
		\begin{align}
			\label{exchange relation of co and S_11}
			[\co(k\vec{c})]*[S_{1,1}]=\sqq[S_{1,1}]*[ \co(k\vec{c})] -(q-1)  [\co(k\vec{c}+\vec{x}_1)].
		\end{align}
		Hence,
		\begin{align*}
			&\widehat{S}(k_1,k_2\mid 0;\star,[1,1])\\
			&=[S_{1,1}]*[\co(k_1\vec{c})]*[\co(k_2\vec{c})] -[2]_\sqq [\co(k_1\vec{c})]*[S_{1,1}]*[\co(k_2\vec{c})]
			+[\co(k_1\vec{c})]*[\co(k_2\vec{c})]*[S_{1,1}]
			\\
			&=\big( \sqq^{-1} [\co(k_1\vec{c})]*[S_{1,1}]+(\sqq-\sqq^{-1})[\co(k_1\vec{c}+\vec{x}_1)]\big)*[\co(k_2\vec{c})]\\
			&\quad-[2]_\sqq [\co(k_1\vec{c})]*[S_{1,1}]*[\co(k_2\vec{c})]
			+[\co(k_1\vec{c})]*\big(\sqq[S_{1,1}]*[\co(k_2\vec{c})]-(q-1)[\co(k_2\vec{c}+\vec{x}_1)]\big)
			\\
			&=(\sqq-\sqq^{-1})  [\co(k_1\vec{c}+\vec{x}_1)]*[\co(k_2\vec{c})]-(q-1) [\co(k_1\vec{c})]* [\co(k_2\vec{c}+\vec{x}_1)].
		\end{align*}
		For $k_1=k_2$, the desired formula \eqref{SerreO11} holds by noting that
		\begin{align*}
			\widehat{\SS}(k_1,k_1\mid 0;\star,[1,1])
			=&\widehat{S}(k_1,k_1\mid 0;\star,[1,1])
			\\
			=&-(q-1)  %[\co(k_1\vec{c}+\vec{x}_1)]*[\co(k_1\vec{c})]-(q-1)
			\big[[\co(k_1\vec{c})], [\co(k_1\vec{c}+\vec{x}_1)]\big]_{\sqq^{-1}}
			\\
			=&-\sqq^{-1}(q-1)^2 [S_{1,1}]*[K_{ \co(k_1\vec{c}) }]
			\\
			\stackrel{\eqref{eqnSSkk1}}{=}&(1-q)^2\widehat{\R}(k_1,k_1\mid 0;\star,[1,1]).
		\end{align*}
		Without loss of generality, we assume $k_1<k_2$ in the following. Then we have
		\begin{align*}
			&\widehat{\SS}(k_1,k_2\mid 0;\star,[1,1])\\
			&=(\sqq-\sqq^{-1})  [\co(k_1\vec{c}+\vec{x}_1)]*[\co(k_2\vec{c})]-(q-1) [\co(k_1\vec{c})]* [\co(k_2\vec{c}+\vec{x}_1)]
			\\
			&\quad+ (\sqq-\sqq^{-1})   [\co(k_2\vec{c}+\vec{x}_1)]*[\co(k_1\vec{c})]-(q-1) [\co(k_2\vec{c})]* [\co(k_1\vec{c}+\vec{x}_1)]
			\\
			&=-(q-1) \Big( \big[[\co(k_2\vec{c})], [\co(k_1\vec{c}+\vec{x}_1)]\big]_{\sqq^{-1}}+\big[[\co(k_1\vec{c})], [\co(k_2\vec{c}+\vec{x}_1)]\big]_{\sqq^{-1}}\Big)
			\\
			&\stackrel{\eqref{eq:commcolcox1}}{=}\sqq^{-2}(q-1)^3\Big(\haTh_{\star,k_2-k_1}^-*[K_{\co(k_1\vec{c}+\vec{x}_1)}]-
			\haTh_{\star,k_2-k_1}^+*[K_{\co(k_1\vec{c})}]\Big).
		\end{align*}
		If $k_2=k_1+1$, then
		\eqref{SerreO11} holds since
		\begin{align*}
			&\widehat{\SS}(k_1,k_1+1\mid 0;\star,[1,1])\\
			&=\sqq^{-2}(q-1)^3\haTh_{\star,1}^-*[K_{\co(k_1\vec{c}+\vec{x}_1)}]-
			\sqq^{-2}(q-1)^3\haTh_{\star,1}^+*[K_{\co(k_1\vec{c})}]\\
			%=&(1+\sqq^{-2})(q-1)^2 [S_{1,0}^{(p_1-1)}]*[K_{\co(k_1\vec{c}+\vec{x}_1)}]
			%-(q-1)^2\big[[S_{1,1}], \haTh_{\star,1}\big]_{\sqq^{-2}}*[K_{\co(k_1\vec{c})}]\\
			&=\sqq^{-1}[2]_\sqq(q-1)^2 [S_{1,0}^{(p_1-1)}]*[K_{\co(k_1\vec{c}+\vec{x}_1)}]
			-(q-1)^2\big[[S_{1,1}], \haTh_{\star,1}\big]_{\sqq^{-2}}*[K_{\co(k_1\vec{c})}]
			\\
			&=(q-1)^2\Big(-\sqq^{-1}[2]_\sqq \haB_{[1,1],-1}*[K_\de] -\big[\haB_{[1,1],0},\haTh_{\star,1}\big]_{\sqq^{-2}}\Big)*[K_{\co(k_1\vec{c})}]
			\\
			&\stackrel{\eqref{eqnSSkk+1}}{=}(1-q)^2\widehat{\R}(k_1,k_1+1\mid 0;\star,[1,1]),
		\end{align*}
		where the second equality follows by Lemma \ref{lem:A1}.
		%by noting that $\haB_{[1,1],-1}=[S_{1,0}^{p_1-1}]*[K_{\alpha_{11}-\de}]$.
		
		If $k_2>k_1+1$,
		then by (\ref{Theta+}) and (\ref{Theta-}), we have
		\begin{align*}
			&\widehat{\SS}(k_1,k_2\mid 0;\star,[1,1])-q\widehat{\SS} (k_1,k_2-2\mid 0;\star,[1,1])*[K_{\delta}]
			\\
			%&\SS(\co(k_1\vec{c}),\co(k_2\vec{c})| S_{1,1})-q\SS(\co(k_1\vec{c}),\co((k_2-2)\vec{c})| S_{1,1})\K_{\delta}\\
			&=\sqq^{-2}(q-1)^3[K_{\co(k_1\vec{c}+\vec{x}_1)}]*\big(\haTh_{\star,k_2-k_1}^--
			q\haTh_{\star,k_2-k_1-2}^-*[K_{\delta}]\big)
			\\
			&\quad-\sqq^{-2}(q-1)^3[K_{\co(k_1\vec{c})}]*\big(\haTh_{\star,k_2-k_1}^+-
			q\haTh_{\star,k_2-k_1-2}^+*[K_{\delta}]\big)\\
			&=(q-1)^2
			\Big(\sqq^{-1}\big[\haTh_{\star, k_2-k_1-1}, [S_{1,0}^{(p_1-1)}]\big]_{\sqq^{-2}}-
			\big[[S_{1,1}],\haTh_{\star, k_2-k_1-2}\big]_{\sqq^{-2}}*[K_{\de-\alpha_{11}}]\Big)*[K_{\co(k_1\vec{c}+\vec{x}_1)}]\\
			&\quad-(q-1)^2\Big(\big[[S_{1,1}],\haTh_{\star, k_2-k_1}\big]_{\sqq^{-2}}-\sqq\big[\haTh_{\star, k_2-k_1-1}, [S_{1,0}^{(p_1-1)}]\big]_{\sqq^{-2}}*[K_{S_{1,1}}]\Big)*[K_{\co(k_1\vec{c})}]\\
			&=(q-1)^2\Big(-[2]_\sqq\big[\haTh_{\star, k_2-k_1-1}, \haB_{[1,1],-1}\big]_{\sqq^{-2}}*[K_{\co((k_1+1)\vec{c})}]
			-\big[\haB_{[1,1],0},\haTh_{\star, k_2-k_1}\big]_{\sqq^{-2}}*[K_{\co(k_1\vec{c})}]\\
			&\quad-\big[\haB_{[1,1],0},\haTh_{\star, k_2-k_1-2}\big]_{\sqq^{-2}}*[K_{\co((k_1+1)\vec{c})}]\Big)
			\\
			&\stackrel{(\ref{Theta and B k_2>k_1+1})}{=}(q-1)^2\Big(\widehat{\R}(k_1,k_2\mid 0;\star,[1,1])-q\widehat{\R} (k_1,k_2-2\mid 0;\star,[1,1])*[K_{\delta}]\Big).
		\end{align*}
		Therefore, the desired formula \eqref{SerreO11} follows by induction.
	\end{proof}

	%%%%%%%%%%%%%%
	\subsection{Proof of Lemma \ref{lem:SerreO111}}
	\label{subsec:proofSerre0111}
	
	In order to prove Lemma \ref{lem:SerreO111}, we need the following lemma.
	%Recall that $\haB_{[1,1],-1}=-[S_{1,0}^{(p_1-1)}]*[K_{\alpha_{11}-\de}]$.
	\begin{lemma}
		\label{lem:HBstar11-1}
		For any $m\geq 1$,
		\begin{align}
			\label{lem:HBstar11-1form-pre}
			&\big[\haTh_{\star,m}, \widehat{B}_{[1,1],-1}\big]+\big[\haTh_{\star,m-2},\widehat{B}_{[1,1],-1}\big]*[K_{\delta}]
			\\\notag
			&=
			\sqq^{-1}\big[\haTh_{\star,m-1},\widehat{B}_{[1,1],0}\big]_{\sqq^{2}}
			+\sqq\big[\haTh_{\star,m-1},\widehat{B}_{[1,1],-2}\big]_{\sqq^{-2}}*[K_{\delta}].
		\end{align}
	\end{lemma}
	
	\begin{proof}
		Recall $\widehat{\Theta}_{\star,m}$ defined in (\ref{def:Theta star}).
		Denote by $\mathcal{L}_k=\{\coker(g)\;|\;g\in\Hom(\co, \co(k\vec{c})), x_1\nmid g\}$ for any $k\geq0$. Then any $M\in\mathcal{L}_k$ has no direct summands from $\scrt_{\bla_1}$, and hence $[M]*[X]=[M\oplus X]=[X]*[M]$ for any $X\in\scrt_{\bla_1}$.
		For any non-zero map $f:\co\to\co(k\vec{c})$, assume $f=x_1^{rp_1}\cdot g$ for some $r\geq 0$ and $x_1\nmid g$.
		Then $\coker (f)\cong S_{1,0}^{(rp_1)}\oplus M_{k-r}$, where $M_{k-r}\in\mathcal{L}_{k-r}$.
		Similar to the proof of Lemma \ref{lem:A1}, it suffices to show for any $r\geq 0$ that
		\begin{align*}
			&\big[[S_{1,0}^{(rp_1)}], \widehat{B}_{[1,1],-1}\big]+q\big[[S_{1,0}^{((r-2)p_1)}],\widehat{B}_{[1,1],-1}\big]*[K_{\delta}]
			\\\notag
			&=
			\big[[S_{1,0}^{((r-1)p_1)}],\widehat{B}_{[1,1],0}\big]_{\sqq^{2}}
			+q\big[[S_{1,0}^{((r-1)p_1)}],\widehat{B}_{[1,1],-2}\big]_{\sqq^{-2}}*[K_{\delta}],
		\end{align*}
		which follows from Lemma \ref{theta and B 4 term in a tube}.
	\end{proof}
	
	Using the equalities $\sqq^{-1}[2]_\sqq[a,b]=[a,b]_{\sqq^{-2}}-[b,a]_{\sqq^{-2}}$ and $[a,b]_{\sqq^{2}}=-\sqq^{2}[b,a]_{\sqq^{-2}}$, \eqref{lem:HBstar11-1form-pre} can be reformulated as follows
	\begin{align}
		\label{lem:HBstar11-1form}
		&\big[\haTh_{\star,m}, \widehat{B}_{[1,1],-1}\big]_{\sqq^{-2}}-[2]_\sqq\big[\haTh_{\star,m-1},\widehat{B}_{[1,1],-2}\big]_{\sqq^{-2}}*[K_{\delta}]
		+\big[\haTh_{\star,m-2},\widehat{B}_{[1,1],-1}\big]_{\sqq^{-2}}*[K_{\delta}]
		\\\notag
		&=\big[\widehat{B}_{[1,1],-1},\haTh_{\star,m}\big]_{\sqq^{-2}}-
		[2]_\sqq\big[\widehat{B}_{[1,1],0},\haTh_{\star,m-1}\big]_{\sqq^{-2}}
		+\big[\widehat{B}_{[1,1],-1},\haTh_{\star,m-2}\big]_{\sqq^{-2}}*[K_{\delta}].
	\end{align}

	\begin{proof}[Proof of Lemma \ref{lem:SerreO111}]
		For any $k\in\Z$, by (\ref{Si0 and co double direction}) we have
		%by direct computations, we get
		%\begin{align*}
		%[S_{1,0}^{(p_1-1)}]*[\co(k\vec{c})]&=[S_{1,0}^{(p_1-1)}\oplus\co(k\vec{c})];
		%\\
		%[\co(k\vec{c})]*[S_{1,0}^{(p_1-1)}]&=\sqq^{-1}\Big([S_{1,0}^{(p_1-1)}\oplus \co(k\vec{c})]+(q-1) [\co((k-1)\vec{c}+\vec{x}_1)]\ast[K_{\de-\alpha_{11}}]\Big).
		%\end{align*}
		\begin{align}\label{exchange relation of co and S_11}
			[S_{1,0}^{(p_1-1)}]*[\co(k\vec{c})]
			=&\sqq[\co(k\vec{c})]*[S_{1,0}^{(p_1-1)}]-(q-1)  [\co((k-1)\vec{c}+\vec{x}_1)]*[K_{\de-\alpha_{11}}].
		\end{align}
		Hence,
		\begin{align*}
			& [S_{1,0}^{(p_1-1)}]*[\co(k_1\vec{c})]*[\co(k_2\vec{c})] -[2]_\sqq [\co(k_1\vec{c})]*[S_{1,0}^{(p_1-1)}]*[\co(k_2\vec{c})]
			\\
			&\quad+[\co(k_1\vec{c})]*[\co(k_2\vec{c})]*[S_{1,0}^{(p_1-1)}]
			\\
			&=\Big(\sqq[\co(k_1\vec{c})]*[S_{1,0}^{(p_1-1)}]-(q-1)[\co((k_1-1)\vec{c}+\vec{x}_1)]*[K_{\de-\alpha_{11}}]\Big)*[\co(k_2\vec{c})]
			\\
			&\quad -[2]_\sqq [\co(k_1\vec{c})]*[S_{1,0}^{(p_1-1)}]*[\co(k_2\vec{c})]
			\\
			&\quad +[\co(k_1\vec{c})]*\Big( \sqq^{-1} [S_{1,0}^{(p_1-1)}]*[\co(k_2\vec{c})] +(\sqq-\sqq^{-1}) [\co((k_2-1)\vec{c}+\vec{x}_1)]  *[K_{\de-\alpha_{11}}]\Big)
			\\
			%=&\sqq[\co(k_1\vec{c})]\ast\Big(\sqq[\co(k_2\vec{c})]*[S_{1,0}^{(p_1-1)}]\ast[K_{\de-\alpha_{11}}]^{-1}
			%-(q-1)  [\co((k_2-1)\vec{c}+\vec{x}_1)]\Big)-(q-1)  [\co((k_1-1)\vec{c}+\vec{x}_1)]*[\co(k_2\vec{c})]\\
			%&-(q+1)[\co(k_1\vec{c})]\ast[\co(k_2\vec{c})]*[S_{1,0}^{(p_1-1)}]\ast[K_{\de-\alpha_{11}}]^{-1}+(v+v^{-1})(q-1) [\co(k_1\vec{c})]\ast [\co((k_2-1)\vec{c}+\vec{x}_1)]\\
			%&+[\co(k_1\vec{c})]*[\co(k_2\vec{c})]*[S_{1,0}^{(p_1-1)}]\ast[K_{\de-\alpha_{11}}]^{-1}\\
			&=(\sqq-\sqq^{-1})  \Big([\co(k_1\vec{c})]*[\co((k_2-1)\vec{c}+\vec{x}_1)]-\sqq [\co((k_1-1)\vec{c}+\vec{x}_1)]* [\co(k_2\vec{c})]\Big)\ast[K_{\de-\alpha_{11}}].
			%=&\blue{S(\co(k_1\vec{c}),\co(k_2\vec{c})| S_{1,0}^{(p_1-1)})(-\vec{x}_1).}
		\end{align*}
		Then
		\begin{align*}
			&\widehat{S}(k_1,k_2\mid -1; \star, [1,1])
			\\	&=-(\sqq-\sqq^{-1})  \Big([\co(k_1\vec{c})]*[\co((k_2-1)\vec{c}+\vec{x}_1)]-\sqq [\co((k_1-1)\vec{c}+\vec{x}_1)]* [\co(k_2\vec{c})]\Big).
		\end{align*}
		For $k_1=k_2$, the desired formula \eqref{SerreO11-1} holds since
		\begin{align*}
			\widehat{\SS}(k_1,k_1\mid -1; \star, [1,1])
			=&\widehat{S}(k_1,k_1\mid -1; \star, [1,1])\\
			=&(q-1)\big[[\co((k_1-1)\vec{c}+\vec{x}_1)],[\co(k_1\vec{c})]\big]_{\sqq^{-1}}
			\\
			=&\sqq^{-1}(q-1)^2 [S_{1,0}^{(p_1-1)}]*[K_{ \co((k_1-1)\vec{c}+\vec{x}_1) }]
			\\
			=&-\sqq^{-1}(q-1)^2 \haB_{[1,1],-1}*[K_{ \co(k_1\vec{c}) }]
			\\
			%=&-\sqq^{-1}(q-1)^2 [\widehat{B}_{j,l}, \haTh _{i,0}]_{\sqq^{-2}}*[K_{k_1\de+\alpha_i}]
			%\\
			\stackrel{\eqref{eqnSSkk1}}{=}&(1-q)^2\widehat{\R}(k_1,k_1\mid -1;\star,[1,1]).
		\end{align*}
		Without loss of generality, we assume $k_1<k_2$ in the following. Then we have
		\begin{align*}
			&\widehat{\SS}(k_1,k_2\mid -1; \star, [1,1]) \\
			&=-(\sqq-\sqq^{-1})  [\co(k_1\vec{c})]*[\co((k_2-1)\vec{c}+\vec{x}_1)]+(q-1) [\co((k_1-1)\vec{c}+\vec{x}_1)]* [\co(k_2\vec{c})]
			\\
			&\quad- (\sqq-\sqq^{-1})  [\co(k_2\vec{c})]*[\co((k_1-1)\vec{c}+\vec{x}_1)]+(q-1) [\co((k_2-1)\vec{c}+\vec{x}_1)]* [\co(k_1\vec{c})]
			\\
			&=(q-1)\big[ [\co((k_2-1)\vec{c}+\vec{x}_1)], [\co(k_1\vec{c})]\big]_{\sqq^{-1}}+(q-1)\big[[\co((k_1-1)\vec{c}+\vec{x}_1)], [\co(k_2\vec{c})] \big]_{\sqq^{-1}}
			\\
			&\stackrel{\eqref{eq:commcolx1co} }{=}-\sqq^{-3}(q-1)^3\haTh_{\star,k_2-k_1-1}^+*[K_{\co(k_1\vec{c})}]+\sqq^{-1}(q-1)^3\haTh_{\star,k_2-k_1+1}^-*[K_{\co((k_1-1)\vec{c}+\vec{x}_1)}].
		\end{align*}
		If $k_2=k_1+1$, then the desired formula \eqref{SerreO11-1}
		holds since
		\begin{align*}
			&\widehat{\SS}(k_1,k_1+1\mid -1; \star, [1,1]) \\
			&=-\sqq^{-3}(q-1)^3\haTh_{\star,0}^+*[K_{\co(k_1\vec{c})}]+\sqq^{-1}(q-1)^3\haTh_{\star,2}^-*[K_{\co((k_1-1)\vec{c}+\vec{x}_1)}]
			\\
			&=-\sqq^{-3}(q-1)^3\cdot\frac{1}{\sqq-\sqq^{-1}}[S_{1,1}]*[K_{\co(k_1\vec{c})}]\\
			&\quad+\sqq^{-1}(q-1)^2
			\cdot\Big(\sqq\big[\haTh_{\star,1}, [S_{1,0}^{(p_1-1)}]\big]_{\sqq^{-2}}-q\big[[S_{1,1}], \haTh_{\star, 0}\big]_{\sqq^{-2}}*[K_{\de-\alpha_{11}}]\Big)*[K_{\co((k_1-1)\vec{c}+\vec{x}_1)}]
			\\&=-\sqq^{-1}[2]_\sqq(q-1)^2[S_{1,1}]*[K_{\co(k_1\vec{c})}]+(q-1)^2
			\big[\haTh_{\star,1}, [S_{1,0}^{(p_1-1)}]\ast[K_{\de-\alpha_{11}}]^{-1} \big]_{\sqq^{-2}}  *[K_{\co(k_1\vec{c})}]\\
			&=(1-q)^2\Big(-[2]_\sqq [\haTh _{\star,0},\haB_{[1,1],0}]_{\sqq^{-2}}- [\haTh _{\star,1},\widehat{B}_{[1,1],-1}]_{\sqq^{-2}}\Big)*[K_{\co(k_1\vec{c})}]
			\\
			&\stackrel{\eqref{lem:HBstar11-1form}}{=}(1-q)^2\Big(-[2]_\sqq [\haTh _{\star,0},\haB_{[1,1],-2}]_{\sqq^{-2}}*[K_\de]- [\widehat{B}_{[1,1],-1}, \haTh _{\star,1}]_{\sqq^{-2}}\Big)*[K_{\co(k_1\vec{c})}]
			\\
			&\stackrel{\eqref{eqnSSkk+1} }{=}(1-q)^2\widehat{\R}(k_1,k_1+1\mid -1; \star, [1,1]).
		\end{align*}
		If $k_2>k_1+1$,
		then by (\ref{Theta+}) and (\ref{Theta-}), we have
		\begin{align*}
			&\frac{1}{(q-1)^2}\Big(\widehat{\SS}(k_1,k_2\mid -1;\star,[1,1])-q\widehat{\SS}(k_1,k_2-2\mid -1;\star,[1,1])*[K_{\delta}]\Big)
			\\
			%&\SS(\co(k_1\vec{c}),\co(k_2\vec{c})| S_{1,1})-q\SS(\co(k_1\vec{c}),\co((k_2-2)\vec{c})| S_{1,1})\K_{\delta}\\
			&=-\sqq^{-3}(q-1)\Big(\haTh_{\star,k_2-k_1-1}^+-
			q\haTh_{\star,k_2-k_1-3}^+*[K_{\delta}]\Big)*[K_{\co(k_1\vec{c})}]\\
			&\quad+\sqq^{-1}(q-1)\Big(\haTh_{\star,k_2-k_1+1}^--
			q\haTh_{\star,k_2-k_1-1}^-*[K_{\delta}]\Big)*[K_{\co((k_1-1)\vec{c}+\vec{x}_1)}]\\
			&=
			-\Big(\sqq^{-1}\big[[S_{1,1}],\haTh_{\star,k_2-k_1-1}\big]_{\sqq^{-2}}-
			\big[\haTh_{\star,k_2-k_1-2}, [S_{1,0}^{(p_1-1)}]\big]_{\sqq^{-2}}*[K_{S_{1,1}}]\Big)*[K_{\co(k_1\vec{c})}]\\
			&\quad
			+\Big(\big[\haTh_{\star,k_2-k_1}, [S_{1,0}^{(p_1-1)}]\big]_{\sqq^{-2}}-
			\sqq\big[[S_{1,1}],\haTh_{\star,k_2-k_1-1}\big]_{\sqq^{-2}}*[K_{\de-\alpha_{11}}]\Big)*[K_{\co((k_1-1)\vec{c}+\vec{x}_1)}]
			\\
			&=-[2]_\sqq\big[[S_{1,1}],\haTh_{\star,k_2-k_1-1}\big]_{\sqq^{-2}}*[K_{\co(k_1\vec{c})}]
			+\big[\haTh_{\star,k_2-k_1-2}, [S_{1,0}^{(p_1-1)}]\big]_{\sqq^{-2}}*[K_{\co(k_1\vec{c}+\vec{x}_1)}]\\
			&\quad+\big[\haTh_{\star,k_2-k_1}, [S_{1,0}^{(p_1-1)}]\big]_{\sqq^{-2}}*[K_{\co((k_1-1)\vec{c}+\vec{x}_1)}]
			\\
			&=-[2]_\sqq\big[\haB_{[1,1],0},\haTh_{\star,k_2-k_1-1}\big]_{\sqq^{-2}}*[K_{\co(k_1\vec{c})}]
			-\big[\haTh_{\star,k_2-k_1-2}, \haB_{[1,1],-1}]\big]_{\sqq^{-2}}*[K_{\co((k_1+1)\vec{c})}]\\
			&\quad-\big[\haTh_{\star,k_2-k_1}, \haB_{[1,1],-1}]\big]_{\sqq^{-2}}*[K_{\co(k_1\vec{c})}].
		\end{align*}
		Therefore,
		\begin{align*}
			&\frac{1}{(q-1)^2}\Big(\widehat{\SS}(k_1,k_2\mid -1;\star,[1,1])-q\widehat{\SS} (k_1,k_2-2\mid -1;\star,[1,1])*[K_{\delta}]\Big)
			\\
			&=\Big(-[2]_\sqq\big[\widehat{B}_{[1,1],0},\haTh_{\star,k_2-k_1-1}\big]_{\sqq^{-2}}
			-\big[\haTh_{\star,k_2-k_1-2}, \widehat{B}_{[1,1],-1}\big]_{\sqq^{-2}}*[K_{\delta}]
			\\
			&\quad-\big[\haTh_{\star,k_2-k_1}, \widehat{B}_{[1,1],-1}\big]_{\sqq^{-2}}\Big)*[K_{\co(k_1\vec{c})}]\\
			&\stackrel{\eqref{lem:HBstar11-1form}}{=}\Big(-[2]_\sqq [\haTh _{\star,k_2-k_1-1},\haB_{[1,1],-2}]_{\sqq^{-2}}*[K_\de]-[\haB_{[1,1],-1},\haTh _{\star,k_2-k_1-2}]_{\sqq^{-2}}* [K_\de]
			\\
			&\quad-[\widehat{B}_{[1,1],-1}, \haTh _{\star,k_2-k_1}]_{\sqq^{-2}}\Big)*[K_{\co(k_1\vec{c})}]\\
			&\stackrel{(\ref{Theta and B k_2>k_1+1})}{=}\widehat{ \R}(k_1,k_2\mid -1;\star,[1,1])-q\widehat{ \R}(k_1,k_2-2\mid -1;\star,[1,1]),
		\end{align*}
		%where $(\dag)$ follows from Sublemma \ref{lem:HBstar11-1}.
		Then \eqref{SerreO11-1} follows by induction and the proof is completed.
	\end{proof}

	%%%%%%%%%%%%%%%%%%%
	\section{Proofs of results in Section \ref{sec:Relationsstari1 II}}
	\label{sec:proofRelatstari1II}

	For a partition $\lambda=(1^{l_1},\dots,n^{l_n})$ and $u>0$, we denote
	
	$\triangleright$ $u\triangleleft\lambda$ -- $u$ is a part of $\lambda$ (i.e., $l_u>0$),

	$\triangleright$ $\lambda\cup \{u\}$ -- the partition $(1^{l_1},\dots, u^{l_u+1},\dots,n^{l_n})$,
	
	$\triangleright$ $\lambda\backslash \{u\}$ -- the partition $(1^{l_1},\dots, u^{l_u-1},\dots,n^{l_n})$ if $u\triangleleft\lambda$.

	\subsection{Proof of Proposition \ref{prop:realroot}}
	\label{subsec:proofrealroot}
	
	Recall definition of $\widehat{H}_{\star,r}$ in \eqref{formula for Hm}.
	Then we have the following results.
	\begin{lemma}
		\label{lem:l=0}
		For any $r>0$, we have
		\begin{align}
			\label{eq:S1H}
			\big[[\![S_{1,1}]\!], \widehat{H}_{\star,r}\big]=&\frac{[r]_\sqq}{r} [M_{r\delta+\alpha_1}]+\frac{[r]_\sqq}{r}[M_{r\delta-\alpha_1}]*[K_{S_{1,1}}],
			\\
			\label{eq:S0H}
			\big[\widehat{H}_{\star,r},[\![S_{1,0}^{(p_1-1)}]\!]\big]=&\frac{[r]_\sqq}{r}[M_{(r+1)\delta-\alpha_{11}}]+\frac{[r]_\sqq}{r} [M_{(r-1)\delta+\alpha_{11}}]*[K_{\de-\alpha_{11}}].
		\end{align}
	\end{lemma}
	
	\begin{proof}
		We only prove \eqref{eq:S1H} since \eqref{eq:S0H} can be proved  dually.
		
		Let $\lambda=(1^{l_1},2^{l_2},\dots,n^{l_n})$ be a partition of $r$. Since $\langle S_{1,0}^{(\lambda)}, S_{1,1}\rangle=0$, we know that  $\Hom_\X(S_{1,0}^{(\lambda)}, S_{1,1})=0=\Ext^1_\X(S_{1,0}^{(\lambda)}, S_{1,1})$. Hence,  $[S_{1,0}^{(\lambda)}]*[S_{1,1}]=[S_{1,0}^{(\lambda)}\oplus S_{1,1}]$, and then
		\begin{align}
			\label{eq:realrootS10S11}
			[\![S_{1,0}^{(\lambda)}]\!]*[\![S_{1,1}]\!]
			=q^{\ell(\lambda)}[\![S_{1,1}\oplus S_{1,0}^{(\lambda)}]\!]
		\end{align}
		by noting that $|\Hom_\X(S_{1,1}, S_{1,0}^{(\lambda)})|=q^{\ell(\lambda)}$.
		Observe that the middle term of any non-split sequence in $\Ext_{\X}^1(S_{1,1}, S_{1,0}^{(\lambda)})$ has the form $S_{1,1}^{(u p_1+1)}\oplus S_{1,0}^{(\lambda\backslash\{ u\})}$ for some $u\triangleleft\lambda$, and $$F_{S_{1,1}, S_{1,0}^{(\lambda)}}^{S_{1,1}^{(u p_1+1)}\oplus S_{1,0}^{(\lambda\backslash\{ u\})}}=F_{S_{1,1}, S_{1,0}^{(u)}}^{S_{1,1}^{(u p_1+1)}}=1.$$
		On the other hand, for any non-zero map $f:S_{1,1}\rightarrow S_{1,0}^{(\lambda)}$, we have $f$ is injective and $\coker(f)=S_{1,0}^{(u p_1-1)}\oplus S_{1,0}^{(\lambda\backslash\{ u\})}$ for some $u\triangleleft\lambda$. Obviously,
		$$\Ext^1_\X(S_{1,0}^{(u p_1-1)}\oplus S_{1,0}^{(\lambda\backslash\{ u\})}, S_{1,1})\cong\Ext_\X^1(S_{1,0}^{(u p_1-1)}, S_{1,1})\cong \mathbf{k}.$$
		Then by Corollary \ref{ihall formula} (1),
		%using $|\Aut(S_{1,1})|=q-1$ and $\Hom(S_{1,0}^{(n)}, S_{1,1})=0$ for $n\geq 1$,
		we get
		\begin{align}
			\label{eq:realrootS11S10}
			&[\![S_{1,1}]\!]*[\![S_{1,0}^{(\lambda)}]\!]
			=[\![S_{1,1}\oplus S_{1,0}^{(\lambda)}]\!]
			+\sum_{u\triangleleft \lambda} [\![S_{1,1}^{(u p_1+1)}\oplus S_{1,0}^{(\lambda\backslash\{ u\})}]\!]
			+\sum_{u\triangleleft \lambda}[\![S_{1,0}^{(u p_1-1)}\oplus S_{1,0}^{(\lambda\backslash\{ u\})}]\!]*[K_{\alpha_{11}}].
		\end{align}
		Observe that $\big[[S_{1,1}], [M]\big]=0$ for any indecomposable torsion sheaf $M\notin\scrt_{\bla_1}$.
		Therefore, it follows  from  \eqref{eq:realrootS10S11}--\eqref{eq:realrootS11S10} that
		\begin{align*}
			&\frac{r}{[r]}  \big[[\![S_{1,1}]\!], \widehat{H}_{\star,r}\big]
			%\notag=&\frac{1}{q-1} [[S_{1,1}], \sum\limits_{|\lambda|=r}\bn(\ell(\lambda)-1)\frac{1}{|\Aut(S_{1,0}^{(\lambda)})|}[S_{1,0}^{(\lambda)}]]\\
			=\sum\limits_{|\lambda|=r}{\bn(\ell(\lambda)-1)}\big[[\![S_{1,1}]\!],[\![S_{1,0}^{(\lambda)}]\!]\big]
			\\
			&=\sum\limits_{|\lambda|=r}\bn(\ell(\lambda)-1)\cdot(1-q^{\ell(\lambda)})[\![S_{1,1}\oplus S_{1,0}^{(\lambda)}]\!]
			+\sum\limits_{|\lambda|=r}\bn(\ell(\lambda)-1)\sum_{u\triangleleft \lambda} [\![S_{1,1}^{(u p_1+1)}\oplus S_{1,0}^{(\lambda\backslash\{ u\})}]\!]\\
			&\quad+\sum\limits_{|\lambda|=r}\bn(\ell(\lambda)-1)\sum_{u\triangleleft \lambda}[\![S_{1,0}^{(u p_1-1)}\oplus S_{1,0}^{(\lambda\backslash\{ u\})}]\!]*[K_{\alpha_{11}}]\\
			&= [M_{r\delta+\alpha_1}]+[M_{r\delta-\alpha_1}]*[K_{\alpha_{11}}],
		\end{align*}
		by using $\bn(\ell(\lambda)-1)\cdot(1-q^{\ell(\lambda)})=\bn(\ell(\lambda))$. The proof is completed.
	\end{proof}
	
	With the help of Lemma \ref{lem:l=0}, we can prove Proposition \ref{prop:realroot} now.
	
	\begin{proof}[Proof of Proposition \ref{prop:realroot}]
		We prove the formulas by induction on $r$.
		For $r=1$,  they follow from \eqref{eq:B11} and Proposition \ref{prop:DrGenA}.
		
		We assume $\haB_{[1,1],l}= (q-1)[M_{l\de+\alpha_{11}}]$ and $\haB_{[1,1],-l}=(1-q) [M_{l\de-\alpha_{11}}]*[K_{-l\de+\alpha_{11}}]$ for any $1\leq l\leq r$.
		Using \eqref{humstarbvl} proved in Proposition \ref{prop:humstarbvl}, we have
		\begin{align*}
			[\widehat{H}_{\star,m},\widehat{\y}_{[1,1],l}]=-\frac{[m]_\sqq}{m} \widehat{\y}_{[1,1],l+m}+\frac{[m]_\sqq}{m} \widehat{\y}_{[1,1],l-m}*[K_{m\de}]
		\end{align*}
		for any $l\in\Z$, $m>0$.
		Recall that $\haB_{[1,1],0}=[S_{1,1}]$, $\haB_{[1,1],-1}=-[S_{1,0}^{(p_1-1)}]*[K_{-\de+\alpha_{11}}]$. Using \eqref{eq:S0H}, one can see that
		$\haB_{[1,1],-(r+1)}=(1-q) [M_{(r+1)\de-\alpha_{11}}]*[K_{-(r+1)\de+\alpha_{11}}]$.
		Together with \eqref{eq:S1H}, we have $\haB_{[1,1],r+1}= (q-1)[M_{(r+1)\de+\alpha_{11}}]$. The proof is completed.
		%\red{By induction.} With the blue colour sign, it coincides with Section 5.
	\end{proof}

	\subsection{Proof of Proposition \ref{prop:imageroot}}
	\label{subsec:proofimroot}

	First, let us give some formulas on the cardinalities of  Hom-spaces and automorphism groups.
	\begin{lemma}
		\label{lem:HomExt}
		Let $\lambda=(1^{l_1},2^{l_2},\dots, n^{l_n})$ be a partition and $u>0$. Then
		\begin{align}
			&|\Aut(S_{1,0}^{(u p_1-1)})|= q^{u-1}(q-1), \qquad   |\Aut(S_{1,1}^{(u p_1+1)})|= q^{u}(q-1),
			\\
			\label{eq:aut}
			&|\Aut(S_{1,0}^{(\lambda)})|=q^{|\lambda|+\sum_i il_i(l_i-1)+2\sum\limits_{i<j}il_il_j}\prod_i(1-q^{-1})(1-q^{-2})\cdots (1-q^{-l_i}),
			\\
			&|\Hom(S_{1,1}^{(u p_1+1)},S_{1,0}^{(\lambda)})|=|\Hom(S_{1,1}^{((u+1)p_1)},S_{1,0}^{(\lambda)})|=q^{\sum_i\min\{i,u+1\}l_i},
			%^{\sum\limits_{i\leq u} il_i +\sum\limits_{i>u}(u+1)l_i },
			\\
			&|\Hom(S_{1,0}^{(\lambda)},S_{1,1}^{(u p_1+1)})|=|\Hom(S_{1,0}^{(\lambda)}, S_{1,0}^{(u p_1)})|=q^{\sum_i\min\{i,u\}l_i},
			%{\sum\limits_{i\leq u} il_i+\sum\limits_{i>u} u l_i},
			\\
			&|\Hom(S_{1,0}^{(u p_1-1)},S_{1,0}^{(\lambda)})|=|\Hom(S_{1,0}^{((u-1)p_1)},S_{1,0}^{(\lambda)})|=q^{\sum_i\min\{i,u-1\}l_i},
			\\
			&|\Hom(S_{1,0}^{(\lambda)},S_{1,0}^{(u p_1-1)})|=|\Hom(S_{1,0}^{(\lambda)}, S_{1,0}^{(u p_1)})|=q^{\sum_i\min\{i,u\}l_i},
			\\
			&|\Aut(S_{1,0}^{(\lambda\backslash\{ u\})})| = \frac{|\Aut(S_{1,0}^{(\lambda)})|\cdot q^{u} }{(1-q^{-l_u})\cdot q^{2\sum_{i}\min\{i,u\}l_i}}
			\\
			&|\Aut(S_{1,0}^{(\lambda\cup\{u\})})| = |\Aut(S_{1,0}^{(\lambda)})|\cdot q^{u}\cdot (1-q^{-l_{u}-1})\cdot q^{2\sum_{i}\min\{i,u\}l_i}.
		\end{align}
	\end{lemma}
	
	\begin{proof}
		See \cite[Lemma 2.8]{Sch12} for \eqref{eq:aut}, and the others follow by direct computations.
	\end{proof}
	
	For any $a,b\geq 1$, we denote by
	\begin{align*}
		\phi_{a,b}=\frac{1}{q-1}\cdot|\Ext^{1}(S_{1,0}^{(ap_1-1)}\oplus S_{1,0}^{(bp_1-1)}, S_{1,1})_{S_{1,0}^{(ap_1-1)}\oplus S_{1,0}^{(bp_1)}}|,
	\end{align*}
	and
	\begin{align*}
		\psi_{a,b}= F_{S_{1,0}^{(\nua  p_1-1)}\oplus S_{1,0}^{((b-1)p)}, S_{1,0}^{(p_1-1)}}^{S_{1,0}^{(\nua  p_1-1)}\oplus S_{1,0}^{(b p_1-1)}}.
	\end{align*}
	A direct calculation shows that
	\begin{align}\label{phiab=psiba}
		\phi_{a,b}=\psi_{b,a}= \begin{cases}
			q,  &\text{ if }a<b,\\
			1,  &\text{ if } a>b, \\
			q+1,  & \text{ if }a=b.
		\end{cases}
	\end{align}

	\begin{lemma}
		\label{lem:image1}
		For any {$r\geq 1$}, we have
		\begin{align*}
			&-\frac{1}{(q-1)^2}\big[\haB_{[1,1],0}, \haB_{[1,1],-r}\big]_{\sqq^{2}}*[K_{r\de-\alpha_{11}}]\\
			&=\sqq^{-1}\sum\limits_{[M]\in\cm_{r\de}}\bn(\ell(M)-1) [\![M]\!]+\frac{\sqq}{q-1}\sum\limits_{|\lambda|=r}\bn(\ell(\lambda))[\![S_{1,0}^{(\lambda)}]\!]\\
			& \quad+{\sqq}\sum\limits_{a+|\nu|=r}\bn(\ell(\nu))\sum_{u\triangleleft\nu}\phi_{a,u}[\![S_{1,0}^{(\nua  p_1-1)}\oplus S_{1,0}^{(u p_1-1)}\oplus S_{1,0}^{(\nu\backslash \{u\})}]\!]*[K_{\alpha_{11}}].
		\end{align*}
		%where $\phi_{u}=|\Ext^1(S_{1,0}^{(\nua  p_1-1)}\oplus S_{1,0}^{(u p_1-1)}\oplus S_{1,0}^{(\nu\backslash \{u\})}, S_{1,1})_{S_{1,0}^{(\nua  p_1-1)}\oplus S_{1,0}^{(\nu )}}|$.
	\end{lemma}
	
	\begin{proof}
		Let $[M]=[S_{1,0}^{(\nua  p_1-1)}\oplus S_{1,0}^{(\nu  )}]\in\cm_{r\delta-\alpha_{11}}$. Then we have $|\nu|+\nua=r$ and $$\langle \widehat{M}, \alpha_{11}\rangle=\langle \alpha_{11}, \widehat{M}\rangle=-\langle \alpha_{11}, \alpha_{11}\rangle=-1.$$
		Since $\Hom(M, S_{1,1})=0$ and $\Ext^1(M, S_{1,1})\cong \mathbf{k}$, we get
		\begin{align*}
			[M]*[S_{1,1}]= {\sqq^{-1}}[M\oplus S_{1,1}]+ {\sqq^{-1}}(q-1)[S_{1,0}^{(\nu\cup\{a\})}].
		\end{align*}
		Let $\nu=(1^{l_1},2^{l_2},\dots,n^{l_n})$. Then we obtain
		\begin{align}
			\label{eq:S10S11}
			[\![M]\!]*[\![S_{1,1}]\!]=\sqq^{-1}q^{\ell(\nu)}[\![M\oplus S_{1,1}]\!]+\sqq^{-1}\frac{q^{l_a+1}-1}{q-1}\cdot q^{\sum\limits_{i>a}l_i}[\![S_{1,0}^{(\nu\cup\{a\})}]\!].
		\end{align}
		Here we use
		$|\Hom(S_{1,1}, M)|=q^{\ell(\nu)}$ and $$\frac{|\Aut(S_{1,0}^{(\nu\cup\{a\})})|}{|\Aut(M)|}=\frac{q^{l_a+1}-1}{q-1}\cdot q^{\sum\limits_{i>a}l_i},$$
		which follows from Lemma \ref{lem:HomExt} by direct computations.
		
		Observe that the middle term of any non-split extension in
		%$0\to S_{1,0}^{(\nua  p_1-1)}\oplus S_{1,0}^{(\nu )}\to L\to S_{1,1}\to0$, we have % $M$ of the
		$\Ext^1(S_{1,1}, M)$ has the form
		$S_{1,1}^{(ap_1)}\oplus S_{1,0}^{(\nu)}$ or $S_{1,0}^{(\nua  p_1-1)}\oplus S_{1,1}^{(u p_1+1)}\oplus S_{1,0}^{(\nu\backslash \{u\})}$ for some $u\triangleleft\nu$, and
		$$F_{S_{1,1}, M}^{S_{1,1}^{(ap_1)}\oplus S_{1,0}^{(\nu)}}=F_{S_{1,1}, S_{1,0}^{(\nua  p_1-1)}}^{S_{1,1}^{(ap_1)}}=1, \text{\qquad} F_{S_{1,1}, M}^{S_{1,0}^{(\nua  p_1-1)}\oplus S_{1,1}^{(u p_1+1)}\oplus S_{1,0}^{(\nu\backslash \{u\})}}=F_{S_{1,1}, S_{1,0}^{(u p_1)}}^{S_{1,1}^{(u p_1+1)}}=1.$$ Moreover, for the split extension we also have
		$F_{S_{1,1}, M}^{S_{1,1}\oplus M}=1$. 
		On the other hand, for any non-zero map $f:S_{1,1}\to M$, we have that $f$ is injective and $\coker(f)=S_{1,0}^{(\nua  p_1-1)}\oplus S_{1,0}^{(u p_1-1)}\oplus S_{1,0}^{(\nu\backslash \{u\})}$ for some $u\triangleleft\nu$, hence $\Hom(\coker(f), S_{1,1})=0$.
		%and $$\Ext^1(S_{1,0}^{(u p_1-1)}\oplus S_{1,0}^{(\lambda\backslash\{ u\})}, S_{1,1})\cong\Ext^1(S_{1,0}^{(u p_1-1)}, S_{1,1})\cong \mathbf{k}.$$
		Then by Corollary \ref{ihall formula} (1), we get
		\begin{align}
			\label{eq:S11S10}
			&[\![S_{1,1}]\!]*[\![M]\!]\\\notag
			&=\sqq^{-1}\Big([\![S_{1,1}\oplus M]\!]+ [\![S_{1,1}^{(\nua  p_1)}\oplus S_{1,0}^{(\nu)}]\!]+ \sum_{u\triangleleft\nu} [\![S_{1,0}^{(\nua  p_1-1)}\oplus S_{1,1}^{(u p_1+1)}\oplus S_{1,0}^{(\nu\backslash \{u\})}]\!]\Big)\\\notag
			&\quad+\sqq\sum_{u\triangleleft\nu}{\phi_{a,u}}[\![S_{1,0}^{(\nua  p_1-1)}\oplus S_{1,0}^{(u p_1-1)}\oplus S_{1,0}^{(\nu\backslash \{u\})}]\!]*[K_{\alpha_{11}}].
		\end{align}
		%where $\phi_{u}=|\Ext^1(S_{1,0}^{(\nua  p_1-1)}\oplus S_{1,0}^{(u p_1-1)}\oplus S_{1,0}^{(\nu\backslash \{u\})}, S_{1,1})_{S_{1,0}^{(\nua  p_1-1)}\oplus S_{1,0}^{(\nu )}}|$.
		Combining \eqref{eq:S10S11}--\eqref{eq:S11S10}, we have
		\begin{align*}
			&\big[[\![S_{1,1}]\!], [\![M]\!]\big]_{\sqq^2}
			\\
			&=\sqq^{-1}\Big((1-q^{\ell(\nu)+1})[\![S_{1,1}\oplus M]\!]{\normalsize {\Large }}+ [\![S_{1,1}^{(\nua  p_1)}\oplus S_{1,0}^{(\nu)}]\!]+ \sum_{u\triangleleft\nu} [\![S_{1,0}^{(\nua  p_1-1)}\oplus S_{1,1}^{(u p_1+1)}\oplus S_{1,0}^{(\nu\backslash \{u\})}]\!]\Big)\\
			&\quad-\sqq\frac{q^{l_a+1}-1}{q-1}\cdot q^{\sum\limits_{i>a}l_i}[\![S_{1,0}^{(\nu\cup\{a\})}]\!]+\sqq\sum_{u\triangleleft\nu}{\phi_{a,u}}[\![S_{1,0}^{(\nua  p_1-1)}\oplus S_{1,0}^{(u p_1-1)}\oplus S_{1,0}^{(\nu\backslash \{u\})}]\!]*[K_{\alpha_{11}}].
		\end{align*}
		Therefore, by definition,
		\begin{align}
			\label{eq:B110B11-r}
			&-\frac{1}{(q-1)^2}\big[\haB_{[1,1],0}, \haB_{[1,1],-r}\big]_{\sqq^{2}}*[K_{r\de-\alpha_{11}}]\\\notag
			&=\sum\limits_{[M]\in\cm_{r\delta-\alpha_{11}}}\bn(\ell(M)-1)\big[[\![S_{1,1}]\!], [\![M]\!]\big]_{\sqq^{2}}\\\notag
			&=\sum\limits_{a+|\nu|=r}\bn(\ell(\nu))\big[[\![S_{1,1}]\!], [\![S_{1,0}^{(\nua  p_1-1)}\oplus S_{1,0}^{(\nu )}]\!]\big]_{\sqq^{2}}\\\notag
			&=\sqq^{-1}\sum\limits_{[M]\in\cm_{r\de}}\bn(\ell(M)-1) [\![M]\!]-\frac{\sqq}{q-1}\sum\limits_{a+|\nu|=r}\bn(\ell(\nu))(q^{l_a+1}-1)\cdot q^{\sum\limits_{i\geq a}l_i}[\![S_{1,0}^{(\nu\cup\{a\})}]\!]\\\notag
			&\quad+{\sqq}\sum\limits_{a+|\nu|=r}\bn(\ell(\nu))\sum_{u\triangleleft\nu}\phi_{a,u}[\![S_{1,0}^{(\nua  p_1-1)}\oplus S_{1,0}^{(u p_1-1)}\oplus S_{1,0}^{(\nu\backslash \{u\})}]\!]*[K_{\alpha_{11}}].
		\end{align}
		
		Let $\lambda=\nu\cup\{a\}=(1^{l_1'},2^{l_2'},\dots,n^{l_n'})$. Then $l_i'=l_i$ for $i\neq a$ and $l_a'=l_a+1$. Hence
		\begin{align}
			\label{eq:coeff1}
			&\sum\limits_{a+|\nu|=r}\bn(\ell(\nu))(q^{l_a+1}-1)\cdot q^{\sum\limits_{i>a}l_i}[\![S_{1,0}^{(\nu\cup\{a\})}]\!]
			\\\notag
			&=\sum\limits_{|\lambda|=r}\bn(\ell(\lambda)-1)\sum_{a\triangleleft\lambda}(q^{l'_a}-1)\cdot q^{\sum\limits_{i>a}l'_i}[\![S_{1,0}^{(\lambda)}]\!]\\\notag
			&=\sum\limits_{|\lambda|=r}\bn(\ell(\lambda)-1)\sum_{a\triangleleft\lambda}\Big(q^{\sum\limits_{i\geq a}l'_i}- q^{\sum\limits_{i>a}l'_i}\Big)[\![S_{1,0}^{(\lambda)}]\!]\\\notag
			&=\sum\limits_{|\lambda|=r}\bn(\ell(\lambda)-1)
			(q^{\ell(\lambda)}-1)[\![S_{1,0}^{(\lambda)}]\!]=-\sum\limits_{|\lambda|=r}\bn(\ell(\lambda))[\![S_{1,0}^{(\lambda)}]\!]
		\end{align}
		by using $\bn(\ell(\lambda)-1)(1-q^{\ell(\lambda)})=\bn(\ell(\lambda))$.
		
		Then the desired formula follows by plugging \eqref{eq:coeff1} into \eqref{eq:B110B11-r}.
	\end{proof}

	\begin{lemma}
		\label{lem:image2}
		For any {$r\geq 3$}, we have
		\begin{align}
			\label{eq:image2}
			&\frac{1}{(q-1)^2}\big[\haB_{[1,1],-r+1}, \haB_{[1,1],-1}\big]_{\sqq^{2}}*[K_{r\de-2\alpha_{11}}]\\\notag
			&={\sqq}\sum\limits_{a+|\nu|=r}\bn(\ell(\nu))\sum_{u\triangleleft\nu}\phi_{a,u}[\![S_{1,0}^{(\nua  p_1-1)}\oplus S_{1,0}^{(u p_1-1)}\oplus S_{1,0}^{(\nu\backslash \{u\})}]\!]
			\\\notag
			&\quad+\sqq^{-3}\sum\limits_{[M]\in\cm_{(r-2)\de}}\bn(\ell(M)-1) [\![M]\!]\ast[K_{\delta-\alpha_{11}}]+
			\frac{\sqq}{q(q-1)}\sum\limits_{|\lambda|=r-2}\bn(\ell(\lambda))[\![S_{1,0}^{(\lambda)}]\!]\ast
			[K_{\delta-\alpha_{11}}].
		\end{align}
	\end{lemma}
	
	\begin{proof}
		Let $[M]:=[S_{1,0}^{(\nua  p_1-1)}\oplus S_{1,0}^{(\mu )}]\in\cm_{(r-1)\delta-\alpha_{11}}$. Then we have $|\mu|+a=r-1$ and $$\langle \delta-\alpha_{11}, \widehat{M}\rangle=\langle \widehat{M}, \delta-\alpha_{11}\rangle=\langle -\alpha_{11}, -\alpha_{11}\rangle=1.$$ Observe that $\Ext^1(S_{1,0}^{(p_1-1)}, M)=0$ and $\Hom(S_{1,0}^{(p_1-1)}, M)\cong \mathbf{k}$. For any non-zero map $f:S_{1,0}^{(p_1-1)}\rightarrow M$, we have $f$ is injective and $\coker(f)=S_{1,0}^{(\mu\cup\{\nua -1\})}$.
		%Hence,
		% \begin{align*}
		%&[S_{1,0}^{(p_1-1)}]*[S_{1,0}^{(\nua  p_1-1)}\oplus S_{1,0}^{(\mu )}]=\sqq^{-1}[S_{1,0}^{(p_1-1)}\oplus S_{1,0}^{(\nua  p_1-1)}\oplus S_{1,0}^{(\mu  )}]+\sqq^{-1}(q-1)[S_{1,0}^{(\mu\bigcup\{\nua -1\})}]*[[K_{\de-\alpha_{11}}]].
		%\end{align*}
		Then
		\begin{align*}
			F_{S_{1,0}^{(\mu\cup\{\nua -1\})}, \ S_{1,0}^{(p_1-1)}}^{M}=1, \quad \quad F_{S_{1,0}^{(p_1-1)},\ M}^{S_{1,0}^{(p_1-1)}\oplus M}= \begin{cases}
				q^{\ell(\mu)+1},  &\text{ if }a\geq 2,\\
				(q+1)q^{\ell(\mu)},  &\text{ if } a=1.
			\end{cases}
		\end{align*}
		Let $\mu=(1^{l_1},2^{l_2},\dots,n^{l_n})$. By Lemma \ref{lem:HomExt}, a direct computation shows that
		\begin{align*}
			\frac{|\Aut(S_{1,0}^{(\mu\cup\{a-1\})})|}{|\Aut(M)|}=\begin{cases}\frac{q^{l_{a-1}+1}-1}{q(q-1)}\cdot q^{-\sum\limits_{i\geq a-1}l_i},  &\text{ if }a\geq 2,\\
				\frac{1}{q-1}\cdot q^{-\ell(\mu)},  &\text{ if } a=1.
			\end{cases}
		\end{align*}
		%Let $\mu=(1^{l_1},2^{l_2},\dots,n^{l_n})$,
		Then by Lemma \ref{ihall formula}, we have for $a\geq 2$,
		\begin{align}
			\label{eq:S10Smua}
			&[\![S_{1,0}^{(p_1-1)}]\!]*[\![M]\!]=\sqq\cdot q^{\ell(\mu)+1}[\![S_{1,0}^{(p_1-1)}\oplus M]\!]
			+\sqq^{-1}\frac{q^{l_{a-1}+1}-1}{q(q-1)}\cdot q^{-\sum\limits_{i\geq a-1}l_i}[\![S_{1,0}^{(\mu\cup\{\nua -1\})}]\!]*[K_{\de-\alpha_{11}}],
		\end{align}
		while for $a=1$,
		\begin{align}
			\label{eq:S10Smu1}
			&[\![S_{1,0}^{(p_1-1)}]\!]*[\![M]\!]=\sqq(q+1)\cdot q^{\ell(\mu)}[\![S_{1,0}^{(p_1-1)}\oplus M]\!]
			+\sqq^{-1}\frac{1}{q-1}\cdot q^{-\ell(\mu)}[\![S_{1,0}^{(\mu)}]\!]*[K_{\de-\alpha_{11}}].
		\end{align}
		%Here we use
		%$\Hom(S_{1,0}^{(\nua  p_1-1)}\oplus S_{1,0}^{(\mu)}, S_{1,0}^{(p_1-1)})=q^{\ell(\mu)+1}$ and $$\frac{\Aut(S_{1,0}^{(\mu\cup\{a-1\})})}{\Aut(S_{1,0}^{(\nua  p_1-1)}\oplus S_{1,0}^{(\mu)})}=\frac{q^{l_{a-1}+1}-1}{q(q-1)}\cdot q^{-\sum\limits_{i\geq a-1}l_i}.$$
		
		On the other hand, the middle term of any non-split extension
		%$0\to S_{1,0}^{(p_1-1)} \to L\to S_{1,0}^{(\nua  p_1-1)}\oplus S_{1,0}^{(\mu )}\to 0$, we have %
		in $\Ext^1(M, S_{1,0}^{(p_1-1)})$ has the form
		$S_{1,0}^{(\nua  p_1-1)}\oplus S_{1,0}^{((u+1) p_1-1)}\oplus S_{1,0}^{(\mu \setminus\{u \})}$ for some $u\triangleleft\mu$, and $$F_{M,S_{1,0}^{(p_1-1)}}^{S_{1,0}^{(\nua  p_1-1)}\oplus S_{1,0}^{((u+1) p_1-1)}\oplus S_{1,0}^{(\mu \setminus\{u \})}}={F_{S_{1,0}^{(\nua  p_1-1)}\oplus S_{1,0}^{(up_1)}, S_{1,0}^{(p_1-1)}}^{S_{1,0}^{(\nua  p_1-1)}\oplus S_{1,0}^{((u+1) p_1-1)}}}=\psi_{a,u+1};$$
		while for the split extension we have
		$$F_{M, S_{1,0}^{(p_1-1)}}^{S_{1,0}^{(p_1-1)}\oplus M}=\psi_{a,1}=\begin{cases}q,  &\text{ if }a\neq 1,\\
			q+1,  &\text{ if } a=1.
		\end{cases}$$
		For any non-zero map $f:M\to S_{1,0}^{(p_1-1)}$, we have $f$ is surjective and $\ker(f)=S_{1,1}^{((\nua -1) p_1)}\oplus S_{1,0}^{(\mu )}$ or $S_{1,0}^{(\nua  p_1-1)}\oplus S_{1,1}^{((u-1) p_1+1)}\oplus S_{1,0}^{(\mu\backslash \{u\})}$ for some $u\triangleleft\mu$. It is easy to see that
		$$\dim_\bfk\Ext^1(S_{1,0}^{(p_1-1)}, \ker(f))=\dim_\bfk\Hom(S_{1,0}^{(p_1-1)}, \ker(f))\leq 1,
		$$
		and equals to $0$ if and only if $\Ker(f)=S_{1,0}^{(\mu)}$.
		%$$\Ext^1(S_{1,0}^{(p_1-1)}, \ker(f))\cong \mathbf{k}\cong \Hom(S_{1,0}^{(p_1-1)}, \ker(f));$$
		%and $$\Ext^1(S_{1,0}^{(u p_1-1)}\oplus S_{1,0}^{(\lambda\backslash\{ u\})}, S_{1,1})\cong\Ext^1(S_{1,0}^{(u p_1-1)}, S_{1,1})\cong \mathbf{k}.$$
		Then by Corollary \ref{ihall formula} (2), we get for $a\geq 2$,
		\begin{align}
			\label{eq:SmuS10a}
			&[\![M]\!]*[\![S_{1,0}^{(p_1-1)}]\!]\\\notag
			&= \sqq\cdot\psi_{a,1} [\![S_{1,0}^{(p_1-1)}\oplus M]\!]+\sqq \sum_{u\triangleleft\mu} \psi_{a,u+1}[\![S_{1,0}^{(\nua  p_1-1)}\oplus S_{1,0}^{((u+1) p_1-1)}\oplus S_{1,0}^{(\mu \setminus\{u \})}]\!]\\\notag
			&\quad+\sqq^{-3} \Big([\![S_{1,1}^{((\nua -1) p_1)}\oplus S_{1,0}^{(\mu )}]\!]+ \sum_{u\triangleleft\mu}[\![S_{1,0}^{(\nua  p_1-1)}\oplus S_{1,1}^{((u-1) p_1+1)}\oplus S_{1,0}^{(\mu\backslash \{u\})}]\!]\Big)\ast[K_{\delta-\alpha_{11}}];
		\end{align}
		while for $a=1$,
		\begin{align}
			\label{eq:SmuS101}
			&[\![M]\!]*[\![S_{1,0}^{(p_1-1)}]\!]\\\notag
			&= \sqq\cdot \psi_{1,1}[\![S_{1,0}^{(p_1-1)}\oplus M]\!]+\sqq\sum_{u\triangleleft\mu} \psi_{1,u+1}[\![S_{1,0}^{(  p_1-1)}\oplus S_{1,0}^{((u+1) p_1-1)}\oplus S_{1,0}^{(\mu \setminus\{u \})}]\!]\\\notag
			&\quad+\sqq^{-1} \frac{1}{q-1}[\![S_{1,0}^{(\mu )}]\!]\ast[K_{\delta-\alpha_{11}}]+ \sqq^{-3}\sum_{u\triangleleft\mu}[\![S_{1,0}^{(p_1-1)}\oplus S_{1,1}^{((u-1) p_1+1)}\oplus S_{1,0}^{(\mu\backslash \{u\})}]\!]\ast[K_{\delta-\alpha_{11}}].
		\end{align}
		Therefore, by combining with \eqref{eq:S10Smua} and \eqref{eq:SmuS10a}, for $a\geq 2$,
		\begin{align}
			\label{ageq2}
			&\big[[\![M]\!], [\![S_{1,0}^{(p_1-1)}]\!]\big]_{\sqq^2}\\
			\notag&= \sqq\cdot\psi_{a,1}(1-q^{\ell(\mu)+1}) [\![S_{1,0}^{(p_1-1)}\oplus M]\!]+\sqq\sum_{u\triangleleft\mu} \psi_{a,u+1}[\![S_{1,0}^{(\nua  p_1-1)}\oplus S_{1,0}^{((u+1) p_1-1)}\oplus S_{1,0}^{(\mu \setminus\{u \})}]\!]\\
			\notag&\quad+\sqq^{-3} \Big([\![S_{1,1}^{((\nua -1) p_1)}\oplus S_{1,0}^{(\mu )}]\!]+ \sum_{u\triangleleft\mu}[\![S_{1,0}^{(\nua  p_1-1)}\oplus S_{1,1}^{((u-1) p_1+1)}\oplus S_{1,0}^{(\mu\backslash \{u\})}]\!]\Big)\ast[K_{\delta-\alpha_{11}}]\\
			\notag&\quad-\frac{\sqq}{q(q-1)}\cdot (q^{l_{a-1}+1}-1)q^{-\sum\limits_{i\geq a-1}l_i}[\![S_{1,0}^{(\mu\bigcup\{\nua -1\})}]\!]*[K_{\delta-\alpha_{11}}];
		\end{align}
		while by combining with \eqref{eq:S10Smu1} and \eqref{eq:SmuS101}, for $a=1$,
		\begin{align}\label{a=1}
			&\big[[\![M]\!], [\![S_{1,0}^{(p_1-1)}]\!]\big]_{\sqq^2}\\
			\notag&= \sqq\cdot\psi_{1,1}(1-q^{\ell(\mu)+1}) [\![S_{1,0}^{(p_1-1)}\oplus M]\!]+\sqq\sum_{u\triangleleft\mu} \psi_{1,u+1}[\![S_{1,0}^{(  p_1-1)}\oplus S_{1,0}^{((u+1) p_1-1)}\oplus S_{1,0}^{(\mu \setminus\{u \})}]\!]\\
			\notag&\quad+\sqq^{-3} \sum_{u\triangleleft\mu}[\![S_{1,0}^{(  p_1-1)}\oplus S_{1,1}^{((u-1) p_1+1)}\oplus S_{1,0}^{(\mu\backslash \{u\})}]\!]\ast[K_{\delta-\alpha_{11}}]\\
			\notag&\quad+\frac{\sqq}{q(q-1)}\cdot(1-q^{-\ell(\mu)+1})[\![S_{1,0}^{(\mu)}]\!]*[K_{\delta-\alpha_{11}}].
		\end{align}
		
		To sum up with \eqref{ageq2} and \eqref{a=1}, we have
		\begin{align}
			\notag
			&\frac{1}{(q-1)^2}\big[\haB_{[1,1],-r+1}, \haB_{[1,1],-1}\big]_{\sqq^{2}}*[K_{r\de-2\alpha_{11}}]\\\notag
			&=\sum\limits_{[M]\in\cm_{(r-1)\delta-\alpha_{11}}}\bn(\ell(M)-1)\big[[\![M]\!], [\![S_{1,0}^{(p_1-1)}]\!]\big]_{\sqq^2}\\\notag
			&=\sum\limits_{a+|\mu|=r-1}\bn(\ell(\mu))\big[[\![S_{1,0}^{(\nua  p_1-1)}\oplus S_{1,0}^{(\mu)}]\!], [\![S_{1,0}^{(p_1-1)}]\!]\big]_{\sqq^{2}}\\
			& =\label{eq:image3term1}\sqq\sum\limits_{a+|\mu|=r-1}\bn(\ell(\mu)+1)\psi_{a,1}[\![S_{1,0}^{( \nua p_1-1)}\oplus S_{1,0}^{(p_1-1)}\oplus S_{1,0}^{(\mu)}]\!]\\
			\label{eq:image3term2}&\quad+\sqq\sum\limits_{a+|\mu|=r-1}\bn(\ell(\mu))
			\sum_{u\triangleleft\mu}\psi_{a,u+1}[\![S_{1,0}^{(ap_1-1)}\oplus S_{1,0}^{((u+1)   p_1-1)}\oplus S_{1,0}^{(\mu\backslash \{u\})}]\!]\\
			\label{eq:image3term3}&\quad+\sqq^{-3}\sum\limits_{[M]\in\cm_{(r-2)\de}}\bn(\ell(M)-1) [\![M]\!]\ast[K_{\delta-\alpha_{11}}]
			\\
			\label{eq:image3term4}&\quad+\sum\limits_{ |\mu|=r-2}\bn(\ell(\mu))\cdot\frac{\sqq}{q(q-1)}\cdot(1-q^{-\ell(\mu)+1})[\![S_{1,0}^{(\mu)}]\!]\ast[K_{\delta-\alpha_{11}}]\\
			\label{eq:image3term5}
			&\quad+\sum\limits_{a\geq 2} \sum\limits_{ a+|\mu|=r-1}\bn(\ell(\mu))\cdot\Big(-\frac{\sqq}{q(q-1)}\cdot (q^{l_{a-1}+1}-1)\cdot q^{-\sum\limits_{i\geq a-1}l_i}[\![S_{1,0}^{(\mu\bigcup\{\nua -1\})}]\!]\Big)\ast[K_{\delta-\alpha_{11}}]
			\\\notag
			&=\text{RHS\,of\,} \eqref{eq:image2},
		\end{align}
		where the last equality follows from the following Lemmas \ref{lem:image3} and \ref{lem:image4}.
	\end{proof}
	
	\begin{lemma}
		\label{lem:image3}
		Keep the notations as above. For any {$r\geq 2$}, we have
		\begin{align}
			\label{eq:image3}
			\eqref{eq:image3term1}+ \eqref{eq:image3term2}=&\sqq\sum\limits_{a+|\nu|=r}\bn(\ell(\nu))\sum_{u\triangleleft\nu}\phi_{a,u}[\![S_{1,0}^{(a p_1-1)}\oplus S_{1,0}^{(u  p_1-1)}\oplus S_{1,0}^{(\nu\backslash \{u\})}]\!].
		\end{align}
	\end{lemma}
	
	\begin{proof} The right-hand side of \eqref{eq:image3} can be divided into two parts according to $a=1$ or not.
		For $a=1$, we have
		\begin{align*}
			&\sum\limits_{|\nu|=r-1}\bn(\ell(\nu))\sum_{u\triangleleft\nu}\phi_{1,u}[\![S_{1,0}^{( p_1-1)}\oplus S_{1,0}^{(up_1-1)}\oplus S_{1,0}^{(\nu\backslash \{u\})}]\!]\\
			\stackrel{\mu=\nu\setminus\{u\}}{=}& \sum\limits_{u+|\mu|=r-1}\bn(\ell(\mu)+1)\psi_{u,1}[\![S_{1,0}^{(u  p_1-1)}\oplus S_{1,0}^{(p_1-1)}\oplus S_{1,0}^{(\mu)}]\!]\\
			\stackrel{u\rightarrow a}{=}&\sum\limits_{a+|\mu|=r-1}\bn(\ell(\mu)+1)\psi_{a,1}[\![S_{1,0}^{(a  p_1-1)}\oplus S_{1,0}^{(p_1-1)}\oplus S_{1,0}^{(\mu)}]\!],
		\end{align*}
		where $u\rightarrow a$ means substituting $u$ by $a$;
		while for $a\geq 2$,
		\begin{align*}
			&\sum\limits_{a\geq 2,\nu;a+|\nu|=r}\bn(\ell(\nu))\sum_{u\triangleleft\nu}\phi_{a,u}[\![S_{1,0}^{(a p_1-1)}\oplus S_{1,0}^{(u  p_1-1)}\oplus  S_{1,0}^{(\nu\backslash \{u\})}]\!]\\
			&\stackrel{\nu\backslash \{u\}=\nu'}{=}\sum\limits_{a\geq 2,u}\sum\limits_{|\nu'|=r-u-a}\bn(\ell(\nu')+1)\phi_{a,u}[\![S_{1,0}^{(a p_1-1)}\oplus S_{1,0}^{(u p_1-1)}\oplus S_{1,0}^{(\nu')}]\!]\\
			&\stackrel{a=b+1}{=}\sum\limits_{b,u}\sum\limits_{|\nu'|=r-u-b-1}\bn(\ell(\nu')+1)\phi_{b+1,u}[\![S_{1,0}^{((b+1) p_1-1)}\oplus S_{1,0}^{(u p_1-1)}\oplus S_{1,0}^{(\nu')}]\!]\\
			&\stackrel{\mu=\nu'\cup\{b\}}{=}\sum\limits_{b,u}\sum\limits_{|\mu|=r-u-1}\bn(\ell(\mu))\psi_{u,b+1}[\![S_{1,0}^{(u p_1-1)}\oplus S_{1,0}^{((b+1)p_1-1)}\oplus S_{1,0}^{(\mu\setminus\{b\})}]\!]\\
			&\stackrel{u\rightarrow a; b\rightarrow u}{=}\sum\limits_{a+|\mu|=r-1}\bn(\ell(\mu))\sum_{u\triangleleft\mu} \psi_{a,u+1}[\![S_{1,0}^{(a p_1-1)}\oplus S_{1,0}^{((u+1) p_1-1)}\oplus S_{1,0}^{(\mu\setminus\{u\})}]\!].
			%=&\sum\limits_{a+|\mu|=r-1}\bn(\ell(\mu))
			%\sum_{u\triangleleft\mu}\psi_{b+1,u}[\![S_{1,0}^{(\nua  p_1-1)}\oplus S_{1,0}^{((u+1) p_1-1)}\oplus S_{1,0}^{(\mu\backslash \{u\})}]\!].
		\end{align*}
		Then the result follows.
	\end{proof}

	\begin{lemma}
		\label{lem:image4}
		Keep the notations as above. For any {$r\geq 2$}, we have
		\begin{align}
			\label{eq:image4}
			\eqref{eq:image3term4}+ \eqref{eq:image3term5}=&\frac{\sqq}{q(q-1)}\sum\limits_{|\lambda|=r-2}
			\bn(\ell(\lambda))[\![S_{1,0}^{(\lambda)}]\!]\ast[K_{\delta-\alpha_{11}}].
		\end{align}
	\end{lemma}
	
	\begin{proof}
		Recall $\mu=(1^{l_1},2^{l_2},\dots,n^{l_n})$.
		For $a\geq 2$, let $\lambda=\mu\cup\{a-1\}=(1^{l_1'},2^{l_2'},\dots,n^{l_n'})$, then $l_i'=l_i$ for $i\neq a-1$ and $l_{a-1}'=l_{a-1}+1$. Hence
		\begin{align*}
			&\sum\limits_{a\geq 2}\sum\limits_{|\mu|=r-a-1}\bn(\ell(\mu))(q^{l_{a-1}+1}-1)\cdot q^{-\sum\limits_{i\geq a-1}l_i}[\![S_{1,0}^{(\mu\cup\{a-1\})}]\!]\\
			&=\sum\limits_{|\lambda|=r-2}\bn(\ell(\lambda)-1)\sum_{a-1\triangleleft\lambda}(q^{l'_{a-1}}-1)\cdot q^{-\sum\limits_{i\geq a-1}l'_i+1}[\![S_{1,0}^{(\lambda)}]\!]\\
			&=\sum\limits_{|\lambda|=r-2}\bn(\ell(\lambda)-1)\sum_{a-1\triangleleft\lambda}q\cdot\Big(q^{-\sum\limits_{i\geq a}l'_i}- q^{-\sum\limits_{i\geq a-1}l'_i}\Big)[\![S_{1,0}^{(\lambda)}]\!]\\
			&=\sum\limits_{|\lambda|=r-2}\bn(\ell(\lambda)-1)\cdot
			q\cdot(1-q^{-\ell(\lambda)})[\![S_{1,0}^{(\lambda)}]\!]\\
			&=-\sum\limits_{|\lambda|=r-2}\bn(\ell(\lambda))\cdot q^{-\ell(\lambda)+1}[\![S_{1,0}^{(\lambda)}]\!].
		\end{align*}
		Then \eqref{eq:image4} follows.
	\end{proof}
	
	With the help of  Lemmas \ref{lem:image1} and \ref{lem:image2}, we can prove  Proposition \ref{prop:imageroot} now.
	
	\begin{proof}[Proof of Proposition \ref{prop:imageroot}]
		First, for $r=1$, it follows from \eqref{eq:Theta11}.
		
		For $r=2$, by Corollary \ref{ihall formula}, we have
		\begin{align*}
			&[\![S_{1,0}^{(p_1-1)}]\!]*[\![S_{1,0}^{(p_1-1)}]\!]=\sqq(q+1)[\![S_{1,0}^{(p_1-1)}\oplus S_{1,0}^{(p_1-1)}]\!]+\frac{\sqq^{-1}}{q-1}[K_{\delta-\alpha_{11}}].
		\end{align*}
		Hence, \begin{align*}
			&\frac{1}{(q-1)^2}\big[\haB_{[1,1],-1}, \haB_{[1,1],-1}\big]_{\sqq^{2}}*[K_{2\delta-\alpha_{11}}]\\
			&=(1-\sqq^2)[\![S_{1,0}^{(p_1-1)}]\!]*[\![S_{1,0}^{(p_1-1)}]\!]*[K_{\alpha_{11}}]\\
			&=\sqq(1-q^2)[\![S_{1,0}^{(p_1-1)}\oplus S_{1,0}^{(p_1-1)}]\!]*[K_{\alpha_{11}}]-\sqq^{-1}[K_{\delta}].
		\end{align*}
		On the other hand, by Lemma \ref{lem:image1}, we have
		\begin{align*}
			&-\frac{1}{(q-1)^2}\big[\haB_{[1,1],0}, \haB_{[1,1],-2}\Big]_{\sqq^{2}}*[K_{2\de-\alpha_{11}}]
			=\sqq^{-1}\sum\limits_{[M]\in\cm_2}\bn(\ell(M)-1) [\![M]\!]\\
			&\quad+\frac{\sqq}{q-1}\sum\limits_{|\lambda|=2}\bn(\ell(\lambda))[\![S_{1,0}^{(\lambda)}]\!]+{\sqq}(1-q^2)[\![S_{1,0}^{(p_1-1)}\oplus S_{1,0}^{(p_1-1)}]\!]*[K_{\alpha_{11}}].
		\end{align*}
		%Here we have used $|\Ext^1(S_{1,0}^{(p_1-1)}\oplus S_{1,0}^{(p_1-1)}, S_{1,1})_{S_{1,0}^{(p_1-1)}\oplus S_{1,0}^{(p_1)}}|=q^2-1$.
		It follows that
		\begin{align*}
			&-\frac{1}{(q-1)^2}\Big(\big[\haB_{[1,1],-1}, \haB_{[1,1],-1}\big]_{\sqq^{2}}+\big[\haB_{[1,1],0}, \haB_{[1,1],-2}\big]_{\sqq^{2}}\Big)*[K_{2\de-\alpha_{11}}]\\
			&=\sqq^{-1}\sum\limits_{[M]\in\cm_{2\de}}\bn(\ell(M)-1) [\![M]\!]+\frac{\sqq}{q-1}\sum\limits_{|\lambda|=2}\bn(\ell(\lambda))[\![S_{1,0}^{(\lambda)}]\!]+\sqq^{-1}[K_{\delta}],
		\end{align*}
		which also equals to $\sqq^{-1}[K_{\delta}]+\haTh_{[1,1],2}$ by (\ref{Bmuk Bmul}). This proves the result for $r=2$.

		For $r\geq 3$, by Lemmas \ref{lem:image1} and \ref{lem:image2},  we have
		
		\begin{align*}
			&\frac{1}{(q-1)^2}\Big([\haB_{[1,1],0}, \haB_{[1,1],-r}]_{\sqq^{2}}+[\haB_{[1,1],-r+1},\haB_{[1,1],-1}]_{\sqq^{2}}\Big)*[K_{r\delta-\alpha_{11}}]\\
			&=-\Big(\sqq^{-1}\sum\limits_{[M]\in\cm_{r\de}}\bn(\ell(M)-1) [\![M]\!]+\frac{\sqq}{q-1}\sum\limits_{|\lambda|=r}\bn(\ell(\lambda))[\![S_{1,0}^{(\lambda)}]\!]\Big)\\
			&\quad+\sqq^{-2}\Big(\sqq^{-1}\sum\limits_{[M]\in\cm_{(r-2)\de}}\bn(\ell(M)-1) [\![M]\!]+
			\frac{\sqq}{q-1}\sum\limits_{|\lambda|=r-2}q^{-\ell(\lambda)}\bn(\ell(\lambda))[\![S_{1,0}^{(\lambda)}]\!]\Big)*[K_\de].
		\end{align*}
		Furthermore, it follows from (\ref{Bmuk Bmul})  that
		\begin{align*}
			&[\haB_{[1,1],0}, \haB_{[1,1],-r}]_{\sqq^{2}}+[\haB_{[1,1],-r+1},\haB_{[1,1],-1}]_{\sqq^{2}}\\
			&=(q-1)^2\Big(-\haTh_{[1,1],r}+\sqq^{-2}\haTh_{[1,1],r-2}*[K_\de]\Big)*[K_{-r\delta+\alpha_{11}}].
		\end{align*}
		Then the result follows by induction. % by noting that $\haTh_{[1,1],0}=\frac{1}{\sqq-\sqq^{-1}}$.
	\end{proof}

	%%%%%%%%%%%%%%
	\subsection{Proof of Lemma \ref{Theta 11 times co}}
	\label{subsec:proofTH11timesO}

	\begin{lemma}\label{middle term for lower case}
		Assume $0<k<r$ and $M,N \in\mathcal{T}_{\bla_1}$.
		\begin{itemize}
			\item[(1)] If $\Ext^1(N, S_{1,0}^{(kp_1)})_M\neq 0$, then $[M]\in\mathcal{M}_{r\de}$ if and only if $[N]\in\mathcal{M}_{(r-k)\de}$.
			\item[(2)] If $\Ext^1(N, S_{1,0}^{(kp_1-1)})_M\neq 0$, then $[M]\in\mathcal{M}_{r\de}$ if and only if $[N]\in\mathcal{M}_{(r-k)\de+\alpha_{11}}$.
		\end{itemize}
	\end{lemma}
	
	\begin{proof}
		Similar to the proof of Corollary \ref{extension between line bundles}, %each element $[M]\in\mathcal{M}_{r\de}$ (resp. $[M]\in\mathcal{M}_{r\de+\alpha_{11}}$) belongs to the image of $F_{\mathbb{X},\mathbb{Y}}$, hence
		we can restrict to consider the weighted projective line of weight type $\bp=(2,1)$. Then the result follows by noting that $S_{1,1}$ is a composition factor of ${\rm{top}}(M)$ with multiplicity one for $[M]\in\cm_{r\de}$.
	\end{proof}
	
	\begin{lemma}
		\label{general formula for sum of lN=lM+1 for hall number 0}
		For any $k>0$ and any $N\in\mathcal{T}_{\bla_1}$, we have
		\begin{align*}
			\sum_{M\in\ct_{\bla_1};\ell(M)=\ell(N)+1}F_{M, \co}^{\co(k\vec{x}_1)\oplus N}=\frac{|\Hom(\co, N)|}{|\Hom(\co, {\rm{top}}(N)|}.
		\end{align*}
	\end{lemma}
	
	\begin{proof}
		Consider the exact sequence in $\coh\X$ of the following forms:
		\begin{align}
			\label{ses1}
			0 \longrightarrow \co \xrightarrow{(\phi, \psi)^t} \co(k\vec{x}_1)\oplus N \longrightarrow M\longrightarrow 0.
		\end{align}
		Then $M\in\ct_{\bla_1}$ if and only if
		$\phi =\kappa x_1^{k}$ for some non-zero scalar $\kappa$.
		Note that the map $(\phi, \psi)^t$ in \eqref{ses1} is injective if and only if $\phi$ is injective; there are $q-1$ many such $\phi$'s. Assume $N$ has a decomposition $N=\bigoplus_{i} N_{i}$ with each $N_i$ indecomposable, and write $\psi=(\psi_i)_{i}$ with $\psi_i:\co\to N_i$.
		Then one can see that $\ell(M)=\ell(N)$ or $\ell(N)+1$. Moreover, $\ell(M)=\ell(N)+1$ if and only if each $\psi_{i}$ is not surjective; there are $\frac{|\Hom(\co, N_i)|}{|\Hom(\co, {\rm{top}}(N_i)|}$ many such $\psi_{i}$'s.
		Therefore,
		\begin{align*}
			\sum_{M\in\ct_{\bla_1};\ell(M)=\ell(N)+1}F_{M, \co}^{\co(k\vec{x}_1)\oplus N}=\frac{q-1}{|\Aut(\co)|}\cdot\prod_i\frac{|\Hom(\co, N_i)|}{|\Hom(\co, {\rm{top}}(N_i)|}=\frac{|\Hom(\co, N)|}{|\Hom(\co, {\rm{top}}(N)|}.
		\end{align*}
	\end{proof}

	\begin{corollary} Assume $0<k<r$. Then for any  $[N_1]\in\cm_{r\de}$, $[N_2]\in\cm_{r\de-\alpha_{11}}$ and any partition $\mu$,
		\begin{align}
			\label{sum formula of N1}
			&\sum\limits_{[M]\in\cm_{(r+k)\de};\ell(M)=\ell(N_1)+1}F_{M, \co}^{\co(k\vec{c})\oplus N_1}=q^{r-\ell(N_1)+1};\\
			\label{sum formula of N2}
			&\sum\limits_{[M]\in\cm_{(r+k)\de};\ell(M)=\ell(N_2)+1}F_{M, \co}^{\co(k\vec{c}+\vec{x}_1)\oplus N_2}=q^{r-\ell(N_2)};\\
			&\sum\limits_{\lambda; \ell(\lambda)=\ell(\mu)+1}F_{S_{1,0}^{(\lambda)},\co}^{\co(k\vec{c})\oplus S_{1,0}^{(\mu)}}=q^{|\mu|-\ell(\mu)}.\label{sum formula of lambda}
		\end{align}
	\end{corollary}
	
	\begin{proof}
		We only prove \eqref{sum formula of N1}, since the proofs for the other two formulas are similar. For any $M\in\ct_{\bla_1}$, $F_{M, \co}^{\co(k\vec{c})\oplus N_1}\neq 0$ if and only if $[M]\in\cm_{(r+k)\de}$ by Lemma \ref{middle term for lower case}. Moreover, $\dim_\bfk\Hom(\co, N_1)=\langle \widehat{\co}, r\de\rangle=r$ and $\dim_\bfk\Hom(\co, {\rm{top}}(N_1))=\ell(N_1)-1$. Then the result follows from Lemma \ref{general formula for sum of lN=lM+1 for hall number 0}.
	\end{proof}

	Now we can prove Lemma \ref{Theta 11 times co}.
	
	\begin{proof}[Proof of Lemma \ref{Theta 11 times co}]
		Recall $\mathcal{S}$ defined in \eqref{def:S}. 
		Recall from \eqref{def:Mr} and \eqref{eq:TH11r} that each term of $\haTh_{[1,1],r}$ comes from ${}^\perp\mathcal{S}$. Hence each extension term in the product $\haTh_{[1,1],r}\ast [\co]$ also comes from ${}^\perp\mathcal{S}$. Similar to the proof of Corollary \ref{extension between line bundles}, by considering the degrees and the ranks,
		%of the extension terms in the product $\haTh_{[1,1],r}\ast [\co]$,
		we see that only the following terms can appear in $\haTh_{[1,1],r}\ast [\co]$:
		\begin{align*}
			[\![\co(k\vec{c})\oplus  S_{1,0}^{(\lambda)}]\!],  \qquad [\![\co(k\vec{c})\oplus M]\!],\qquad[\![\co(k\vec{c}+\vec{x}_1)\oplus N]\!],
		\end{align*}
		where $k>0$, $\lambda$ is a partition of $r-k$, $[M]\in\cm_{(r-k)\de}$ and $[N]\in\cm_{(r-k)\de-\alpha_{11}}$.

		Recall from \eqref{eq:TH11r} that
		\begin{align*}
			\haTh_{[1,1],r}=
			&\frac{\sqq}{q-1} \sum\limits_{|\lambda|=r}\bn(\ell(\lambda))  [\![S_{1,0}^{(\lambda)}]\!]
			+\sqq^{-1} \sum\limits_{[M]\in\cm_{r\de}}\bn(\ell(M)-1) [\![M]\!].
		\end{align*}
		
		For any partition $\lambda$ of $r$, we have $\langle \widehat{S_{1,0}^{(\lambda)}}, \widehat{\co} \rangle=-r$. Hence,
		\begin{align}
			\label{eq:SmultO}
			&\frac{\sqq}{q-1}\sum\limits_{|\lambda|=r}\bn(\ell(\lambda))   [\![S_{1,0}^{(\lambda)}]\!]\ast [\![\co]\!]\\\notag
			&=\frac{\sqq}{q-1}\sum\limits_{|\lambda|=r}\sum\limits_{k\geq 0}\sum\limits_{|\mu|=r-k}\bn(\ell(\lambda))\cdot\sqq^{-r} F_{S_{1,0}^{(\lambda)},\co}^{\co(k\vec{c})\oplus S_{1,0}^{(\mu)}}[\![\co(k\vec{c})\oplus S_{1,0}^{(\mu)} ]\!]\\\notag
			&=\frac{\sqq^{-r+1}}{q-1}\sum\limits_{k\geq 0}\sum\limits_{|\mu|=r-k}\sum\limits_{|\lambda|=r}\bn(\ell(\lambda)) \cdot F_{S_{1,0}^{(\lambda)},\co}^{\co(k\vec{c})\oplus S_{1,0}^{(\mu)}}[\![\co(k\vec{c})\oplus S_{1,0}^{(\mu)} ]\!]\\\notag
			&=\frac{\sqq^{-r+1}}{q-1}\sum\limits_{|\mu|=r}\bn(\ell(\mu))\cdot F_{S_{1,0}^{(\mu)},\co}^{\co\oplus S_{1,0}^{(\mu)}}[\![\co\oplus S_{1,0}^{(\mu)} ]\!]+\frac{\sqq^{-r+1}}{q-1}\sum\limits_{k>0}
			\sum\limits_{|\mu|=r-k}\\\notag
			&\quad \Big(\bn(\ell(\mu)+1)\sum\limits_{\lambda; \ell(\lambda)=\ell(\mu)+1}F_{S_{1,0}^{(\lambda)},\co}^{\co(k\vec{c})\oplus S_{1,0}^{(\mu)}}+\bn(\ell(\mu))\sum\limits_{\lambda; \ell(\lambda)=\ell(\mu)}F_{S_{1,0}^{(\lambda)},\co}^{\co(k\vec{c})\oplus S_{1,0}^{(\mu)}}\Big)[\![\co(k\vec{c})\oplus S_{1,0}^{(\mu)} ]\!]\\\notag
			&\stackrel{\eqref{sum formula of lambda}}{=}\frac{\sqq^{-r+1}}{q-1}\sum\limits_{|\mu|=r}\bn(\ell(\mu))\cdot q^{r}[\![\co\oplus S_{1,0}^{(\mu)} ]\!]+\frac{\sqq^{-r+1}}{q-1}\sum\limits_{k>0}
			\sum\limits_{|\mu|=r-k}\\\notag
			&\quad \Big(\bn(\ell(\mu)+1)\cdot q^{|\mu|-\ell(\mu)}+\bn(\ell(\mu))\cdot\big(q^{|\mu|}-q^{|\mu|-\ell(\mu)}\big)\Big)[\![\co(k\vec{c})\oplus S_{1,0}^{(\mu)} ]\!]\\\notag
			&=\frac{\sqq^{r+1}}{q-1}\sum\limits_{|\mu|=r}\bn(\ell(\mu))[\![\co\oplus S_{1,0}^{(\mu)} ]\!]-\sum\limits_{k>0}
			\sum\limits_{|\mu|=r-k}\bn(\ell(\mu))\cdot \sqq^{r-2k+1}[\![\co(k\vec{c})\oplus S_{1,0}^{(\mu)} ]\!].
		\end{align}
		
		For any $[M]\in\cm_{r\de}$, we have $\langle \widehat{M}, \widehat{\co} \rangle=-r$ and then by Lemma \ref{middle term for lower case},
		\begin{align}
			\label{M*O}[\![M]\!]\ast[\![\co]\!]  =&
			%\sqq^{\langle M, \co\rangle}\sum\limits_{[X]}F_{M,\co}^{X}[\![X]\!]=
			\sqq^{-r}\sum\limits_{k\geq 0}\sum\limits_{[N_1]\in\cm_{(r-k)\de}}F_{M,\co}^{\co(k\vec{c})\oplus N_1}[\![\co(k\vec{c})\oplus N_1]\!]\\\notag
			&+  \sqq^{-r}\sum\limits_{k\geq 0}\sum\limits_{[N_2]\in\cm_{(r-k)\de-\alpha_{11}}}F_{M,\co}^{\co(k\vec{c}+\vec{x}_1)\oplus N_2}[\![\co(k\vec{c}+\vec{x}_1)\oplus N_2]\!].
		\end{align}
		
		For any $k\geq 0$ and $[N_2]\in\cm_{(r-k)\de-\alpha_{11}}$, the coefficient of $[\![\co(k\vec{c}+\vec{x}_1)\oplus N_2]\!]$ in $\haTh_{[1,1],r}*[\![\co]\!]$ equals to (up to a power of $\sqq$)
		\begin{align}\label{zero coefficient}&\sum\limits_{[M]\in\cm_{r\de}}\bn(\ell(M)-1)\cdot F_{M,\co}^{\co(k\vec{c}+\vec{x}_1)\oplus N_2}\\
			\notag&=\bn(\ell(N_2))\sum\limits_{[M]\in\cm_{r\de},\ell(M)=\ell(N_2)+1}F_{M,\co}^{\co(k\vec{c}+\vec{x}_1)\oplus N_2}+\bn(\ell(N_2)-1)\sum\limits_{[M]\in\cm_{r\de},\ell(M)=\ell(N_2)}F_{M,\co}^{\co(k\vec{c}+\vec{x}_1)\oplus N_2}\\
			\notag&\stackrel{\eqref{sum formula of N2}}{=}\bn(\ell(N_2))\cdot q^{r-k-\ell(N_2)}+\bn(\ell(N_2)-1)\cdot (q^{r-k}-q^{r-k-\ell(N_2)})=0.
		\end{align}
		Similarly, for $k>0$ and $[N_1]\in\cm_{(r-k)\de}$ we have
		%$$\sum\limits_{[M]\in\cm_{r\de},\ell(M)=\ell(N_1)+1}F_{M,\co}^{\co(k\vec{c})\oplus N_1}=q^{r-k-\ell(N_1)+1},$$
		%hence
		\begin{align}\label{k>0M*co}&\sum\limits_{[M]\in\cm_{r\de}}\bn(\ell(M)-1)\cdot F_{M,\co}^{\co(k\vec{c})\oplus N_1}\\
			\notag&=\bn(\ell(N_1))\sum\limits_{[M]\in\cm_{r\de},\ell(M)=\ell(N_1)+1}F_{M,\co}^{\co(k\vec{c})\oplus N_1}+\bn(\ell(N_1)-1)\sum\limits_{[M]\in\cm_{r\de},\ell(M)=\ell(N_1)}F_{M,\co}^{\co(k\vec{c})\oplus N_1}\\
			\notag&\stackrel{\eqref{sum formula of N1}}{=}\bn(\ell(N_1))\cdot q^{r-k-\ell(N_1)+1}+\bn(\ell(N_1)-1)\cdot(q^{r-k}-q^{r-k-\ell(N_1)+1})\\
			\notag&=\bn(\ell(N_1)-1)\cdot q^{r-k-\ell(N_1)+1}(1-q^{\ell(N_1)}+q^{\ell(N_1)-1}-1)\\
			\notag&=\bn(\ell(N_1)-1)\cdot q^{r-k}(1-q);
		\end{align}
		while for $k=0$ and $[N_1]\in\cm_{r\de}$,
		\begin{align}\label{k=0M=N}
			\sum\limits_{[M]\in\cm_{r\de}}\bn(\ell(M)-1)\cdot F_{M,\co}^{\co\oplus N_1}
			=\bn(\ell(N_1)-1)\cdot F_{N_1,\co}^{\co\oplus N_1}
			=\bn(\ell(N_1)-1)\cdot q^{r}.
		\end{align}
		
		Combining with \eqref{M*O}--\eqref{k=0M=N}, we obtain
		\begin{align}
			\label{eq:MO}
			&\sqq^{-1} \sum\limits_{[M]\in\cm_{r\de}}\bn(\ell(M)-1)[\![M]\!]\ast [\![\co]\!]\\\notag
			&=\sqq^{-r-1}\sum\limits_{k\geq 0}\sum\limits_{[N]\in\cm_{(r-k)\de}}\Big(\sum\limits_{[M]\in\cm_{r\de}}\bn(\ell(M)-1) F_{M,\co}^{\co(k\vec{c})\oplus N}\Big)[\![\co(k\vec{c})\oplus N]\!]
			\\\notag
			&=\sqq^{r-1}\sum\limits_{[N]\in\cm_{r\de}}\bn(\ell(N)-1)[\![\co\oplus N]\!] \\
			&\quad+\sum\limits_{k>0}\sum\limits_{[N]\in\cm_{(r-k)\de}}\bn(\ell(N)-1)\cdot(1-q)\sqq^{r-2k-1}[\![\co(k\vec{c})\oplus N]\!].
		\end{align}

		Combining \eqref{eq:MO} and \eqref{eq:SmultO}, we finish the proof of Lemma \ref{Theta 11 times co}. %Using \eqref{eq:TH11r}, the desired formula follows by combining  \eqref{eq:MO} and \eqref{eq:SmultO}.
	\end{proof}

	%%%%%%%%%%%%%%%%5
	\subsection{Proof of Lemma \ref{co times Theta 11}}
	\label{subsec:proofOtimesTH11}

	For any $[M]\in\cm_{r\de}$, we have $\langle\widehat{\co}, \widehat{M} \rangle=r$, $\Ext^1(\co, M)=0$ and $\Hom(\co, M)\cong\bfk^r$. For any non-zero morphism $f:\co\rightarrow M$, we have $\Im(f)\cong S_{1,0}^{(kp_1)}$ or $S_{1,0}^{(k  p_1-1)}$ for some $k>0$.
	
	If $\Im(f)=S_{1,0}^{(kp_1)}$, then $\ker(f)=\co(-k\vec{c})$ and $[\coker(f)]\in\cm_{(r-k)\de}$ by Lemma \ref{middle term for lower case}. For any $[N]\in\cm_{(r-k)\de}$, we have
	\begin{align*}&|\{f:\co\rightarrow M\,|\,\ker(f)\cong \co(-k\vec{c}), \coker(f)\cong N\}|\\&=|\{g:S_{1,0}^{(kp_1)}\rightarrowtail M\,|\, \coker(g)\cong N\}|\\
		&=F_{N,S_{1,0}^{(kp_1)}}^{M}\cdot|\Aut(S_{1,0}^{(kp_1)})|.
	\end{align*}
	Similarly, if $\Im(f)=S_{1,0}^{(kp_1-1)}$, then $\ker(f)=\co(-k\vec{c}+\vec{x}_1)$ and $[\coker(f)]\in\cm_{(r-k)\de+\alpha_{11}}$. For any $[N']\in\cm_{(r-k)\de+\alpha_{11}}$,
	\begin{align*}&|\{f:\co\rightarrow M\,|\,\ker(f)\cong \co(-k\vec{c}+\vec{x}_1), \coker(f)\cong N'\}|\\&=|\{g:S_{1,0}^{(kp_1-1)}\rightarrowtail M\,|\, \coker(g)\cong N'\}|\\
		&=F_{N',S_{1,0}^{(kp_1-1)}}^{M}\cdot|\Aut(S_{1,0}^{(kp_1-1)})|.
	\end{align*}
	Then by \eqref{Hallmult1}, we have
	\begin{align*}
		&[\co]\ast[M]
		=\sqq^{-r}\sum\limits_{k\geq 0}\sum\limits_{[N]\in\cm_{(r-k)\de}} F_{N,S_{1,0}^{(kp_1)}}^{M}|\Aut(S_{1,0}^{(kp_1)})|[\co(-k\vec{c})\oplus N]\ast [K_{k\delta}]\\
		&\quad+\sqq^{-r}\sum\limits_{k\geq 1}\sum\limits_{[N']\in\cm_{(r-k)\de+\alpha_{11}}} F_{N',S_{1,0}^{(kp_1-1)}}^{M}|\Aut(S_{1,0}^{(kp_1-1)})|[\co(-k\vec{c}+\vec{x}_1)\oplus N']\ast [K_{k\delta-\alpha_{11}}];
		% &=\sqq^{-r}\sum\limits_{k\geq 0}\sum\limits_{[N]} \frac{|\Ext^1(N, S_{1,0}^{(kp_1)})_{M}|\cdot|\Aut(M)|}{|\Hom(N, S_{1,0}^{(kp_1)})|\cdot|\Aut(N)|}\cdot[\co(-k\vec{c})\oplus N]\ast [K_{k\delta}]\\
		%  &\quad+\sqq^{-r}\sum\limits_{k\geq 1}\sum\limits_{[N']} \frac{|\Ext^1(N', S_{1,0}^{(kp_1-1)})_{M}|\cdot|\Aut(M)|}{|\Hom(N', S_{1,0}^{(kp_1-1)})|\cdot|\Aut(N')|}\cdot[\co(-k\vec{c}+\vec{x}_1)\oplus N']\ast [K_{k\delta-\alpha_{11}}].\\
	\end{align*}
	or equivalently, by using the Riedtmann-Peng formula, we have
	\begin{align}\label{formula of co times M}
		&[\![\co]\!]\ast[\![M]\!]
		=\sqq^{-r}\sum\limits_{k\geq 0}\sum\limits_{[N]\in\cm_{(r-k)\de}}q^{r-k}\frac{|\Ext^1(N, S_{1,0}^{(kp_1)})_{M}|}{|\Hom(N, S_{1,0}^{(kp_1)})|}[\![\co(-k\vec{c})\oplus N]\!]*[K_{k\delta}] \\\notag
		&\quad+\sqq^{-r}\sum\limits_{k\geq 1}\sum\limits_{[N']\in\cm_{(r-k)\de+\alpha_{11}}}q^{r-k+1}\frac{|\Ext^1(N', S_{1,0}^{(kp_1-1)})_{M}|}{|\Hom(N', S_{1,0}^{(kp_1-1)})|}[\![\co(-k\vec{c}+\vec{x}_1)\oplus N']\!]* [K_{k\delta-\alpha_{11}}].
	\end{align}
	
	Similar arguments show that
	\begin{align}\label{the formula of co times S lambda}
		&[\![\co]\!]\ast  [\![S_{1,0}^{(\lambda)}]\!]= \sqq^{-r}\sum\limits_{k\geq 0}\sum\limits_{|\mu|=r-k} q^{r-k}\frac{|\Ext^1(S_{1,0}^{(\mu)}, S_{1,0}^{(kp_1)})_{S_{1,0}^{(\lambda)}}|}{|\Hom(S_{1,0}^{(\mu)}, S_{1,0}^{(kp_1)})|}\cdot[\![\co(-k\vec{c})\oplus S_{1,0}^{(\mu)}]\!]\ast [K_{k\delta}].
	\end{align}
	
	Now we can prove Lemma \ref{co times Theta 11}.

	\begin{proof}[Proof of Lemma \ref{co times Theta 11}]
		For any $k>0$ and $[L]\in\cm_{(r-k)\de+\alpha_{11}}$, we have $\langle L,S_{1,0}^{(kp_1-1)}\rangle=-1$ and $\dim_\bfk \Hom(L,S_{1,0}^{(p_1-1)})=\ell(L)-1$. Then the coefficient of $[\![\co(-k\vec{c}+\vec{x}_1)\oplus L]\!]\ast [K_{k\delta-\alpha_{11}}]$ in the product $[\![\co]\!]\ast \haTh_{[1,1],r}$ equals to (up to a power of $\sqq$)
		\begin{align}
			\label{zero coefficient for imf=S10(kp-1)}&\sum\limits_{[M]\in\cm_{r\de}}\bn(\ell(M)-1)\cdot|\Ext^1(L, S_{1,0}^{(kp_1-1)})_{M}|\\\notag
			%=&\bn(\ell(L))\sum\limits_{[M]\in\cm_{r\de},\ell(M)=\ell(N)+1}|\Ext^1(L, S_{1,0}^{(kp_1-1)})_{M}|\\
			%&\qquad+\bn(l(L)-1)\sum\limits_{[M]\in\cm_{r\de},\ell(M)=\ell(L)}|\Ext^1(L, S_{1,0}^{(kp_1-1)})_{M}|\\
			&=\bn(\ell(L))\sum\limits_{[M];\ell(M)=\ell(L)+1}|\Ext^1(L, S_{1,0}^{(kp_1-1)})_{M}|+\bn(\ell(L)-1)\sum\limits_{[M];\ell(M)=\ell(L)}|\Ext^1(L, S_{1,0}^{(kp_1-1)})_{M}|\\\notag
			&=\bn(\ell(L))\frac{|\Hom(L,S_{1,0}^{(kp_1-1)})|}{|\Hom(L,S_{1,0}^{(p_1-1)})|}
			+\bn(\ell(L)-1)\Big(|\Ext^1(L, S_{1,0}^{(kp_1-1)})|-\frac{|\Hom(L,S_{1,0}^{(kp_1-1)})|}{|\Hom(L,S_{1,0}^{(p_1-1)})|}\Big)\\\notag
			&=\bn(\ell(L)-1) \Big(\frac{|\Hom(L,S_{1,0}^{(kp_1-1)})|}{|\Hom(L,S_{1,0}^{(p_1-1)})|}(1-q^{\ell(L)})
			+\frac{|\Hom(L, S_{1,0}^{(kp_1-1)})|}{q^{\langle L, S_{1,0}^{(kp_1-1)}\rangle}}-\frac{|\Hom(L,S_{1,0}^{(kp_1-1)})|}{|\Hom(L,S_{1,0}^{(p_1-1)})|}\Big)\\\notag
			&=\bn(\ell(L)-1)\cdot |\Hom(L,S_{1,0}^{(kp_1-1)})|\cdot \Big(\frac{1-q^{\ell(L)}}{q^{\ell(L)-1}}
			+q-\frac{1}{q^{\ell(L)-1}}\Big)
			=0.
		\end{align}
		Here the second equality follows from Lemma \ref{sum of l(N)=l(M)+1 for p-1}.

		Similarly, for any $k>0$ and $[N]\in\cm_{(r-k)\de}$,  we have $\dim_\bfk \Ext^1(N,S_{1,0}^{(p_1-1)})=\ell(N)-1$. By using Lemma \ref{sum of l(N)=l(M)+1 for p},
		\begin{align}\label{coeff for imf=kdelta}&\sum\limits_{[M]\in\cm_{r\de}}\bn(\ell(M)-1)\cdot|\Ext^1(N, S_{1,0}^{(kp_1)})_{M}|\\\notag
			%=&\bn(l(N))\sum\limits_{[M]\in\cm_{r\de},\ell(M)=\ell(N)+1}|\Ext^1(N, S_{1,0}^{(kp_1)})_{M}|
			%\\
			%&+\bn(l(N)-1)\sum\limits_{[M]\in\cm_{r\de},\ell(M)=\ell(N)}|\Ext^1(N, S_{1,0}^{(kp_1)})_{M}|\\
			&=\bn(l(N))\sum\limits_{\ell(M)=\ell(N)+1}|\Ext^1(N, S_{1,0}^{(kp_1)})_{M}|
			+\bn(l(N)-1)\sum\limits_{\ell(M)=\ell(N)}|\Ext^1(N, S_{1,0}^{(kp_1)})_{M}|\\&\notag
			{=}\bn(\ell(N)) \frac{|\Hom(N,S_{1,0}^{(kp_1)})|}{q^{\ell(N)-1}}+\bn(l(N)-1)\Big(|\Ext^1(N, S_{1,0}^{(kp_1)})|-\frac{|\Hom(N,S_{1,0}^{(kp_1)})|}{q^{\ell(N)-1}})|\Big)\\\notag
			%=&\bn(\ell(N)-1)\Big( \frac{1-q^{\ell(N)}}{q^{\ell(N)-1}}|\Hom(N,S_{1,0}^{(kp_1)})|+|\Hom(N, S_{1,0}^{(kp_1)})|-\frac{|\Hom(N,S_{1,0}^{(kp_1)})|}{q^{\ell(N)-1}}\Big)\\
			&=\bn(\ell(N)-1)\cdot(1-q)\cdot|\Hom(N,S_{1,0}^{(kp_1)})|;
		\end{align}
		while for $k=0$,
		\begin{align}\label{coeff for imf=0}
			\sum\limits_{[M]\in\cm_{r\de}}\bn(\ell(M)-1)\cdot|\Ext^1(N, S_{1,0}^{(kp_1)})_{M}|
			=\bn(l(N)-1).
		\end{align}

		Therefore, combining \eqref{formula of co times M}--\eqref{coeff for imf=0}, we obtain
		
		\begin{align}\label{the formula of co times M}
			&\sqq^{-1} [\![\co]\!]\ast \sum\limits_{[M]\in\cm_{r\de}}\bn(\ell(M)-1) [\![M]\!]\\\notag
			%=& \sum\limits_{[M]\in\cm_{r\de}}\bn(\ell(M)-1) \sqq^{-1}[\![\co]\!]\ast [\![M]\!]\\
			%\stackrel{\eqref{coeff of cokc oplus N}}{=}
			&= \sqq^{-(r+1)}\sum\limits_{k\geq 0}\sum\limits_{[N]\in\cm_{(r-k)\de}}\sum\limits_{[M]\in\cm_{r\de}}\bn(\ell(M)-1) q^{r-k}\frac{|\Ext^1(N, S_{1,0}^{(kp_1)})_{M}|}{|\Hom(N, S_{1,0}^{(kp_1)})|}[\![\co(-k\vec{c})\oplus N]\!]\ast [K_{k\delta}]\\\notag
			% & +\sum\limits_{k\geq 1}\sum\limits_{[N']}\frac{\sqq^{-(r+1)}\cdot|\Hom(\co(-k\vec{c}+\vec{x}_1),N')|}{|\Hom(N', S_{1,0}^{(kp_1-1)})|}\cdot\\
			% &\quad\sum\limits_{[M]\in\cm_{r\de}}\bn(\ell(M)-1)\cdot|\Ext^1(N', S_{1,0}^{(kp_1-1)})_{M}|\cdot[\![\co(-k\vec{c}+\vec{x}_1)\oplus N']\!]\ast [K_{k\delta-\alpha_{11}}]\\
			%\stackrel{\eqref{coeff for imf=kdelta}}{=}
			&=\sum\limits_{[N]\in\cm_{r\de}} \bn(\ell(N)-1)\sqq^{r-1}[\![\co\oplus N]\!]\\\notag
			&\quad+\sum\limits_{k> 0}\sum\limits_{[N]\in\cm_{(r-k)\de}}  (1-q) \sqq^{r-2k-1} \bn(\ell(N)-1)[\![\co(-k\vec{c})\oplus N]\!]\ast [K_{k\delta}].
		\end{align}

		On the other hand, for any $k>0$ and $|\mu|=r-k$, by using Lemma \ref{sum of l(N)=l(M)+1 for p} we obtain
		%\red{Using Lemma \ref{sum of l(N)=l(M)+1 for p}, by similar arguments as above, for $k>0$ and $|\mu|=r-k$, we have}
		\begin{align}\label{coeff of co and lambda for k not 0}
			&\sum\limits_{|\lambda|=r}\bn(\ell(\lambda))\cdot|\Ext^1(S_{1,0}^{(\mu)}, S_{1,0}^{(kp_1)})_{S_{1,0}^{(\lambda)}}|\\\notag
			&=\bn(\ell(\mu)+1)\sum\limits_{\ell(\lambda)=\ell(\mu)+1}|\Ext^1(S_{1,0}^{(\mu)}, S_{1,0}^{(kp_1)})_{S_{1,0}^{(\lambda)}}|+\bn(\ell(\mu))\sum\limits_{\ell(\lambda)=\ell(\mu)}|\Ext^1(S_{1,0}^{(\mu)}, S_{1,0}^{(kp_1)})_{S_{1,0}^{(\lambda)}}|\\\notag
			&=\bn(\ell(\mu)+1)\cdot \frac{|\Hom(S_{1,0}^{(\mu)}, S_{1,0}^{(kp_1)})|}{|\Hom(S_{1,0}^{(\mu)}, S_{1,0}^{(p_1)})|}
			+\bn(\ell(\mu))\Big(|\Ext^1(S_{1,0}^{(\mu)}, S_{1,0}^{(kp_1)})|-\frac{|\Hom(S_{1,0}^{(\mu)}, S_{1,0}^{(kp_1)})|}{|\Hom(S_{1,0}^{(\mu)}, S_{1,0}^{(p_1)})|}\Big)\\\notag
			&=\bn(\ell(\mu))\cdot|\Hom(S_{1,0}^{(\mu)}, S_{1,0}^{(kp_1)})|\cdot \Big(\frac{1-q^{\ell(\mu)+1}}{q^{\ell(\mu)}}+1-\frac{1}{q^{\ell(\mu)}}\Big)\\\notag
			&=(1-q)\cdot \bn(\ell(\mu))\cdot |\Hom(S_{1,0}^{(\mu)}, S_{1,0}^{(kp_1)})|,\end{align}
		while for $k=0$ and $|\mu|=r$, \begin{align}\label{coeff of co and lambda for k=0}
			\sum\limits_{|\lambda|=r}\bn(\ell(\lambda))\cdot|\Ext^1(S_{1,0}^{(\mu)}, S_{1,0}^{(kp_1)})_{S_{1,0}^{(\lambda)}}|=\bn(\ell(\mu)).
		\end{align}
		Combining \eqref{the formula of co times S lambda}, \eqref{coeff of co and lambda for k not 0} and \eqref{coeff of co and lambda for k=0}, we obtain
		\begin{align}\label{the formula of sum of co times S lambda}
			&\frac{\sqq}{q-1}[\![\co]\!]\ast \sum\limits_{|\lambda|=r}\bn(\ell(\lambda))  [\![S_{1,0}^{(\lambda)}]\!]\\\notag
			% =&  \sum\limits_{|\lambda|=r}\bn(\ell(\lambda)) \frac{\sqq}{q-1} [\![\co]\!]\ast[\![S_{1,0}^{(\lambda)}]\!]\\
			&= \frac{\sqq^{-r+1}}{q-1}\sum\limits_{k\geq 0}\sum\limits_{|\mu|=r-k}\sum\limits_{|\lambda|=r}\bn(\ell(\lambda))\cdot q^{r-k}\frac{|\Ext^1(S_{1,0}^{(\mu)}, S_{1,0}^{(kp_1)})_{S_{1,0}^{(\lambda)}}|}{|\Hom(S_{1,0}^{(\mu)}, S_{1,0}^{(kp_1)})|}\cdot[\![\co(-k\vec{c})\oplus S_{1,0}^{(\mu)}]\!]\ast [K_{k\delta}]\\\notag
			&=\frac{\sqq^{r+1}}{q-1}\sum\limits_{|\mu|={r}}  \bn(\ell(\mu))\cdot[\![\co\oplus S_{1,0}^{(\mu)}]\!]
			-\sum\limits_{k>0}\sum\limits_{|\mu|={r-k}}\sqq^{r-2k+1}\cdot \bn(\ell(\mu))\cdot[\![\co(-k\vec{c})\oplus S_{1,0}^{(\mu)}]\!]\ast [K_{k\delta}].
		\end{align}
		
		Combining  \eqref{the formula of co times M} and \eqref{the formula of sum of co times S lambda}, we finish the proof of Lemma \ref{co times Theta 11}.
	\end{proof}

	%%%%%%%
	%%%%%%%

\end{document}